FACHBEREICH
MATHEMATIK

# Courcelle's Theorem

## A Self-Contained Proof and a Path-Width Variant

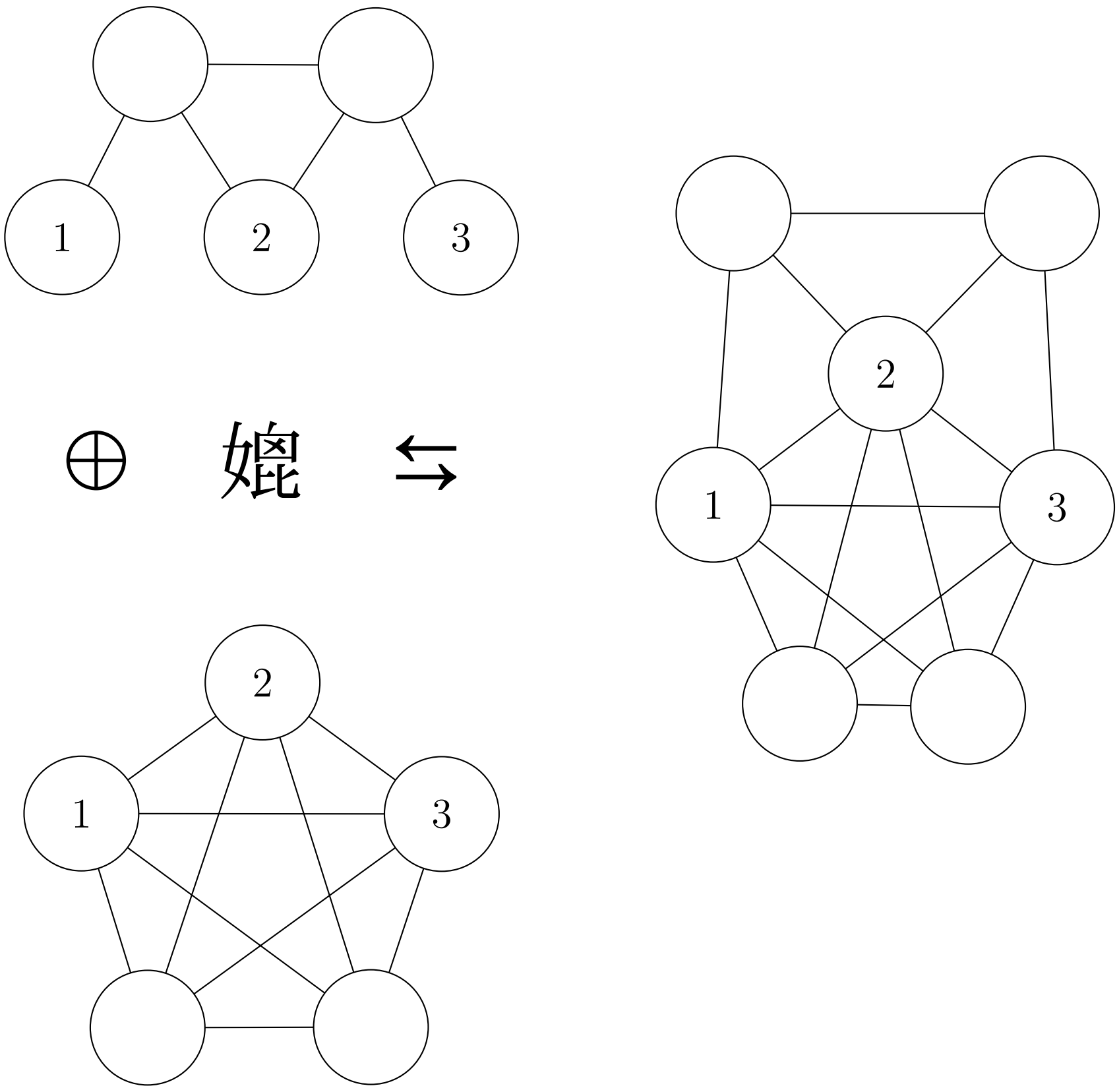

Master's Thesis

Adrian Rettich

September 16th, 2022

supervised by:
Prof. Dr. Sven Oliver Krumke
and M. Sc. Oliver Bachtler

TECHNISCHE UNIVERSITÄT
KAISERSLAUTERN

# Contents

# Chapter 1
# Introduction

Courcelle's Theorem ([Cou90]) is an important result in graph theory, proving the existence of linear-time algorithms for many decision problems on graphs whose tree-width is bounded by a constant. The result has been extended ([CM93]) to cover even more problems (namely optimisation problems) and more graph classes (namely all those that can be built using a hyperedge replacement grammar). Due to its generality, however, many a graph theorist might shy away from trying to apply it to their special case because of the "black box" nature of the theorem. Existing literature largely fails to provide an accessible insight into the workings of the theorem, requiring the reader to be already familiar with hypergraph grammars, typed algebras, or other constructions not traditionally taught together with graph theory.

The purpose of the current text is twofold: to provide an explanation and step-by-step proof of Courcelle's Theorem as applied to graphs of tree-width bounded by a constant, and to show explicitly (on the example of path-width[1]) how to apply the same principles to other graph classes.

We present the proof of Courcelle's Theorem in a way that does not assume any particular knowledge on the part of the reader except a basic understanding of mathematics and possibly the fundamentals of graph theory. A more advanced reader may therefore find it desirable to skip certain chapters so as not to become the victim of boredom. Such a reader is advised to peruse the section “How to Read this Thesis” on page 10.

[1] A note for the reader already familiar with Courcelle's Theorem: the traditional formulation does allow to encode the path-width of a graph, making this redundant in a certain sense. However, encoding the path-width by monadic second-order logic is done by encoding the (finitely many) forbidden minors for the desired path-width, which are at the time of this writing known only up to a path-width of 2.

There exist in the literature two basic techniques to prove Courcelle’s Theorem: the one used by Courcelle himself utilises algebraic methods to show that certain classes of graphs have properties which then let one construct an automaton that solves the membership problem of these classes in linear time. An alternative approach is seen in [KL09], where Kneis and Langer use Hintikka games to arrive at a practical implementation of Courcelle’s Theorem with better constants in the runtime, but slightly restricted scope.

We have opted for a variation on the former approach. We introduce the algebraic and automata-theoretic basics without assuming any prior knowledge on the part of the reader and show a theoretical proof of Courcelle’s Theorem. We then construct explicitly the corollaries for the special cases of classes of graphs of tree-width bounded by a constant before doing the same with respect to path-width. We have chosen to omit from our construction certain parts of the main statement, namely the extension of monadic second-order logic to counting monadic second-order logic and the case where our graphs are not free of loops. Each of these topics has been exiled into an appendix, where only the most tireless of readers shall find the formal constructions needed to integrate them into the framework of the main chapters.

It should be noted that our goal is to provide a theoretical approach to Courcelle’s Theorem. Algorithms and runtime considerations will be provided where they are of interest, but questions of computability and concrete implementation are not discussed. Furthermore, there are certain generalisations and extensions which are, for the sake of brevity, not included in this text; these are discussed in the closing remarks on page 213.

We hope that the present work can provide a gentle introduction to the topic for the interested graph theorist who has only in passing heard of this “monadic second-order logic” and wishes to learn more without having to compile the pertinent material themselves, and that the examples and constructions given turn Courcelle’s Theorem from a black box into a useful implement in the reader’s personal tool shed.

# Chapter 2 Overview

Before we dive into the deep theory, the reader might want to get their feet wet by splashing around the shallow end for a bit. This chapter shall serve as a primer on what it is we are trying to prove, followed by a road map for the rest of the thesis.

The reader already intrinsically motivated may of course skip this chapter to get right to the nitty-gritty mathematics.

## 1. What Is Courcelle's Theorem?

Courcelle's Theorem ([Cou90]) states the following.[1]

**Courcelle's Theorem**

> Every definable subset of $FG(A)_k$ is an effectively recognizable set of graphs.

Even knowing that $FG(A)_k$ is the set of isomorphism classes of directed hypergraphs of type $k$ over a finite set of edge labels $A$, this statement might not be immediately helpful. Let us dissect it.

The finite set $A$ is just a collection of edge labels. Each edge label $a \in A$ is also assigned a *type*, which is just a natural number[2]. A label $a$ can only go on edges that connect the same amount of vertices[3] as the type of $a$. For example, a label of type 2 can only go on an edge with exactly two end points.

[1] The exact statement in [Cou90] is "Every definable subset of $FG(A)_k$ is an effectively *given* recognizable set of graphs", which appears to be a typo since Courcelle's proof explicitly states "effectively recognisable".

[2] In this and forever hence, the natural numbers shall contain the number 0.

We shall see in this work that the set $A$ does not actually have to be finite. Indeed, it can safely be omitted without affecting the truth of the statement.[4]

The natural number $k$ denotes the number of "terminal" vertices in our graph and is called the graph's *type*.[5] These vertices are not to be confused with, for example, the sources and sinks of a network flow. Rather, they are special only in the sense that we designate them as points where we may "glue" two graphs together (in a sense to be made precise later). A graph of type 0 is a "normal" graph as the reader surely knows it (where all vertices are equally important and worthy of love).

Graph homomorphisms between such hypergraphs are defined exactly as one would expect and respect terminal vertices; in particular, two graphs of different type can never be isomorphic.

We hence fix a natural number $k$ and only look at hypergraphs of this type, up to isomorphism (in layman's terms, up to renaming of vertices and edges). For the reader uncomfortable with the notion of "graph type", they may simply prescribe $k = 0$ to stay in familiar waters[6]. For the reader uncomfortable with hypergraphs, not much is lost by restricting to the case where all edges are of type 2 (that is, have exactly two end points).[7]

A "definable" subset of this set of (isomorphism classes of) graphs is a set that is defined by a formula in *counting monadic second-order logic*, which is a certain restricted version of the second-order logic the reader may or may not know. For now, it suffices to know that counting monadic second-order logic defines rules by which one may build logical formulas, and a set $X$ of graphs is called definable if and only if it there is a counting monadic second-order logic formula $\varphi$ such that

$$X = \{\, G \in FG(A)_k : G \text{ satisfies } \varphi \,\}.$$

---

[3] Remember that we are talking about hypergraphs, so an edge may have zero, one, two, or finitely many end points instead of just one or two.

[4] However, $A$ may as well be chosen to be finite, since a larger set does not affect any of the practical implications of the theorem.

[5] The reader will at some point notice that several concepts in the general vicinity of Courcelle's Theorem are called the "type" of something. Unless otherwise stated, it should not be assumed that these share anything besides their name – an edge's "type" is its number of end points, while a graph's "type" is its type in a certain typed algebra as defined in chapter 6.

[6] We are committed to this nautical metaphor now and shall ride it until the bitter end. Probably when we fall off the edge of the world in our little logic boat.

[7] For the reader uncomfortable with mathematics in general, we recommend [Pip20].

In other words, a set of graphs is definable if it is the set of all graphs having a certain property (like "all graphs which admit a vertex cover of size 8"), and that property is definable in counting monadic second-order logic (which the aforementioned vertex cover property happens to be).

Courcelle now asserts that any such set of graphs is *effectively recognisable.*

"Recognisable" means that the set is precisely the set of accepted inputs of a certain deterministic finite state automaton.[8]

"Effectively" recognisable, then, means that this automaton not only exists, but can actually be computed, without invoking trickery like the axiom of choice or applying uncomputable functions.[9]

Since deterministic finite state automata always halt, the latter property also implies that the set $X$ of graphs is *decidable* – that there exists an algorithm that, given any graph, tells us in finite time whether or not that graph belongs to the set. In fact, it even tells us something about the runtime of this algorithm, as we will see later.

Picking up the pieces we have shaved off the theorem so far, we get the following slightly less confusing corollary.

[8] Technically, *recognisable* is defined as being the preimage of a homomorphism from a certain locally finite typed algebra. We shall see later that these concepts are interchangeable for certain classes of graphs, for example the class of graphs with tree-width at most $k$ for a fixed $k \in \mathbb{N}$.

[9] We do not discuss efficacy or computability in this thesis. The constructive nature of the provided proofs might already convince the reader that all steps must be computable. Otherwise, a look at [Cou90] and the other papers in the same series should prove enlightening once one has understood the constructions introduced in the present work.

**Courcelle's Theorem**

Let $k \in \mathbb{N}$. Let $\varphi$ be a property of graphs that is expressible in counting monadic second-order logic.

Then there is an algorithm which, given any graph of tree-width at most $k$, decides in finite time whether that graph satisfies $\varphi$.

Actually, we shall show that the algorithm mentioned above runs in *linear* time.

We also prove several extensions of this result, for example to other classes of graphs, and present a generalised framework in which to generate similar results.

# 2. How to Read this Thesis

The chapters of this thesis are mostly designed for linear consumption. Some can be read independently of each other, but since this thesis is geared towards graph theorists, most examples will involve graphs in one form or another.

Some generalisations have been omitted from the main corpus of this work in order to make the first-time reader's life easier. These concepts have instead been put into appendices.

The reader should start with chapter 3, which introduces some simple preliminaries and explains our notation, and then make use of the following road map.

- Chapter 4 introduces the necessary graph-theoretic concepts. They are not needed for all the subsequent chapters, but they are needed to understand all examples.
- Chapter 5 introduces the concepts of formal logic in a way accessible to the average mathematician with no such theoretical background. One can read most of this chapter without any advanced knowledge of graph theory, if desired.
- Chapter 6 introduces typed algebras. We show that graphs form such

a typed algebra, whence we can translate many theoretical results about algebras into statements about finite graphs.

- Chapter 7 represents the payoff of our previous work and uses the results to show the theoretical formulation of Courcelle's Theorem. It says nothing about algorithms or complexity yet. We show Courcelle's Theorem for a subset of counting monadic second-order logic that suffices for many practical applications to cut down on notation. The proof of the full statement is deferred to appendix A.
- Chapter 8 introduces tree automata, a variation on deterministic finite state automata. It assumes that the reader has read at least the section on trees from chapter 4.
- Chapter 9 extends the proof from chapter 7 to show that on certain classes of graphs, for example those of tree-width bounded by a constant, Courcelle's Theorem delivers a linear-time algorithm to check whether any given graph satisfies a certain logical property. We focus in this chapter on graphs which may have parallel edges, but no loops, the latter being delayed until appendix B. We discuss some interesting special cases of the results from the previous chapters.
- Appendix A expands upon the concepts from chapter 5 and introduces the "counting" part of counting monadic second-order logic. We then re-prove all previous statements in this slightly more general setting.
- Appendix B shows how we can extend everything we have shown in chapter 9 to work on graphs with loops.

The acyclic graph on the next page represents the dependencies between the chapters. Dashed edges represent a "soft" dependency, where omission of the previous chapter still allows the reader to understand most of the new chapter, minus some examples or possibly a clearly-labelled section that deals with connections between the two concepts.

Of course, reading all chapters in ascending order (putting appendices last) is a valid option and yields in particular a topological sorting for said graph.

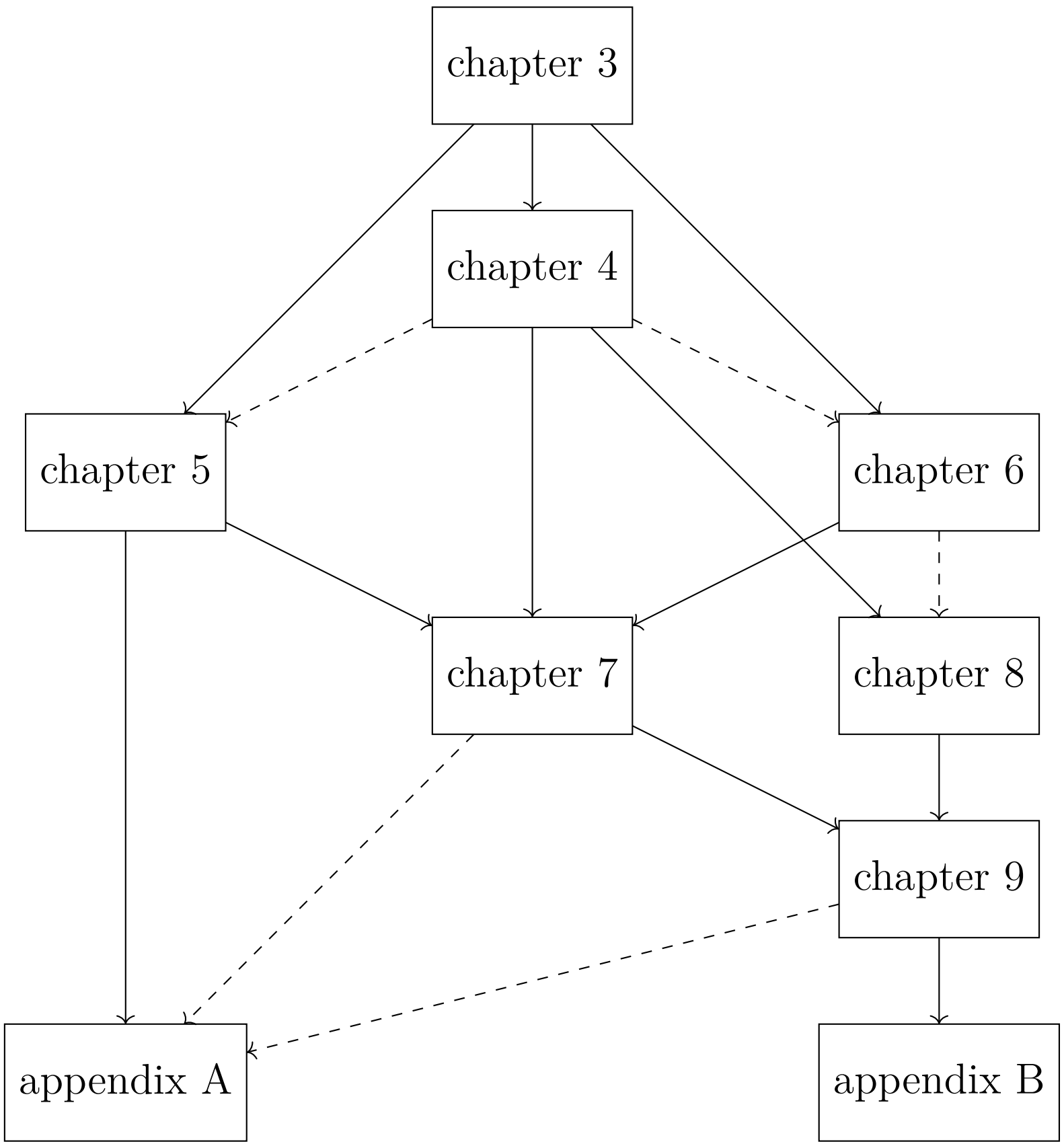

# Chapter 3 Preliminaries

The journey to speaking any language begins with understanding its writing system.[1] We therefore introduce in this chapter the constituent particles from which all words and sentences of logic will be built. We further collect some miscellaneous notation that may or may not be standard. The reader may safely skip section 3.2 and only use it as a glossary.

## 1. Alphabets

The most important thing to remember is that the languages we define are, on their own, devoid of meaning. They are purely syntactical constructs, with semantics only added later on. For example, the average reader would presume the sentence

$$1 = 2$$

to be "false" in some Boolean interpretation. This, however, is presumptuous in the sense that it assumes that the language and meanings used are those known to the reader. A mathematician of some unknown extraterrestrial race may well have defined the *symbol* "1" to mean "the number three", the symbol "2" to mean "the natural numbers", and the symbol "=" to mean "is an element of". Of course, no earth mathematician in their right mind would define such a language, rather sticking to symbols like

$$3 \in \mathbb{N},$$

[1] This is, of course, a blatant lie (think of how one learns one's mother tongue), but in the case of formal logic, a bottom-up approach is generally favourable.

but this exaggerated example shall serve to remind the reader to treat as separate *syntax* and *semantics* of formal logic.

**Definition 3.1.1**

An *alphabet* is a set of symbols.

This is a somewhat vacuous definition (why not just say "set"?), but we shall use the word *alphabet* to emphasise that the elements of this set are syntactical building blocks with no meaning or interactions attached to them.

**Definition 3.1.2**

Let $\Sigma$ be an alphabet. A *word* over $\Sigma$ is a finite sequence of symbols $\alpha_1, \alpha_2, \dots, \alpha_n \in \Sigma$, written like

$$\alpha_1 \alpha_2 \dots \alpha_n.$$

The number $n$ is called the *length* of the word, denoted $|\alpha_1 \alpha_2 \dots \alpha_n|$.

We denote by $\varepsilon_\Sigma$ the word of length zero, called the *empty word*.

The set of all words over $\Sigma$ is denoted by the so-called *Kleene star* $\Sigma^*$.

We usually drop the index from $\varepsilon_\Sigma$, writing simply $\varepsilon$.

From time to time, we might want to map words to other words. This is achieved in a canonical way.

**Definition 3.1.3**

Let $\Sigma, \Theta$ be alphabets and let $f\colon \Sigma \to \Theta$ be a function. The function $f^*\colon \Sigma^* \to \Theta^*$ which maps $\alpha_1 \dots \alpha_n$ to $f(\alpha_1) \dots f(\alpha_n)$ is called the *Kleene extension* of $f$.

One less common operation that we shall need is the permutation of words, that is, the notion that the words "$\alpha\alpha\beta\beta$" and "$\alpha\beta\alpha\beta$" are somehow related in a way that "$\alpha\alpha\beta\beta$" and "$\alpha\alpha\alpha\alpha$" are not.

**Definition 3.1.4**

Let $n \in \mathbb{N}$. We denote by

$$\mathfrak{B}_n := \{\, f\colon \{\, 1, \dots, n \,\} \to \{\, 1, \dots, n \,\} : f \text{ is bijective} \,\}$$

the set of permutations on $\{\, 1, \dots, n \,\}$.

**Definition 3.1.5**

Let $\Sigma$ be an alphabet, $\omega = \alpha_1 \dots \alpha_n \in \Sigma^*$, and $\pi \in \mathfrak{B}_n$. We denote by

$$\pi(\omega) := \alpha_{\pi(1)} \dots \alpha_{\pi(n)}$$

the word $\omega$ reordered according to $\pi$. If $\omega, \eta$ are words in $\Sigma^*$ of equal length $n$ such that there exists a permutation $\pi \in \mathfrak{B}_n$ with $\pi(\omega) = \eta$, we call $\eta$ a *reordering* of $\omega$ and write $\omega \approx \eta$.

# 2. Miscellaneous Notation

As far as possible, we try in this work to stick to preexisting notation. At times, however, this would result in clashes. Rather than have two objects with the same name, we opt to change the notation in this case.

For the reader's convenience, the type of letter used correlates to the type of object it references.

- Standard antiqua letters refer to graph theory – $G$ is a graph, $v \in V$ is a vertex in a set of vertices, and so on.
- Calligraphic letters refer to typed algebras and their signatures – $\mathscr{A}$ is an algebra, $\mathscr{f} \in \mathscr{F}$ is a function symbol in a set of function symbols, and so on.
- Greek letters refer to formal logic – $\varphi$ is a well-formed formula, $\nu_0$ a formal variable, $\Gamma$ is a language, and so on. In this context, we consider general strings and alphabets a part of the broader topic of formal logic.
- Fraktur letters are used for constants rather than variables – $\mathfrak{G}$ is the

> algebra of finite graphs, for example, and we have already seen $\mathfrak{B}_n$. This means that whenever the reader re-encounters a particular fraktur letter, it will refer to exactly the same object as it did when they first met it.

For tree automata, which we need only in few places, we have reserved letters from the ancient Phoenician alphabet:

| | | | | |
|---|---|---|---|---|
| 𐤀 (alef) | 𐤁 (bet) | 𐤂 (gaml) | 𐤃 (delt) | 𐤄 (he) |
| 𐤅 (wau) | 𐤆 (zai) | 𐤇 (het) | 𐤈 (tet) | 𐤉 (yod) |
| 𐤊 (kaf) | 𐤋 (lamd) | 𐤌 (mem) | 𐤍 (nun) | 𐤎 (semk) |
| 𐤏 (ain) | 𐤐 (pe) | 𐤑 (sade) | 𐤒 (qof) | 𐤓 (rosh) |
| 𐤔 (shin) | 𐤕 (tau) | | | |

𐤄 (he) represents a deterministic bottom-up finite tree automaton.

In case the reader is, against expectations, not fluent in ancient Phoenician, a helpful margin note such as the one on this page shall remind them of the symbol's name.

## 2.1. Sets and Tuples

We need in several places the power set of a given set.

**Definition 3.2.1**

> Let $X$ be a set. We denote by $\mathbb{2}^X$ the power set of $X$, that is, the set of all subsets of $X$.

Since our definitions are built iteratively on top of each other, there will be many nested tuples. For instance, an undirected graph will be a triple $(V, E, \llparenthesis_\rrparenthesis)$, while a directed graph will be a pair $(G, \star)$, where $G$ is an undirected graph. This means that writing out a directed graph looks like

$$((V, E, \llparenthesis_\rrparenthesis), \star),$$

which quickly gets ugly as the nested parentheses accumulate.

We allow ourselves to drop the additional parentheses where no confusion can occur.

**Notation 3.2.2**

Given a nested tuple

$$((x_1, \dots, x_n), y),$$

we write

$$(x_1, \dots, x_n, y)$$

to denote that same nested tuple.

## 2.2. Functions with More or Less than One Parameter

We shall often deal with families of functions of the form $\{f_n\}_{n \in \mathbb{N}}$, where $f_n$ is a function from $X^n$ to some set $X'$. We understand in this case that $X^2, X^3, \dots$ contain tuples and hence $f_2, f_3, \dots$ take multiple input parameters. The set $X^1$ should naturally contain 1-tuples, but is usually identified with the set $X$.

A function $X^1 \to X'$ is called "unary", a function $X^2 \to X'$ is called "binary", and so on, after the number of input parameters it takes, even though in reality, a function $X^2 \to X'$ takes just *one* 2-tuple from $X \times X$.

No matter what $X$ actually is, we set the convention that $X^0 = \{ () \}$, a one-element set. A function $X^0 \to X'$ is called "nullary" because it takes no real parameters. In reality, since it is a function from a one-element set (and not from the empty set), it does of course take one parameter, it just happens to always be $()$. Thus, a nullary function is always constant. When defining such a function mapping to some element $x' \in X'$, we write

$$f \colon X^0 \to X', () \mapsto x'.$$

## 2.3. Indices

In most cases, we shall use both upper and lower indices for our objects. In some cases, in order to avoid triple and quadruple indices, we also use *left* indices, like in the twine of two graphs $G$ and $G'$:

$$G\ {}^{n}_{m}{\otimes}^{k}_{K}\, G'.$$

To save the reader from having to remember what every index means, we adopt the following conventions: a *left* index (upper or lower) always indicates in some way the input type of an object – for example, a function which is denoted as ${}_{\mathbb{N}}f$ can be expected to take a natural number as input, while ${}_{\Sigma}f$ will take words from some alphabet $\Sigma$. *Right* indices relate to things other than the input type.

In most cases, the more "important" index will be the lower (left or right) one, while "auxiliary" properties are upper indices. For example, what others would write as $f_{x,y}$, we write as $f_x^y$, whereas $f_{y,x}$ becomes $f_y^x$.

If we need to indicate the *output* type of a function, we shall always choose the upper right index to do so.

## 2.4. Formal Logic

Since we shall be dealing with logical formulas quite a bit, and also with formulas which in turn contain other formulas, we adopt the convention that the symbols of a formal language shall be printed in red. For example, the formula

$$\lambda_1(\nu_1, \nu_2) \Leftrightarrow \nu_1 = \nu_2$$

means not "in the formal language, $\lambda_1(\nu_1, \nu_2)$ and $\nu_1 = \nu_2$ are the same formula", but rather "$\lambda_1(\nu_1, \nu_2)$ is true precisely when the symbols $\nu_1$ and $\nu_2$ are the same", so perhaps it is the definition of an equality-checking symbol $\lambda_1$.

## 2.5. Algorithmic Complexity

For the algorithmic specifics in chapter 9, we assume that the reader is familiar with the very basics of complexity theory.

**Notation 3.2.3**

We write **P** for the class of problems solvable in polynomial time and **NP** for the class of problems verifiable in polynomial time. We write $\mathcal{O}(f)$ for the class of functions asymptotically bounded by $f$.

# Chapter 4
# Graphs

We begin our journey where it will end: with finite graphs. Before the impatient graph theorist jumps to the next chapter, we briefly touch on some peculiarities of our notation.

As our ultimate goal is to relate graph properties and formal logic, we must take special care in our definitions. When talking about formal logic, it can be rather inconvenient when the structures one examines are too big to fit into a set. If the reader recalls the standard definition of a graph, they will notice that usually the vertex set $V$ is simply “a set”, meaning the vertices $v \in V$ can be whatever the author pleases. This means that the vertices themselves could again be sets – in particular, for every set $X$, there is a (finite!) graph $G$ which contains $X$ as a vertex. Hence the collection of all graphs includes all possible sets and consequently cannot be itself a set[1], and this remains true even if one requires all graphs to be finite.

The standard construction to pare down the class of graphs to a manageable size is to define what a graph isomorphism is, then consider isomorphism classes of graphs. If the reader prefers this, they may pretend that this is what we have done and pretend that every vertex mentioned in the proofs starting from chapter 5 is actually a vertex of a representative of such an isomorphism class, that graph morphisms are actually classes of morphisms between classes of graphs, and so forth – this requires no more notational fudgery than is usually involved in the isomorphism class construction, as virtually no author states their theorems in the form “let $G$

[1] Suppose there were a set $X$ containing all sets. In Zermelo-Fraenkel set theory, we know (for example by Cantor diagonalisation) that $|\mathbb{2}^X| > |X|$ for any set. But $X$ contains all sets, in particular the elements of $\mathbb{2}^X$, hence $\mathbb{2}^X \subseteq X$ – a contradiction.

$\mathbb{2}^X$ denotes the power set of $X$.

be a representative of an isomorphism class of cubic graphs …", even if they later treat the class of cubic graphs as if it were a set.

We present here an alternative approach that integrates nicely with our usage of formal languages and which allows to work with a set (instead of a proper class) of finite graphs while staying formally sound throughout. The price we pay for this is that most of what the reader would consider "graphs" are no longer graphs in our definition – we do not allow the vertices and edges to be arbitrary objects anymore.

Is this an absurd restriction? It is indeed not, because for any "graph" the reader might invent, there is a graph in our definition that is isomorphic to it. We simply prohibit cases where the vertices of the graph are sets, or classes, or oranges, and only ever allow them to be the symbols $v_1, v_2, \ldots$.

A "graph" as the reader knows it, where the vertices and edges can be arbitrary objects, is instead demoted to the status of "pseudograph".

**A note on hypergraphs**

When we talk about graphs in this thesis, we actually mean hypergraphs, that is, graphs where a single edge can connect more than two vertices. If the reader feels uncomfortable with this, they may simply pretend that edges can only ever connect two vertices – whenever we mention the number of end points of an edge, they may mentally substitute the number 2.
In case further clarification seems in order, a box not entirely unlike this one shall help the reader out.

# 1. What is a Graph?

All of our graphs are hypergraphs where an edge can connect arbitrarily many vertices and can also contain the same vertex more than once. We omit the prefix "hyper" in most cases.

**Definition 4.1.1**

An *undirected pseudograph* is a tuple $(V, E, \langle\!\langle _ \rangle\!\rangle)$ with the following properties.

- $V$ is a nonempty set, called the set of *vertices*.

- $E$ is a set, called the set of *edges*.
- $\langle\!|_|\!\rangle$ is a function $E \to V^*$, called the *end points* of the edges.
- $V$ and $E$ are disjoint.

For an edge $e \in E$, the natural number $|\langle\!| e |\!\rangle|$ is called the *edge type* of $e$, denoted $|e|$.

To illustrate, let us consider the (not very practical) undirected pseudograph $G = (V, E, \langle\!|_|\!\rangle)$ with

- $V := \{ v_1, \mathbb{N}, f,\ 3.7 \}$,
- $E := \{ 2, V, G, \{ 2 \} \}$,
- $\langle\!| e |\!\rangle := \begin{cases} v_1\mathbb{N} & \text{if } e = \ 2 \\ v_1v_1 & \text{if } e = V \\ \mathbb{N}v_1f & \text{if } e = \{ 2 \} \\ fff3.7 & \text{if } e = G. \end{cases}$

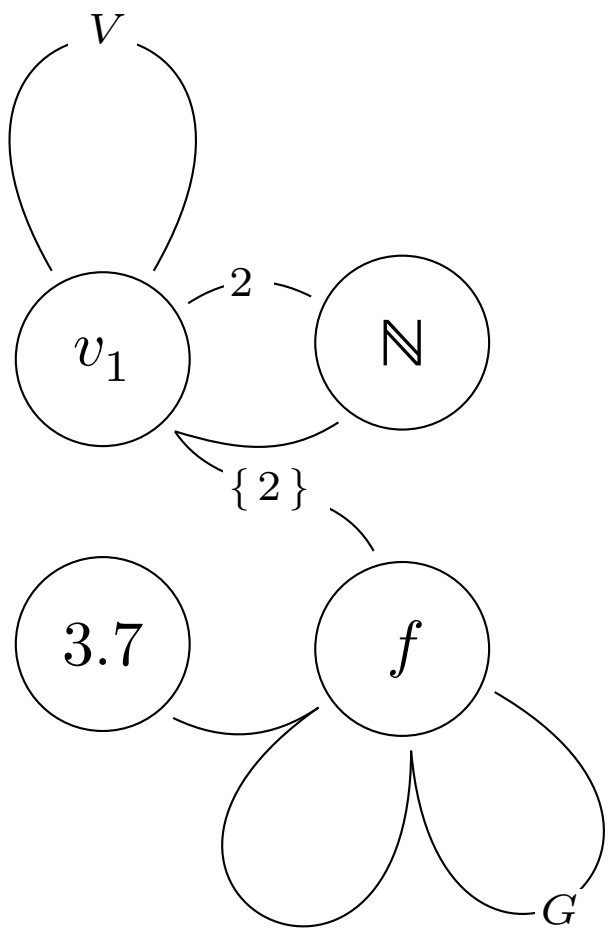

The edge 2 is an ordinary edge, going from $\mathbb{N}$ to $v_1$ (or vice versa, since the graph is not directed). The edge $V$ is not too uncommon either, it is simply a loop. The edge $\{ 2 \}$ connects three vertices, while the edge $G$ connects to 3.7 once and to $f$ thrice, forming a multiloop.

**What about the empty graph?**

By definition, all of our graphs have at least one vertex. This is a deliberate choice that makes several proofs in later chapters easier to write down. Should it be the reader's deepest wish to allow for empty graphs in their algorithms, they may simply add that as a special case without changing any of the runtimes, since any property is easy to check on the empty graph.

Apart from showcasing the kind of edges that a hypergraph admits, $G$ illustrates our point about pseudographs: allowing arbitrary sets for $V$ and $E$ leads to strange objects which are still graphs, like this one containing itself as an edge. We shall soon disallow these constructions.

Since these graphs are supposed to be undirected, it makes sense to regard graphs where the end points of edges are simply reordered as "the same".[2]

### Definition 4.1.2

Let $G = (V, E, \langle_\rangle)$, $G' = (V', E', \langle\!\langle_\rangle\!\rangle)$ be undirected pseudographs. An *undirected pseudograph morphism* from $G$ to $G'$ is a pair $(g, h)$ with the following properties.

- $g$ is a function $V \to V'$.
- $h$ is a function $E \to E'$.
- for every edge $e \in E$, we have $g^* \langle e \rangle \approx \langle\!\langle he \rangle\!\rangle$.

$\approx$ denotes two strings which are identical up to reordering. (def. 3.1.5, p. 15)

Stated in prose, a pseudograph morphism maps vertices to vertices in any way it likes, then maps edges to edges such that the number of end points of each edge (with multiplicity) is preserved and such that if $v$ was an end point of $e$ with multiplicity $k$, $gv$ is an end point of $he$ with multiplicity *at least* $k$. The "look" of $e$ can still change, as $g$ might identify several of its end points with the same vertex $v'$:

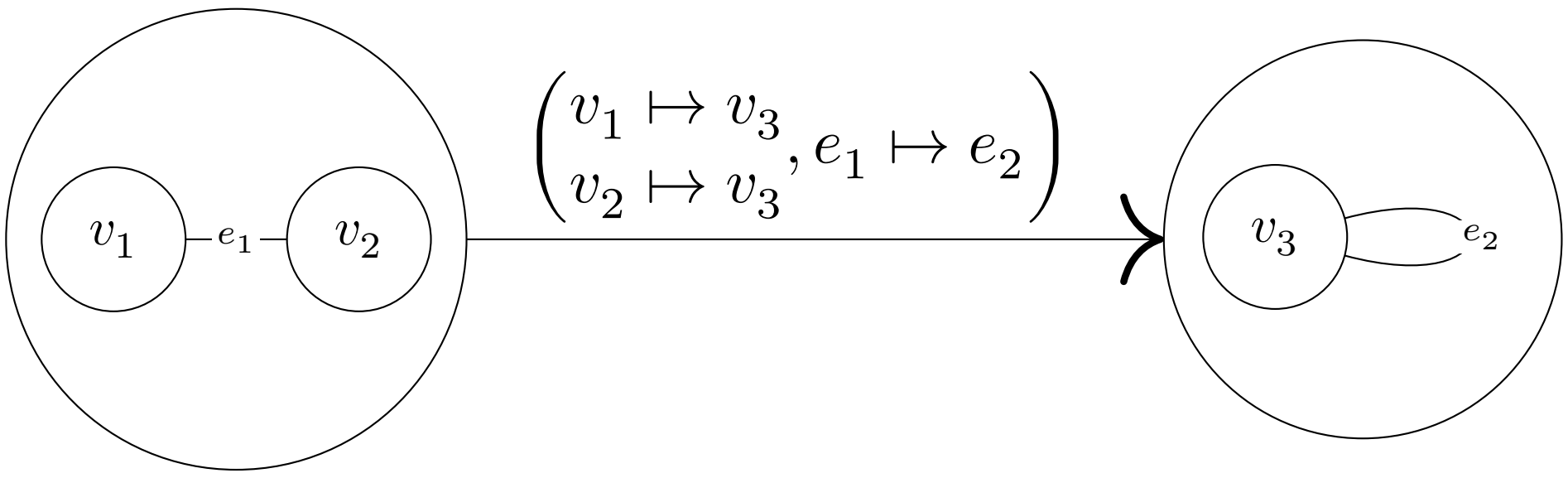

The reader may easily verify that composition of pseudograph morphisms yields a pseudograph morphism and that the identity on a pseudograph

[2] This could alternatively be achieved by prescribing $\langle e \rangle$ to be a multiset rather than a word. This makes the notion of graph homomorphism simpler, at the cost of more cluttered notation and making the definition of directed graphs harder. The reader may simply pretend that we are using their preferred version, since nothing depends on the actual word-ness of $\langle_\rangle$.

(denoted by id) is a pseudograph morphism, thus justifying the name and the following definition.

### Definition 4.1.3

An undirected pseudograph morphism $(g, h)\colon G \to G'$ is called an *undirected pseudograph isomorphism* if there is an undirected pseudograph morphism $(g', h')$ such that $(g, h) \circ (g', h') = \mathrm{id}_{G'}$ and such that $(g', h') \circ (g, h) = \mathrm{id}_G$.

Two pseudographs $G$ and $G'$ are called *isomorphic* if there exists a pseudograph isomorphism $G \to G'$. We write $G \cong G'$.

The following lemma will be well-known to the reader.

### Lemma 4.1.4

Let $f = (g, h)$ be a pseudograph morphism. Then $f$ is a pseudograph isomorphism if and only if $g$ and $h$ are both bijective. In that case, we have $f^{-1} = (g^{-1}, h^{-1})$.

**Proof.** Let $G = (V, E, \Lbag_\Rbag), G' = (V', E', \llparenthesis_\rrparenthesis)$ be pseudographs.

First, let $(g, h)\colon G \to G'$ be a pseudograph isomorphism. Then by definition there is a pseudograph morphism $(g', h')\colon G' \to G$ such that we have both $(g \circ g', h \circ h') = (\mathrm{id}_{V'}, \mathrm{id}_{E'})$ and $(g' \circ g, h' \circ h) = (\mathrm{id}_V, \mathrm{id}_E)$, proving the claim.

Let conversely $(g, h)\colon G \to G'$ be a pseudograph morphism with $g, h$ bijective. It suffices to show that $(g^{-1}, h^{-1})$ is a pseudograph morphism from $G'$ to $G$, that is,

$$\forall e' \in E' \colon (g^{-1})^*(\llparenthesis e' \rrparenthesis) \approx \Lbag h^{-1}(e') \Rbag.$$

$\approx$ denotes two strings which are identical up to reordering. (def. 3.1.5, p. 15)

Let to this end $e' \in E'$. Since $h$ is bijective, we can find $e \in E$ such that $e' = h(e)$ and immediately recognise that

$$\begin{aligned}(g^{-1})^*(\langle\!\langle e' \rangle\!\rangle) &= (g^{-1})^*(\langle\!\langle h(e) \rangle\!\rangle)\\ &\approx (g^{-1})^*(g^*(\langle e \rangle))\\ &= (g^*)^{-1}(g^*(\langle e \rangle))\\ &= \langle e \rangle\\ &= \langle h^{-1}(h(e)) \rangle\\ &= \langle h^{-1}(e') \rangle,\end{aligned}$$

$\approx$ denotes two strings which are identical up to reordering. (def. 3.1.5, p. 15)

finishing the proof.

□

We are now ready to introduce "real" graphs. We start by fixing a set of meaningless symbols[3].

**Definition 4.1.5**

A *graph language* is a pair $(\Sigma_V, \Sigma_E)$ of two disjoint alphabets. The elements of $\Sigma_V$ are called *vertex symbols*, while the elements of $\Sigma_E$ are called *edge symbols*.

The alphabets $\Sigma_V$ and $\Sigma_E$ may be proper classes if one really wants to construct all possible graphs. We shall, however, soon restrict them to be quite small for our purposes.

A graph in a language is now simply a pseudograph that only uses the allowed symbols.

**Definition 4.1.6**

Let $\Gamma = (\Sigma_V, \Sigma_E)$ be a graph language. An *undirected $\Gamma$-hypergraph* is a pseudograph $(V, E, \langle _ \rangle)$ such that $V \subseteq \Sigma_V$ and $E \subseteq \Sigma_E$.

As long as both $\Sigma_V$ and $\Sigma_E$ are sets, this ensures that the collection of all $\Gamma$-hypergraphs also forms a set.

The notions of morphism and isomorphism carry over without any modification. They are reproduced here for completeness only.

[3] The best kind of symbols.

### Definition 4.1.7

Let $\Gamma$ be a graph language, and let $G = (V, E, \langle\!|_|\!\rangle)$, $G' = (V', E', (\!|_|\!))$ be undirected $\Gamma$-hypergraphs. An *undirected $\Gamma$-hypergraph morphism* from $G$ to $G'$ is a pair $(g, h)$ with the following properties.

- $g$ is a function $V \to V'$.
- $h$ is a function $E \to E'$.
- for every edge $e \in E$, we have $g^* \langle\!| e |\!\rangle \approx (\!| he |\!)$.

$\approx$ denotes two strings which are identical up to reordering. (def. 3.1.5, p. 15)

### Definition 4.1.8

Let $\Gamma$ be a graph language, and let $G$, $G'$ be undirected $\Gamma$-hypergraphs. An undirected $\Gamma$-hypergraph morphism $(g, h) \colon G \to G'$ is called an *undirected $\Gamma$-hypergraph isomorphism* if there is an undirected $\Gamma$-hypergraph morphism $(g', h')$ such that $(g, h) \circ (g', h') = \mathrm{id}_{G'}$ and $(g', h') \circ (g, h) = \mathrm{id}_G$.

Two $\Gamma$-hypergraphs $G$ and $G'$ are called *isomorphic* if there exists a $\Gamma$-hypergraph isomorphism $G \to G'$.

In particular, every $\Gamma$-hypergraph morphism is also a pseudograph morphism.

We formalise the idea that a graph from $\Gamma$ is "just as good" as a graph where vertices and edges are arbitrary objects, provided that the language is large enough.

### Lemma 4.1.9

Let $\Gamma = (\Sigma_V, \Sigma_E)$ be a graph language, and let $G = (V, E, \langle\!|_|\!\rangle)$ be a pseudograph with $|V| \leq |\Sigma_V|$ and $|E| \leq |\Sigma_E|$. Then there exists a $\Gamma$-hypergraph which is isomorphic to $G$.

**Proof.** Say the assumptions hold. The fact that $|V| \leq |\Sigma_V|$ implies that there exists an injective function $g \colon V \to \Sigma_V$. Similarly, there is an injective function $h \colon E \to \Sigma_E$. Since $h$ is injective, it admits a retraction $h^{-1} \colon \Sigma_E \to E$ with $h^{-1} \circ h = \mathrm{id}_E$.

We set $(\!|_|\!) := g^* \circ \langle\!|_|\!\rangle \circ h^{-1}|_{h(E)}$, a function $h \colon E \to \Sigma_V^*$.

We claim that the $\Gamma$-hypergraph $G' = (g(V), h(E), \llparenthesis_\rrparenthesis)$ is isomorphic to $G$. The reader may first wish to convince themselves of the intuitive sensibility of our choice by noting that the following diagram commutes.

$$\begin{array}{ccc} E & \xhookrightarrow{h} & \Sigma_E \\ {\scriptstyle \lbag_\rbag}\downarrow & \circlearrowleft & \downarrow{\scriptstyle g^*\lbag h^{-1}\rbag} \\ V^* & \xhookrightarrow[g^*]{} & \Sigma_V^* \end{array}$$

The isomorphism in question is the pair $(g, h)$. We need to show that this is a morphism, and that it admits an inverse.

By definition of $\llparenthesis_\rrparenthesis$, we have for every edge $e \in E$

$$\begin{aligned} \llparenthesis h(e)\rrparenthesis &= g^*(\lbag h^{-1}|_{h(E)}(h(e))\rbag) \\ &= g^*(\lbag e\rbag) \\ &\approx g^*(\lbag e\rbag), \end{aligned}$$

$\approx$ denotes two strings which are identical up to reordering. (def. 3.1.5, p. 15)

whence $g, h$ is a pseudograph morphism.

The functions $g$ and $h$ are injective, and of course they are surjective onto their respective images. Hence by lemma 4.1.4, $(g, h)$ is an isomorphism.

$\square$

We allow the following shorthand.

### Notation 4.1.10

Let $\Gamma$ be a graph language, and let $G = (V, E, \lbag_\rbag)$ be an undirected $\Gamma$-hypergraph. By $v \in G$ we mean $v \in V$, and by $|G|$ we mean $|V|$.

As our main concern is finding algorithms that run in finite time, we are only ever interested in graphs with finitely many vertices and edges. For this reason, we fix now and forever a reasonably small language for our graphs to live in.[4]

[4] The reader may recall that fraktur letters are reserved for global constants, which means the meaning of these symbols shall not change until the end of the thesis.

**Definition 4.1.11**

We set $\mathfrak{F} = (\mathfrak{V}, \mathfrak{E})$ with two countably infinite sets of symbols

$$\mathfrak{V} := \{ v_0, v_1, \dots \}$$

and

$$\mathfrak{E} := \{ e_0, e_1, \dots \}$$

and call this the *language of finite graphs*.

Given any pseudograph $G = (V, E, \langle\!\!\!\langle _ \rangle\!\!\!\rangle)$ with finite vertex and edge sets in any reasonable encoding, one can then compute in linear time an $\mathfrak{F}$-hypergraph $G' = (V', E', (\!|_|\!))$ which is isomorphic to $G$: the "reasonable" encoding should allow one to enumerate the elements of $V \cup E$ in some order. Traverse this enumeration. When a vertex comes up, assign to it the lowest-index unassigned vertex symbol from $\mathfrak{V}$. When an edge comes up, assign to it the lowest-index unassigned edge symbol from $\mathfrak{E}$. When the enumeration is finished, we have produced a bijective mapping $f$ from $V \cup E$ to some finite subset of $\mathfrak{V} \cup \mathfrak{E}$. For an edge $e' \in E'$, computing $(\!|e'|\!)$ is as simple as computing $f^{-1}e'$ and applying $f$ to its end points.

Given an encoding that does not itself already blow up the time to compute elements of $G$, the above conversion procedure uses time linear in $|V| + |E|$. We shall hence assume that any reasonably encoded finite graph is really an $\mathfrak{F}$-graph, since adding a linear preprocessing step does not slow down any of our algorithms.

**Notation 4.1.12**

By *graph*, we always mean an undirected $\mathfrak{F}$-hypergraph with finitely many vertices and edges.

Since $\mathfrak{F}$ admits only countably many hypergraphs, we have thus achieved our goal of paring down the class of graphs to a set.

# 2. Directed Graphs

For the sake of completeness, we also introduce formally the notion of *directed* hypergraph. All of our constructions in later chapters will work equally well on directed as on undirected graphs. The reader who wishes to avoid the technical construction may well skip this section and pretend that there is some intuitive notion of edges going “from” one vertex “to” another.

For our formal approach, we have chosen to define directed graphs in a way that takes an undirected graph and “enhances” it by adding an orientation to each edge. This provides the added bonus that a directed graph becomes undirected by simply forgetting its orientation function.

**Definition 4.2.1**

The reader can easily memorise this notation because ☆ is a *star* and indicates the *start* points of an edge.

Let $G = (V, E, \langle\!\langle _ \rangle\!\rangle)$ be an undirected graph. An *orientation* for $G$ is a map $☆\colon E \to \mathbb{N}$ with

$$\forall e \in E\colon 0 < ☆(e) < |\langle\!\langle e \rangle\!\rangle|.$$

For an edge $e \in E$ with $\langle\!\langle e \rangle\!\rangle = v_1 \dots v_n$, we call $v_1 \dots v_{☆(e)}$ the *start points* of $e$ and $v_{☆(e)+1} \dots v_n$ its *end points.*

A *directed graph* is a pair $(G, ☆)$, where $G$ is a graph and ☆ is an orientation for $G$.

For ordinary (non-hyper) graphs, this coincides with the classical notion of directedness. For hypergraphs, note that an edge of type 3 cannot go “from $v_1$ to $v_2$ to $v_3$”. It only has start and end points, with no particular distinction between starts. Hence it can only go “from $v_1$ and $v_2$ to $v_3$” or “from $v_1$ to $v_2$ and $v_3$” (or other permutations on the vertex order, of course).

We also prohibit edges with no start points or no end points; in particular, a directed hypergraph cannot contain edges of type 1. If the reader disagrees with this convention, it should be straightforward to adapt the definitions, taking some care to fix the definition of a tree in section 4.4.

**Definition 4.2.2**

Let $G = ((V, E, \langle\!| _ |\!\rangle), \star)$ and $G' = ((V', E', \llparenthesis _ \rrparenthesis), \bigstar)$ be directed graphs. A *directed graph morphism* from $G$ to $G'$ is a graph morphism $(g, h)\colon (V, E, \langle\!| _ |\!\rangle) \to (V', E', \llparenthesis _ \rrparenthesis)$ such that for all $e \in E$, we have $\bigstar(he) = \star(e)$ and such that for every $e \in E$ with $\langle\!| e |\!\rangle = v_1 \dots v_n$ and $\llparenthesis he \rrparenthesis = w_1 \dots w_n$, we have $w_1 \dots w_{\star(e)} \approx g^*(v_1 \dots v_{\star(e)})$.

$\approx$ denotes two strings which are identical up to reordering. (def. 3.1.5, p. 15)

Note that this immediately implies that

$$w_{\star(e)+1} \dots w_n \approx g^*(v_{\star(e)+1} \dots v_n).$$

In other words, a directed graph morphism is simply a graph morphism which leaves start points as start points and end points as end points (both up to reordering).

For an ordinary directed graph as we know it (where an edge goes from one vertex to another), definition 4.2.2 coincides with the classical definition (since there is exactly one start and one end for each edge).

**Definition 4.2.3**

A directed graph morphism whose underlying hypergraph morphism is a hypergraph isomorphism is called a *directed graph isomorphism*.

The reader is reminded that, where no confusion is likely to occur, we abuse notation to let $(V, E, \langle\!| _ |\!\rangle, \star)$ mean the same thing as $((V, E, \langle\!| _ |\!\rangle), \star)$.

## 3. More Definitions

We now introduce notation for several useful concepts.

**Definition 4.3.1**

Let $G = (V, E, \langle\!| _ |\!\rangle)$ be a graph.

For an edge $e \in E$ and a vertex $v \in V$, the number of times that the symbol $v$ occurs in the string $\langle\!| e |\!\rangle$ is called the *incidence of $v$ in $e$*, denoted by $\operatorname{inc}_e v$.

For a vertex $v \in V$, the number $\deg v := \sum_{e \in E} \operatorname{inc}_e v$ is called the *degree* of $v$.

Two vertices $v, w \in V, v \neq w$ are called *adjacent* if there is an edge $e$ such that both $v$ and $w$ occur in $\Lbag e \Rbag$.

For a vertex $v \in V$, the set of all vertices adjacent to $v$ is called its *neighbourhood*, denoted $\mathrm{N}v$.

### Definition 4.3.2

Let $G = (V, E, \Lbag _ \Rbag, \star)$ be a directed graph.

For a vertex $v \in V$ and an edge $e \in E$, the number of times that the symbol $v$ occurs in the start points of $e$ is called the *out-incidence of $v$ in $e$*, denoted $\mathrm{inc}_e^{\mathrm{out}} v$. The number of times that the symbol $v$ occurs in the end points of $e$ is called the *in-incidence of $v$ in $e$*, denoted $\mathrm{inc}_e^{\mathrm{in}} v$.

For a vertex $v \in V$, the number $\deg^{\mathrm{in}} v := \sum_{e \in E} \mathrm{inc}_e^{\mathrm{in}} v$ is called the *in-degree* of $v$. The number $\deg^{\mathrm{out}} v := \sum_{e \in E} \mathrm{inc}_e^{\mathrm{out}} v$ is called the *out-degree* of $v$.

Let $v, w \in V, v \neq w$. If there is an edge $e \in E$ with $\mathrm{inc}_e^{\mathrm{out}} v > 0$ and $\mathrm{inc}_e^{\mathrm{in}} w > 0$, we say that $v$ is a *predecessor* of $w$ and that $w$ is a *successor* of $v$.

For $v \in V$, the set of all its successors is called its *out-neighbourhood*, denoted $\mathrm{N}^{\mathrm{out}} v$. The set of all its predecessors is called its *in-neighbourhood*, denoted $\mathrm{N}^{\mathrm{in}} v$.

We can now look at paths through a hypergraph.

### Definition 4.3.3

Let $G = (V, E, \Lbag _ \Rbag)$ be a graph. An *undirected path in $G$ from $v \in V$ to $w \in V$* is a word $P \in (V \cup E)^*$ fulfilling the following conditions.

- $P$ alternates between symbols from $V$ and symbols from $E$.
- The first symbol of $P$ is $v$, its last symbol is $w$. In particular, they both lie in $V$.
- Let $e$ be an edge symbol in $P$ such that $P = \dots xey \dots$. Then $x$ and $y$ must be incident to $e$, that is, $\mathrm{inc}_e x > 0$ and $\mathrm{inc}_e y > 0$.

The *length* of $P$, denoted $|P|$, is the number of edge symbols in $P$, or equivalently, $|P| := \frac{|P|-1}{2}$, where the $|_|$ used on the right hand is the ordinary string length.

If the first and last symbol of a path $P$ of length at least 1 are identical, $P$ is called a *cycle*.

A path $P$ is called *simple* if no vertex symbol and no edge symbol occurs more than once, except that the first and last symbol may be the same.

The graph $G$ is called *acyclic* if it admits no simple cycles.

It is called *connected* if for any $v, w \in V$, there is a path in $G$ from $v$ to $w$.

### Definition 4.3.4

Let $G = (V, E, \langle\!\!|_|\!\!\rangle, \star)$ be a directed graph. A *directed path in $G$ from $v \in V$ to $w \in V$* is a word $P \in (V \cup E)^*$ fulfilling the following conditions.

- $P$ alternates between symbols from $V$ and symbols from $E$.
- The first symbol of $P$ is $v$, its last symbol is $w$. In particular, they both lie in $V$.
- Let $e$ be an edge symbol in $P$ such that $P = \ldots xey \ldots$. Then we have $\mathrm{inc}_e^{\mathrm{out}} x > 0$ and $\mathrm{inc}_e^{\mathrm{in}} y > 0$.

The *length* of $P$, denoted $|P|$, is the number of edge symbols in $P$, or equivalently, $|P| := \frac{|P|-1}{2}$, where the $|_|$ used on the right hand is the ordinary string length.

If the first and last symbol are identical and $|P| \geq 1$, then $P$ is called a *cycle*.

It is called *simple* if no vertex symbol and no edge symbol occurs more than once, except that the first and last symbol may be the same.

The graph $G$ is called *acyclic* if it admits no simple cycles.

It is called *strongly acyclic* if $(V, E, \langle\!\!|_|\!\!\rangle)$ is acyclic.

It is called *connected* (or *weakly connected*) if for any $v, w \in V$, there is an undirected path from $v$ to $w$ in $(V, E, \langle\!\!|_|\!\!\rangle)$.

It is called *strongly connected* if for any $v, w \in V$, there is a directed path in $G$ from $v$ to $w$.

**Definition 4.3.5**

Let $G = (V, E, \langle\!\langle_\rangle\!\rangle)$ be a graph. A graph $G' = (V', E', (\!|_|\!))$ is called a *subgraph* of $G$ if $V' \subseteq V$, $E' \subseteq E$, and $(\!|_|\!) = \langle\!\langle_\rangle\!\rangle|_{E'}$.

It is called a *full* or *induced subgraph* of $G$ if it is maximal among subgraphs with vertex set $V'$, that is, if it contains all edges which are not incident to a vertex in $V \setminus V'$. In this case, we say that it is the (full) subgraph *induced* by $V'$ and write $G[V'] := G'$.

For any $V' \subseteq V$, we write $G - V' := G[V \setminus V']$.

For any $E' \subseteq E$, we write $G - E' := (V, E \setminus E', \langle\!\langle_\rangle\!\rangle|_{E \setminus E'})$.

**Definition 4.3.6**

Let $G$ be a directed or undirected graph, $v$ a vertex of $G$. The *connected component* of $v$ is the largest connected full subgraph of $G$ containing $v$, that is, the full subgraph induced by all vertices (weakly) connected to $v$.

**Definition 4.3.7**

Le $k \in \mathbb{N}$. We say that a directed or undirected graph $G$ is *$k$-uniform* if every edge of $G$ has type $k$.

In particular, a non-hyper graph is the same thing as a 2-uniform graph.

# 4. Trees

For the study of trees, we restrict ourselves to 2-uniform graphs, that is, graphs where every edge has exactly two end points (both of which may be the same vertex).

**Definition 4.4.1**

A 2-uniform graph $G$ is called a *forest* if it is acyclic.

A connected forest is called a *tree*.

This is a tree:

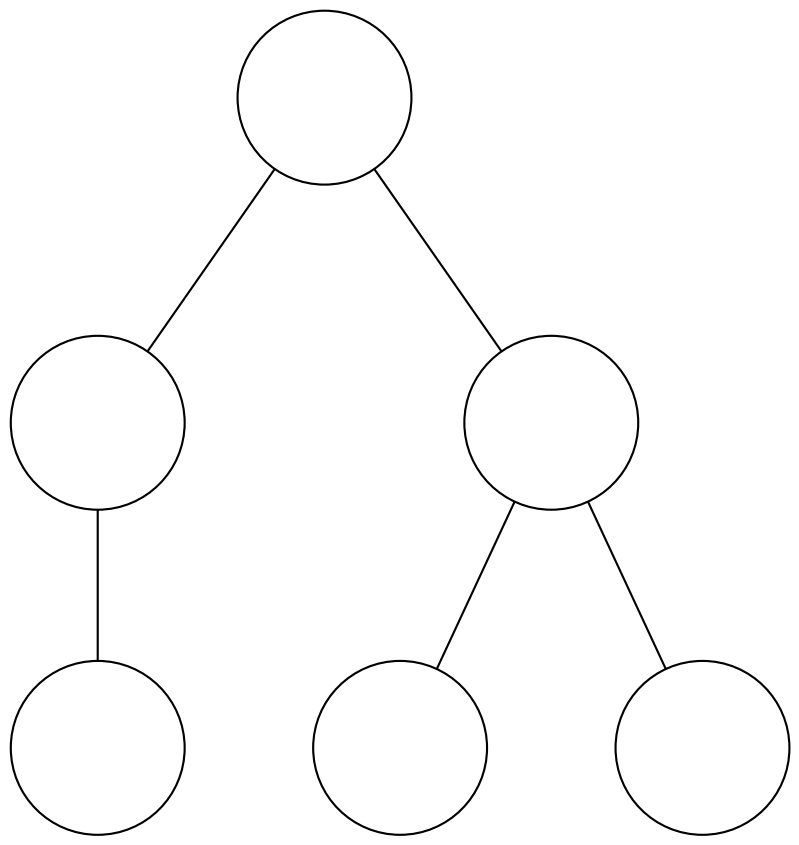

**Definition 4.4.2**

A 2-uniform directed graph $G$ is called a *directed tree* (or *rooted tree*) if it is strongly acyclic and there exists a $v \in V$ such that for every $w \in V$ there is a directed path from $v$ to $w$.

We call $v$ the *root* of $G$, denoted $\sqrt{G} := v$.

A vertex without successors in a rooted tree is called a *leaf*.

A directed graph $G$ is called a *directed forest* if each of its connected components is a rooted tree.

This is a directed tree with root $v$:

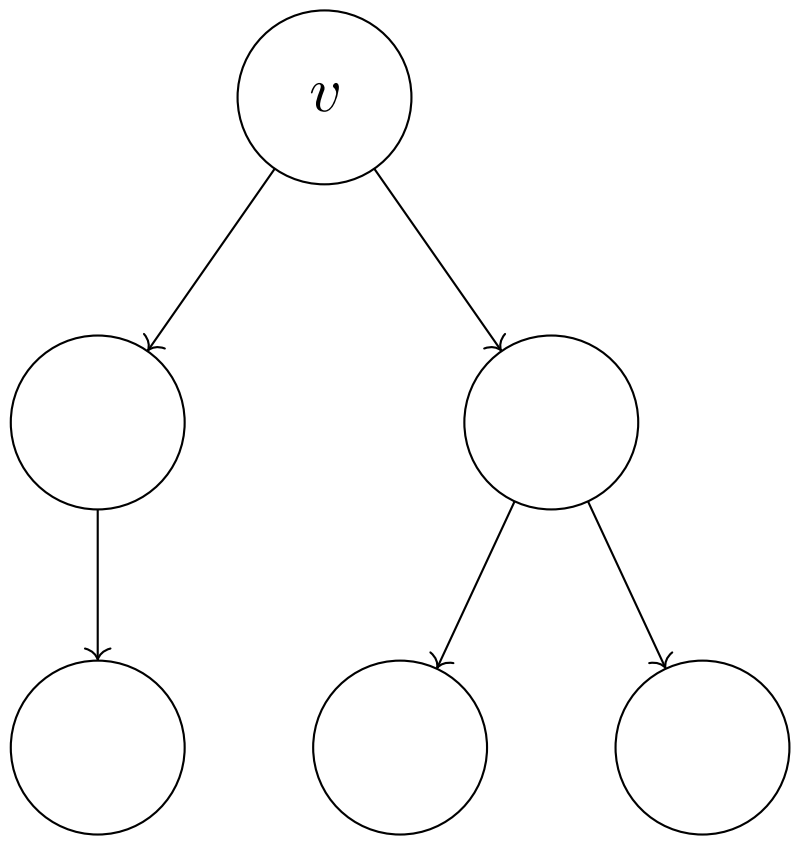

In particular, the underlying undirected graph of a rooted tree is always a tree, whereas a directed graph whose underlying undirected graph is a tree need not be a rooted tree itself.

**Definition 4.4.3**

Let $T$ be a rooted tree, $v \in T$. The largest subtree of $T$ with root $v$ is called the *induced subtree* of $v$, denoted $T[v]$.

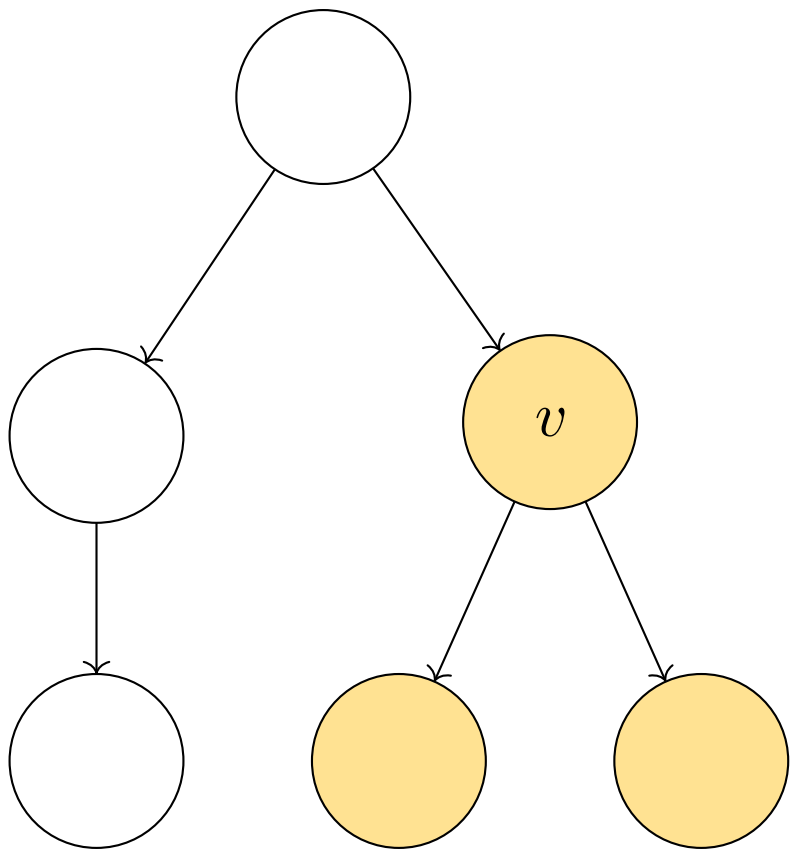

Note that, unlike most of our definitions, this notion makes sense only for directed graphs because for an undirected tree, any vertex can serve as a "root" for the entire tree.

The following notion is another we only need for directed trees.

**Definition 4.4.4**

The *height* of a rooted tree $T$ is one plus the length of the longest directed path in $T$, that is, the number of vertices in that path.

In addition to these classical mathematical notions of trees, we need a way to impose an ordering on the successors of a vertex.

**Definition 4.4.5**

A *traversal tree* is a pair $(T, \preccurlyeq)$ where $T = (V, E, \langle\!|_|\!\rangle, \star)$ is a rooted tree and $\preccurlyeq = \{\preccurlyeq_v\}_{v \in V}$ is a set of relations such that for all $v \in V$, $\preccurlyeq_v$ is a total ordering of $\mathrm{N}^{\mathrm{out}}v$.

One can visualise a traversal tree as follows.

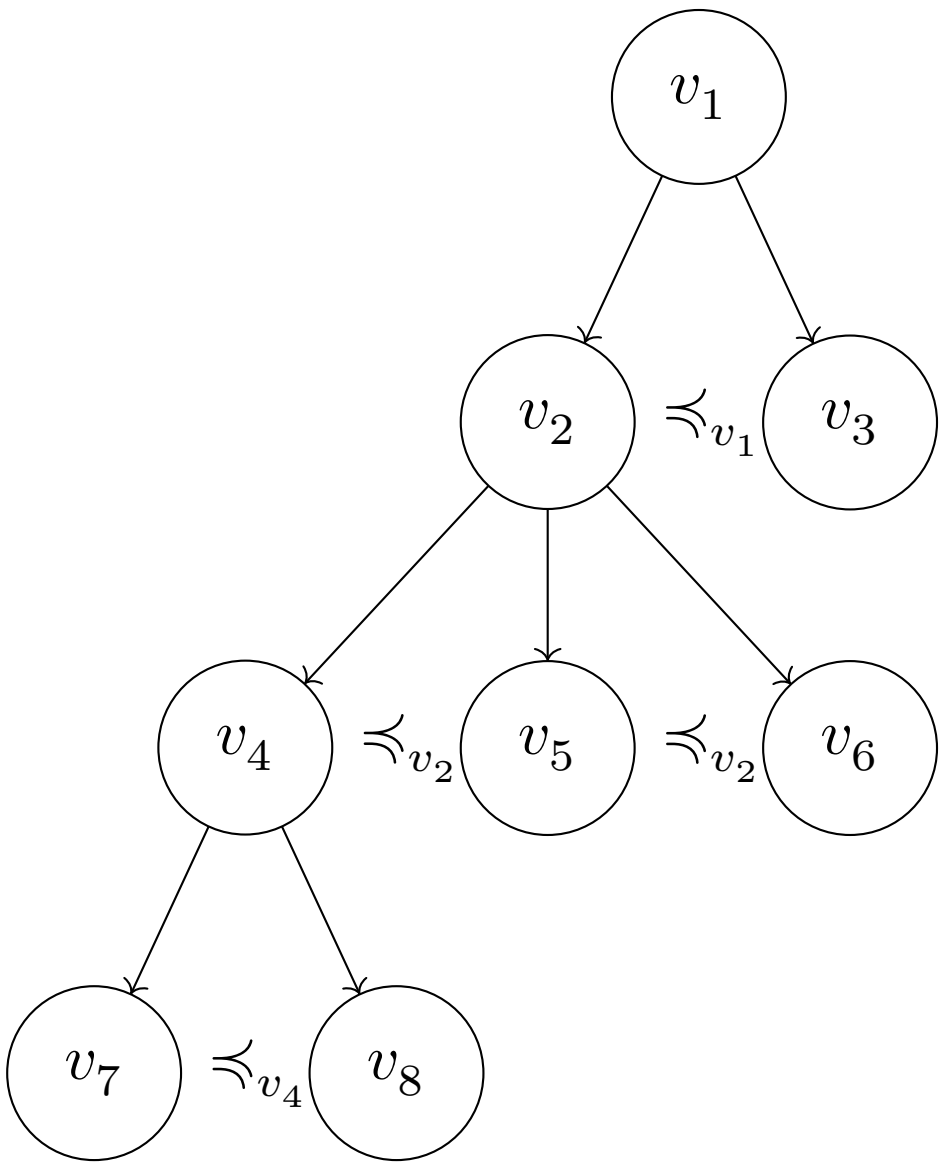

The reason this is called a "traversal tree" is that, given a traversal strategy such as *inorder*, *preorder* or *postorder*[5], the orderings given by $\preccurlyeq$ suffice to uniquely determine the traversal order of $T$.

In the future, we adopt the (hopefully intuitive) convention that in a traversal tree ...

- ... edges are always directed from top to bottom.
- ... successors are always ordered from left to right according to $\preccurlyeq$.

Our example above hence becomes simply the following.

[5] For more on this, check [Knu68]. We shall not actually need these definitions though.

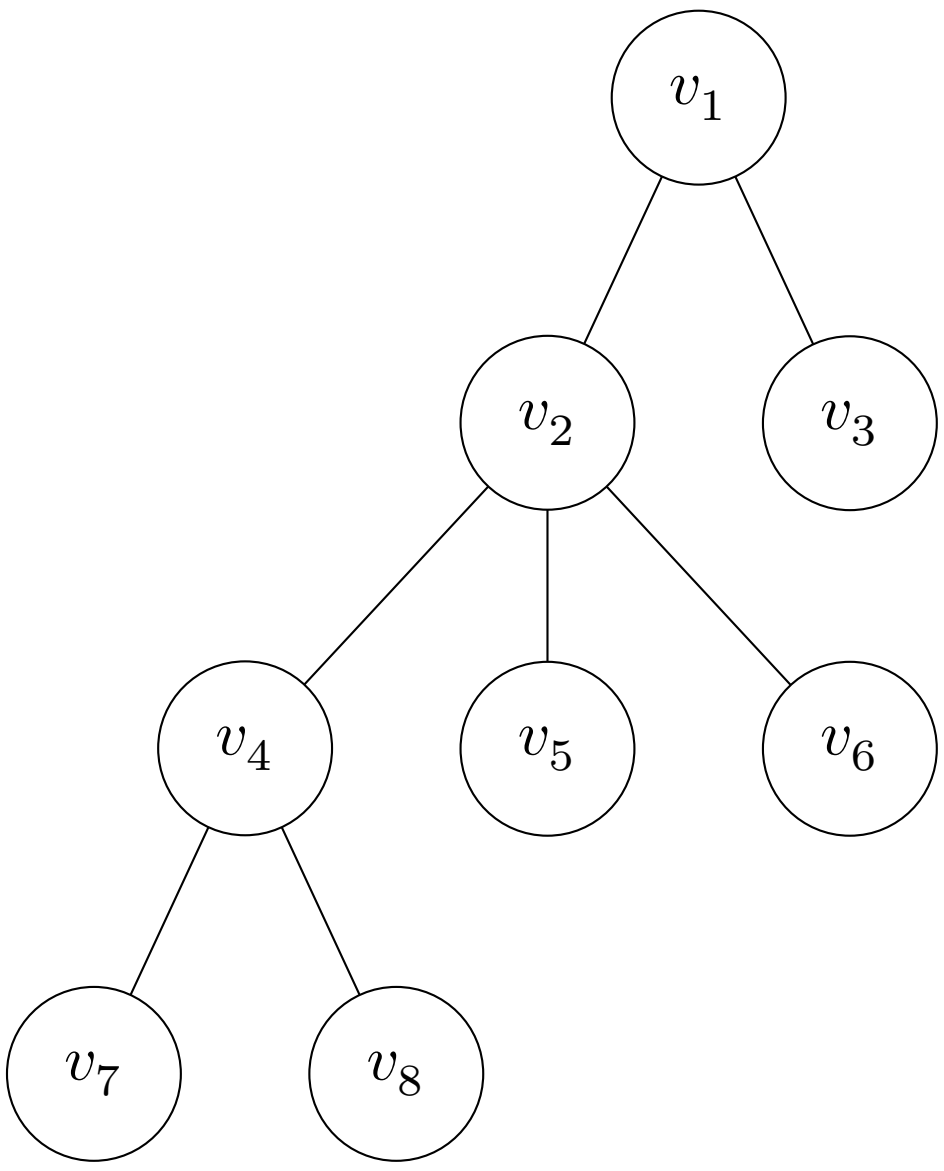

Whenever we talk about traversal trees, we shall make use of the following shorthand.

**Notation 4.4.6**

Let $T = (V, E, \langle\!|_|\!\rangle, \star, \preccurlyeq)$ be a traversal tree, and let $v \in V$. We write

$$\mathrm{N}^{\mathrm{out}} v = [v_1, \dots, v_n]$$

to mean that $\mathrm{N}^{\mathrm{out}} v = \{ v_1, \dots, v_n \}$ and $v_1 \preccurlyeq_v v_2 \preccurlyeq_v \dots \preccurlyeq_v v_n$.

In the example above, we would write

$$\mathrm{N}^{\mathrm{out}} v_1 = [v_2, v_3], \quad \mathrm{N}^{\mathrm{out}} v_2 = [v_4, v_5, v_6],$$

and so forth.

**Definition 4.4.7**

Let $(T, \preccurlyeq)$ and $(S, \trianglelefteq)$ be traversal trees. A *traversal tree morphism* from $(T, \preccurlyeq)$ to $(S, \trianglelefteq)$ is a directed hypergraph morphism $(g, h) \colon T \to S$ such that

$$\forall v \in T \colon \forall x, y \in \mathrm{N}^{\mathrm{out}} v \colon x \preccurlyeq_v y \Rightarrow gx \trianglelefteq_{gv} gy.$$

In other words, we require morphisms between traversal trees to be compatible with the imposed ordering. Note that we need not worry about the fact that hypergraph morphisms allow reordering of end points since we are in the realm of directed trees, where every edge has exactly one start point and one end point.

# 5. Typed Graphs

Later on, we shall want to build every graph from a few simple constituents. In short, we would like the following graph:

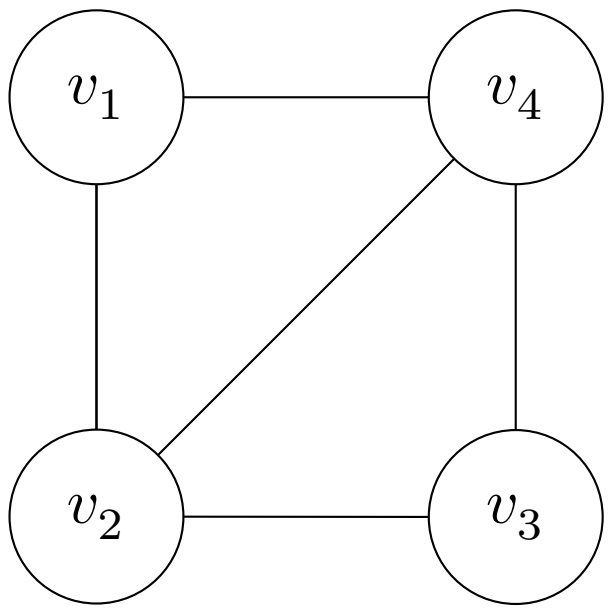

to be constructed from the two smaller graphs

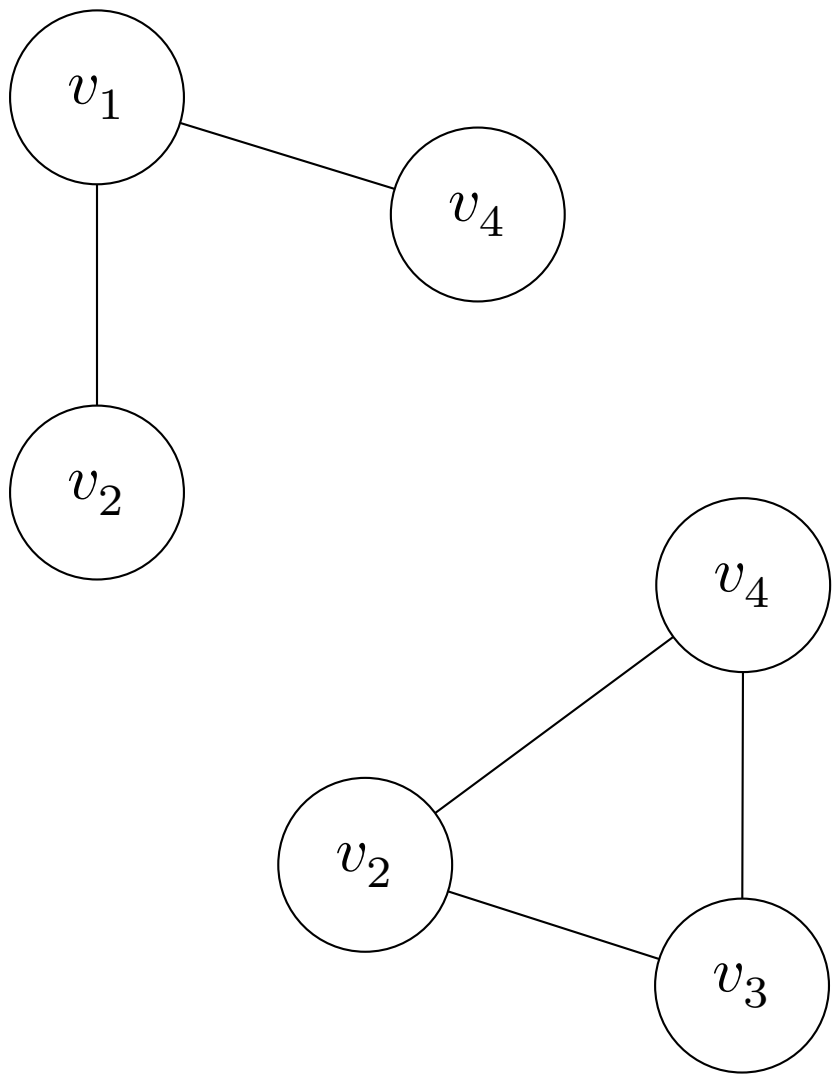

by somehow "gluing" the same-label vertices together:

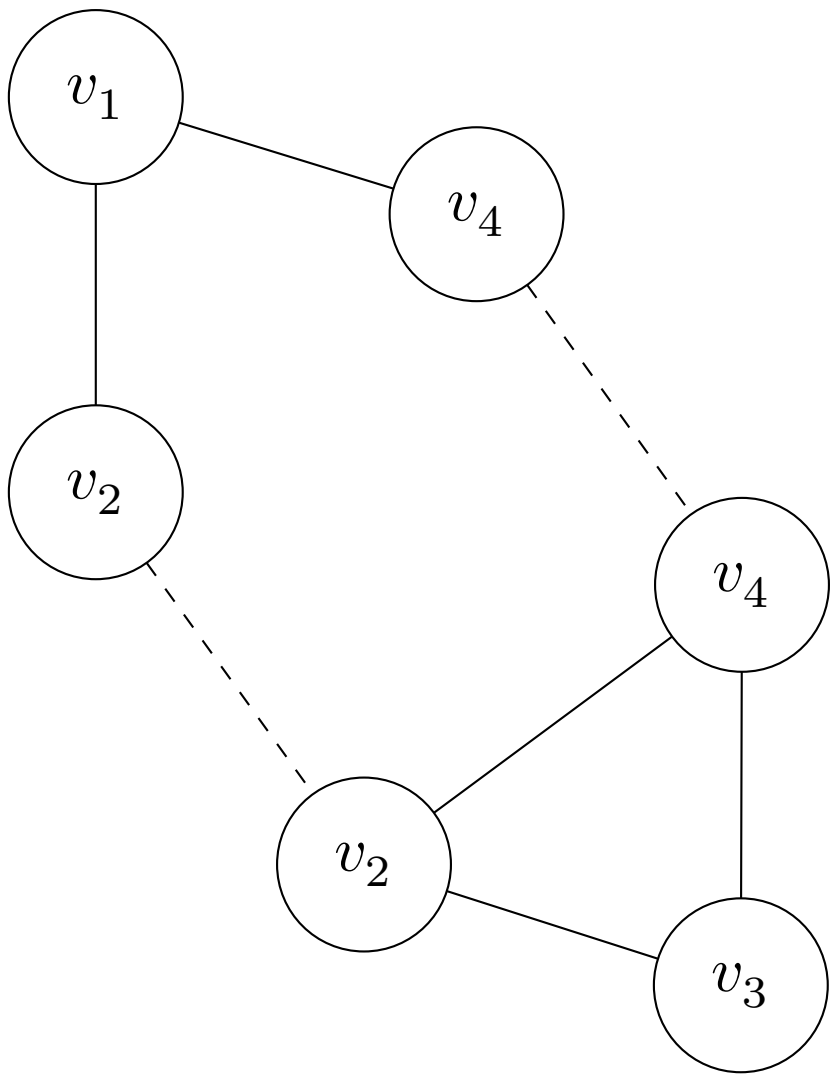

The reason we desire this is that in the proof of Courcelle's Theorem, we can then use the trick that if a graph $G$ is obtained by gluing together graphs $G_1$ and $G_2$, then certain properties of $G$ can be checked by only looking at $G_1$ and $G_2$, which are smaller than the original graph – effectively a divide and conquer approach to the proof.

The reader can easily imagine that, given arbitrarily many copies of the graph with one edge and two vertices and the ability to glue vertices together, they can build any finite 2-uniform graph they desire.[6] However, the structures we develop in chapter 6 to prove Courcelle's Theorem only allow gluing at "some" vertices in order to preserve the properties we need.

To prepare for this restriction, we slightly extend our definition of "graph". A *typed* graph is a graph as we know it, with the additional provision that some of its vertices are marked as special. Courcelle calls these vertices "sources", but the reader should be advised that these bear no relation to the sources of a network flow. We choose the less overloaded denomination "terminals" instead.

The only special property of these terminals is that they mark the vertices at which we are allowed to "glue". All other vertices are glue-resistant.

[6] This is proven in a more general form in theorem 6.4.2.

All definitions in this section work equally well for directed as for undirected graphs. We give the undirected version and trust the reader to be able to extend the concepts as necessary.

**Definition 4.5.1**

Let $n \in \mathbb{N}$. A *typed graph of type n* is a pair $(G, t)$ where $G = (V, E, \langle\!\!\langle _ \rangle\!\!\rangle)$ is an undirected graph and $t \colon \{1, \ldots, n\} \to V$.

The vertices in the range of $t$ are called *terminals*.

For ease of notation, we denote the $i$-th terminal vertex by $t_i := t(i)$.

In our examples, we shall colour the terminal vertices, and since non-terminal vertices are inherently uninteresting, we usually omit their labels and only label the terminal vertices (by their terminal index or indices). In the example above, we would have $v_2 = t(1)$ and $v_4 = t(2)$.

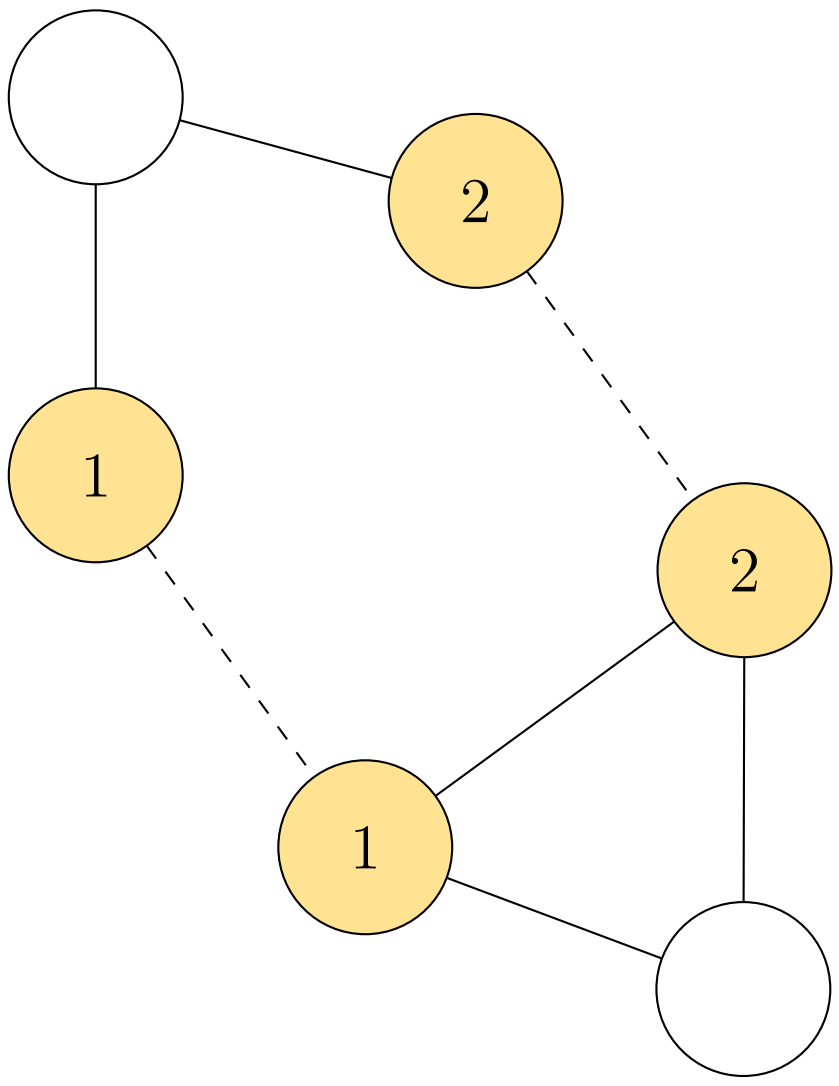

We give a name to the collection of graphs of a certain type.

**Definition 4.5.2**

We denote the set of all finite graphs of type $n \in \mathbb{N}$ by $\mathfrak{G}_n$.

We now define procedures that operate on graphs. Whenever we are given finitely many graphs, we assume without loss of generality that their sets of vertex and edge symbols are disjoint.

For each construction in this section, we give an intuitive explanation and then a formal definition.

The first construction, the disjoint sum of two graphs, will be familiar to most readers. We simply take two graphs, draw them side by side, and take this as our new graph.

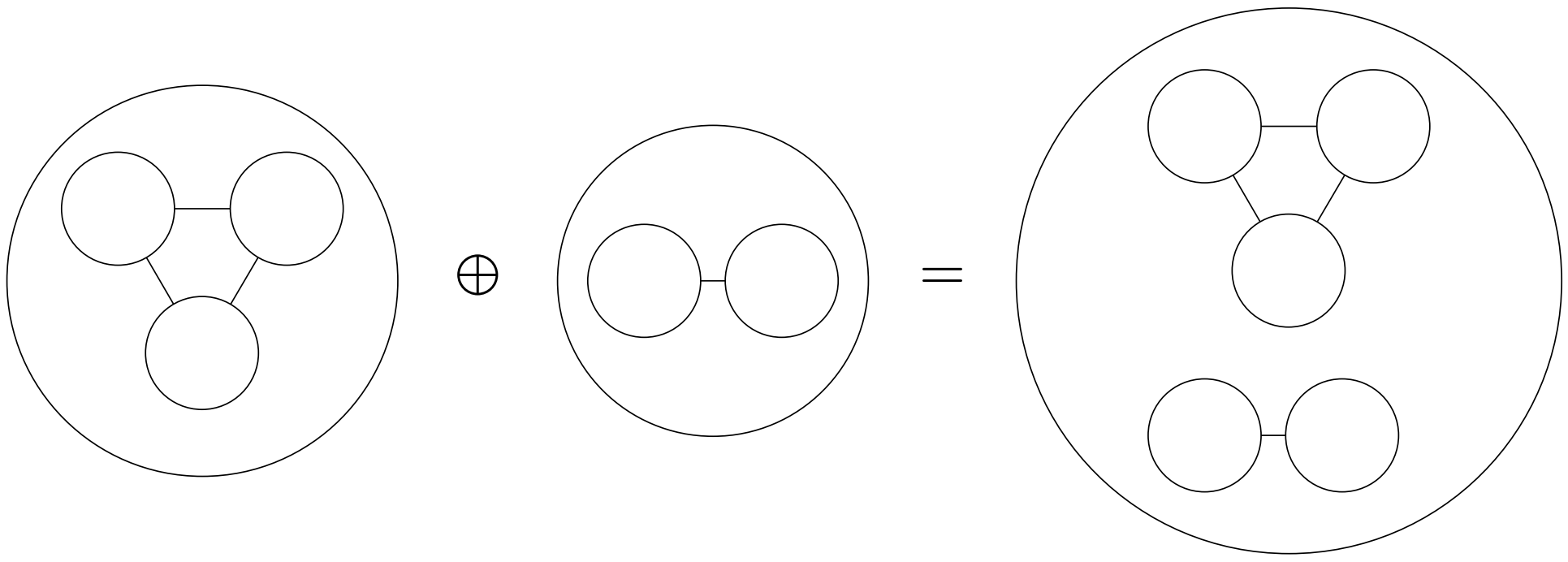

### Definition 4.5.3

Let $G = (V, E, \langle\!|_|\!\rangle, t)$, $G' = (V', E', (\!|_|\!), t')$ be graphs of type $n$ and $m$, respectively. Then $G \oplus G'$, called the *disjoint sum* of $G$ and $G'$, is the typed graph $(V'', E'', [\!|_|\!], t'')$ of type $n + m$ with

- $V'' := V \cup V'$,
- $E'' := E \cup E''$,
- $[\!|_|\!] \colon e \mapsto \begin{cases} \langle\!|e|\!\rangle & \text{if } e \in E \\ (\!|e|\!) & \text{if } e \in E', \end{cases}$
- $t'' \colon \{1, \dots, n + m\} \to V'', i \mapsto \begin{cases} t_i & \text{if } i \leq n \\ t'_{i-n} & \text{if } i > n. \end{cases}$

Next we introduce a way to switch the order of terminals or forget about one or more terminal vertices, essentially changing the labels under which the terminals are known.

Say we are given the following type 3 graph.

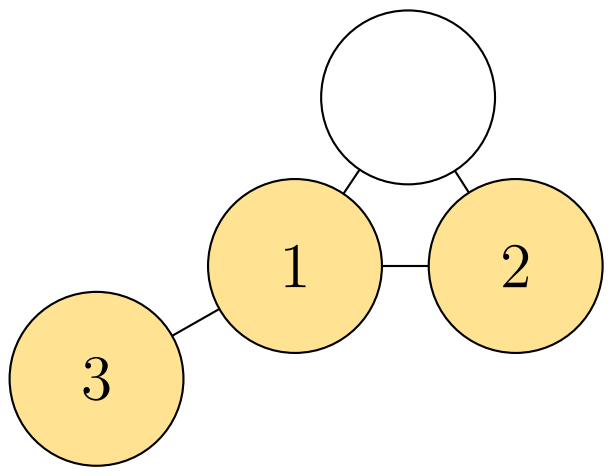

We would like a type 2 graph instead, dropping the first terminal. Consequently, the other vertices must be relabeled. We can model this by expressing the "new" terminals in terms of the old ones: terminal 1 should be the old terminal 2, while terminal 2 should be the old terminal 3.

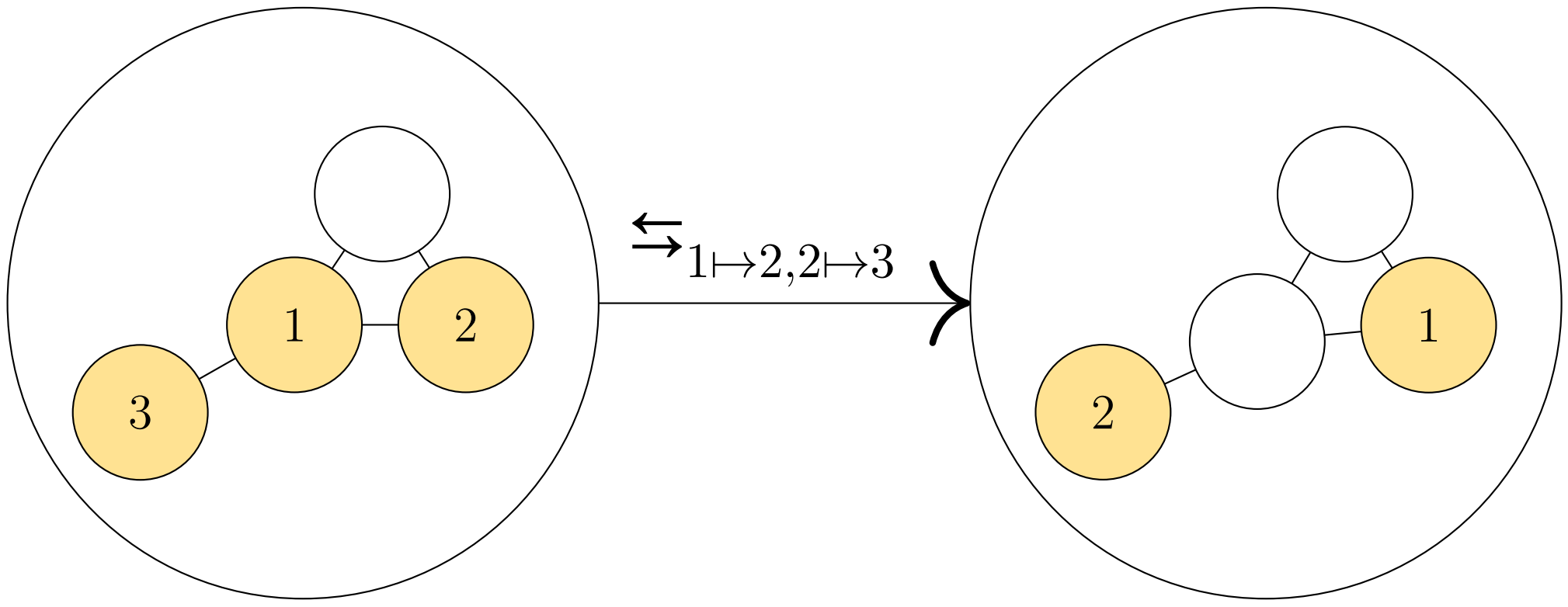

### Definition 4.5.4

Let $n \in \mathbb{N}$ and let $G = (V, E, \langle\!|_|\!\rangle, t)$ be a graph of type $n$. Let further $k \in \mathbb{N}$ and let $\sigma\colon \{1, \dots, k\} \to \{1, \dots, n\}$ be an arbitrary map.

Then $\leftrightarrows_\sigma G$ is the typed graph $(V', E', \llparenthesis_\rrparenthesis, t')$ of type $k$ with

- $V' := V$,
- $E' := E$,
- $\llparenthesis_\rrparenthesis := \langle\!|_|\!\rangle$,
- $t' := t \circ \sigma$.

We call this construction *terminal redefinition*.

Terminal redefinition makes a graph of type $n$ into a graph of type $k$. If $k = 0$, this results in a graph without terminal vertices. If $k > n$, the resulting graph will necessarily have at least one vertex that occurs more than once in the sequence of terminal vertices, because source redefinition cannot promote a non-terminal vertex to a terminal vertex.

The latter is an important fact to remember – once a vertex loses its terminal status, it can never regain it.

Last in our list of constructions is a way to glue exactly two terminal vertices together.

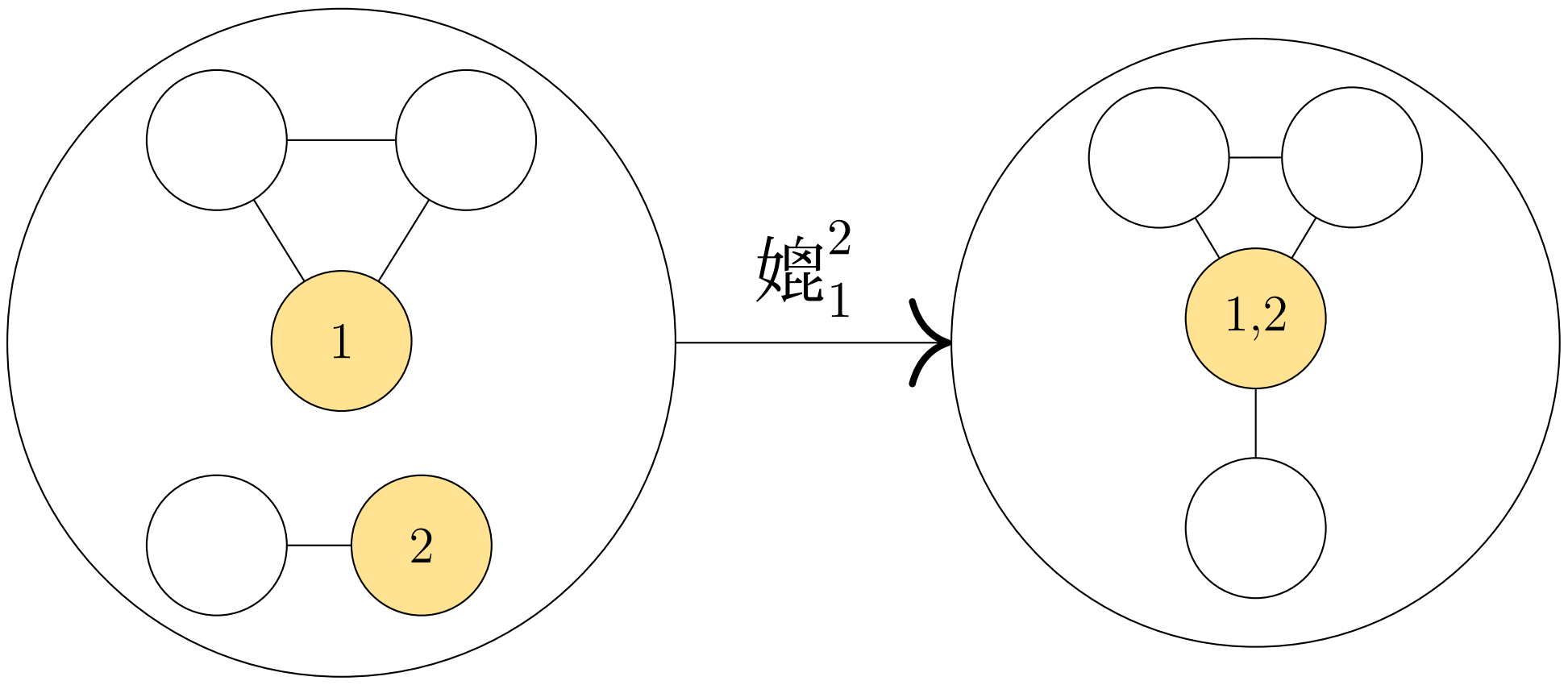

## Definition 4.5.5

Let $n \in \mathbb{N}$ and let $G = (V, E, \langle _ \rangle, t)$ be a graph of type $n$. Let now $a, b \in \{1, \dots, n\}$.

媲 (tsureai) is a Japanese character meaning to marry, to pair off.

For $t_a \neq t_b$, we set 媲$_a^b G$ to be the type $n$ graph $G' = (V', E', \llparenthesis _ \rrparenthesis, t')$, called the *fusion* of $G$ with respect to $a$ and $b$, with

- $V' := V \setminus \{t(b)\}$,
- $E' := E$,
- $\llparenthesis _ \rrparenthesis \colon e \mapsto \zeta^*(\langle e \rangle)$, where $\zeta \colon V \to V', v \mapsto \begin{cases} v & v \neq t(b) \\ t(a) & v = t(b), \end{cases}$

- $t' \colon \{1, \dots, n\} \to V', i \mapsto \begin{cases} t(i) & t(i) \neq t(b) \\ t(a) & t(i) = t(b). \end{cases}$.

For $t_a = t_b$, we set $\text{媲}_a^b G := G$.

In other words, fusion identifies two terminal vertices with each other, regluing edges as needed. The multiplicity of the new terminal vertex grows accordingly such that the graph's type remains unchanged.

Last but not least we introduce the basic building blocks from which all other graphs will eventually be constructed using the above tools.

**Definition 4.5.6**

We denote by $\mathfrak{v}$ the type 1 graph with one vertex and no edges.

For $n \in \mathbb{N}_{>0}$, we denote by $\mathfrak{e}_n$ the type $n$ graph $(V, E, \langle\!| _ |\!\rangle, t)$ with $n$ vertices $v_1, \dots, v_n$, one edge $e$ connecting all the vertices of $V$ (that is, $\langle\!| e |\!\rangle = v_1 \dots v_n$) and $t \colon i \mapsto v_i$.

We draw the first few graphs from the above definition. Note how $\mathfrak{e}_1$ is the peculiar hypergraph with an edge with only one end point, which is particularly hard to draw.

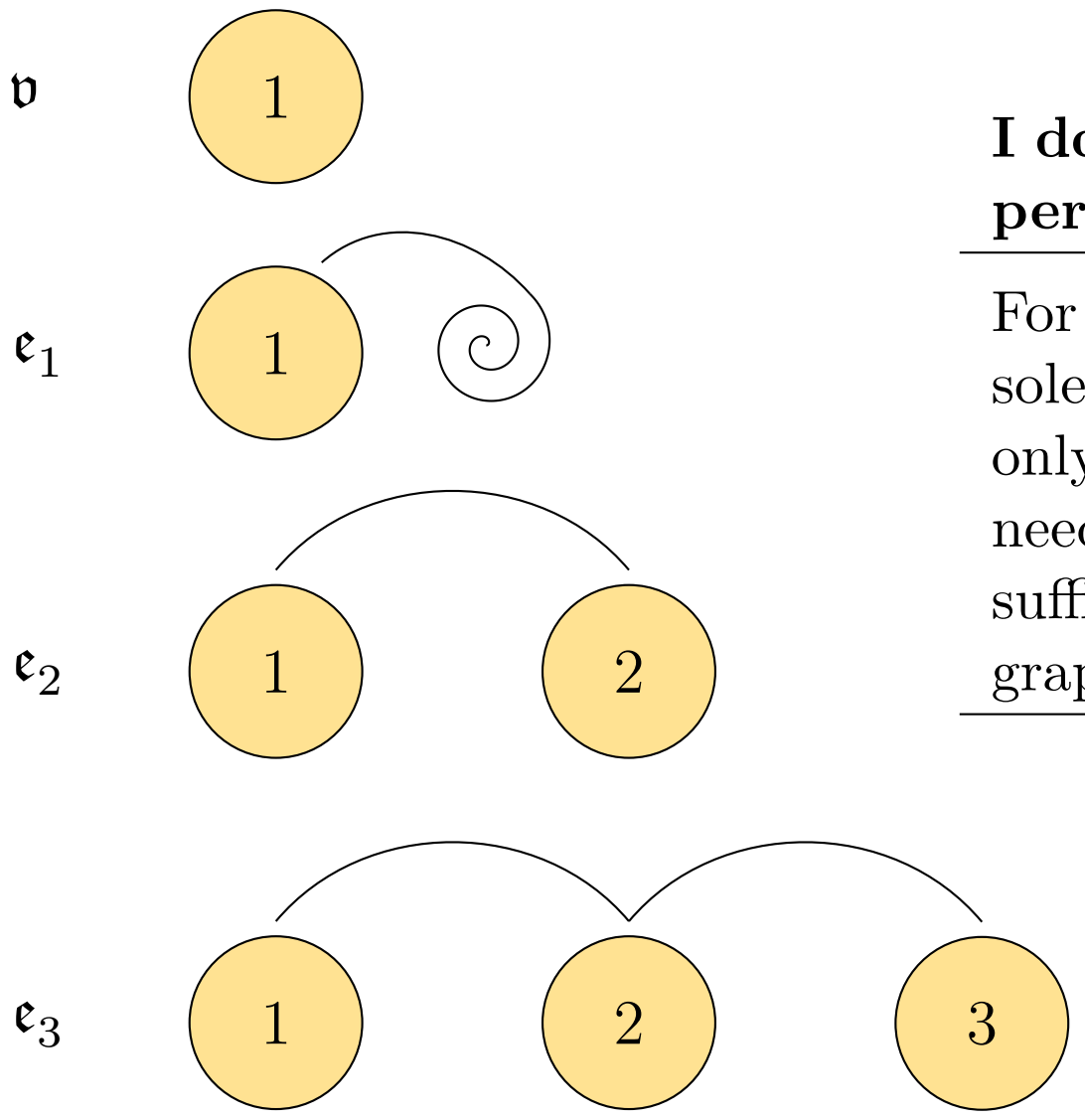

**I do not want to build hypergraphs!**

For the reader interested solely in 2-uniform graphs, only two trivial graphs are needed: $\mathfrak{v}$ and $\mathfrak{e}_2$. These suffice to build all "normal" graphs.

We call these building blocks the "trivial graphs". They suffice to build from scratch any finite (typed) hypergraph in finitely many steps, as we prove in theorem 6.4.2.

# 6. Tree- and Path-Decompositions

In chapter 9, we show how to apply the theoretical result of Courcelle's Theorem to obtain information on graphs of a certain tree- or path-width. We present here the relevant basics; additional (but entirely optional) information can be found in [Die05].

**I want to build more hypergraphs!**

The classical definition of hypergraphs prohibits having an edge $e$ with $\langle e \rangle = \varepsilon$. If the reader wishes to allow such edges, however, it is trivial to add an additional step to each of the algorithms presented later on that simply adds the requisite amount of end point-less edges to the hypergraph.

### Definition 4.6.1

Let $G = (V, E, \langle _ \rangle)$ be graph. A *tree-decomposition* of $G$ is a tuple $(T, X, b)$ with the following properties.

- $X$ is a family of subsets of $V$, called the *bags* of the decomposition.
- $\bigcup_{x \in X} x = V$, that is, every vertex of $G$ is contained in at least one bag.
- $\forall e \in E \colon \langle e \rangle = v_1 \dots v_{|e|} \Rightarrow \exists x \in X \colon v_1, \dots, v_{|e|} \in x$, that is, the collection of vertices of a given edge always shares at least one bag.[7]
- $T = (V', E', \langle\!\langle _ \rangle\!\rangle)$ is a tree.
- $b$ is a function $V' \to X$, that is, it associates with each node of $T$ a bag of vertices of $G$.
- $b$ is bijective.
- For any $v \in V$, the subgraph of $T$ induced by $\{ w \in T : v \in b(w) \}$ is connected, that is, if $v$ is contained in two bags $x, y$, it must also be contained in all bags on the path from $b^{-1}(x)$ to $b^{-1}(y)$.

The integer

$$\max_{x \in X} |x| - 1$$

is called the *width* of the decomposition, denoted $|(T, X, b)|$.

Note well that $X$ is a family, not a set – more than one node of $T$ can be associated the same set of vertices of $G$, which we still consider as different "bags". An alternative way to write this would be to index $X$ with the nodes of $T$, writing $X_v$ instead of $b(v)$.

Tree-decompositions allow us to define an important graph property, namely the smallest possible width of a tree-decomposition of $G$.

**Node-ice something?**

When talking about graphs, whether directed or undirected, we always talk about "vertices" and "edges", never "nodes" or "arcs". The only exception are tree- and path-decompositions, where we shall refer to the vertices of the underlying tree exclusively as "nodes" in order to make it easier for the reader to distinguish between vertices of the decomposed graph and nodes of the decomposition itself.

### Definition 4.6.2

Let $G$ be a graph. The number

$$\min_{(T,X,b) \text{ is a tree-decomposition of } G} |(T, X, b)|$$

is called the *tree-width* of $G$, denoted $\operatorname{tw}(G)$.

Tree-decompositions lend themselves nicely to dynamic programming. We make them even nicer in order to run algorithms on them later on. Since these algorithms need to decide in which order to descend into a tree, we usually assign an arbitrary orientation to the underlying tree.

### Definition 4.6.3

Let $G$ be a graph. A tree-decomposition $(T, X, b)$ for $G$ is called *rooted* if $T$ is a directed tree. The root of $T$ is then also called the root of the tree-decomposition, and we write $\sqrt{(T, X, b)} := \sqrt{T}$.

[7] For 2-uniform graphs, this condition simplifies to $\langle\!| e |\!\rangle = v_1 v_2 \Rightarrow \exists x \in X \colon v_1, v_2 \in x$.

**Definition 4.6.4**

Let $G = (V, E)$ be a graph. A rooted tree-decomposition $(T, X, b)$ for $G$ is called *nice* if every node $v \in T$ is of one of the following types.

- $v$ is a leaf, and $|b(v)| = 1$.
- $v$ has exactly one successor $w$, and there is a vertex $v' \in b(v)$ such that $b(w) = b(v) \setminus \{ v' \}$. We then call $v$ an *introduce node.*
- $v$ has exactly one successor $w$, and there is a vertex $v' \in b(w)$ such that $b(v) = b(w) \setminus \{ v' \}$. We then call $v$ a *forget node.*
- $v$ has exactly two successors $w_1$ and $w_2$, and they contain the same vertices, that is, $b(v) = b(w_1) = b(w_2)$. We then call $v$ a *join node.*

It is well-known that once a tree-decomposition is found, one can turn it into a nice tree-decomposition.

**Theorem 4.6.5**

Let $G$ be a graph, and let $Z = (T, X, b)$ be a rooted tree-decomposition for $G$. Then there exists also a nice tree-decomposition $Z' = (S, Y, c)$ for $G$ with $|Z'| \leq |Z|$ and $c(\sqrt{Z'}) = b(\sqrt{Z})$, and such a decomposition can be computed in time $\mathcal{O}(|Z| \cdot |T|)$.

**Proof.** While the base result is well-known, an algorithm which also preserves the root bag (an additional property we need in chapter 9 for some constructions which explicitly invoke the root of a given tree-decomposition) can be found in [KN12, pp. 359–362]. The runtime of their algorithm is straightforward to see.

□

In [RS83], Robertson and Seymour introduce in addition the notion of *path-decomposition.*

**Definition 4.6.6**

A *path-decomposition* of a graph $G$ is a tree-decomposition $(T, X, b)$ of $G$ such that $T$ is a path graph, that is, there is a simple path in $T$ such that all nodes of $T$ lie on that path.

The integer

$$\min_{(T,X,b) \text{ is a path-decomposition of } G} |(T, X, b)|$$

is called the *path-width* of $G$, denoted $\mathrm{pw}(G)$.

The decomposition $(T, X, b)$ is called *rooted* if $T$ is, and the root of $T$ is then also called the root of the decomposition, written $\sqrt{(T, X, b)} := \sqrt{T}$.

A rooted path-decomposition is called *nice* if it is nice as a tree-decomposition.

Note that the root of a rooted path-decomposition is not required to be one of the end points of the path. Thus, a nice path-decomposition can have up to one join node.

Just as with tree-decompositions, we can always find a particularly nice path-decomposition.

**Theorem 4.6.7**

Let $G$ be a graph, and let $Z = (T, X, b)$ be a rooted path-decomposition of width $k$ for $G$. Then there exists a nice path-decomposition $Z' = (S, Y, c)$ for $G$ of width at most $k$ and with $c(\sqrt{Z'}) = b(\sqrt{Z})$, and such a decomposition can be computed in time $\mathcal{O}(k \cdot |T|)$.

**Proof.** Let $G$ be a graph and let $(T, X, b)$ be a rooted path-decomposition for $G$ of width $k$ with $\sqrt{T} = v$. Since a path-decomposition is in particular a tree-decomposition, theorem 4.6.5 yields a nice tree-decomposition $(S, Y, c)$ of the same or smaller width with $\sqrt{S} = v'$ and $c(v') = b(v)$. But the only time the proof of theorem 4.6.5 creates a join node is when it encounters a node in $T$ with more than two successors, of which there are none.

Hence $(S, Y, c)$ is nice with at most one join node and thus actually a path-decomposition.

□

# Chapter 5
# Logic

When giving a proof for Courcelle's Theorem, one can hardly avoid talking at least a little bit about formal logic. After all, the statement itself already talks about monadic second-order logic.

On the other hand, once one starts digging into the formal foundations of logic, one is often compelled to keep digging until one emerges somewhere in New Zealand, exhausted and covered in dirt.

We shall try an intermediate approach where we introduce only the parts of formal logic vital to our explanations.

This chapter introduces the formal foundations in a way that assumes no prior knowledge on the part of the reader. A reader already familiar with mathematical logic should feel free to skim only the graph-related sections.

Any introduction we give here must necessarily be incomplete and in parts informal due to space constraints. For a more complete picture, the reader is welcome to check [End72], after which much of our notation is modeled.

Sections 5.1 and 5.2 should be read as a mostly intuitive (as opposed to formally sound) introduction. Formal definitions are given in section 5.3. Sections 5.4 to 5.6 apply our newly acquired knowledge to the realm of finite graphs.

# 1. Propositional Logic

By the principles of induction, the journey to understanding $k$-th-order logic is begun by understanding 0th-order logic and then incrementing $k$. *Propositional logic* is what can conceivably be called 0th-order logic, and is something with which any first-year undergrad is familiar. Syntactically, the language of propositional logic consists of so-called *well-formed formulas* built from the four symbols $($, $)$, $\neg$, $\wedge$ and countably many *proposition symbols* $\lambda_0, \lambda_1, \ldots$ (remember how we use colour to visually distinguish the symbols of formal languages). It is easy to define formally what "well-formed" actually means, but there is little to be gained from such an excursion in the current context. We instead trust that the reader possesses mathematical intuition enough to see that the formula

$$(\lambda_0 \wedge \lambda_1) \wedge \neg \lambda_2$$

is well-formed, while

$$)))\lambda_0 \wedge \neg \wedge$$

is not so.

Traditionally, one also allows symbols such as $\vee$, $\rightarrow$, or even $\dot{\vee}$[1] as shortcuts for commonly used constructions. An example for those confused by our omission of these symbols from the formal base language: the expression

$$\lambda_0 \vee \lambda_1$$

is tautologically equivalent to the expression

$$\neg(\neg \lambda_0 \wedge \neg \lambda_1),$$

which is to say, no matter which combination of truth values[2] we assign to $\lambda_0$ and $\lambda_1$, the formulas yield the same boolean result.

The decision to use $\wedge$ as a base symbol is somewhat arbitrary; all results in logic follow just as well if we use $\vee$ (or even $\rightarrow$) as the base and demote $\wedge$ to the status of shortcut.

[1] This being one of the many symbols logicians over the years have proposed for the exclusive or.

[2] We have indeed not defined truth values yet. The impatient reader may skip ahead about three sentences to see what they are about.

Semantically, the proposition symbols represent statements that are either true or false, such as

$$\lambda_0 = \text{“Sokrates loves chocolate.”}$$

In the coldhearted world of propositional logic, Sokrates either loves chocolate or does not love it; there is no middle ground. Formally, in order to do anything useful with a well-formed formula, we need to assign truth values to all its propositional symbols.[3] The truth values are the elements of the set $\{ \top, \bot \}$, where the verum $\top$ is interpreted as “true” (Sokrates does love chocolate) and the falsum $\bot$ is interpreted as “false” (Sokrates’ feelings toward chocolate are not of an adoring nature).

A *truth assignment* for a set $\Lambda$ of proposition symbols is then a function $\tau\colon \Lambda \to \{ \top, \bot \}$. One can show[4] that any such truth assignment uniquely extends to an assignment of truth values to the well-formed formulas using only proposition symbols in $\Lambda$, whence it is reasonable to speak of the truth or falsity of not just propositional variables, but also of formulas.

Given a propositional formula like, for example, $\varphi = \lambda_0 \wedge \lambda_1$ and a truth assignment $\tau\colon \{ \lambda_0, \lambda_1 \} \to \{ \top, \bot \}$, we write $\varphi[\tau]$ to mean “the formula $\varphi$ with each proposition symbol replaced by its image under $\tau$”. In particular, we note again the difference between syntax and semantics: while $\varphi$ is an abstract string of symbols, $\varphi[\tau]$ can be assigned a truth value – in this example, it is true if $\tau(\lambda_0) = \tau(\lambda_1) = \top$ and false otherwise.

As a shorthand, we allow writing the values assigned to the proposition symbols in $\varphi$ as a list, like $\varphi[\top, \bot]$, where the list is sorted by the index of the variables.

In a certain sense, this covers (informally) everything there is to know about propositional logic, at least for the purposes of this work. As a result of this simplicity, propositional logic is of severely limited expressiveness.

Enter the first induction step.

[3] Remember our example at the beginning of section 3.1 where we talked about the difference between syntax and semantics!

[4] This statement follows from the recursion and unique readability theorems, both of which can be found in [End72].

# 2. Predicate Logic

Whilst perusing the previous section, the reader might well have wondered about the lack of mention of the quantifier symbols $\exists$ and $\forall$. The reason for this is now revealed: propositional logic has no notion of "universe", that is, no notion of a certain space where variables can live. Consequently, it cannot talk about "all" variables, and indeed not even assert that a variable "exists" satisfying some property, because it would have no idea where to look for this object.

This oversight is fixed by first-order logic, often called *predicate logic* for its introduction of *predicate variables*. Of course, the added expressiveness (indeed, predicate logic is the tool that a mathematician uses, perhaps unwittingly, in most of their proofs) comes at the price of greater complexity.

## 2.1. Syntactical Building Blocks

As with propositional logic, our predicate syntax admits the atomic symbols $($, $)$, $\neg$, and $\wedge$. However, the propositional variables are replaced by the following. We give a short intuition of what these symbols are going to represent, even though right now they are purely syntactical.

First, we introduce countably many *variable symbols* $\nu_0$, $\nu_1$, and so forth. These are going to be just that: variables, to facilitate expressions like "$\forall \nu_0 \colon \nu_0$ likes chocolate".

Next are the fabled *predicates*. For every positive integer $n$, we introduce a set ${}_n\mathit{\Lambda}$ of $n$*-place predicate symbols*. The set ${}_n\mathit{\Lambda}$ may be empty, finite, or infinite, and need not be countable. These will be relations between the constants of our universe. For example, by introducing the 1-place predicate symbol ${}_1\mathit{\Lambda}_0$ with the meaning "loves chocolate", our previous assertion about Sokrates becomes "${}_1\mathit{\Lambda}_0(\text{Sokrates})$". If our universe happens to contain more foodstuffs than just cocoa derivatives, we might instead introduce a 2-place predicate symbol ${}_2\mathit{\Lambda}_0(\nu_0, \nu_1)$ with the meaning "$\nu_0$ loves $\nu_1$" and reformulate the assertion as "${}_2\mathit{\Lambda}_0(\text{Sokrates, chocolate})$". Of course, there is nothing preventing us from introducing both ${}_1\mathit{\Lambda}_0$ *and* ${}_2\mathit{\Lambda}_0$ into our syntax. The intuitive notion that the predicates should agree (that is, it cannot be that ${}_1\Lambda_0(\text{Sokrates})$ is true, but ${}_2\Lambda_0(\text{Sokrates, chocolate})$ is false) is a purely semantic one and hence not relevant at this stage. In the language of numbers, a well-known 2-place predicate symbol is the symbol $<$.

Now, for every *non-negative* integer $n$ (note well that this is one integer more than in the previous paragraph) we introduce a set ${}_n\Delta$ of $n$*-place function symbols*, again empty, finite, or even infinite. These will, as is the wont of functions, provide ways to turn $n$ elements of our universe into a different element of that same universe. In the language of numbers, addition is a 2-place function symbol.

Finally, we have enticed the reader with the promise of quantifiers. We deliver half of that promise by introducing the *universal quantifier* $\forall$. Given a well-formed formula $\varphi$, we decree that also the formula $\forall\nu_n(\varphi)$ shall be considered well-formed for all $n \in \mathbb{N}$. (Note that quantification is thus only allowed over variable symbols, not over predicates or functions.)

Any collection of symbols such as these defines a viable *first-order language*. As before, the language itself is a purely syntactical construct with no inherent meaning. Semantics are defined only once we choose an appropriate *structure*, which will be taking the place of the truth assignments from the previous section.

The smallest possible first-order language in our definition is the alphabet

$$\{\, (, ), \wedge, \neg, \forall, \nu_0, \nu_1, \dots \},$$

which contains no predicate symbols and no function symbols. This is also the least exciting first-order language, since it contains exactly zero well-formed formulas.

## 2.2. Shorthands

As before, it is convenient to agree upon several shorthand notations. For example, the reader might have been expecting an existential quantifier in the previous section. This is however not an atomic building block, as for a well-formed formula $\varphi$ depending on a variable symbol $\nu_0$, the formula

$$\exists\nu_0(\varphi)$$

is tautologically equivalent to the formula

$$\neg\forall\nu_0(\neg\varphi).$$

Likewise, we have not defined any symbols to refer to the constants of our universe. For example, if we were working in the universe of real numbers (that is, variables are understood to range over $(-\infty, \infty)$), the formula

$$\exists\nu_0(\nu_0 < 3)$$

would not actually be valid, because the symbol "3" has been assigned no meaning. However, a 0-place function symbol readily serves this function: it takes no parameters, then outputs an element of our universe. In our example, it makes sense to name the 0-place function symbols after the number they represent, that is, the constant function that outputs the number 3 would itself be called $3$. Then the expression

$$\exists\nu_0(\nu_0 < 3())$$

is well-formed and has exactly the meaning we intended. As a shorthand, it is only natural to drop the empty parentheses after the symbol $3$ and arrive at the expression we had originally planned on using.

Another shorthand that we have already snuck in is the use of *infix notation*. Syntactically, in our language of real numbers, the assertion that one and one makes two would use the 2-place function symbol $+$ and the 2-place predicate symbol $=$, hence reading

$$=(+(1,1),2).$$

This is quite convenient for formally defining which formulas are well-formed, but the average human might be more comfortable reading the formula

$$1+1=2.$$

Where no confusion is likely to occur, we silently drop prefix notation, as well as many, many "obvious" parentheses. We also sneak the symbol "$:$" into our formulas, which serves no semantic purpose, but can make statements such as

$$\forall\nu_0 : \forall\nu_1 : {}_2\Lambda_0(\nu_0, \nu_1)$$

more readable. Similarly, we use commas to separate symbols in expressions like

$$ {}_3\Lambda_0(\nu_0, \nu_1, \nu_2) $$

even though the correct formula would read

$$ {}_3\Lambda_0(\nu_0\nu_1\nu_2). $$

## 2.3. Semantics of Predicate Logic

With predicates and functions, a simple truth assignment will no longer do. With the concepts involved being less widely-known, we give a formal definition this time.

**Definition 5.2.1**

Let $\Gamma$ be a first-order language.

A *structure* for $\Gamma$ is a function $\Omega$ that satisfies the following properties:

- The domain of $\Omega$ is the union of $\{\forall\}$ with all sets of predicates and all sets of functions.
- The *universe* $|\Omega| := \Omega(\forall)$ of $\Gamma$ is a nonempty set.
- Let $\Lambda$ be an $n$-place predicate symbol of $\Gamma$. Then $\Omega(\Lambda) \subseteq |\Omega|^n$.
- Let $\Delta$ be an $n$-place function symbol of $\Gamma$. Then $\Omega(\Delta)$ is a function $|\Omega|^n \to |\Omega|$.[5]

Intuitively, these conditions mean the following:

- The universe is the collection of "all things about which we want to say something". When we write "$\forall\nu_0$", we mean "for all elements $\nu_0$ of the universe".
- Every $n$-place predicate is associated with an $n$-ary relation on the universe. Thus, for a 3-place predicate symbol $\Lambda$, we mean that $\Lambda(\nu_0, \nu_1, \nu_2)$ is true if and only if $(\nu_0, \nu_1, \nu_2) \in \Omega(\Lambda)$.
- Every $n$-place function symbol $\Delta$ is associated with an actual *function* (as distinct from the function *symbol*) that takes $n$ elements from the universe and outputs a new element from the universe. In particular, this fits with what we have said before about constants in $\mathit{\Gamma}$ being defined by 0-place function symbols.

[5] Recall that for any set $X$, we have $X^0 = \{\,()\,\}$, a one-element set.

If the reader has ever dabbled in set theory, they should by now be saying: this cannot be right, the above definition does not even let me define a structure in which to state the Zermelo-Fraenkel axioms, since the universe of such a structure would have to contain all sets and therefore cannot be a set itself. Indeed, some authors allow for universes that are not sets. However, since everything described in this work fits into a set, we have opted to circumnavigate this additional rock in our already stormy logical sea. For the reader not satisfied with this constraint, it should be straightforward to adapt the definitions to include proper classes.

We can now use the structure $\mathit{\Omega}$ to talk about the validity of formulas.

**Definition 5.2.2**

A well-formed first-order formula with no free variable symbols is called a *sentence*.

As always, a “free” variable symbol is one that is not “captured” by a preceding quantifier, that is, the expression

$$\forall \nu_0 (\mathit{\Lambda}(\nu_0) \wedge \mathit{\Lambda}(\nu_1))$$

contains $\nu_1$ as a free variable symbol and $\nu_0$ as a non-free variable symbol.

By assigning to each logical symbol its usual meaning, this means that given a structure $\mathit{\Omega}$ and a sentence $\varphi$, either $\varphi$ is true in that structure, or it is not.

**Definition 5.2.3**

Let $\mathit{\Gamma}$ be a first-order language, $\mathit{\Omega}$ a structure for $\mathit{\Gamma}$, and let $\varphi$ be a sentence of $\mathit{\Gamma}$. If $\varphi$, interpreted as an ordinary logical formula[6] in $\mathit{\Omega}$, is true, we say that $\mathit{\Omega}$ is a *model* for $\varphi$ and write $\vDash_{\mathit{\Omega}} \varphi$.

[6] This is of course informal, but we trust that the reader knows exactly how to deal with concrete logical formulas. If an even more rigorous definition is desired, the reader may want to read the entirety of [End72].

For formulas which are not sentences, this notion makes no sense: consider again the language of numbers, endowed with the structure of the real numbers. While the two sentences

$$1 + 1 = 2, \quad \exists \nu_0 (1 + \nu_0 = 2)$$

are both true and the sentences

$$1 + 2 = 2, \quad \forall \nu_0 (1 + \nu_0 = 2)$$

are false, the well-formed formula

$$1 + \nu_0 = 2$$

is neither true nor false, since it still contains a free variable with no assigned meaning.

# 3. Second-Order Logic

Predicate logic can express much more than its propositional cousin, but there is still room for improvement. Consider, for example, the following first-order statement for a 1-place predicate symbol $\Lambda$.

$$(\forall \nu_0 \Lambda(\nu_0)) \to (\exists \nu_0 \Lambda(\nu_0)).$$

This is a sentence as defined in definition 5.2.2, and since we have required the universe to be nonempty, it is clearly true.

Does it matter what exactly $\Lambda$ is? Forsooth, it does not, since the truth of this sentence does not actually depend on $\Lambda$ – if all $\nu_0$ satisfy $\Lambda$, then the nonempty universe contains at least one $\nu_0$ satisfying $\Lambda$, and if not all $\nu_0$ satisfy $\Lambda$, then the implication is vacuously true.

Consequently, anyone with an intuitive understanding of the concepts would agree that the "sentence"

$$\forall \Lambda \left( (\forall \nu_0 \Lambda(\nu_0)) \to (\exists \nu_0 \Lambda(\nu_0)) \right)$$

should be true. However, this is not even a well-formed formula in predicate logic, because we can quantify only over variable symbols, not over predicate or function symbols.

This is where second-order logic comes in.

## 3.1. Syntax of Second-Order Logic

A *second-order* language is defined just as a first-order language, except that we drop the variables $\nu_0, \nu_1, \ldots$ and instead introduce the following:

- For each positive integer $n$, countably many $n$-place *predicate variables* ${}_n\lambda_0, {}_n\lambda_1, \ldots$. These will allow quantification over predicates, as seen in the motivating example above.
- For each non-negative integer $n$, countably many $n$-place *function variables* ${}_n\delta_0, {}_n\delta_1, \ldots$. These, of course, allow quantification over functions.

Note that the role of variables previously played by $\nu_0, \nu_1, \ldots$ is now covered by the 0-place function variables ${}_0\delta_0, {}_0\delta_1, \ldots$ – not to be confused with the 0-place function *symbols*: a 0-place function symbol represents a fixed constant in our universe, while a 0-place function variable represents a single non-fixed element in our universe over which we can quantify.

Consequently, we widen the definition of a well-formed formula to include expressions built by quantification over the new symbols and specify that a sentence is a well-formed formula that contains no free predicate or function variables.

## 3.2. Second-Order Structures

Surprisingly, the transition to second-order brings nothing new in the semantics department. A structure is defined just as in definition 5.2.1 (recall that the definition makes no use of the variables). The truth or falsity of a sentence now takes into consideration quantification over the new variables.

## 3.3. Formal Definitions

We take the concepts introduced intuitively in the previous sections and roll them into a series of formal definitions for later reference. The fact that we assign names and notation to every little subset of a language might seem over-designed at first, but it will come in handy in chapter 6 when we carry out inductions over the structure of a well-formed formula.

**Definition 5.3.1**

A *second-order language* is an alphabet

$$\Gamma = \gamma_\Gamma \cup \bigcup_{i \in \mathbb{N}_{>0}} {}_i\lambda_\Gamma \cup \bigcup_{i \in \mathbb{N}} {}_i\delta_\Gamma \cup \bigcup_{i \in \mathbb{N}_{>0}} {}_i\Lambda_\Gamma \cup \bigcup_{i \in \mathbb{N}} {}_i\Delta_\Gamma$$

consisting of the following symbols.

- The set $\gamma_\Gamma := \{\, (, ), \wedge, \neg \,\} \cup \{\, {}^{\lambda}_{n}\forall : n \in \mathbb{N}_{>0} \,\} \cup \{\, {}^{\delta}_{n}\forall : n \in \mathbb{N} \,\}$ of *atomic symbols*.
- For every $n \in \mathbb{N}_{>0}$, the set ${}_n\lambda_\Gamma := \{\, {}_n\lambda_0, {}_n\lambda_1, \ldots \,\}$ of *n-place predicate variables*.
- For every $n \in \mathbb{N}$, the set ${}_n\delta_\Gamma := \{\, {}_n\delta_0, {}_n\delta_1, \ldots \,\}$ of *n-place function variables*.
- For every $n \in \mathbb{N}_{>0}$, a set ${}_n\Lambda_\Gamma$ of *n-place predicate symbols*. Any of these may be empty.
- For every $n \in \mathbb{N}$, a set ${}_n\Delta_\Gamma$ of *n-place function symbols*. Any of these may be empty.

All of these sets are pairwise disjoint.

The separation of the universal quantifiers into different symbols is purely for the benefit of notation in some later definitions. In virtually all circumstances we shall simply write $\forall$ and trust that it is clear from context which universal quantifier is meant.

As before, we allow shortcuts and clarifying notation such as commas and colons as long as their meaning seems clear.

Whenever we talk about variable symbols, we use the symbols defined above if we know what kind of variable symbol is meant, or the symbols $\mu_0, \mu_1, \ldots$ if we know not whether the variable symbol is a predicate or a function.

Note that, since the first three bullet points are fixed, defining a new second-order language only means specifying the sets ${}_n\Lambda_\Gamma$ and ${}_n\Delta_\Gamma$.

We now define recursively exactly which words of $\Gamma^*$ we consider to be well-formed. The elementary building blocks are the 0-place function symbols and 0-place function variables: like a magician from his hat, they fetch an element out of our universe $|\Omega|$ without needing any input.

**Definition 5.3.2**

Let $\Gamma$ be a second-order language. An *elementary term* over $\Gamma$ is a string of the form

$$\chi$$

for some $\chi \in {}_0\delta_\Gamma \cup {}_0\Delta_\Gamma$.

The set of all elementary terms over $\Gamma$ is denoted $\Gamma_{\text{Elem}}$.

${}_0\delta_\Gamma$ is the set of 0-place function variables of $\Gamma$.

${}_0\Delta_\Gamma$ is the set of 0-place function symbols of $\Gamma$.

Now that we have, from thin air, procured at least one element of $|\Omega|$, we can use it as input for functions of higher arity.

**Definition 5.3.3**

Let $\Gamma$ be a second-order language. We set

$$\Gamma^0_{\text{Term}} := \Gamma_{\text{Elem}}$$

and for $n \in \mathbb{N}_{>0}$

$$\Gamma^n_{\text{Term}} := \{\, \zeta(\chi_1, \dots, \chi_k) : k \in \mathbb{N}_{>0}, \zeta \in {}_k\delta_\Gamma \cup {}_k\Delta_\Gamma, \chi_1, \dots, \chi_k \in \Gamma^{n-1}_{\text{Term}} \,\}$$

and call

$$\Gamma_{\text{Term}} := \bigcup_{i \in \mathbb{N}} \Gamma^i_{\text{Term}}$$

the *set of terms over* $\Gamma$.

This allows us to write all combinations of constants and functions in $|\Omega|$. Next, we want to apply predicates to these objects to arrive at a truth value.

**Definition 5.3.4**

Let $\Gamma$ be a second-order language. An *atomic formula* over $\Gamma$ is a string of the form

$$\rho(\chi_1, \dots, \chi_n)$$

for some $n \in \mathbb{N}_{>0}$, some $\rho \in {}_n\lambda_\Gamma \cup {}_n\Lambda_\Gamma$ and some $\chi_1, \dots, \chi_n \in \Gamma_{\text{Term}}$. The set of all atomic formulas over $\Gamma$ is denoted $\Gamma_{\text{Atom}}$.

${}_n\lambda_\Gamma$ is the set of $n$-place predicate variables of $\Gamma$.

${}_n\Lambda_\Gamma$ is the set of $n$-place predicate symbols of $\Gamma$.

Note that once we have put a term or terms into a predicate, we have converted it from an element of the universe (for example, a number) to a truth value in $\{\top, \bot\}$. This process is irreversible – all information about the element of $|\Omega|$ has been lost, we only know whether it has the property encoded by $\rho$. For this reason, the inductive definition is one-way – we do not provide a way to go from atomic formulas back to terms. Instead, everything we do from this point on deals with truth values, meaning that we can introduce logical symbols like $\wedge$ which make sense between two truth values, but not for example between two numbers.

**Definition 5.3.5**

Let $\Gamma$ be a second-order language. We set

$$|\Gamma|^{0,0} := \Gamma_{\text{Atom}}$$

and for $h, w \in \mathbb{N}$ with $h + w > 0$

$$\begin{aligned} |\Gamma|^{h,w} := & \bigcup_{i,j\in\mathbb{N}, i<h, j<w} |\Gamma|^{i,j} \\ & \cup \{ \neg(\varphi) : \varphi \in |\Gamma|^{h-1,w} \} \\ & \cup \{ (\varphi)\wedge(\psi) : \varphi, \psi \in |\Gamma|^{h-1,w} \} \\ & \cup \{ {}^{\lambda}_{n}\forall \lambda(\varphi) : \varphi \in |\Gamma|^{h,w-1}, n \in \mathbb{N}_{>0}, \lambda \in {}_n\lambda_\Gamma \} \\ & \cup \{ {}^{\delta}_{n}\forall \delta(\varphi) : \varphi \in |\Gamma|^{h,w-1}, n \in \mathbb{N}, \delta \in {}_n\delta_\Gamma \}. \end{aligned}$$

${}_n\delta_\Gamma$ is the set of $n$-place function variables of $\Gamma$.

The elements of the set

$$|\Gamma| := \bigcup_{h,w\in\mathbb{N}} |\Gamma|^{h,w}$$

are called the *well-formed formulas* of $\Gamma$.

For a well-formed formula $\varphi$, the integer

$$\min\{ h \in \mathbb{N} : \exists w \in \mathbb{N} : \varphi \in |\Gamma|^{h,w} \}$$

is called the *height* of $\varphi$, while the integer

$$\min\{\, w \in \mathbb{N} : \exists h \in \mathbb{N} \colon \varphi \in |\Gamma|^{h,w} \,\}$$

is called its *width*.

Intuitively, the height of a formula denotes the height of its logical tree, while the width of a formula tells us the maximum number of nested quantifiers.

We have thus provided a way to construct all well-formed formulas over $\Gamma$. If a string in $\Gamma^*$ cannot be constructed by the rules above, then it is not well-formed.

We already know from section 5.2 what a free variable is.

**Definition 5.3.6**

$\Gamma_{\text{Term}}$ is the set of all terms of $\Gamma$. (def. 5.3.3, p. 60)

$|\Gamma|$ is the set of all well-formed formulas over $\Gamma$. (def. 5.3.5, p. 61)

Let $\Gamma$ be a second-order language, and let $\chi \in \Gamma_{\text{Term}}$. We denote by $\overset{\leftharpoonup}{\chi}$ the set of all variable symbols occurring in $\chi$.

Let now $\varphi \in |\Gamma|$. We denote by $\overset{\leftharpoonup}{\varphi}$ the set of all free variable symbols occurring in $\varphi$.[7]

The formula $\varphi$ is called a *sentence* if it contains no free variables, that is, if $\overset{\leftharpoonup}{\varphi} = \varnothing$.

The set of all sentences over $\Gamma$ is denoted by $\|\Gamma\| := \{\, \varphi \in |\Gamma| : \overset{\leftharpoonup}{\varphi} = \varnothing \,\}$.

For the sake of completeness, we also note again what a structure looks like in this case.

[7] The definitions for terms and well-formed formulas are compatible, but in a term, all variables are automatically free (since a term cannot contain quantifiers), whereas a well-formed formula may contain non-free variables.

**Definition 5.3.7**

Let $\Gamma$ be a second-order language. A *structure* for $\Gamma$ is a function $\Omega$ that satisfies the following conditions.

- The domain of $\Omega$ is the set

$$\{\, {}^{\lambda}_{n}\forall : n \in \mathbb{N}_{>0} \,\} \cup \{\, {}^{\delta}_{n}\forall : n \in \mathbb{N} \,\} \cup \bigcup_{n\in\mathbb{N}_{>0}} {}_{n}\Lambda_{\Gamma} \cup \bigcup_{n\in\mathbb{N}} {}_{n}\Delta_{\Gamma}.$$

- The *universe* $|\Omega| := {}^{\delta}_{0}|\Omega| := \Omega({}^{\delta}_{0}\forall)$ is a nonempty set.
- For $n \in \mathbb{N}_{>0}$, the *$n$-place predicate universe* ${}^{\lambda}_{n}|\Omega| := \Omega({}^{\lambda}_{n}\forall)$ is a set of subsets of $|\Omega|^n$.
- For $n \in \mathbb{N}_{>0}$, the *$n$-place function universe* ${}^{\delta}_{n}|\Omega| := \Omega({}^{\delta}_{n}\forall)$ is a set of functions $|\Omega|^n \to |\Omega|$.
- Let $n \in \mathbb{N}$, $\Lambda \in {}_{n}\Lambda_{\Gamma}$. Then $\Omega(\Lambda) \subseteq |\Omega|^n$.
- Let $n \in \mathbb{N}$, $\Delta \in {}_{n}\Delta_{\Gamma}$. Then $\Omega(\Delta)$ is a function $|\Omega|^n \to |\Omega|$.

${}_{n}\Lambda_{\Gamma}$ is the set of $n$-place predicate symbols of $\Gamma$.

${}_{n}\Delta_{\Gamma}$ is the set of $n$-place function symbols of $\Gamma$.

The universe tells us over which variables we can quantify, the predicate and function universes allow quantification over a chosen subset of functions and predicates, and the existing function and predicate symbols in formulas are assigned meaning.

Note that the $n$-place predicate universe ${}^{\lambda}_{n}|\Omega|$ need not contain all subsets that are available for $n$-place predicate *symbols*, or it may contain entirely different ones. The same goes for functions. This means that the universal quantifier might “see” fewer (or more) predicates than there are available for building well-formed formulas.

As before, truth values are only assigned to sentences, that is, well-formed formulas without free variables. However, we shall need to deal with free variables in several proofs.

### Definition 5.3.8

Let $\Gamma$ be a second-order language, and let $\Omega$ be a structure for $\Gamma$. A *variable assignment* in $\Omega$ is a function $\tau$ whose domain is a subset of all variable symbols of $\Gamma$, that is, a set

$${}_n\lambda_\Gamma$ is the set of $n$-place predicate variables of $\Gamma$.

${}_n\delta_\Gamma$ is the set of $n$-place function variables of $\Gamma$.

${}^{\lambda}_{n}|\Omega|$ denotes the $n$-place predicate universe. (def. 5.3.7, p. 63)

${}^{\delta}_{n}|\Omega|$ denotes the $n$-place function universe. (def. 5.3.7, p. 63)

$$X \subseteq \bigcup_{n \in \mathbb{N}_{>0}} {}_n\lambda_\Gamma \cup \bigcup_{n \in \mathbb{N}} {}_n\delta_\Gamma$$

and which fulfils

- $\forall n \in \mathbb{N}_{>0}\colon \tau(X \cap {}_n\lambda_\Gamma) \subseteq {}^{\lambda}_{n}|\Omega|$,
- $\forall n \in \mathbb{N}\colon \tau(X \cap {}_n\delta_\Gamma) \subseteq {}^{\delta}_{n}|\Omega|$.

We write $\varepsilon$ for the variable assignment with domain $\varnothing$ and call this the *empty assignment*.

In other words, a variable assignment takes some (but not necessarily all) variable symbols and assigns to them elements from the universe of the "correct" kind – predicates to predicate variables, functions to function variables, and of the arity the symbol expects.

### Definition 5.3.9

Let $\Gamma$ be a second-order language, $\Omega$ a structure for $\Gamma$, $\tau$ a variable assignment in $\Omega$, and $\varphi$ a term or a well-formed formula of $\Gamma$.

We say that $\tau$ is a *full* variable assignment for $\varphi$ if every free variable in $\varphi$ is in the domain of $\tau$.

We can now define recursively what "truth" means.[8] First, a variable assignment takes all *terms* and transforms them from meaningless symbols into elements of our universe.

[8] Not to be confused with the *circular* definition of truth, which we shall leave to the attorneys.

**Definition 5.3.10**

Let $\mathit{\Gamma}$ be a second-order language, $\mathit{\Omega}$ a structure for $\mathit{\Gamma}$, $\chi$ a term of $\mathit{\Gamma}$, and $\tau$ a full variable assignment for $\chi$ in $\mathit{\Omega}$.

We write $\chi[\tau]$ to mean the following.

- If $\chi \in \mathit{\Gamma}_{\text{Elem}}$:

$$\chi[\tau] := \begin{cases} \mathit{\Omega}(\mathit{\Delta})() & \text{if } \chi = \mathit{\Delta} \text{ for some } \mathit{\Delta} \in {}_0\mathit{\Delta}_{\mathit{\Gamma}} \\ \tau(\delta)() & \text{if } \chi = \delta \text{ for some } \delta \in {}_0\delta_{\mathit{\Gamma}}. \end{cases}$$

- If $\chi \in \mathit{\Gamma}^n_{\text{Term}}$ for some $n \in \mathbb{N}_{>0}$, $\chi = \zeta(\chi_1, \dots, \chi_k)$:

$$\chi[\tau] := \begin{cases} \mathit{\Omega}(\mathit{\Delta})(\chi_1[\tau], \dots, \chi_k[\tau]) & \text{if } \zeta = \mathit{\Delta} \text{ for some } \mathit{\Delta} \in {}_k\mathit{\Delta}_{\mathit{\Gamma}} \\ \tau(\delta)(\chi_1[\tau], \dots, \chi_k[\tau]) & \text{if } \zeta = \delta \text{ for some } \delta \in {}_k\delta_{\mathit{\Gamma}}. \end{cases}$$

$\mathit{\Gamma}_{\text{Elem}}$ is the set of elementary terms over $\mathit{\Gamma}$. (def. 5.3.2, p. 60)

${}_0\mathit{\Delta}_{\mathit{\Gamma}}$ is the set of 0-place function symbols of $\mathit{\Gamma}$.

${}_0\delta_{\mathit{\Gamma}}$ is the set of 0-place function variables of $\mathit{\Gamma}$.

$\mathit{\Gamma}^n_{\text{Term}}$ is the set of terms of depth at most $n$. (def. 5.3.3, p. 60)

We can then take a well-formed formula and see whether the elements yielded by the variable assignment fulfil it.

**Definition 5.3.11**

Let $\mathit{\Gamma}$ be a second-order language, $\mathit{\Omega}$ a structure for $\mathit{\Gamma}$, $\varphi$ a well-formed formula of $\mathit{\Gamma}$, and $\tau$ a full variable assignment for $\varphi$ in $\mathit{\Omega}$.

We write $\varphi[\tau] \leftrightarrow \top$, pronounced "$\varphi$ is true under $\tau$", to mean the following.

- If $\varphi$ is atomic, say $\varphi = \rho(\chi_1, \dots, \chi_n)$:

$$\varphi[\tau] \leftrightarrow \top :\Leftrightarrow \begin{cases} (\chi_1[\tau], \dots, \chi_n[\tau]) \in \mathit{\Omega}(\mathit{\Lambda}) & \text{if } \rho = \mathit{\Lambda} \text{ for some } \mathit{\Lambda} \in {}_n\mathit{\Lambda}_{\mathit{\Gamma}} \\ (\chi_1[\tau], \dots, \chi_n[\tau]) \in \tau(\lambda) & \text{if } \rho = \lambda \text{ for some } \lambda \in {}_n\lambda_{\mathit{\Gamma}}. \end{cases}$$

- If $\varphi = \neg(\psi)$:

$$\varphi[\tau] \leftrightarrow \top :\Leftrightarrow \psi[\tau] \not\leftrightarrow \top.$$

- If $\varphi = (\psi) \wedge (\zeta)$:

$$\varphi[\tau] \leftrightarrow \top :\Leftrightarrow \psi[\tau] \leftrightarrow \top \text{ and } \zeta[\tau] \leftrightarrow \top.$$

${}_n\mathit{\Lambda}_{\mathit{\Gamma}}$ is the set of $n$-place predicate symbols of $\mathit{\Gamma}$.

${}_n\lambda_{\mathit{\Gamma}}$ is the set of $n$-place predicate variables of $\mathit{\Gamma}$.

- If $\varphi = {}^{\lambda}_{n}\forall\lambda(\psi)$:

$$\varphi[\tau] \leftrightarrow \top :\Leftrightarrow \psi[\kappa] \leftrightarrow \top \text{ for every full variable assignment } \kappa \text{ for } \psi \text{ in } \Omega \text{ with } \kappa|_{\overline{\varphi}} = \tau|_{\overline{\varphi}}.$$

$\overline{\varphi}$ denotes the set of free variables of $\varphi$. (def. 5.3.6, p. 62)

- If $\varphi = {}^{\delta}_{n}\forall\delta(\psi)$:

$$\varphi[\tau] \leftrightarrow \top :\Leftrightarrow \psi[\kappa] \leftrightarrow \top \text{ for every full variable assignment } \kappa \text{ for } \psi \text{ in } \Omega \text{ with } \kappa|_{\overline{\varphi}} = \tau|_{\overline{\varphi}}.$$

We write $\varphi[\tau] \leftrightarrow \bot$ to mean $\varphi[\tau] \not\leftrightarrow \top$, pronounced "$\varphi$ is false under $\tau$", and given a second formula $\psi$ such that $\tau$ is also full for $\psi$, we write $\varphi[\tau] \leftrightarrow \psi[\tau]$ to mean that either both formulas are true or both formulas are false under $\tau$.

Note again that a truth value can only be assigned if every free variable is taken care of by a variable assignment.

Since a sentence has no free variables, the empty assignment is a full variable assignment for every sentence. We restate the definition of a sentence's truth value for future reference.

### Definition 5.3.12

$\|\Gamma\|$ denotes the set of sentences over $\Gamma$. (def. 5.3.6, p. 62)

Let $\Gamma$ be a second-order language, $\varphi \in \|\Gamma\|$, and let $\Omega$ be a structure for $\Gamma$. We say that $\Omega$ is a *model* for $\varphi$, written

$$\models_{\Omega} \varphi,$$

if $\varphi[\varepsilon] \leftrightarrow \top$.

Of course, since there are no free variables to replace, any other variable assignment on a sentence will produce the same truth value.

We use the same notation to denote samety of truth values.

### Notation 5.3.13

$|\Gamma|$ is the set of all well-formed formulas over $\Gamma$. (def. 5.3.5, p. 61)

Let $\Gamma$ be a second-order language, $\Omega$ a structure for $\Gamma$, $\varphi, \psi \in |\Gamma|$, and

let $\tau$ be a full variable assignment for $\varphi$ in $\Omega$ and $\kappa$ a full variable assignment for $\psi$ in $\Omega$. We write

$$\varphi[\tau] \leftrightarrow \psi[\kappa]$$

to mean

$$(\varphi[\tau] \leftrightarrow \top \wedge \psi[\kappa] \leftrightarrow \top) \vee (\varphi[\tau] \leftrightarrow \bot \wedge \psi[\kappa] \leftrightarrow \bot)\,.$$

This notation differs subtly from the final sentence of definition 5.3.11: there, the variable assignment was the same on both sides.

This finishes all that we need to know in order to define a second-order language and fill it with meaning.

## 3.4. A Note on Typed Universes

Seeing as how our interest lies not in the sublime complexity of logic in all its generality, but in graph theory, the inhabitants of our universe will necessarily be vertices and edges.[9]

This, however, presents a problem. Consider a 2-place predicate $\Lambda$ which checks whether its parameters are adjacent vertices, that is, $\Lambda(\delta_0, \delta_1)$ returns $\top$ if $\delta_0, \delta_1$ are vertices and are adjacent and $\bot$ if $\delta_0, \delta_1$ are not adjacent or if at least one of them is an edge.

Suppose now that we want to express the property "$G$ is a complete graph". This is equivalent to all vertices being neighbours, so the sentence

$$\forall \delta_0 (\forall \delta_1 (\Lambda(\delta_0, \delta_1)))$$

should satisfy our needs.

Except that it does not, because the universal quantifier ranges over the entire universe, including edges, and we know that for an edge $e$, the

[9] Remember that while we want to know whether a given *graph* fulfils a property, the graph is not the logical object of our language. The formulas we consider will always be something like "all vertices have at most two neighbours", which is a statement about the constituent particles of the graph.

statement $\Lambda(e, \delta_1)$ is always false.

One approach, used for example in Courcelle's original paper [Cou90], is to introduce a notion of "type" into our universe. Every element of the universe $|\Omega|$ is assigned exactly one type (in our case, either "vertex" or "edge"), and rather than introducing predicate variables ${}_n\lambda_0, {}_n\lambda_1, \ldots$ and function variables ${}_n\delta_0, {}_n\delta_1, \ldots$, we introduce the *typed* predicate and function variables ${}^V_n\lambda_0, {}^V_n\lambda_1, \ldots, {}^V_n\delta_0, {}^V_n\delta_1, \ldots$ for vertices and ${}^E_n\lambda_0, {}^E_n\lambda_1, \ldots, {}^E_n\delta_0, {}^E_n\delta_1, \ldots$ for edges and prescribe that the vertex variables can only range over vertices (respectively predicates or functions that take only vertices as their arguments) and the edge variables can only range over edges (respectively predicates or functions that take only edges as their arguments).

This may at first seem like a radical change in our approach, but it is in fact equivalent to non-sorted second-order logic with one small adjustment. Consider for instance the typed sentence

$$\forall {}^V_0\delta_0 (\forall {}^V_0\delta_1 (\Lambda({}^V_0\delta_0, {}^V_0\delta_1)))$$

(where $\Lambda$ is as before), which now correctly expresses our thought "the graph is complete".

Suppose we were to introduce a 1-place predicate $\Lambda_V$ that returns $\top$ if and only if its argument is a vertex. We could then achieve the same statement in an untyped second-order language by means of the sentence

$$\forall {}_0\delta_0 (\forall {}_0\delta_1 (\underbrace{(\Lambda_V({}_0\delta_0) \wedge \Lambda_V({}_0\delta_1))}_{\text{type check}} \to \Lambda({}_0\delta_0, {}_0\delta_1)))$$

or, without the shorthands,

$$\forall {}_0\delta_0 (\forall {}_0\delta_1 (\neg(\Lambda_V({}_0\delta_0) \wedge \Lambda_V({}_0\delta_1) \wedge \neg\Lambda({}_0\delta_0, {}_0\delta_1))))$$

(or indeed with more intuitive shorthands, $\forall \delta_0, \delta_1 : \delta_0, \delta_1 \in V \to \Lambda(\delta_0, \delta_1)$).

By analogous means, any second-order language with finitely many types can be expressed in an untyped second-order language by

- introducing finitely many 1-place predicates and
- extending formula length by a constant factor.

A more rigorous treatment of the topic may be found in [End72].

Since the additional complexity comes in the form of a constant factor (which, additionally, is relatively low since we only have two types) and will thus not increase the runtimes discussed in chapter 9, we shall ignore the distinction for the sake of ease of notation. We continue to work in the framework of untyped second-order logic as introduced in this chapter, but whenever it is convenient, we allow ourselves to make our formulas easier to read by omitting the type check and instead telling the reader which type of variable is expected.

## 3.5. Different Logical Frameworks

Structures as we have introduced them are quite powerful, but also difficult to handle in their generality. It will at times prove useful to restrict the type of universe we want to consider. For example, we may want to allow only universes that contain a certain predicate, or indeed ones that do *not* contain a certain construction that would hamper our progress.

**Definition 5.3.14**

A *logical framework* is a pair $\Psi = (\Gamma, \Xi)$, where $\Gamma$ is a language and $\Xi$ is a collection of structures for $\Gamma$, called the *multiverse* of $\Psi$. The framework $\Psi$ is called *first-order* if $\Gamma$ is a first-order language and every $\Omega \in \Xi$ is a first-order structure. It is called *second-order* if $\Gamma$ is a second-order language and every $\Omega \in \Xi$ is a second-order structure.

## 3.6. Monadic Second-Order Logic

In the previous sections, we have seen that second-order logic is more expressive than first-order logic. Regrettably, this additional expressiveness is too powerful even for Courcelle's Theorem. In order to rescue as strong a statement as possible, we restrict to a subset of second-order logic which is still a bit stronger than first-order logic.

*Monadic* second-order logic is the fragment of second-order logic in which quantification is allowed only over 0-place function variables (which already

subsumes first-order logic) and over 1-place predicate variables (hence the name "monadic").

Because a monadic predicate is equivalent to a set[10], some authors ([CE12]) like to call this "quantification over sets" rather than over 1-place predicate variable symbols.

### Definition 5.3.15

A second-order logical framework $(\Gamma, \Xi)$ is called *monadic* if every structure $\Omega \in \Xi$ of its multiverse satisfies the following conditions.

- For $n \in \mathbb{N}_{>0}$, the $n$-place function universe is empty: ${}^{\delta}_{n}|\Omega| = \varnothing$.
- For $n \in \mathbb{N}_{>1}$, the $n$-place predicate universe is empty: ${}^{\lambda}_{n}|\Omega| = \varnothing$.

${}^{\delta}_{n}|\Omega|$ denotes the $n$-place function universe. (def. 5.3.7, p. 63)

${}^{\lambda}_{n}|\Omega|$ denotes the $n$-place predicate universe. (def. 5.3.7, p. 63)

Since this renders the higher-arity universal quantifiers moot (if there are no 3-place functions, the sentence ${}^{\delta}_{3}\forall\, {}_{3}\delta_0 ({}_{3}\delta_0(\delta_0, \delta_1, \delta_2))$ is vacuously true), we allow ourselves to treat the language as if these symbols did not exist, thus cutting down on notation overhead.

Courcelle's Theorem actually works on a slight extension of monadic second-order logic called *counting monadic second-order logic* where one can ascertain whether a given set's cardinality is equal to $n \in \mathbb{N}$ modulo some $m \in \mathbb{N}_{>0}$. Expressibility-wise, counting monadic second-order logic is stronger than monadic second-order logic and weaker than second-order logic, but it introduces much additional notation. Since none of our examples use it, we have opted to restrict ourselves to monadic second-order logic for the main part of this thesis. The definition of counting monadic second-order logic, as well as formal proofs that all of our constructions also work for this larger framework, are instead included separately in appendix A. We advise the reader to completely understand the non-counting proof of Courcelle's Theorem before attempting to move on to said appendix.

[10] Every element of our universe either satisfies the predicate or not; equivalently, every element of our universe is either in the set (of elements which satisfy the predicate) or not.

# 4. The Monadic Second-Order Language of Graphs

In order to prove the result that we can decide whether a given graph fulfils a logical formula, we must first define what "fulfilling" a logical formula even means.

Consider once again our toy formula to check whether a given 2-uniform graph is complete:

$$\forall_0\delta_0(\forall_0\delta_1((\Lambda_V({}_0\delta_0) \wedge \Lambda_V({}_0\delta_1)) \Rightarrow \Lambda({}_0\delta_0, {}_0\delta_1))).$$

Remember that $\Lambda_V$ checks whether its argument is a vertex, while $\Lambda$ checks whether its arguments are neighbours. This means that given a graph $G$ which we want to check, the "inputs" of our formula are the vertices and edges of $G$. What are the inputs of a second-order formula? They are the elements of some universe $|\Omega|$! Consequently, we must model our logic on graphs in a way that ensures that the vertices and edges of every graph are somehow contained in the universe of some structure.

One approach one might consider would be to use as our universe the graph language $\mathfrak{F}$, such that all vertex and edge symbols are contained in it. One would then have to define a structure $\Omega$ assigning to the predicate symbol $\Lambda_V$ the set of all vertex symbols and to the predicate symbol $\Lambda$ … what, exactly? The set $\Omega(\Lambda)$ must contain all pairs of vertices $(v, v')$ such that $v$ and $v'$ are adjacent. But which pairs of vertex symbols are adjacent depends on $G$ itself! Do we now introduce for every graph $G$ its own predicate symbol $\Lambda_G$ meaning "$v$ and $v'$ are adjacent *in* $G$"? But then we would have to adapt our formula to each graph by changing the index $G$.

$\mathfrak{F}$ denotes the language of finite graphs. (def. 4.1.11, p. 27)

This motivates a different approach: rather than look for one "big" structure $\Omega$ and then somehow specify on which graph we are, we let each graph $G$ induce its own "small" universe containing just the vertices and edges of $G$. This way, we can fix the above formula to mean "our graph is complete" and then evaluate it on any given graph $G$ by checking whether the structure induced by $G$ is a *model* for this formula.

We begin by specifying the valid symbols of our language. The reader may want to cross-reference notation with definition 5.3.1.

## Definition 5.4.1

We denote by $\mathfrak{L}$ the second-order language with

- ${}_1\Lambda_{\mathfrak{L}} := \{ \lambda^0_{\text{conn}}, \lambda_{\text{vert}}, \lambda_{\text{edge}} \}$,
- ${}_2\Lambda_{\mathfrak{L}} := \{ \lambda^1_{\text{conn}}, \equiv \}$,
- $\forall n \in \mathbb{N}_{>2} \colon {}_n\Lambda_{\mathfrak{L}} := \{\lambda^{n-1}_{\text{conn}}\}$,
- $\forall n \in \mathbb{N} \colon {}_n\Delta_{\mathfrak{L}} := \varnothing$

and call this the *direct monadic second-order language of graphs.*

${}_1\Lambda_{\mathfrak{L}}$ is the set of 1-place predicate symbols of $\mathfrak{L}$.

${}_n\Delta_{\mathfrak{L}}$ is the set of $n$-place function symbols of $\mathfrak{L}$.

These symbols will shortly be assigned the following meanings:

The 2-place predicate $\equiv$ checks vertices or edges for equality, while the $(n+1)$-place predicate $\lambda^n_{\text{conn}}$ takes an edge $e$ and $n$ vertices, say $v_1, \dots, v_n$, and checks whether the vertex sequence of $e$ is precisely $v_1 \dots v_n$. The 0-place function variables will represent vertices and edges, and 1-place predicate variables will allow quantification over sets of vertices or edges. Lastly, the 1-place predicates $\lambda_{\text{vert}}$ and $\lambda_{\text{edge}}$ check whether their argument is a vertex or an edge, respectively.[11]

**The Direct Monadic Second-Order Language of 2-Uniform Graphs**

For 2-uniform graphs and using a typed universe as discussed in section 5.3.4, the number of symbols in definition 5.4.1 shrinks considerably:

- ${}_1\Lambda_{\mathfrak{L}} := \varnothing$,
- ${}_2\Lambda_{\mathfrak{L}} := \{ \equiv \}$,
- ${}_3\Lambda_{\mathfrak{L}} := \{ \lambda^2_{\text{conn}} \}$,
- $\forall n \in \mathbb{N}_{>3} \colon {}_n\Lambda_{\mathfrak{L}} := \varnothing$,
- $\forall n \in \mathbb{N} \colon {}_n\Delta_{\mathfrak{L}} := \varnothing$.

The reason that this is called the “direct” language of graphs is that it is an immediate formalisation of our intuitive requirements. We shall later see a second logical language on graphs which will be less “direct”.

In order to assign a truth value to a formula “in” a graph, we view the graph as a structure for $\mathfrak{L}$. The reader is welcome to cross-reference with definition 5.3.7.

[11] As discussed in section 5.3.4, the reader may pretend that these predicates do not exist and that instead our variable symbols are typed.

### Definition 5.4.2

Let $G = (V, E, \langle\!|_|\!\rangle)$ be a graph. The *induced direct second-order structure* of $G$, denoted $\Omega(G)$, is the following second-order structure on $\mathfrak{L}$.

- The universe of $\Omega(G)$ is $|\Omega(G)| := {}^{\delta}_{0}|\Omega(G)| := V \cup E$.
- The 1-place predicate universe of $\Omega(G)$ is ${}^{\lambda}_{1}|\Omega(G)| := \mathbb{2}^V \cup \mathbb{2}^E$.
- For $n \in \mathbb{N}_{>1}$, ${}^{\lambda}_{n}|\Omega(G)| := \varnothing$.
- For $n \in \mathbb{N}_{>0}$, ${}^{\delta}_{n}|\Omega(G)| := \varnothing$.
- $\Omega(G)(\lambda_{\text{vert}}) := V \subseteq {}^{\delta}_{0}|\Omega(G)|$.
- $\Omega(G)(\lambda_{\text{edge}}) := E \subseteq {}^{\delta}_{0}|\Omega(G)|$.
- $\Omega(G)(\equiv) := \{ (x, x) : x \in {}^{\delta}_{0}|\Omega(G)| \} \subseteq \left({}^{\delta}_{0}|\Omega(G)|\right)^2$.
- For every $n \in \mathbb{N}$,

$$\begin{aligned}\Omega(G)(\lambda^{n}_{\text{conn}}) &:= \{ (e, v_1, \dots, v_n) : e \in E, \langle\!|e|\!\rangle = v_1 \dots v_n \} \\ &\subseteq \left({}^{\delta}_{0}|\Omega(G)|\right)^{n+1}.\end{aligned}$$

$\mathbb{2}^V$ denotes the power set of $V$.

This definition finally enables us to give meaning to the intuitive notion of a "graph property": a graph property is any sentence $\varphi$ in $\mathfrak{L}$, and a graph $G$ fulfils this property if $\vDash_{\Omega(G)} \varphi$.

The reader should take a minute to convince themselves that the definitions above adequately represent the intuitive meaning of the language.

We gather up all of these structures to form the logical framework for our proof.

### Definition 5.4.3

We call

$$\mathfrak{M} := \bigcup_{n \in \mathbb{N}} \{ \Omega(G) : G \in \mathfrak{G}_n \}$$

the *direct multiverse of finite graphs* and

$$\mathfrak{X} := (\mathfrak{L}, \mathfrak{M})$$

the *direct logical framework of finite graphs.*

$\mathfrak{G}_n$ denotes the set of all graphs of type $n$. (def. 4.5.1, p. 39)

Note that, in contrast to "most" multiverses one might encounter, the multiverse of graphs is actually small enough to be a set.

Note also that $\mathfrak{X}$, by definition, is a second-order logical framework, and due to the way we have defined induced structures, it is monadic.

# 5. Staging an Example

We use this opportunity to take a break from all the definitions and come up for breath.[12] Now that we know, in theory, how to encode graph-theoretic properties in logical formulas, we should have a look at some practical examples. The reader who cares not about such applied mathematics may of course skip to the next section for even more definitions and, eventually, even a theorem.

## 5.1. A Question of Colour

One question that is famously hard is whether or not a (2-uniform, loop-free) graph is vertex-colourable with at most $k$ colours for some $k \in \mathbb{N}$, $k > 2$.

To turn this into a decision problem, we fix $k$. Since 1-colourability is not too interesting, let us first ask "is our graph 2-colourable?"

In order to turn this question into a monadic second-order formula, we formulate it in terms of sets: a 2-colouring partitions the graph's vertices into sets $V, V'$ such that no two vertices of $V$ are adjacent, nor are any two vertices of $V'$.

[12] Remember the nautical metaphor we had at one time?

In mathematical terms:

$$\begin{aligned}
\varphi := \exists\lambda_0\exists\lambda_1 :& \\
&\forall\delta_0 : \lambda_0(\delta_0) \rightarrow \lambda_{\text{vert}}(\delta_0) \\
&\wedge\forall\delta_0 : \lambda_1(\delta_0) \rightarrow \lambda_{\text{vert}}(\delta_0) \\
&\wedge\neg(\exists\delta_0 : \lambda_0(\delta_0) \wedge \lambda_1(\delta_0)) \\
&\wedge\forall\delta_0 : \lambda_{\text{vert}}(\delta_0) \rightarrow (\lambda_0(\delta_0) \vee \lambda_1(\delta_0)) \\
&\wedge\neg(\exists\delta_0\exists\delta_1\exists\delta_2 : \lambda_0(\delta_0) \wedge \lambda_0(\delta_1) \wedge \lambda^2_{\text{conn}}(\delta_2, \delta_0, \delta_1)) \\
&\wedge\neg(\exists\delta_0\exists\delta_1\exists\delta_2 : \lambda_1(\delta_0) \wedge \lambda_1(\delta_1) \wedge \lambda^2_{\text{conn}}(\delta_2, \delta_0, \delta_1)).
\end{aligned}$$

$\lambda_{\text{vert}}$ checks whether its argument is a vertex.

Once the pattern is detected, it is easily expanded to detect $k$-colourability for any constant $k \in \mathbb{N}$:

$$\begin{aligned}
\varphi := \exists\lambda_1 \dots \exists\lambda_k :& \\
&\forall\delta_0 : \lambda_1(\delta_0) \rightarrow \lambda_{\text{vert}}(\delta_0) \\
&\wedge \dots \\
&\wedge\forall\delta_0 : \lambda_k(\delta_0) \rightarrow \lambda_{\text{vert}}(\delta_0) \\
&\wedge\forall\delta_0 : \lambda_1(\delta_0) \rightarrow (\neg\lambda_2(\delta_0)\wedge \dots \wedge\neg\lambda_k(\delta_0)) \\
&\wedge \dots \\
&\wedge\forall\delta_0 : \lambda_k(\delta_0) \rightarrow (\neg\lambda_1(\delta_0)\wedge \dots \wedge\neg\lambda_{k-1}(\delta_0)) \\
&\wedge\forall\delta_0 : \lambda_{\text{vert}}(\delta_0) \rightarrow (\lambda_1(\delta_0)\vee \dots \vee\lambda_k(\delta_0)) \\
&\wedge\neg(\exists\delta_0\exists\delta_1\exists\delta_2 : \lambda_1(\delta_0) \wedge \lambda_1(\delta_1) \wedge \lambda^2_{\text{conn}}(\delta_2, \delta_0, \delta_1)) \\
&\wedge \dots \\
&\wedge\neg(\exists\delta_0\exists\delta_1\exists\delta_2 : \lambda_k(\delta_0) \wedge \lambda_k(\delta_1) \wedge \lambda^2_{\text{conn}}(\delta_2, \delta_0, \delta_1)).
\end{aligned}$$

One should note how here, as in other cases, the formula differs for different $k$ – since our language has no symbols for numbers, we cannot encode "$G$ is $k$-colourable" with both $G$ and $k$ as part of the input.

The example also illustrates the importance of second-order structures – in a first-order language, we could not quantify over the two sets of vertices in the first line.

# 6. The Other Language of Graphs

The logical framework we have just introduced is very nice for encoding graph properties "by hand", because it breaks down to simple questions of equality and adjacency. Regrettably, it is not well-suited for the proof technique we use to show Courcelle's Theorem. We introduce a *second* second-order logical framework of graphs that is uglier to use, but much nicer to work with in our proofs.[13] We then show that any formula expressible in our first framework has an equivalent in the second, meaning that we "lose" nothing by using the second framework – given a graph property expressible in the first framework, we can simply convert its formula into the equivalent formula in the second framework and work with that. The result we thus show will actually be slightly *stronger* than the one we need, as the framework we construct is strictly more expressive than $\mathfrak{X}$.

$\mathfrak{X}$ denotes the direct logical framework of graphs. (def. 5.4.3, p. 73)

In our new language, there are no vertices or edges. All variables ${}_0\delta_n$ will be *sets* of vertices or edges. For this reason, we can dispense with the 1-place predicates we have used in $\mathfrak{L}$ to check for sets.

$\mathfrak{L}$ denotes the second-order language of graphs. (def. 5.4.1, p. 72)

[13] The first-time reader shall erstwhile have to trust us that this second framework somehow has nicer properties than the first one. The repeat customer might remember that we need to prove inductiveness of our formulas, a proof which relies in crucial places on the second framework's ability to distinguish terminal vertices. The reader is invited to try their hand at a proof using only the first framework and, in case of success, send it to the author for inclusion as an addendum.

**Definition 5.6.1**

We denote by $\mathring{\mathfrak{L}}$ the second-order language with

- ${}_1\Lambda_{\mathring{\mathfrak{L}}} := \{\, \lambda^0_{\text{conn}}, \lambda_{\text{sgl}}, \lambda_{\text{vert}}, \lambda_{\text{edge}} \,\}$,
- ${}_2\Lambda_{\mathring{\mathfrak{L}}} := \{\, \lambda^1_{\text{conn}}, \sqsubseteq \,\}$,
- $\forall n \in \mathbb{N}_{>2}\colon\ {}_n\Lambda_{\mathring{\mathfrak{L}}} := \{\lambda^{n-1}_{\text{conn}}\}$,
- ${}_0\Delta_{\mathring{\mathfrak{L}}} := \{\varnothing\}$,
- ${}_1\Delta_{\mathring{\mathfrak{L}}} := \bigcup_{k \in \mathbb{N}_{>0}} \{\, \delta^K_{\text{term}} : K \in 2^{\{1,\dots,k\}} \,\}$,
- $\forall n \in \mathbb{N}_{>1}\colon\ {}_n\Delta_{\mathring{\mathfrak{L}}} := \varnothing$.

and call this the *circuitous monadic second-order language of graphs.*

${}_1\Lambda_{\mathring{\mathfrak{L}}}$ is the set of 1-place predicate symbols of $\mathring{\mathfrak{L}}$.

${}_0\Delta_{\mathring{\mathfrak{L}}}$ is the set of 0-place function symbols of $\mathring{\mathfrak{L}}$.

$2^{\{1,\dots,k\}}$ denotes the power set of $\{1,\dots,k\}$.

The predicate symbols will soon be assigned the following meanings.

- $x \sqsubseteq y$ is true if and only if $x$ is a subset of $y$ (keep in mind that all variables are sets).
- $\lambda_{\text{sgl}}(x)$ is true if and only if $|x| = 1$, that is, if $x$ is a "singleton" set.
- $\lambda_{\text{vert}}(x)$ is true if and only if all elements of $x$ are vertices, and $\lambda_{\text{edge}}(x)$ is true if and only if all elements of $x$ are edges.
- $\lambda^n_{\text{conn}}(E, V_1, \dots, V_n)$ is true if and only if there is at least one element $e \in E$ which is an edge and elements $v_1 \in V_1, \dots, v_n \in V_n$ which are vertices such that $\langle\!\langle e \rangle\!\rangle = v_1 \dots v_n$.

The function symbols will mean the following.

- $\varnothing$ yields the empty set.
- For a graph of type $n \in \mathbb{N}$ with terminals $t(1), \dots, t(n)$, a set $V$ containing only vertices, and a set $K \subseteq \{\, 1, \dots, n \,\}$, the term $\delta^K_{\text{term}}(V)$ denotes the set $V \cup \{\, t(i) : i \in K \,\}$.

**The Circuitous Monadic Second-Order Language of 2-Uniform Graphs**

For 2-uniform graphs and using a typed universe, the number of symbols in definition 5.6.1 shrinks again:

$_1\Lambda_{\mathring{\mathfrak{L}}}$ is the set of 1-place predicate symbols of $\mathring{\mathfrak{L}}$.

- $_1\Lambda_{\mathring{\mathfrak{L}}} := \{ \lambda_{\text{sgl}} \}$,
- $_2\Lambda_{\mathring{\mathfrak{L}}} := \{ \sqsubseteq \}$,
- $_3\Lambda_{\mathring{\mathfrak{L}}} := \{ \lambda^2_{\text{conn}} \}$,
- $\forall n \in \mathbb{N}_{>3}\colon\ _n\Lambda_{\mathring{\mathfrak{L}}} := \varnothing$,

$_0\Delta_{\mathring{\mathfrak{L}}}$ is the set of 0-place function symbols of $\mathring{\mathfrak{L}}$.

$2^{\{1,\dots,k\}}$ denotes the power set of $\{1,\dots,k\}$.

- $_0\Delta_{\mathring{\mathfrak{L}}} := \{ \varnothing \}$,
- $_1\Delta_{\mathring{\mathfrak{L}}} := \bigcup_{k\in\mathbb{N}_{>0}} \{ \delta^K_{\text{term}} : K \in 2^{\{1,\dots,k\}} \}$,
- $\forall n \in \mathbb{N}_{>1}\colon\ _n\Delta_{\mathring{\mathfrak{L}}} := \varnothing$.

Completely analogously to definition 5.4.2, we define an induced structure for this language.

### Definition 5.6.2

Let $G = (V, E, \langle_\rangle, t)$ be a graph of type $k$. The *induced circuitous second-order structure* of $G$, denoted $\mathring{\Omega}(G)$, is the following second-order structure on $\mathring{\mathfrak{L}}$.

$\mathring{\mathfrak{L}}$ is the circuitous language of graphs. (def. 5.6.1, p. 77)

- The universe of $\mathring{\Omega}(G)$ is $|\mathring{\Omega}(G)| := {}^{\delta}_{0}|\mathring{\Omega}(G)| := 2^V \cup 2^E$.
- For $n \in \mathbb{N}_{>0}$, ${}^{\lambda}_{n}|\mathring{\Omega}(G)| := \varnothing$.
- For $n \in \mathbb{N}_{>0}$, ${}^{\delta}_{n}|\mathring{\Omega}(G)| := \varnothing$.
- $\mathring{\Omega}(G)(\lambda_{\text{sgl}}) := \{ x : x \in |\mathring{\Omega}(G)|, |x| = 1 \} \subseteq |\mathring{\Omega}(G)|$.
- $\mathring{\Omega}(G)(\lambda_{\text{vert}}) := 2^V \subseteq |\mathring{\Omega}(G)|$.
- $\mathring{\Omega}(G)(\lambda_{\text{edge}}) := 2^E \subseteq |\mathring{\Omega}(G)|$.
- $\mathring{\Omega}(G)(\sqsubseteq) := \{ (x, y) \in |\mathring{\Omega}(G)|^2 : x \subseteq y \}$.

$\lambda_{\text{vert}}$ checks whether its argument is a set of vertices.

$\lambda_{\text{edge}}$ checks whether its argument is a set of edges.

- For every $n \in \mathbb{N}$,

$$\mathring{\Omega}(G)(\lambda^n_{\text{conn}}) := \{ (E', V_1, \dots, V_n) : \exists e \in E' : \exists v_1 \in V_1, \dots, \exists v_n \in V_n : e \in E, \langle\!\langle e \rangle\!\rangle = v_1 \dots v_n \} \subseteq |\mathring{\Omega}(G)|^{n+1}.$$

- $\mathring{\Omega}(G)(\varnothing) \colon |\mathring{\Omega}(G)|^0 \to |\mathring{\Omega}(G)|, () \mapsto \varnothing$.
- For $n \in \mathbb{N}_{>0}$, for $K \subseteq \{ 1, \dots, n \}$,

$$\mathring{\Omega}(G)(\delta^K_{\text{term}}) \colon |\mathring{\Omega}(G)| \to |\mathring{\Omega}(G)|,$$
$$x \mapsto \begin{cases} x \cup \{ t(i) : i \in K, i \leq k \} & \text{if } x \subseteq V \\ x & \text{otherwise.} \end{cases}$$

As always, the reader is at their leisure to ignore the predicates $\lambda_{\text{vert}}$ and $\lambda_{\text{edge}}$ and instead assume a typed universe as discussed in section 5.3.4.

$\lambda_{\text{vert}}$ checks whether its argument is a set of vertices.

$\lambda_{\text{edge}}$ checks whether its argument is a set of edges.

### Definition 5.6.3

We call

$$\mathring{\mathfrak{M}} := \bigcup_{n \in \mathbb{N}} \{ \mathring{\Omega}(G) : G \in \mathfrak{G}_n \}$$

the *circuitous multiverse of finite graphs* and

$$\mathring{\mathfrak{X}} := \left(\mathring{\mathfrak{L}}, \mathring{\mathfrak{M}}\right)$$

the *circuitous logical framework of finite graphs.*

$\mathfrak{G}_n$ denotes the set of all graphs of type $n$. (def. 4.5.1, p. 39)

We now show that our new language is at least as expressive as the old one. To this end, we adopt the following convention.

## Notation 5.6.4

$\mathring{\mathfrak{L}}$ is the circuitous language of graphs. (def. 5.6.1, p. 77)

We use as 0-place variable symbols for $\mathring{\mathfrak{L}}$ the set

$$\{\, \delta_0, \delta_1, \ldots \,\} \cup \{\, \lambda_0, \lambda_1, \ldots \,\}.$$

While this "doubles" the number of available variable symbols in some sense, it is clear that this is purely notational, since the number of symbols is still countable. We can safely reuse the symbols $\lambda_0, \ldots$ since $\mathring{\mathfrak{L}}$ has no predicate variables.

The only reason for this notational trick is that it makes the proof of theorem 5.6.6 easier to read.

Recall that we use $\mu_0$ to denote a symbol that could be either $\delta_0$ or $\lambda_0$.

## Definition 5.6.5

$\Omega(G)$ is the induced structure of $G$. (def. 5.4.2, p. 73)

Let $G$ be a graph, and let $\tau$ be a variable assignment in $\Omega(G)$ with domain $X = \{\, \mu_1, \ldots, \mu_n \,\}$ for some $n \in \mathbb{N}$.

The *induced circuitous assignment* of $\tau$ is the variable assignment $\mathring{\tau}$ in $\Omega(G)$ with domain $X$ and

$$\mathring{\tau} \colon \mu_i \mapsto \begin{cases} \{\tau(\mu_i)\} & \text{if } \mu_i \in {}_0\delta_{\mathfrak{L}} \\ \tau(\mu_i) & \text{if } \mu_i \in {}_1\lambda_{\mathfrak{L}} \end{cases}.$$

Note that, while $\tau(\mu_i)$ is either an edge or a vertex of $G$ or a 1-place predicate (also known as a set), the values of $\mathring{\tau}$ are all sets.

$\mathfrak{L}$ denotes the second-order language of graphs. (def. 5.4.1, p. 72)

We can now formulate what we mean by "at least as expressive".

## Theorem 5.6.6

$|\mathfrak{L}|$ is the set of all well-formed formulas over $\mathfrak{L}$. (def. 5.3.5, p. 61)

$\overleftarrow{\varphi}$ denotes the set of free variables of $\varphi$. (def. 5.3.6, p. 62)

Let $\varphi \in |\mathfrak{L}|$. Then there exists a formula $\mathring{\varphi} \in |\mathring{\mathfrak{L}}|$ with $\overleftarrow{\mathring{\varphi}} = \overleftarrow{\varphi}$ such that for every graph $G$ and for every full variable assignment $\tau$ in $\Omega(G)$ for $\varphi$, we have

$$\varphi[\tau] \leftrightarrow \top \Leftrightarrow \mathring{\varphi}[\mathring{\tau}] \leftrightarrow \top.$$

**Proof.** We proceed by induction on the structure of a well-formed formula in $\mathfrak{L}$.

Let $\varphi \in |\mathfrak{L}|$. We construct recursively a formula $\mathring{\varphi} \in |\mathring{\mathfrak{L}}|$ and show two things: that the two formulas have the same set of free variable symbols (albeit in different languages), and that the equivalence from the statement holds.

$\mathfrak{L}$ denotes the second-order language of graphs. (def. 5.4.1, p. 72)

$|\mathfrak{L}|$ is the set of all well-formed formulas over $\mathfrak{L}$. (def. 5.3.5, p. 61)

Note that the former property implies that if $\tau$ is a full variable assignment for $\varphi$, then $\mathring{\tau}$ is full for $\mathring{\varphi}$, which is necessary for the latter property to even make sense.

$\mathring{\mathfrak{L}}$ is the circuitous language of graphs. (def. 5.6.1, p. 77)

The reader may want to consult definitions 5.3.2 to 5.3.5 to remind themselves of the structure over which we run our induction.

Note first that since $\mathfrak{L}$ has no function symbols, all terms over the language are elementary (definition 5.3.2).

We run over the various cases from definition 5.3.5.

*Case 1: $\varphi$ is atomic.*

*Case 1.1: $\varphi = \lambda^n_{\text{conn}}(\delta_0, \delta_1, \dots, \delta_n)$ for some $n \in \mathbb{N}$.* This is an easy case to warm up our proof-writing fingers. We set $\mathring{\varphi} := \lambda^n_{\text{conn}}(\delta_0, \delta_1, \dots, \delta_n)$. This preserves the set of free variables, and if $\tau$ is a full variable assignment for $\varphi$ in $\Omega(G)$ for some graph $G$, we immediately see that

$\Omega(G)$ is the induced structure of $G$. (def. 5.4.2, p. 73)

$$\begin{aligned}
\varphi[\tau] \leftrightarrow \top &\Leftrightarrow (\tau(\delta_0), \tau(\delta_1), \dots, \tau(\delta_n)) \in \Omega(G)(\lambda^n_{\text{conn}}) \\
&\Leftrightarrow (\{\tau(\delta_0)\}, \{\tau(\delta_1)\}, \dots, \{\tau(\delta_n)\}) \in \mathring{\Omega}(G)(\lambda^n_{\text{conn}}) \\
&\Leftrightarrow (\mathring{\tau}(\delta_0), \mathring{\tau}(\delta_1), \dots, \mathring{\tau}(\delta_n)) \in \mathring{\Omega}(G)(\lambda^n_{\text{conn}}) \\
&\Leftrightarrow \mathring{\varphi}[\mathring{\tau}] \leftrightarrow \top,
\end{aligned}$$

$\mathring{\Omega}(G)$ is the circuitous induced structure of $G$. (def. 5.6.2, p. 78)

which is exactly what we wanted.

Since most cases are rather similar, we shall omit easy computations like the one above. We shall also not mention the sets of free variables again when it is obvious that they are the same.

*Case 1.2: $\varphi = \delta_0 \equiv \delta_1$.* We set $\mathring{\varphi} := \delta_0 \sqsubseteq \delta_1 \wedge \delta_1 \sqsubseteq \delta_0$. Equivalence is immediate.

*Case 1.3: $\varphi = \lambda_0(\delta_0)$.* Recall that all predicate variables are unary, where-

fore this is the only case of this type that we need to check. Recall also that while $\tau(\lambda_0)$ is a unary predicate and $\tau(\delta_0)$ is a vertex or an edge, $\mathring{\tau}(\lambda_0)$ and $\mathring{\tau}(\delta_0)$ are both sets.

We set $\mathring{\varphi} := \delta_0 \sqsubseteq \lambda_0$ and get

$$\begin{aligned}\varphi[\tau] \leftrightarrow \top &\Leftrightarrow \tau(\delta_0) \in \tau(\lambda_0)\\ &\Leftrightarrow \{\tau(\delta_0)\} \subseteq \tau(\lambda_0)\\ &\Leftrightarrow \mathring{\tau}(\delta_0) \subseteq \mathring{\tau}(\lambda_0)\\ &\Leftrightarrow \mathring{\varphi}[\mathring{\tau}] \leftrightarrow \top,\end{aligned}$$

as expected.

$\lambda_{\text{vert}}$ checks whether its argument is a vertex.

*Case 1.4:* $\varphi = \lambda_{\text{vert}}(\delta_0)$. As discussed before, the reader is welcome to ignore this predicate and instead assume a typed universe. Otherwise, they should set $\mathring{\varphi} := \lambda_{\text{vert}}(\delta_0)$.

$\lambda_{\text{vert}}$ checks whether its argument is a set of vertices.

*Case 1.5:* $\varphi = \lambda_{\text{edge}}(\delta_0)$. This is analogous to the previous case.

$\lambda_{\text{edge}}$ checks whether its argument is an edge.

This proves the claims for all atomic formulas. Any other formula we encounter must be a composite, where we assume that the claims hold for its constituent parts.

*Case 2:* $\varphi = \neg\psi$. We set, unsurprisingly, $\mathring{\varphi} := \neg\mathring{\psi}$ and get

$$\begin{aligned}\varphi[\tau] \leftrightarrow \top &\Leftrightarrow \psi[\tau] \leftrightarrow \bot\\ &\Leftrightarrow \mathring{\psi}[\mathring{\tau}] \leftrightarrow \bot\\ &\Leftrightarrow \mathring{\varphi}[\mathring{\tau}] \leftrightarrow \top\end{aligned}$$

by induction.

*Case 3:* $\varphi = \psi \wedge \zeta$. We set $\mathring{\varphi} := \mathring{\psi} \wedge \mathring{\zeta}$ and get equivalence by induction.

*Case 4:* $\varphi = \forall \mu_0 \psi$. We set

$$\mathring{\varphi} := \forall \mu_0 (\lambda_{\text{sgl}}(\mu_0) \rightarrow \mathring{\psi})$$

and get[14]

$$\begin{aligned}\varphi[\tau] \leftrightarrow \top &\Leftrightarrow (\forall \mu_0 \psi)[\tau] \leftrightarrow \top \\ &\Leftrightarrow \forall x \in |\Omega(G)| \colon \psi[\tau, \mu_0 \mapsto x] \leftrightarrow \top \\ &\Leftrightarrow \forall x \in |\Omega(G)| \colon \mathring{\psi}[\mathring{\tau}, \mu_0 \mapsto \{\, x \,\}] \leftrightarrow \top \\ &\Leftrightarrow \forall x \in |\mathring{\Omega}(G)| \colon |x| = 1 \Rightarrow \mathring{\psi}[\mathring{\tau}, \mu_0 \mapsto x] \leftrightarrow \top \\ &\Leftrightarrow \forall x \in |\mathring{\Omega}(G)| \colon (\lambda_{\mathrm{sgl}}(x) \to \mathring{\psi})[\mathring{\tau}, \mu_0 \mapsto x] \leftrightarrow \top \\ &\Leftrightarrow (\forall \mu_0 (\lambda_{\mathrm{sgl}}(\mu_0) \to \mathring{\psi}))[\mathring{\tau}] \leftrightarrow \top \\ &\Leftrightarrow \mathring{\varphi}[\mathring{\tau}] \leftrightarrow \top,\end{aligned}$$

$\Omega(G)$ is the induced structure of $G$. (def. 5.4.2, p. 73)

$\mathring{\Omega}(G)$ is the circuitous induced structure of $G$. (def. 5.6.2, p. 78)

finishing this case.

The two claims hence hold for any well-formed formula.

□

The reader might notice that the formula constructed in the proof, while longer than the original, is only larger by a constant factor, and the growth caused by translating into $\mathring{\mathfrak{L}}$ is thus linear. This is nice, but not actually needed – later on, the formula will not be part of our input, so its length is irrelevant when looking at the asymptotic running time of our algorithms.

$\mathring{\mathfrak{L}}$ is the circuitous language of graphs. (def. 5.6.1, p. 77)

We get the following nice corollary.

[14] Note that the obvious choice, $\forall \mu_0 \mathring{\psi}$, does not work: while the universal quantifier of $\Omega(G)$ ranges over the vertices and edges of $G$, the universal quantifier of $\mathring{\Omega}(G)$ ranges over *sets* of vertices and edges. Consider the formula $\varphi = \forall \mu_0 \forall \mu_1 \colon \neg(\mu_0 \equiv \mu_1) \to \neg(\exists \mu_2 \colon \lambda^2_{\mathrm{conn}}(\mu_2, \mu_0, \mu_1))$ expressing that the 2-uniform graph $G$ contains only loops, but no edges between different vertices. Consider further a graph $G$ with three vertice $v_1, v_2, v_3$ and an edge $e$ with $\langle e \rangle = v_1 v_1$. While $\varphi \leftrightarrow \top$ for this graph, the circuitous formula $\mathring{\varphi} = \forall \mu_0 \forall \mu_1 \colon \neg(\mu_0 \sqsubseteq \mu_1 \wedge \mu_1 \sqsubseteq \mu_0) \to \neg(\exists \mu_2 \colon \lambda^2_{\mathrm{conn}}(\mu_2, \mu_0, \mu_1))$ is false: the sets $X_1 := \{v_1, v_2\}$ and $X_2 := \{v_1, v_3\}$ are not equal, but the statement $\lambda^2_{\mathrm{conn}}(\{e\}, X_1, X_2)$ is nonetheless true because the circuitous predicate $\lambda^2_{\mathrm{conn}}$ only cares whether there is any vertex in $X_1$ connected to any vertex in $X_2$ by any edge in $\{e\}$.

## Corollary 5.6.7

$\mathfrak{L}$ denotes the second-order language of graphs. (def. 5.4.1, p. 72)

Let $\varphi \in \|\mathfrak{L}\|$. Then there exists a sentence $\mathring{\varphi} \in \|\mathring{\mathfrak{L}}\|$ such that for every graph $G$, we have

$$\vDash_{\Omega(G)} \varphi \Leftrightarrow \vDash_{\mathring{\Omega}(G)} \mathring{\varphi}\,.$$

$\|\mathfrak{L}\|$ denotes the set of sentences over $\mathfrak{L}$. (def. 5.3.6, p. 62)

The converse of theorem 5.6.6 is not true.

## Theorem 5.6.8

$\mathring{\mathfrak{L}}$ is the circuitous language of graphs. (def. 5.6.1, p. 77)

There exists a sentence $\mathring{\varphi} \in \|\mathring{\mathfrak{L}}\|$ for which there is no sentence $\varphi \in \|\mathfrak{L}\|$ which satisfies

$\Omega(G)$ is the induced structure of $G$. (def. 5.4.2, p. 73)

$$\vDash_{\Omega(G)} \varphi \Leftrightarrow \vDash_{\mathring{\Omega}(G)} \mathring{\varphi}$$

for every graph $G$.

$\mathring{\Omega}(G)$ is the circuitous induced structure of $G$. (def. 5.6.2, p. 78)

**Proof.** This is actually easy to see, since $\mathring{\mathfrak{L}}$ can detect terminal vertices, whereas $\mathfrak{L}$ cannot. Consider the two type 2 graphs $(G, t) := \mathfrak{e}_2$ and $(G, t') := \leftrightarrows_{1,2\mapsto 1} G$.

$\mathfrak{e}_2$ is the type 2 graph with 2 vertices and one edge. (def. 4.5.6, p. 43)

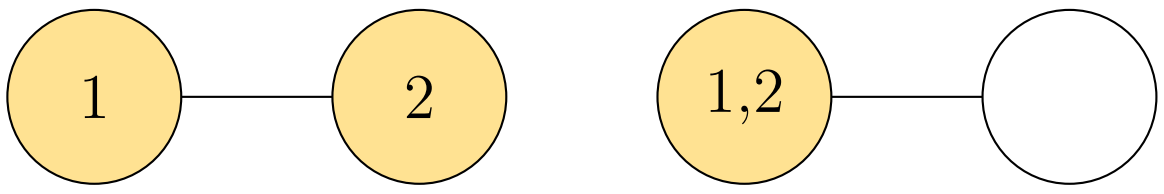


$\leftrightarrows$ denotes the terminal redefinition. (def. 4.5.4, p. 41)

These have the same underlying graph $G$. A quick look back at definition 5.4.2 reveals that $\Omega(G)$ is defined purely on this underlying graph, without accounting for terminals. Therefore, for any sentence $\varphi \in \|\mathfrak{L}\|$, we have

$$\vDash_{\Omega((G,t))} \varphi \Leftrightarrow \vDash_{\Omega((G,t'))} \varphi\,.$$

Consider now the circuitous sentence

$$\mathring{\varphi} := \exists \delta_0 \colon \lambda^2_{\text{conn}}(\delta_0, \delta^{\{1\}}_{\text{term}}(\varnothing), \delta^{\{2\}}_{\text{term}}(\varnothing)).$$

We have $\vDash_{\Omega((G,t))} \mathring{\varphi}$, but $\vDash_{\Omega((G,t))} \neg\mathring{\varphi}$, whence no sentence of $\mathfrak{L}$ can have the same truth values even on just the two graphs shown here.

□

The only problem are, of course, the terminal vertices, so the languages are indeed equivalent if evaluated only on type 0 graphs. Either way, this is not relevant to the proof of Courcelle's Theorem, since we shall always start with a sentence from $\mathfrak{L}$ which we then convert to an equivalent sentence of $\mathring{\mathfrak{L}}$ using theorem 5.6.6.

# Chapter 6

# Typed Algebras

In our journey to answering the question "which graphs fulfil a certain property $\varphi$", we have now reached the point where we are able to properly formulate the question.

In order to also find the answer (without seven and a half million years of computation), we need to construct a framework in which the question whether a given graph fulfils a certain property can be answered efficiently, preferably in time linear in the number of vertices and edges.

To this end, we show that certain classes of graphs form a so-called "typed algebra", and that any typed algebra with certain niceness properties admits a linear-time algorithm. The latter result is the topic of chapter 8, while the former is discussed right here.

## 1. Basic Definitions

In analogy to our treatment of formal languages in chapter 5, we begin by specifying a purely syntactical basis, called a *signature*. It provides a list of "sorts" for variables and a list of function symbols. Each function symbol is assigned a "type", that is, it knows which sorts of variables go into it and what sort of variable comes out. For languages, this was the arity of a function symbol. This time around, we need not only the *number* of variables that fit into a function symbol, but also the *sort* of variables it

takes.[1]

**Definition 6.1.1**

A *signature* is a triple $(\mathscr{T}, \mathscr{F}, \langle_\rangle)$ consisting of the following data:

- a set $\mathscr{T}$ of *sorts* or *types*,
- a set $\mathscr{F}$ of *function symbols*,
- and a *typing* $\langle_\rangle\colon \mathscr{F} \to \mathscr{T}^* \times \mathscr{T}$.

For any $f \in \mathscr{F}$, we call $\langle f\rangle$ the *type* of $f$.

For $\langle f\rangle = (t_1 \dots t_n, t)$, we call $\langle f\rangle^{\text{in}} := (t_1, \dots, t_n)$ the *in-type* of $f$, we call $\langle f\rangle^{\text{out}} := t$ the *out-type* of $f$, and we call $|f| := n$ the *arity* of $f$.

The reader familiar with functional programming has seen signatures and algebras before: one simple example for a signature would be the tuple $\mathscr{S} = (\mathscr{T}, \mathscr{F}, \langle_\rangle)$ with

$$\begin{aligned} \mathscr{T} &= \{ n, w \} \\ \mathscr{F} &= \{ +, -, \mathbf{Str} \} \\ \langle + \rangle &= (n\,n, n) \\ \langle - \rangle &= (n\,n, n) \\ \langle \mathbf{Str} \rangle &= (n, w). \end{aligned}$$

This, as with all of our syntactic constructs, means nothing on its own. However, one could fill it with meaning in the following way:

- Prescribe $n$ to be the set of natural numbers.
- Let $w$ be all number words of the English language.
- Identify $+$ with the addition function $_+_\colon \mathbb{N} \times \mathbb{N} \to \mathbb{N}$.
- Identify $-$ with the subtraction $\mathbb{N} \times \mathbb{N} \to \mathbb{N}$, capped at a minimum of 0, that is, $5 - 7 = 0$.
- And let $\mathbf{Str}\colon \mathbb{N} \to w$ take a natural number and output its English name, that is, $\mathbf{Str}(5) =$ “five” and so forth.

[1] Since standard literature has found no consensus, we shall be using the words “sort” and “type” interchangeably for the elements of $\mathscr{T}$. Everything else is always called a “type”, so if we say “sort” in this thesis, it is clear that we mean an element of $\mathscr{T}$.

Let us call this structure, for now, $\mathbb{N}_{\text{English}}$.

One can see now why $\mathscr{S}$ is called a "signature": it prescribes for each function symbol a number of input types and an output type – addition takes two numbers and outputs another number, while stringification takes only one number and outputs a string.

Of course, $\mathbb{N}_{\text{English}}$ is not the only way to assign semantics to our toy signature. If one dislikes throwing away negative numbers, one could easily use the following instead, which we shall call $\mathbb{Q}_{\text{English}}$:

- Prescribe $n$ to be the set of rational numbers.
- Let $w$ be strings of the form "$\omega$ divided by $\eta$", where $\omega$ and $\eta$ are English number words.
- Identify $+$ with the addition function $_+_ \colon \mathbb{Q} \times \mathbb{Q} \to \mathbb{Q}$.
- Identify $-$ with the subtraction $\mathbb{Q} \times \mathbb{Q} \to \mathbb{Q}$.
- And let $\mathbf{Str} \colon \mathbb{Q} \to w$ take a rational number and output English name of its reduced form, that is, $\mathbf{Str}(\frac{3}{6}) = \mathbf{Str}(\frac{1}{2}) =$ "one divided by two" and so forth.

Both of the above, as realisations of the signature $\mathscr{S}$, are valid $\mathscr{S}$-*algebras*[2]. They are to signatures what structures were to languages.

### Definition 6.1.2

Let $\mathscr{S} = (\mathscr{T}, \mathscr{F}, \langle _ \rangle)$ be a signature. An $\mathscr{S}$-*algebra* is a pair $(\mathscr{C}, \mathscr{O})$ with the following data.

- A family of disjoint *carrier sets* $\mathscr{C} = \{\mathscr{C}_t\}_{t \in \mathscr{T}}$.
- A family of *functions* or *operations* $\mathscr{O} = \{\mathscr{O}_f\}_{f \in \mathscr{F}}$ such that

$$\forall f \in \mathscr{F} : \langle f \rangle = (t_1 \dots t_n, t) \Rightarrow \mathscr{O}_f : \mathscr{C}_{t_1} \times \dots \times \mathscr{C}_{t_n} \to \mathscr{C}_t.$$

An algebra defines for every *sort* of variable a set of constants of that sort and for every function *symbol* an actual function that takes constants from

[2] Courcelle [Cou90] calls the construction a *typed magma*, but this is not how it is known in standard literature.

the sets indicated by the function symbol's type and outputs a constant of the appropriate sort.

We allow a shorthand similar to that which we use on graphs.

### Notation 6.1.3

Let $\mathscr{S} = (\mathscr{T}, \mathscr{F}, \langle_\rangle)$ be a signature, and let $\mathscr{A} = (\mathscr{C}, \mathscr{O})$ be an $\mathscr{S}$-algebra. We write $c \in \mathscr{A}$ to mean $\exists t \in \mathscr{T} : c \in \mathscr{C}_t$.

We also need morphisms between such structures. Consider again our toy signature above, and consider a mathematician from the planet Kronos, who might want the following algebra instead, called $\mathbb{N}_{\text{Klingon}}$:

- Prescribe $n$ to be the set of natural numbers.
- Let $w$ be all number words of the *Klingon* language.[3]
- Identify $+$ with the addition function $_ + _ : \mathbb{N} \times \mathbb{N} \to \mathbb{N}$.
- Identify $-$ with the subtraction $\mathbb{N} \times \mathbb{N} \to \mathbb{N}$, capped at a minimum of 0, that is, $5 - 7 = 0$.
- And let $\mathbf{Str} : \mathbb{N} \to w$ take a natural number and output its *Klingon* name, that is, $\mathbf{Str}(5) =$ "vagh" and so forth.

This is almost the same thing as the first algebra we defined, and we can easily imagine a "translation function" $\mathbb{N}_{\text{English}} \to \mathbb{N}_{\text{Klingon}}$ which leaves the numbers alone and maps strings to their translations, that is, "zero" $\mapsto$ "pagh", "one" $\mapsto$ "wa'", and so forth.

Of course, our translation also needs to map the functions: it leaves addition and subtraction alone (since they depend not on the language), and it maps English stringification to Klingon stringification.

Even if we translate also the number representations into a different system, we can expect a reasonable translation $h$ to be consistent in the sense that

$$\mathbf{Str}(h(1 + 1)) = \mathbf{Str}(h(1) + h(1)).$$

In mathematical lingo, our translation should commute with the functions involved.

[3] Klingon cardinal numbers start with "pagh", "wa'", "cha'", "wej", "loS", "vagh", "jav", "Soch", "chorgh", "Hut", "wa'maH". Pronunciations can be found in [Okr92].

**Definition 6.1.4**

Let $\mathscr{S} = (\mathscr{T}, \mathscr{F}, \langle _ \rangle)$ be a signature, and let $\mathscr{A} = (\mathscr{C}, \mathscr{O})$, $\mathscr{B} = (\mathscr{D}, \mathscr{Q})$ be $\mathscr{S}$-algebras. An *$\mathscr{S}$-algebra morphism* from $\mathscr{A}$ to $\mathscr{B}$ is a family of maps $\{ h_t : \mathscr{C}_t \to \mathscr{D}_t \}_{t \in \mathscr{T}}$ such that

$$\begin{aligned}
&\forall f \in \mathscr{F} : \langle f \rangle^{\text{in}} = (t_1, \dots, t_n) \Rightarrow \\
&\forall c_1 \in \mathscr{C}_{t_1} \dots \forall c_n \in \mathscr{C}_{t_n} : \\
&h_t(\mathscr{O}_f(c_1, \dots, c_n)) = \mathscr{Q}_f(h_{t_1}(c_1), \dots, h_{t_n}(c_n)).
\end{aligned}$$

Equivalently, one may require that all diagrams of the following form commute.

$$\begin{array}{ccc}
\mathscr{C}_{t_1} \times \dots \times \mathscr{C}_{t_n} & \xrightarrow{\mathscr{O}_f} & \mathscr{C}_t \\
{\scriptstyle h_{t_1} \times \dots \times h_{t_n}} \big\downarrow & \circlearrowleft & \big\downarrow {\scriptstyle h_t} \\
\mathscr{D}_{t_1} \times \dots \times \mathscr{D}_{t_n} & \xrightarrow[\mathscr{Q}_f]{} & \mathscr{D}_t
\end{array}$$

We introduce the first niceness property for later use: a signature or algebra is called "locally finite" if something is finite for every sort. For signatures, this "something" are the function symbols outputting that sort, while for algebras, we can require that there should only be finitely many elements of that sort. The reason we want these properties is that later, we often construct something on "all function symbols of a certain type" or similar, and we would like those constructions to stay finite.

**Definition 6.1.5**

A signature $(\mathscr{T}, \mathscr{F}, \langle _ \rangle)$ is called *locally finite* if for every $t \in \mathscr{T}$, the set

$$\{ f \in \mathscr{F} : \langle f \rangle^{\text{out}} = t \}$$

is finite.

**Definition 6.1.6**

Let $\mathscr{S} = (\mathscr{T}, \mathscr{F}, \langle _ \rangle)$ be a signature. An $\mathscr{S}$-algebra $(\mathscr{C}, \mathscr{O})$ is called

*locally finite* if for every $t \in \mathscr{T}$, the set $\mathscr{C}_t$ is finite.

# 2. The Algebra of Graphs

In section 4.5, we have seen some functions which map one or two typed graphs to another typed graph: the disjoint sum, the terminal redefinition, and the terminal fusion. We now interpret these functions as symbols of a signature and show that the set of all typed graphs is a typed algebra, where the set of sorts is simply $\mathbb{N}$ – the sort of a graph is its number of terminal vertices.

We first define the set of function symbols. As always, these are for now purely syntactical.

**Definition 6.2.1**

We set

$\oplus$ is the disjoint sum. (def. 4.5.3, p. 40)

- $\mathfrak{F}_{\oplus} := \{ {}^m_n\oplus : n, m \in \mathbb{N} \}$,
- ${}^m_n\mathfrak{F}_{\leftrightarrows} := \{ {}_n\leftrightarrows_\sigma : \sigma\colon \{1, \dots, n\} \to \{1, \dots, m\} \}$,
- ${}_n\mathfrak{F}_{\text{媲}} := \{ {}_n\text{媲}^b_a : a, b \in \{1, \dots, n\} \}$,
- $\mathfrak{F}_{\mathrm{triv}} := \{ \mathfrak{v}, \mathfrak{e}_n : n \in \mathbb{N}_{>0} \}$,
- and $\mathfrak{F} := \mathfrak{F}_{\oplus} \cup \left( \bigcup_{i \in \mathbb{N}} \bigcup_{j \in \mathbb{N}} {}^j_i\mathfrak{F}_{\leftrightarrows} \right) \cup \left( \bigcup_{i \in \mathbb{N}} {}_i\mathfrak{F}_{\text{媲}} \right) \cup \mathfrak{F}_{\mathrm{triv}}$.

$\mathfrak{v}$ is the type 1 graph with one vertex. (def. 4.5.6, p. 43)

$\mathfrak{e}_n$ is the type $n$ graph with $n$ vertices and one edge. (def. 4.5.6, p. 43)

We then define $\langle _ \rangle\colon \mathfrak{F} \to \mathbb{N}^* \times \mathbb{N}$,

$$\langle \mathfrak{f} \rangle := \begin{cases} (nm, n+m)\text{[4]} & \mathfrak{f} = {}^m_n\oplus \\ (m, n) & \mathfrak{f} \in {}^m_n\mathfrak{F}_{\leftrightarrows} \\ (n, n) & \mathfrak{f} \in {}_n\mathfrak{F}_{\text{媲}} \\ (\varepsilon, 1) & \mathfrak{f} = \mathfrak{v} \\ (\varepsilon, n) & \mathfrak{f} = \mathfrak{e}_n \end{cases}$$

and set $\mathfrak{S} := (\mathbb{N}, \mathfrak{F}, \langle _ \rangle)$.

[4] Here, $nm$ denotes the word $n$ followed by the word $m$, not a product.

With this definition, the triple $(\mathbb{N}, \mathfrak{F}, \langle _ \rangle)$ is a valid signature.

Note how we have added a left subscript to the function symbols to ensure that the sets of symbols for different input types are disjoint. We will of course silently drop these left subscripts whenever the context allows.

We now have to fill these symbols with meaning, and we do so in the canonical way.

**Definition 6.2.2**

We define an $(\mathbb{N}, \mathfrak{F}, \langle _ \rangle)$-algebra $\mathfrak{G} = (\mathfrak{C}, \mathfrak{O})$ by setting

$$\mathfrak{C}_n := \mathfrak{G}_n$$

and by letting $\mathfrak{O}$ assign to each function symbol the corresponding graph construction from definitions 4.5.3 to 4.5.5 respectively the corresponding graph from definition 4.5.6.

$\mathfrak{G}_n$ denotes the set of all graphs of type $n$. (def. 4.5.1, p. 39)

The functions $\mathfrak{v}$ and $\mathfrak{e}_n$ are nullary functions.

Of course, we chose $\langle _ \rangle$ in just such a way that this makes $\mathfrak{G}$ into a valid $(\mathbb{N}, \mathfrak{F}, \langle _ \rangle)$-algebra. The reader is invited to turn back to definition definition 4.5.3 and its colleagues to verify this.

$\mathfrak{v}$ is the type 1 graph with one vertex. (def. 4.5.6, p. 43)

$\mathfrak{e}_n$ is the type $n$ graph with $n$ vertices and one edge. (def. 4.5.6, p. 43)

$\mathfrak{G}$ is the algebra of graphs. (def. 6.2.2, p. 93)

# 3. Expressions

We have mentioned before that the trivial graphs $\mathfrak{v}$ and $\mathfrak{e}_n$, $n \in \mathbb{N}$ together with disjoint sum, redefinition, and fusion suffice to build all finite hypergraphs, and that $\mathfrak{v}$ and $\mathfrak{e}_2$ together with those functions can build all 2-uniform graphs.

In this section, we make this statement precise.

Consider the following 2-uniform graph of type 1, where the only terminal vertex is the middle one. We shall call it $G$.

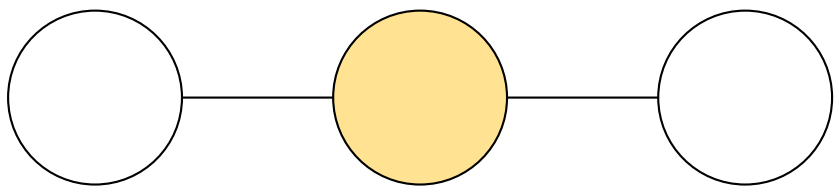

$\mathfrak{e}_2$ is the type 2 graph with 2 vertices and one edge. (def. 4.5.6, p. 43)

It is obvious how to build this graph using our basic building blocks: take two copies of $\mathfrak{e}_2$, the graph with two vertices and one edge, and glue them together.

媲 (tsureai, Japanese for *to marry*) denotes the source fusion. (def. 4.5.5, p. 42)

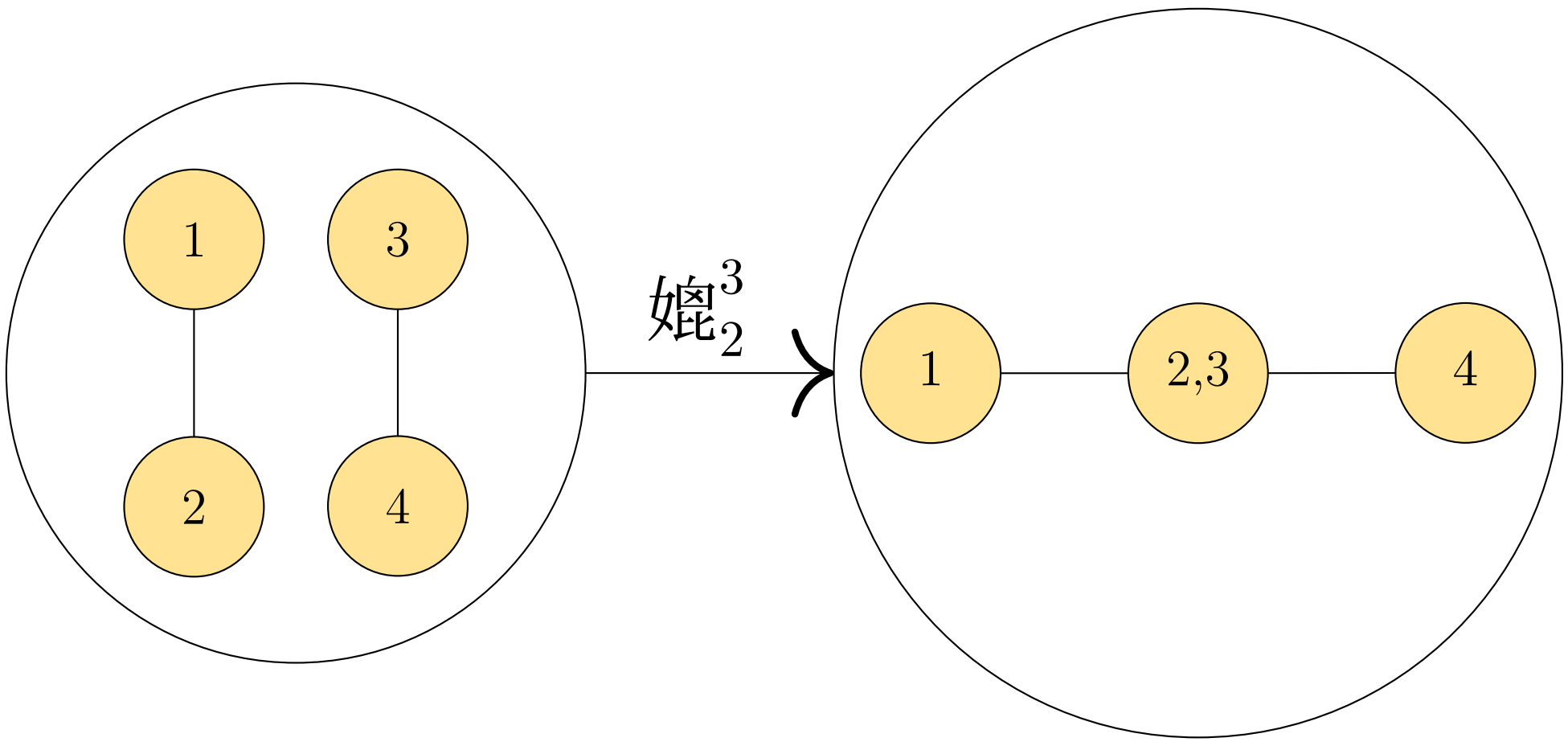


One can then discard the superfluous terminal labels using terminal redefinition.

$\leftrightarrows$ denotes the terminal redefinition. (def. 4.5.4, p. 41)

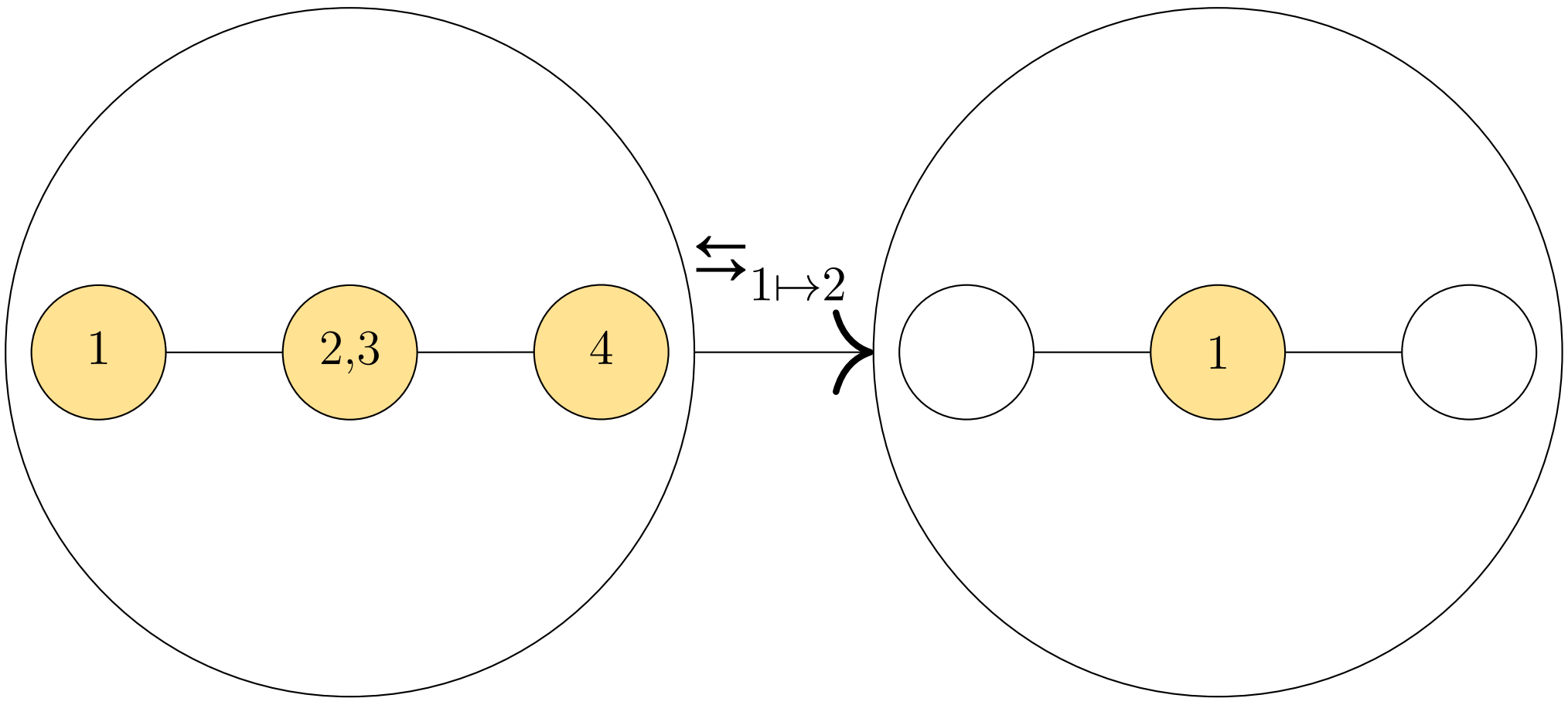


$\oplus$ is the disjoint sum. (def. 4.5.3, p. 40)

Mathematically, we could write

$$G = \leftrightarrows_{1\mapsto 2} \text{媲}^3_2 \left(\mathfrak{e}_2 \oplus \mathfrak{e}_2\right),$$

which is reasonably parseable, but this approach becomes unreadable quickly for larger graphs – the complete 2-uniform graph on four vertices is

$$\leftrightarrows_{\varnothing} \text{媲}_6^9 \text{媲}_4^{11} \text{媲}_2^{10} \text{媲}_1^{12} \left( \text{媲}_4^8 \text{媲}_1^5 \left( \text{媲}_2^3 (\mathfrak{e}_2 \oplus \mathfrak{e}_2) \oplus \text{媲}_6^7 (\mathfrak{e}_2 \oplus \mathfrak{e}_2) \right) \oplus \mathfrak{e}_2 \oplus \mathfrak{e}_2 \right).$$

$\mathfrak{e}_2$ is the type 2 graph with 2 vertices and one edge. (def. 4.5.6, p. 43)

$\oplus$ is the disjoint sum. (def. 4.5.3, p. 40)

媲 (tsureai, Japanese for *to marry*) denotes the source fusion. (def. 4.5.5, p. 42)

$\leftrightarrows$ denotes the terminal redefinition. (def. 4.5.4, p. 41)

We can instead visualise the building instructions pictorially:

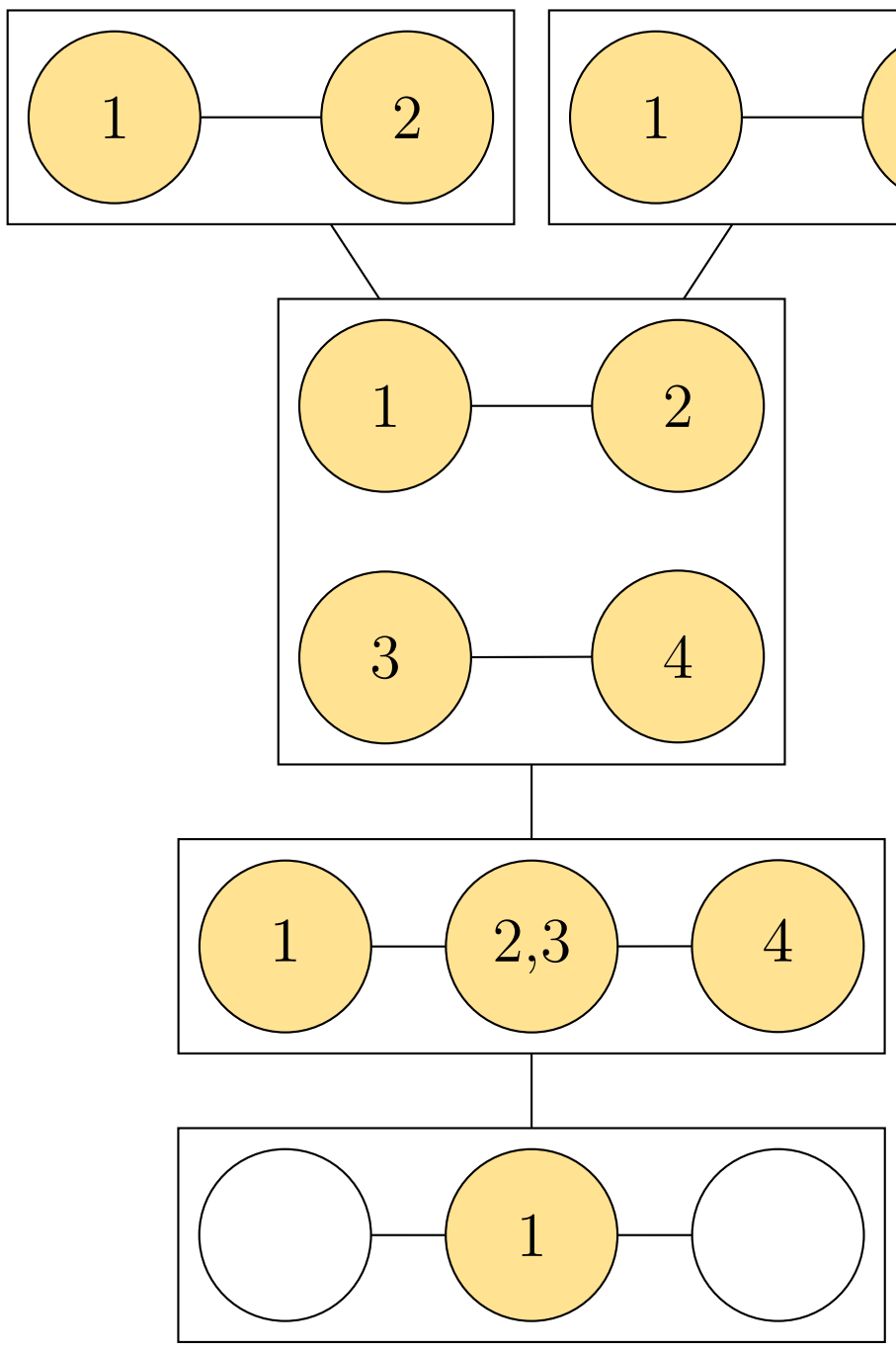


Since the graphs obtained in the intermediate steps are uniquely determined by the functions used, it is actually more informative to only give the function symbols:

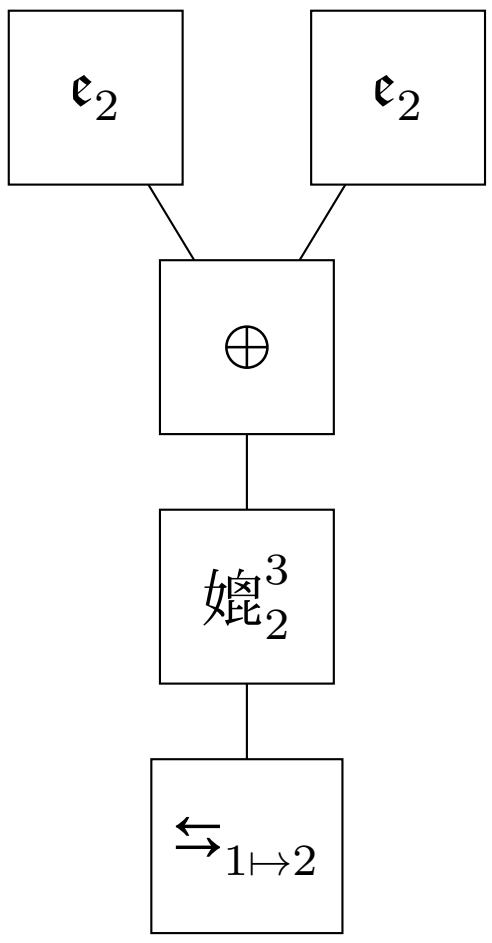


$\mathfrak{e}_2$ is the type 2 graph with 2 vertices and one edge. (def. 4.5.6, p. 43)

$\oplus$ is the disjoint sum. (def. 4.5.3, p. 40)

媲 (tsureai, Japanese for *to marry*) denotes the source fusion. (def. 4.5.5, p. 42)

$\leftrightarrows$ denotes the terminal redefinition. (def. 4.5.4, p. 41)

Flipping this picture upside down yields something that looks suspiciously akin to a rooted tree.

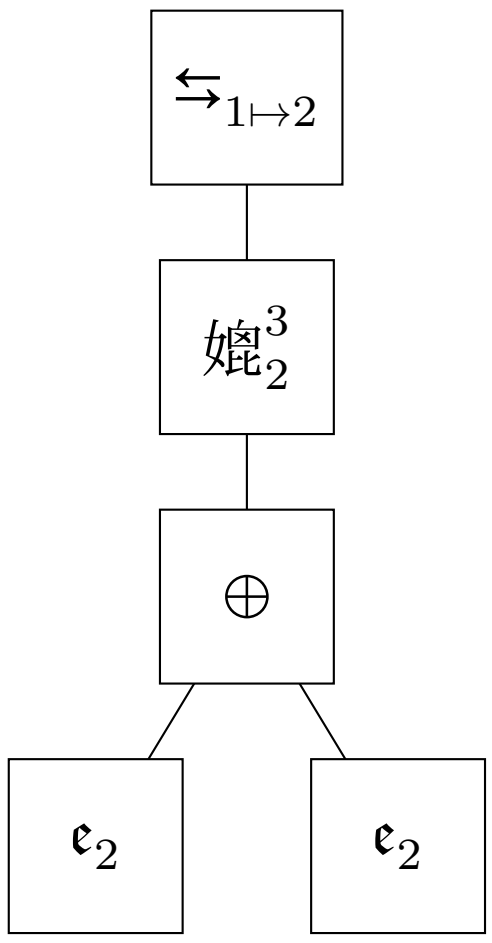


$\mathfrak{G}$ is the algebra of graphs. (def. 6.2.2, p. 93)

The vertices of this tree are labeled with function symbols in the algebra $\mathfrak{G}$, and the inputs of each function are the values of the labels of the successor vertices. The leaves are nullary functions yielding the graph $\mathfrak{e}_2$.

A tree like this is called an *expression* in $\mathfrak{G}$. Returning to our other favourite example, an expression could look something like this:

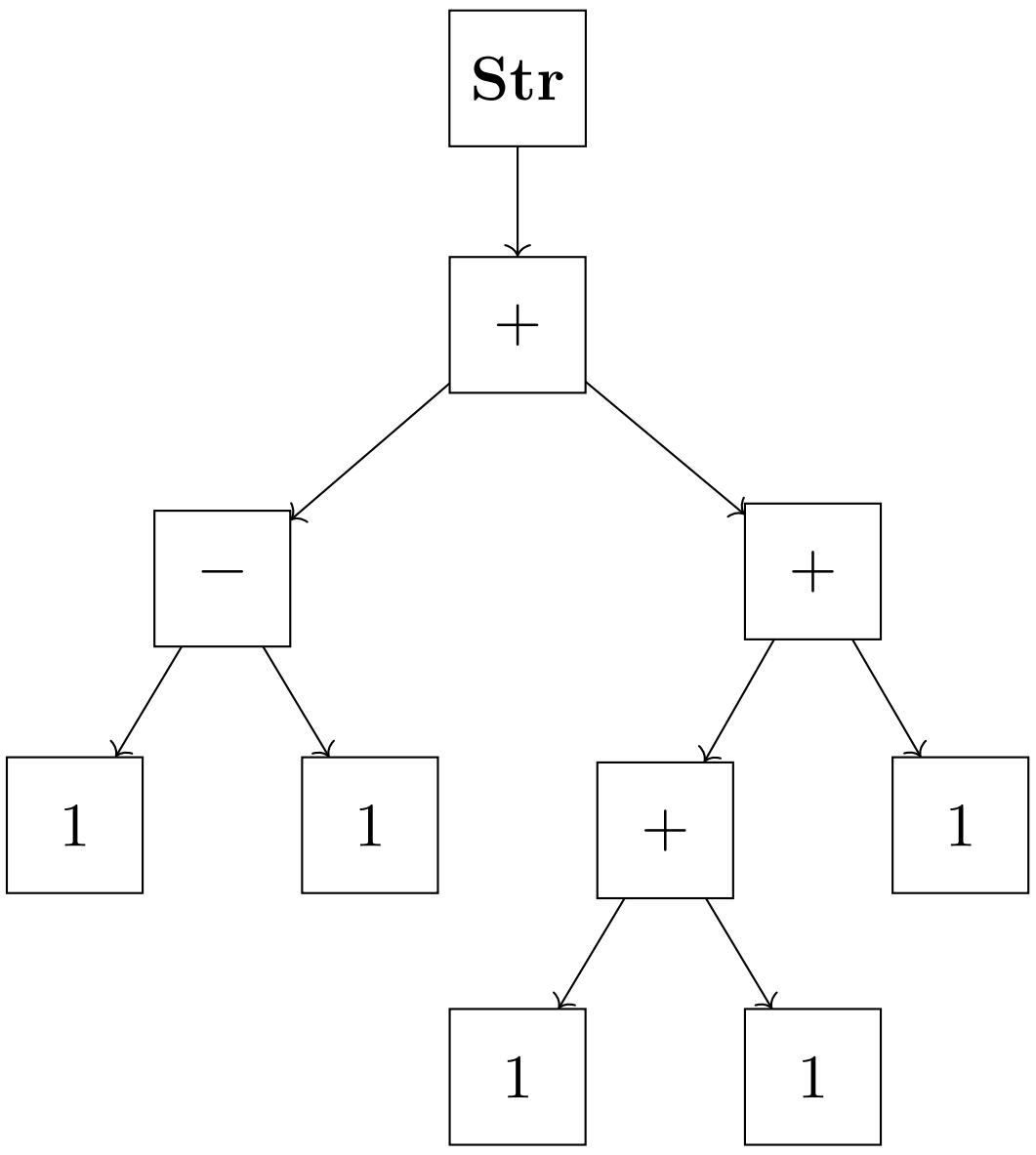


In $\mathbb{N}_{\text{English}}$, this should evaluate to the word “three”.

What if we want to plug in different numbers than just all ones? Do we change the entire expression?

Of course not. Rather, we cleverly cut the leaves off our tree to arrive at something like the following:

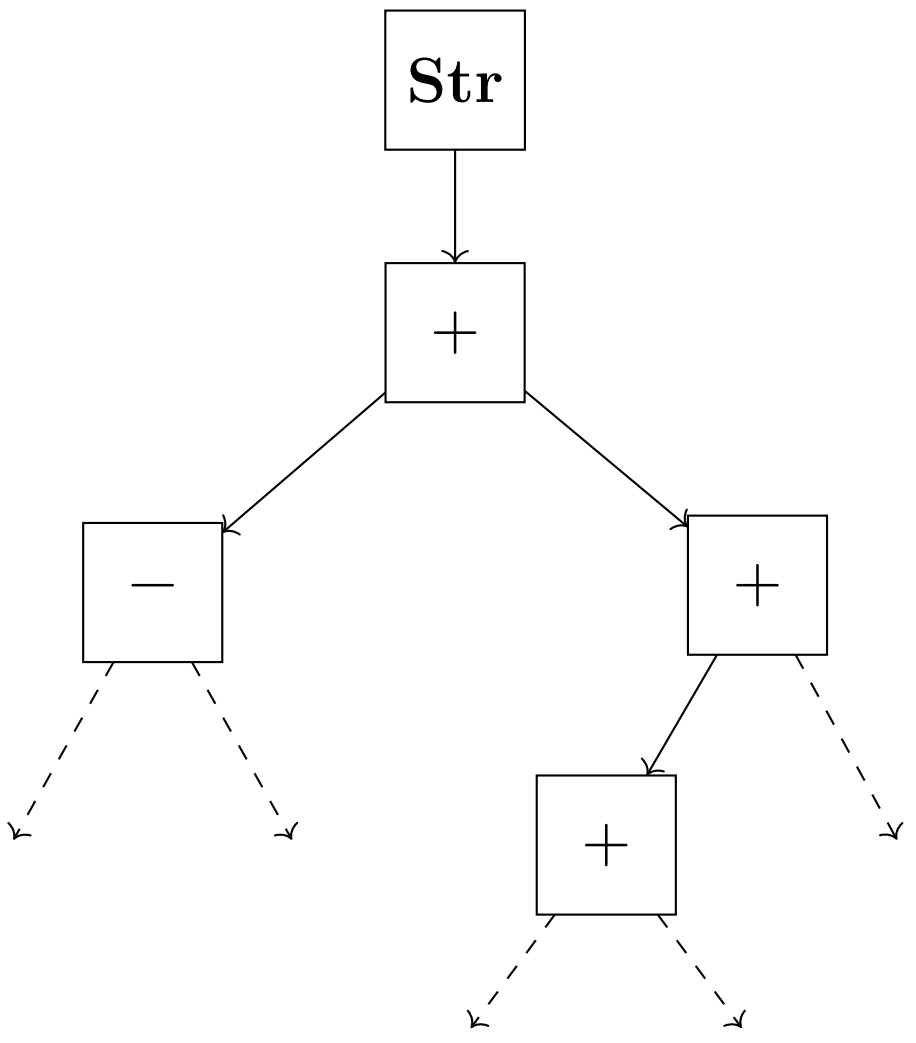

Here, the dashed arrows indicate places where we can “plug in” values from our algebra.

Any values? No: our algebra contains numbers (which are valid inputs for + and −) and English or Klingon words, which are not valid inputs.

Since the type of input can be checked at signature level (independent of, for example, the language of the number words), it makes sense to define expressions not for algebras, but already on signatures.

We collect the properties that we have discovered so far into a formal definition.

### Definition 6.3.1

Let $\mathscr{S} = (\mathscr{T}, \mathscr{F}, \langle_\rangle)$ be a signature. A *pre-expression* over $\mathscr{S}$ is a pair $e = (T, \langle_\rangle)$ fulfilling the following conditions.

- $T = (V, E, \lhd_\rhd, \star, \preccurlyeq)$ is a traversal tree.
- $\langle_\rangle \colon V \to \mathscr{F}$ is a function assigning to each vertex a function symbol.
- $\forall v \in V\colon \deg^{\text{out}} v = 0 \lor \deg^{\text{out}} v = |\langle v \rangle|$.
- For every $v \in V$ with $\langle\langle v \rangle\rangle^{\text{in}} = (t_1, \dots, t_n)$ and $\mathrm{N}^{\text{out}} v = [v_1, \dots, v_n]$ we have
$$\forall i \in \{1, \dots, n\} \colon \langle\langle v_i \rangle\rangle^{\text{out}} = t_i.$$

We call $\langle e \rangle^{\text{out}} := \langle\langle \sqrt{T} \rangle\rangle^{\text{out}}$ the *output sort* of $e$.

If the leaves of $T$, in traversal order, are $v_1, \dots, v_m$ for some $m \in \mathbb{N}$ with $\langle\langle v_i \rangle\rangle^{\text{in}} = \left(t^i_1, \dots, t^i_{k_i}\right)$, we call

$$\langle e \rangle^{\text{in}} := \left(t^1_1, t^1_2, \dots, t^1_{k_1}, t^2_1, \dots, t^2_{k_2}, \dots, t^n_1, \dots, t^n_{k_n}\right)$$

the *input sort* of $e$.

The height of a tree is defined in definition 4.4.4.

The *height* of $e$ is the height of the underlying tree.

We write

$$\langle e \rangle := \left(t^1_1 t^1_2 \dots t^1_{k_1} t^2_1 \dots t^2_{k_2} \dots t^n_1 \dots t^n_{k_n}, \langle e \rangle^{\text{out}}\right).$$

The set of all pre-expressions over $\mathscr{S}$ is denoted by $|\mathscr{S}|$.

Note that we require a *traversal* tree (which orders the children of a vertex, see definition 4.4.5) because some operations may not be commutative – it is relevant whether we compute $1-3$ or $3-1$.

In our example signature, the following should be valid pre-expressions:

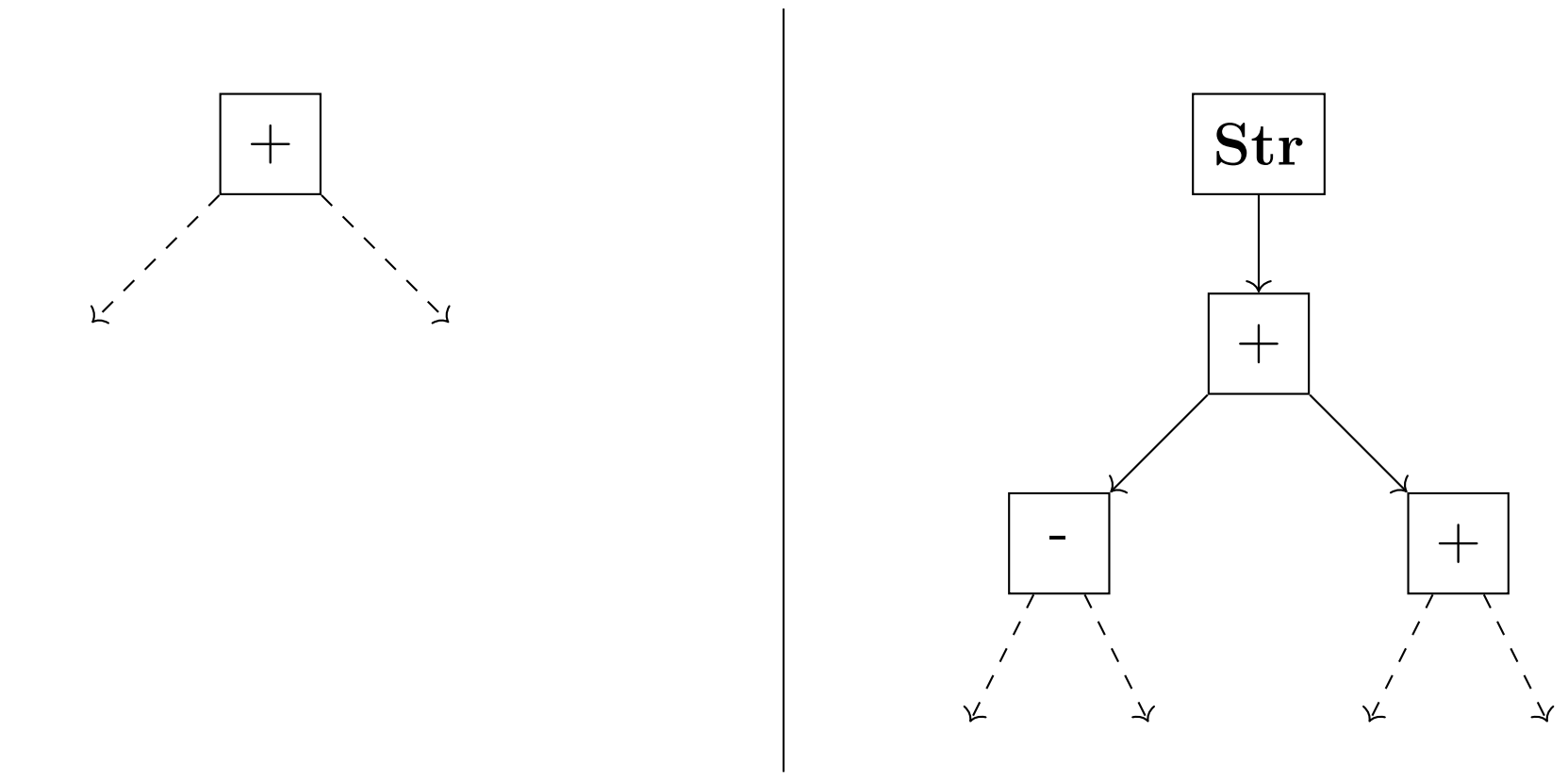


Again, we have indicated the missing inputs by dashed arrows. The input sort of the left-hand pre-expression is $n\,n$, and its output sort is $n$. The input sort of the right-hand pre-expression is $n\,n\,n\,n$, while its output sort is $w$.

Meanwhile, the following should *not* be valid pre-expressions – one violates the type constraint, while the other tries to cram two numbers into one string function.

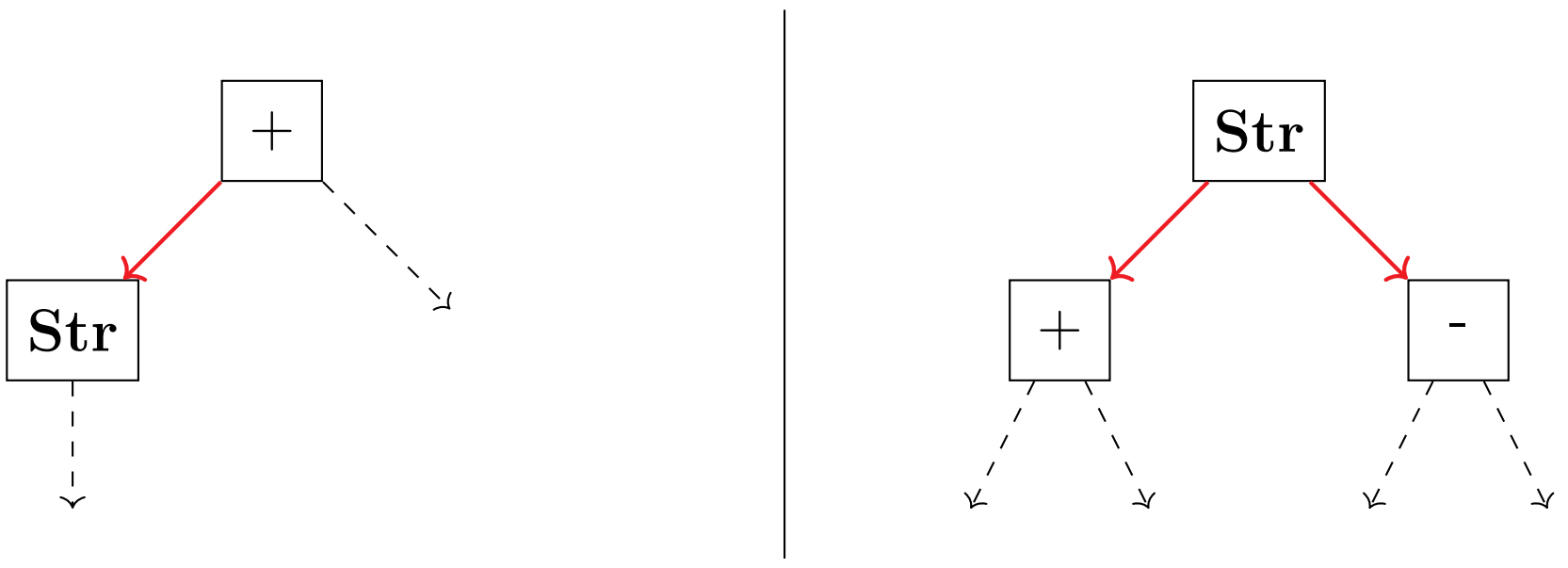

Note that the following is also not a pre-expression as by our definition, because the marked vertex expects an input, but is not a leaf.

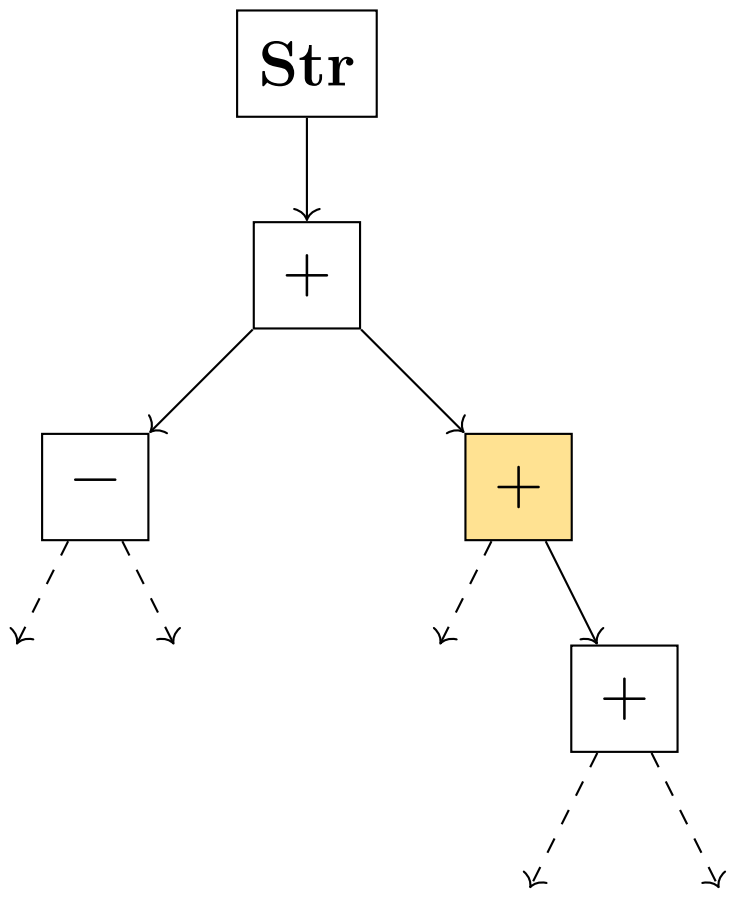


However, if our algebra contains an identity function id for every type, a functionally equivalent expression is easily constructed:

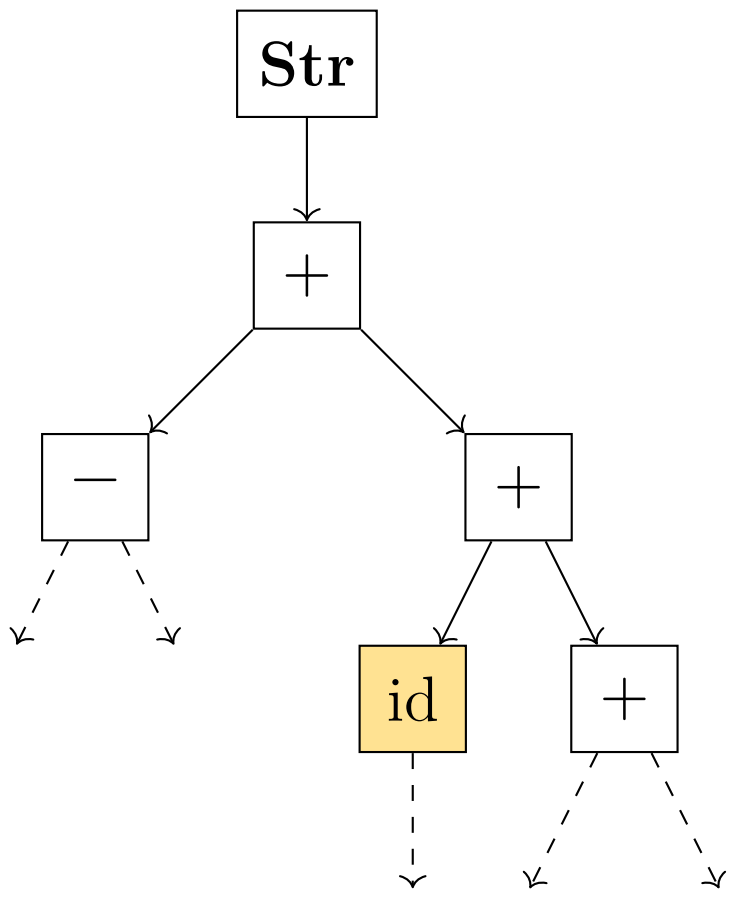


We therefore adopt the convention that we treat the two pictures as the same – the algebras with which we are dealing always contain the identity[5], which we omit from pictorial representations.

[5] The identity function for graphs of type $n \in \mathbb{N}$ is the same function as ${}_{n}\leftrightarrows_{i \mapsto i}$.

As always, it is now time to fill this syntactical construct with meaning. Essentially, we take elements from the carrier sets of our algebra which have the correct types to "fit into" the leaves and evaluate the tree from the bottom up. Again, a pre-expression over a signature is purely abstract. One needs to choose an algebra (a concrete realisation of the signature) in order to "do" something with the pre-expression. In this sense, algebras are to pre-expressions as structures were to well-formed formulas.

We give first the definition, followed by an example.

**Definition 6.3.2**

Let $\mathscr{S} = (\mathscr{T}, \mathscr{F}, \langle_\rangle)$ be a signature and let $\mathscr{A} = (\mathscr{C}, \mathscr{O})$ be a $\mathscr{S}$-algebra. Let further $e = (V, E, \langle\!\lvert_\rvert\!\rangle, \star, \prec, \langle_\rangle)$ be a pre-expression in $\mathscr{S}$ with leaves $w_1, \dots, w_n$ and $\langle\langle w_i \rangle\rangle^{\text{in}} = \left(t^i_1, \dots, t^i_{k_i}\right)$.

Let finally for $i \in \{1, \dots, n\}$ and for $j \in \{1, \dots, k_i\}$ elements $c^i_j \in \mathscr{C}_{t^i_j}$ be given.

We recursively set for $v \in V$ with $\mathrm{N}^{\text{out}} v = [v_1, \dots, v_m]$

$$e[v]_{\mathscr{A}}(c^1_1, \dots, c^n_{k_n}) := \begin{cases} \mathscr{O}_{\langle w_i \rangle}(c^i_1, \dots, c^i_{k_i}) & \text{if } v = w_i \\ \mathscr{O}_{\langle v \rangle}(e[v_1]_{\mathscr{A}}(c^1_1, \dots, c^n_{k_n}), \dots, \\ \qquad e[v_m]_{\mathscr{A}}(c^1_1, \dots, c^n_{k_n})) & \text{otherwise} \end{cases}$$

and call

$$e_{\mathscr{A}}(c^1_1, c^1_2, \dots, c^1_{k_1}, c^2_1, \dots, c^n_{k_n}) := e[\sqrt{e}]_{\mathscr{A}}(c^1_1, c^1_2, \dots, c^1_{k_1}, c^2_1, \dots, c^n_{k_n})$$

the *value* of $e$ at $c^1_1, c^1_2, \dots, c^1_{k_1}, c^2_1, \dots, c^n_{k_n}$.

We apply this definition to the example on the previous page, silently extending our toy signature and algebra by the identity functions. The input type of our example pre-expression is $\mathfrak{nnnnn}$, so we choose five natural numbers, say $5, 2, 1, 0, 3$. We note to the left the values if the pre-expression is evaluated in $\mathbb{Q}_{\text{English}}$ (where 5 is actually $\frac{5}{1}$ and so forth), to the right the evaluation in $\mathbb{N}_{\text{Klingon}}$.

We represent the inputs (which are not vertices of our traversal tree) as circles.

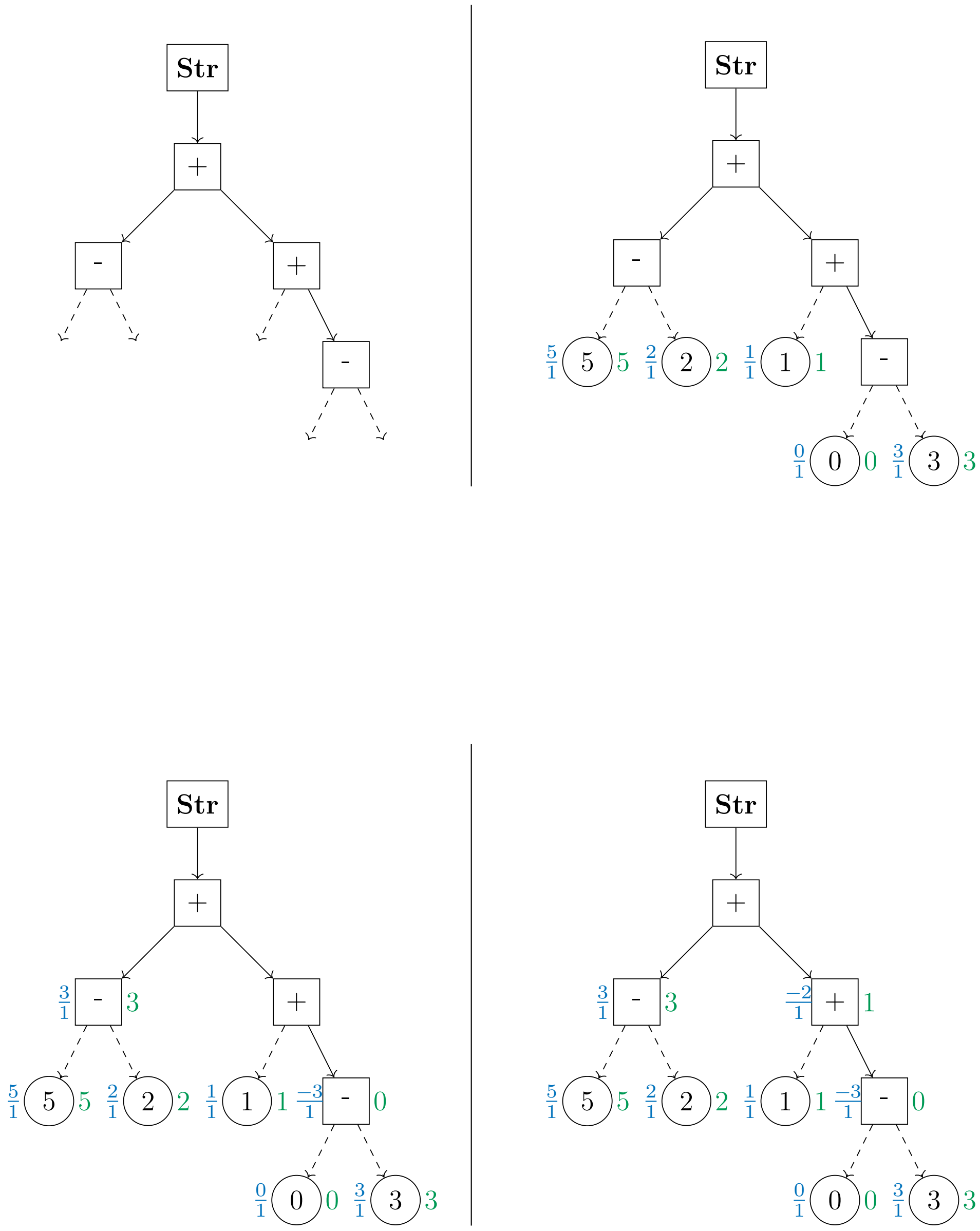

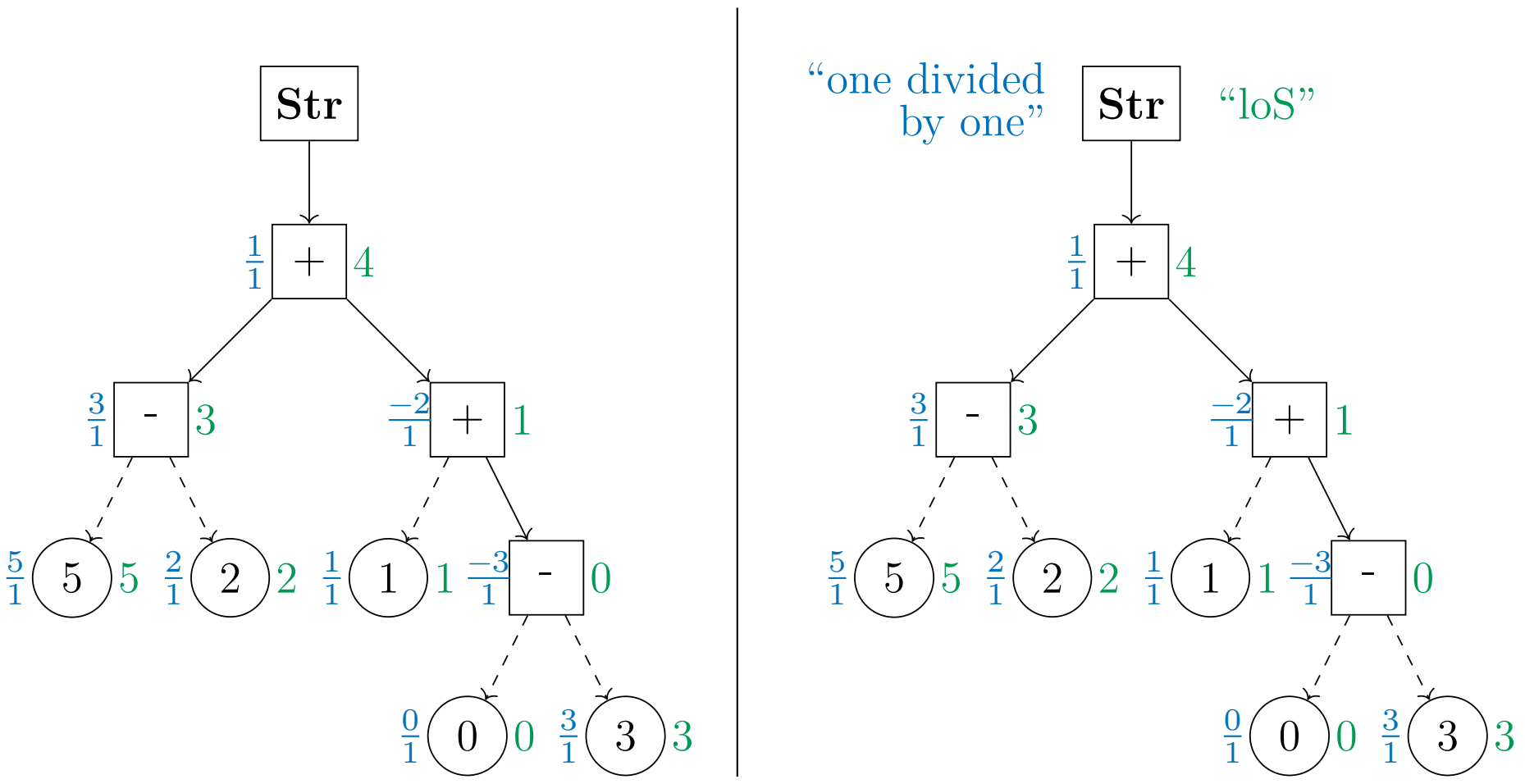


The value of the pre-expression at $(5, 2, 1, 0, 3)$ is then the label assigned to the root vertex that is, “one divided by one” in $\mathbb{Q}_{\text{English}}$ and “loS” in $\mathbb{N}_{\text{Klingon}}$.

The reader should keep in mind this intuitive understanding of how a pre-expression is evaluated.

If there are pre-expressions, there should probably be expressions as well. Indeed, an expression shall be the special case of a pre-expression that requires no input.

### Definition 6.3.3

Let $\mathscr{S} = (\mathscr{T}, \mathscr{F}, \langle _ \rangle)$ be a signature. An *expression* over $\mathscr{S}$ is a pre-expression over $\mathscr{S}$ with input type $\varepsilon$.

The set of all $\mathscr{S}$-expressions is denoted $\|\mathscr{S}\|$.

For an $\mathscr{S}$-algebra $\mathscr{A}$ and an $\mathscr{S}$-expression $e$, we call

$$\operatorname{val}_{\mathscr{A}} e := e_{\mathscr{A}}()$$

the *value* of $e$ in $\mathscr{A}$.

We call $\langle e \rangle := \langle e \rangle^{\text{out}}$ the *sort* of $e$.

$e_{\mathscr{A}}(\dots)$ denotes the result of $e$ when evaluated in $\mathscr{A}$. (def. 6.3.2, p. 101)

Since pre-expressions cannot be empty, an expression will always "end" with nullary function symbols on the leaves, which replace the input we had before.

We note a quick observation about the interaction between algebra morphisms and expressions.

**Lemma 6.3.4**

Let $\mathscr{S} = (\mathscr{T}, \mathscr{F}, \langle_\rangle)$ be a signature, $\mathscr{A} = (\mathscr{C}, \mathscr{O})$ and $\mathscr{B} = (\mathscr{D}, \mathscr{Q})$ two $\mathscr{S}$-algebras, and let $\{h_t\}_{t\in\mathscr{T}} : \mathscr{A} \to \mathscr{B}$ be an $\mathscr{S}$-algebra morphism.

Then for every expression $e \in \|\mathscr{S}\|$ of sort $t \in \mathscr{T}$, we have

$$h_t(\mathrm{val}_{\mathscr{A}} e) = \mathrm{val}_{\mathscr{B}} e.$$

$\|\mathscr{S}\|$ denotes the set of expressions over $\mathscr{S}$. (def. 6.3.3, p. 103)

$\mathrm{val}_{\mathscr{A}} e$ denotes the result of $e$ when evaluated in $\mathscr{A}$. (def. 6.3.3, p. 103)

**Proof.** This is a simple observation from induction. Suppose we have an expression $e = (V, E, \langle\!\langle_\rangle\!\rangle, \star, \preccurlyeq, \langle_\rangle)$ of sort $t$ as in the assumptions above, say of height $h \in \mathbb{N}_{>0}$.

For $h = 1$, $V$ consists of a single vertex $v$ with no children, so we must have $\langle\!\langle v\rangle\!\rangle = (\varepsilon, t)$. By definition, this means

$$h_t \mathrm{val}_{\mathscr{A}} e = h_t \mathscr{O}_{\langle v\rangle}() = \mathscr{Q}_{\langle v\rangle}() = \mathrm{val}_{\mathscr{B}} e.$$

For $h > 1$, assume the claim has been proven for expressions of height $h-1$. Let $v$ be the root of $e$, $\langle\!\langle v\rangle\!\rangle = (t_1 \dots t_n, t)$, and $\mathrm{N}^{\mathrm{out}} v = [v_1, \dots, v_n]$. But then again by definition, we have

$$\begin{aligned}
h_t(\mathrm{val}_{\mathscr{A}} e) &= h_t(\mathrm{val}_{\mathscr{A}} e[v]) \\
&= h_t(\mathscr{O}_{\langle v\rangle}(\mathrm{val}_{\mathscr{A}} e[v_1], \dots, \mathrm{val}_{\mathscr{A}} e[v_n])) \\
&= \mathscr{Q}_{\langle v\rangle}(h_{t_1}(\mathrm{val}_{\mathscr{A}} e[v_1]), \dots, h_{t_n}(\mathrm{val}_{\mathscr{A}} e[v_n])) \\
&= \mathscr{Q}_{\langle v\rangle}(\mathrm{val}_{\mathscr{B}} e[v_1], \dots, \mathrm{val}_{\mathscr{B}} e[v_n]) \\
&= \mathrm{val}_{\mathscr{B}} e[v] \\
&= \mathrm{val}_{\mathscr{B}} e,
\end{aligned}$$

proving the claim.

□

For practical applications, some finiteness assumptions will be needed. We give a name to expressions where not only inputs and output, but even all intermediary values lie in a certain subset of sorts. Algorithmically, this will be useful because when explicitly computing the value of such an expression, we never exceed a certain bound on the values kept in memory.

### Definition 6.3.5

Let $\mathscr{S} = (\mathscr{T}, \mathscr{F}, \langle_\rangle)$ be a signature, and let $\mathscr{Y} \subseteq \mathscr{T}$. An expression $e = (V, E, \langle\!\!|_|\!\!\rangle, \star, \preccurlyeq, \langle_\rangle) \in \|\mathscr{S}\|$ is called $\mathscr{Y}$*-local* if

$$\forall v \in V\colon \langle\langle v \rangle\rangle^{\text{out}} \in \mathscr{Y}.$$

$\|\mathscr{S}\|$ denotes the set of expressions over $\mathscr{S}$.
(def. 6.3.3, p. 103)

Note that due to the recursive structure of expressions, definition 6.3.5 automatically implies that also the input sorts of every vertex are in $\mathscr{Y}$.

We extend this nomenclature to sets of expressions.

### Definition 6.3.6

Let $\mathscr{S} = (\mathscr{T}, \mathscr{F}, \langle_\rangle)$ be a signature, $\mathscr{Y} \subseteq \mathscr{T}$. A set of $\mathscr{S}$-expressions is called $\mathscr{Y}$-local if all of its elements are.

A set $\mathscr{K}$ of $\mathscr{S}$-expressions is called *finitely typed* if there exists a set $\mathscr{Y} \subseteq \mathscr{T}$ with $|\mathscr{Y}| < \infty$ such that $\mathscr{K}$ is $\mathscr{Y}$-local.

As we have stated more than once, we want to build graphs from smaller graphs, such that our graphs can be seen as the value of an expression in a certain algebra. This necessarily means that there should in the first place *exist* an expression that evaluates to the desired graph.

### Definition 6.3.7

Let $\mathscr{S} = (\mathscr{T}, \mathscr{F}, \langle_\rangle)$ be a signature, $\mathscr{A} = (\mathscr{C}, \mathscr{O})$ an $\mathscr{S}$-algebra, let $t \in \mathscr{T}$, and let $\mathscr{Y} \subseteq \mathscr{T}$. A set $\mathscr{L} \subseteq \mathscr{C}_t$ is called $\mathscr{Y}$*-expressible* if for every $c \in \mathscr{L}$, there exists a $\mathscr{Y}$-local expression $e \in \|\mathscr{S}\|$ such that $c = \operatorname{val}_{\mathscr{A}} e$.

$\operatorname{val}_{\mathscr{A}} e$ denotes the result of $e$ when evaluated in $\mathscr{A}$.
(def. 6.3.3, p. 103)

A set $\mathscr{L} \subseteq \mathscr{C}_t$ is called *expressible* if it is $\mathscr{T}$-expressible.

A set $\mathscr{L} \subseteq \mathscr{C}_t$ is called *finitely expressible* if it is $\mathscr{Y}$-expressible for some

finite set $\mathscr{Y} \subseteq \mathscr{T}$.

An algebra is called $\mathscr{Y}$-expressible respectively expressible respectively finitely expressible if all of its carrier sets are.

It is important to note that for arbitrary algebras, being expressible is not a trivial requirement. Even if our set of function symbols admits enough constants (that is, nullary functions) to build useful expressions at all, not every element of a carrier set needs to be the output of a function, let alone of an expression tree. Even given an expressible algebra, simply adding a new element to one of the carrier sets and letting all functions map it to something other than itself makes that element inexpressible while still leading to a valid algebra.

However, finitely expressible algebras have properties that we shall need in chapter 7.

# 4. Building Graphs

As promised, we can now formally prove that every finite graph can be built with the tools from section 4.5. We first need a small technical lemma.

**Lemma 6.4.1**

$\mathfrak{G}$ is the algebra of graphs. (def. 6.2.2, p. 93)

Let $G = (V, E, \langle\!\lbrace _ \rbrace\!\rangle)$ be an undirected hypergraph with $V = \{ v_1, \dots, v_n \}$. If the type $n$ graph $(G, i \mapsto v_i)$ is the value of some expression in $\mathfrak{G}$, then any graph of the form $(G, t)$, $t\colon \{ 1, \dots, k \} \to V$ for any $k \in \mathbb{N}$ is the value of some expression in $\mathfrak{G}$.

**Proof.** Suppose we are given an undirected hypergraph $G = (V, E, \langle\!\lbrace _ \rbrace\!\rangle)$ with $V = \{ v_1, \dots, v_n \}$ such that the type $n$ graph $(G, i \mapsto v_i)$ can be obtained from a graph expression. Let $k \in \mathbb{N}$, $t\colon \{ 1, \dots, k \} \to V$. It suffices now to show that we can obtain the type $k$ graph $(G, t)$ from $(G, i \mapsto v_i)$ by finite application of direct sum, source redefinition, and source fusion.

$\leftrightarrows$ denotes the terminal redefinition. (def. 4.5.4, p. 41)

But one can now simply find a map $\sigma\colon \{ 1, \dots, k \} \to \{ 1, \dots, n \}$ such that $\leftrightarrows_\sigma(G, i \mapsto v_i) = (G, t)$.

□

## Theorem 6.4.2

For every expression $e \in \|\mathfrak{S}\|$, its value $\mathrm{val}_{\mathfrak{G}}\, e$ is a finite graph. Conversely, every finite graph is the value of an $\mathfrak{S}$-expression in $\mathfrak{G}$.

$\mathfrak{S}$ is the signature of $\mathfrak{G}$. (def. 6.2.1, p. 92)

$\|\mathfrak{S}\|$ denotes the set of expressions over $\mathfrak{S}$. (def. 6.3.3, p. 103)

**Proof.** For $e \in \|\mathfrak{S}\|$, the value $\mathrm{val}_{\mathfrak{G}}\, e$ is by definition an element of $\mathfrak{G}_n$ for some $n \in \mathbb{N}$.

We prove the converse by induction on the number of edges of the desired hypergraph.

$\mathrm{val}_{\mathfrak{G}}\, e$ denotes the result of $e$ when evaluated in $\mathfrak{G}$. (def. 6.3.3, p. 103)

By lemma 6.4.1, it suffices to show that for any given graph $G = (V, E, \langle\!|_|\!\rangle)$ with $V = \{ v_1, \dots, v_n \}$, we can construct the type $n$ graph $(G, i \mapsto v_i)$.

$\mathfrak{G}$ is the algebra of graphs. (def. 6.2.2, p. 93)

For $|E| = 0$, this is easily achieved by taking the disjoint sum over $n$ copies of the graph $\mathfrak{v}$.

$\mathfrak{G}_n$ denotes the set of all graphs of type $n$. (def. 4.5.1, p. 39)

Suppose now that we have proven our statement for graphs with up to $m-1$ edges, $m \in \mathbb{N}_{>0}$. Given $G = (V, E, \langle\!|_|\!\rangle)$ with $|E| = m$, we pick an arbitrary edge $e \in E$. We denote by $G' = (V, E', (\!|_|\!))$ the graph obtained by removing this edge, that is, $E' = E \setminus \{e\}$, $(\!|_|\!) = \langle\!|_|\!\rangle|_{E'}$.

$\mathfrak{v}$ is the type 1 graph with one vertex. (def. 4.5.6, p. 43)

By induction hypothesis, the type $n$ graph $X := (G', i \mapsto v_i)$ can be constructed.

Set $Y = (V'', E'', (\!|_|\!), t'') := X \oplus \mathfrak{e}_{|\langle\!| e |\!\rangle|}$.

$\oplus$ is the disjoint sum. (def. 4.5.3, p. 40)

We assume without loss of generality[6] that the end points of the (removed) edge $e$ are $t''(1), \dots, t''(|\langle\!| e |\!\rangle|)$ and that the vertices of the newly added $\mathfrak{e}_{|\langle\!| e |\!\rangle|}$ are $t''(|V''| + 1), \dots, t''(|V''| + |\langle\!| e |\!\rangle|)$.

$\mathfrak{e}_{|\langle\!| e |\!\rangle|}$ is the type $|\langle\!| e |\!\rangle|$ graph with $|\langle\!| e |\!\rangle|$ vertices and one edge. (def. 4.5.6, p. 43)

But now

$$嫐_1^{|V''|+1} 嫐_2^{|V''|+2} \dots 嫐_{|\langle\!| e |\!\rangle|}^{|V''|+|\langle\!| e |\!\rangle|} Y$$

is the desired graph up to terminal redefinition.

□

嫐 (tsureai, Japanese for *to marry*) denotes the source fusion. (def. 4.5.5, p. 42)

Thus, every finite graph can be thought of as a tree in $\|\mathfrak{S}\|$ and vice versa. This is good news for us, because we shall show that we can detect certain

[6] The terminals can be reordered to satisfy this assumption by means of a single terminal redefinition.

subsets of such trees in linear time, without explicitly computing the graph to which the expression evaluates.

$\mathfrak{G}$ is the algebra of graphs. (def. 6.2.2, p. 93)

$\mathfrak{G}_n$ denotes the set of all graphs of type $n$. (def. 4.5.1, p. 39)

$\mathfrak{v}$ is the type 1 graph with one vertex. (def. 4.5.6, p. 43)

A caveat, however, is that this naïve construction only shows that $\mathfrak{G}$ is *expressible*, but not *finitely* so. For $n \in \mathbb{N}$, constructing a graph $G \in \mathfrak{G}_n$ as in the proof above can cause an arbitrary blowup of the types of input graphs needed: for the graph with no edges and $k \in \mathbb{N}$ vertices, we begin by taking the disjoint sum over $k$ copies of $\mathfrak{v}$, ending up with a type $k$ graph. As $k$ was not bounded, this means that we can put no upper bound on the type of graphs from which we draw during our construction. However, if the reader refers forward to definition 6.3.7, they will discover that for a set of graphs to be *finitely* expressible (which would in turn enable us to apply theorem 8.6.2), the *input* graph types would have to be bounded from above.

We shall later refine our approach to get *finitely* expressible algebras of graphs. For now, we table this observation and focus on the theoretical result. Chapter 8 will pick up where we left off.

## 5. Recognisable Sets

This section is entirely based on [Cou90].

We want to show that certain subsets of the set of all graphs are “recognisable”, which will then in certain cases allow the construction of a linear-time detection algorithm. Intuitively, we now know that all graphs form an expressible algebra, and we want to show that the set of graphs fulfilling a certain monadic second-order property $\varphi$ is a “nice” subset of this algebra.

We now define what “nice” means.

### Definition 6.5.1

Let $\mathscr{S} = (\mathscr{T}, \mathscr{F}, \langle_\rangle)$ be a signature, $\mathscr{A} = (\mathscr{C}, \mathscr{O})$ an $\mathscr{S}$-algebra, $t \in \mathscr{T}$. A subset $\mathscr{L} \subseteq \mathscr{C}_t$ is called *$\mathscr{A}$-recognisable* if there exist an $\mathscr{S}$-algebra $\mathscr{B} = (\mathscr{D}, \mathscr{Q})$ and an $\mathscr{S}$-algebra morphism $\{ h_t \}_{t \in \mathscr{T}}$ fulfilling the following conditions.

- $\mathscr{B}$ is locally finite.

- There is a subset $\mathscr{M} \subseteq \mathscr{D}_t$ such that $\mathscr{L} = h_t^{-1}\mathscr{M}$.

Recall that a locally finite algebra is one where every sort admits only finitely many elements. In particular, $\mathscr{L}$ is the preimage of a finite set.

# 6. Congruences

The following will, in a minute, turn out to yield an equivalent characterisation of recognisability.

**Definition 6.6.1**

Let $\mathscr{S} = (\mathscr{T}, \mathscr{F}, \langle_\rangle)$ be a signature, $\mathscr{A} = (\mathscr{C}, \mathscr{O})$ an $\mathscr{S}$-algebra. A *congruence* on $\mathscr{A}$ is a family $\{\sim_t\}_{t\in\mathscr{T}}$ such that for every $t \in \mathscr{T}$, $\sim_t$ is an equivalence relation on $\mathscr{C}_t$ and compatible with function symbols in the following sense:

$$\forall f \in \mathscr{F} \text{ with } \langle f\rangle^{\text{in}} = (t_1, \dots, t_n) \colon \forall c_1, d_1 \in \mathscr{C}_{t_1}, \dots, \forall c_n, d_n \in \mathscr{C}_{t_n} \colon$$
$$c_1 \sim_{t_1} d_1 \wedge \dots \wedge c_n \sim_{t_n} d_n \Rightarrow \mathscr{O}_f(c_1, \dots, c_n) \sim_t \mathscr{O}_f(d_1, \dots, d_n).$$

In other words, a congruence is a family of equivalence relations that commutes with all operations.

It is easy to see that this definition ensures that the intuitive notion of a quotient algebra $\mathscr{A}/\{\sim_t\}_{t\in\mathscr{T}}$ is well-defined.

**Definition 6.6.2**

Let $\mathscr{S} = (\mathscr{T}, \mathscr{F}, \langle_\rangle)$ be a signature, $\mathscr{A} = (\mathscr{C}, \mathscr{O})$ an $\mathscr{S}$-algebra, and let $\{\sim_t\}_{t\in\mathscr{T}}$ be a congruence on $\mathscr{A}$. The *quotient algebra* $\mathscr{A}/\{\sim_t\}_{t\in\mathscr{T}}$ is the $\mathscr{S}$-algebra $\left(\{\mathscr{D}_t\}_{t\in\mathscr{T}}, \{\mathscr{Q}_f\}_{f\in\mathscr{F}}\right)$ with

$$\forall t \in \mathscr{T} : \mathscr{D}_t := \mathscr{C}_t/\sim_t$$

and

$$\forall f \in \mathscr{F} \text{ with } \langle f \rangle^{\text{in}} = (t_1, \dots, t_n) \colon \forall c_1 \in \mathscr{C}_{t_1}, \dots, c_n \in \mathscr{C}_{t_n} \colon$$
$$\mathscr{Q}_f([c_1], \dots, [c_n]) := [\mathscr{O}_f(c_1, \dots, c_n)],$$

where $[c]$ denotes the equivalence class of $c \in \mathscr{C}_t$ under $\sim_t$.

In order to relate quotient algebras to recognisability, we need local finiteness.

**Definition 6.6.3**

Let $\mathscr{S} = (\mathscr{T}, \mathscr{F}, \langle _ \rangle)$ be a signature, $\mathscr{A} = (\mathscr{C}, \mathscr{O})$ an $\mathscr{S}$-algebra, and let $\{\sim_t\}_{t \in \mathscr{T}}$ be a congruence on $\mathscr{A}$. We call $\{\sim_t\}_{t \in \mathscr{T}}$ *locally finite* if for every $t \in \mathscr{T}$, the equivalence relation $\sim_t$ has only finitely many equivalence classes.

A set $\mathscr{L}$ that is recognisable will (in the section 6.8) turn out to be exactly a set which admits a locally finite congruence such that $\mathscr{L}$ is a union of equivalence classes under that congruence. We give a name to this.

**Definition 6.6.4**

Let $\mathscr{S} = (\mathscr{T}, \mathscr{F}, \langle _ \rangle)$ be a signature, $\mathscr{A} = (\mathscr{C}, \mathscr{O})$ an $\mathscr{S}$-algebra, and $\sim = \{\sim_t\}_{t \in \mathscr{T}}$ a congruence on $\mathscr{A}$. Let $t \in \mathscr{T}$. A subset $\mathscr{L} \subseteq \mathscr{C}_t$ is called $\sim$*-saturated* if for every equivalence class, either all representatives of the class are in $\mathscr{L}$ or none of them, that is,

$$\forall c, d \in \mathscr{C}_t \colon c \in \mathscr{L} \wedge c \sim_t d \Rightarrow d \in \mathscr{L}.$$

# 7. Inductive Sets

Recall that in monadic second-order logic, any predicate can be seen as "set membership", or as dividing the universe into "true" and "false" objects. We now introduce predicates whose universe is an $\mathscr{S}$-algebra in order to work with formal logic within such an algebra. The predicates with which we shall be primarily concerned are those which map a graph $G$ to $\top$ if $G$

fulfils a certain property (for example, if it is 2-colourable) and to $\bot$ if it does not.

**Definition 6.7.1**

Let $X$ be a set. A *predicate* on $X$ is a map $X \to \{\top, \bot\}$.

The extension to typed algebras looks as follows.

**Definition 6.7.2**

Let $\mathscr{S} = (\mathscr{T}, \mathscr{F}, \langle_\rangle)$ be a signature, and let $\mathscr{A} = (\mathscr{C}, \mathscr{O})$ be an $\mathscr{S}$-algebra. A *family of predicates on $\mathscr{A}$* is a pair $(\mathscr{P}, \langle_\rangle)$ such that

- $\mathscr{P}$ is a set of functions,
- $\langle_\rangle \colon \mathscr{P} \to \mathscr{T}$,
- and each $p \in \mathscr{P}$ is a predicate on $\mathscr{C}_{\langle p \rangle}$.

That is, a family of predicates is simply a collection of predicates together with the information to which carrier set each predicate can be applied.

We are later interested in the set of all graphs fulfilling a certain predicate, say "the set of all 3-colourable graphs". So interested are we in this, we introduce a notation for it.

**Definition 6.7.3**

Let $\mathscr{S} = (\mathscr{T}, \mathscr{F}, \langle_\rangle)$ be a signature, let $\mathscr{A} = (\mathscr{C}, \mathscr{O})$ be an $\mathscr{S}$-algebra, and let $(\mathscr{P}, \langle_\rangle)$ be a family of predicates on $\mathscr{A}$. For $p \in \mathscr{P}$, we write $\lceil p \rceil := \{ c \in \mathscr{C}_{\langle p \rangle} : p(c) = \top \} \subseteq \mathscr{C}_{\langle p \rangle}$.

The reader can remember this notation by the notion that $\top$ is the "ceiling" of some poset.

Our main concern will be to show that certain sets of the form $\lceil p \rceil$ are recognisable. To this end, we meet a special class of predicate families.

## 7.1. Intuition

Say we are given a 2-uniform graph $G$, and we want to know whether its vertices are of bounded degree, say whether no vertex has more than $k$ neighbours. To this end, we are given a family of predicates $p_0, p_1, \ldots$ where $p_i(G) = \top$ if and only if no vertex of $G$ has more than $i$ neighbours.

Of course, we can check this in linear time by iterating over the vertices of $G$ and counting the neighbours. But what if we have already checked the degrees for some smaller graphs? Can we use this to speed up our computation?

Assume we know that our graph $G$ is the direct sum of two smaller graphs, say $G = G' \oplus G''$.

$\oplus$ is the disjoint sum. (def. 4.5.3, p. 40)

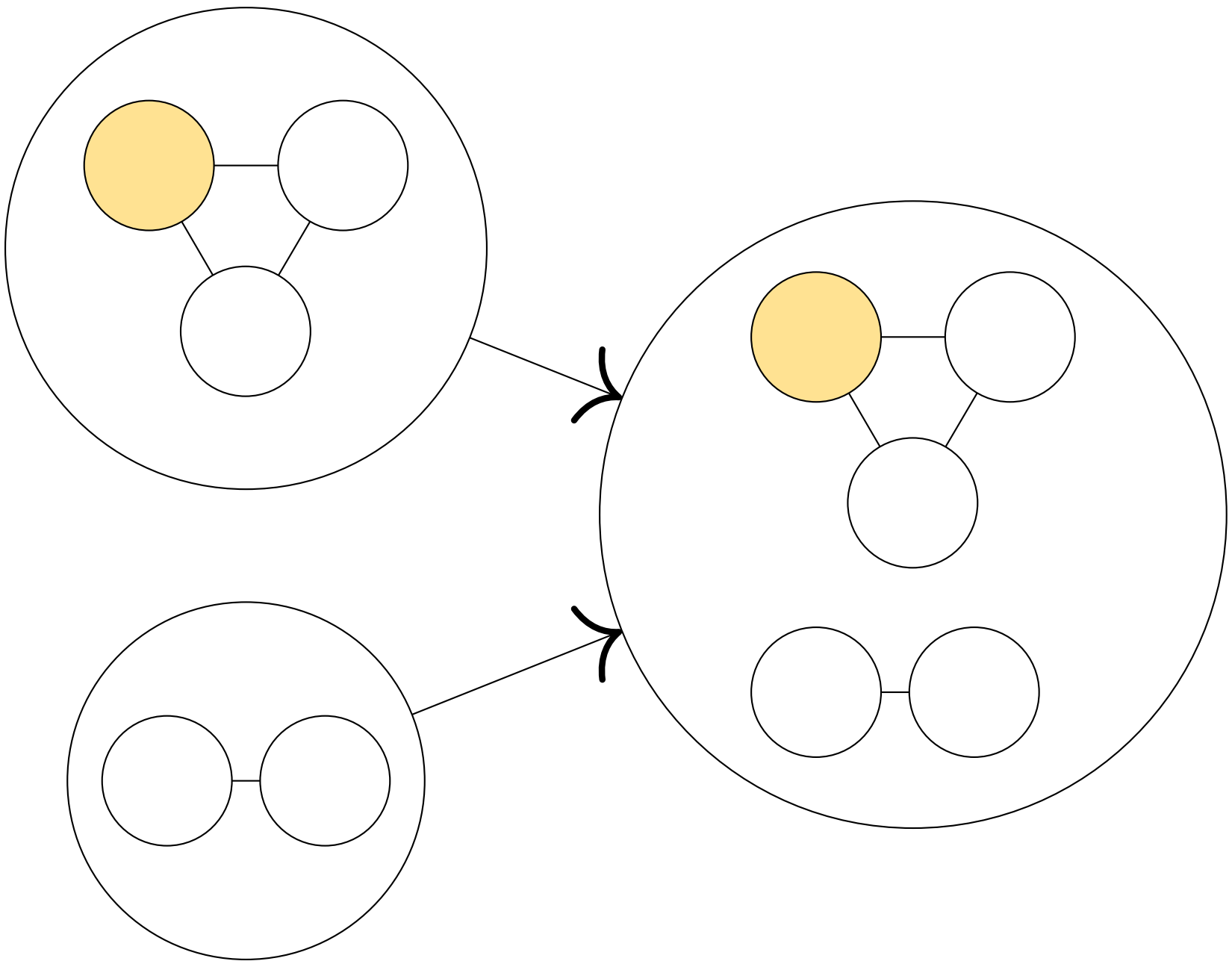

Does the number of neighbours of the marked vertex change? Of course not, since disjoint sum adds no edges.

Hence, if we have already computed in an earlier step that $\wp_2(G') = \top$ and $\wp_2(G'') = \top$, then we know instantly that $\wp_2(G) = \top$!

More generally, we can "decompose" any $\wp_i$ on a graph $G = G' \oplus G''$ into

$$\wp_i(G) = \begin{cases} \top & \text{if } \wp_i(G') = \wp_i(G'') = \top \\ \bot & \text{otherwise.} \end{cases}$$

What about terminal fusion?

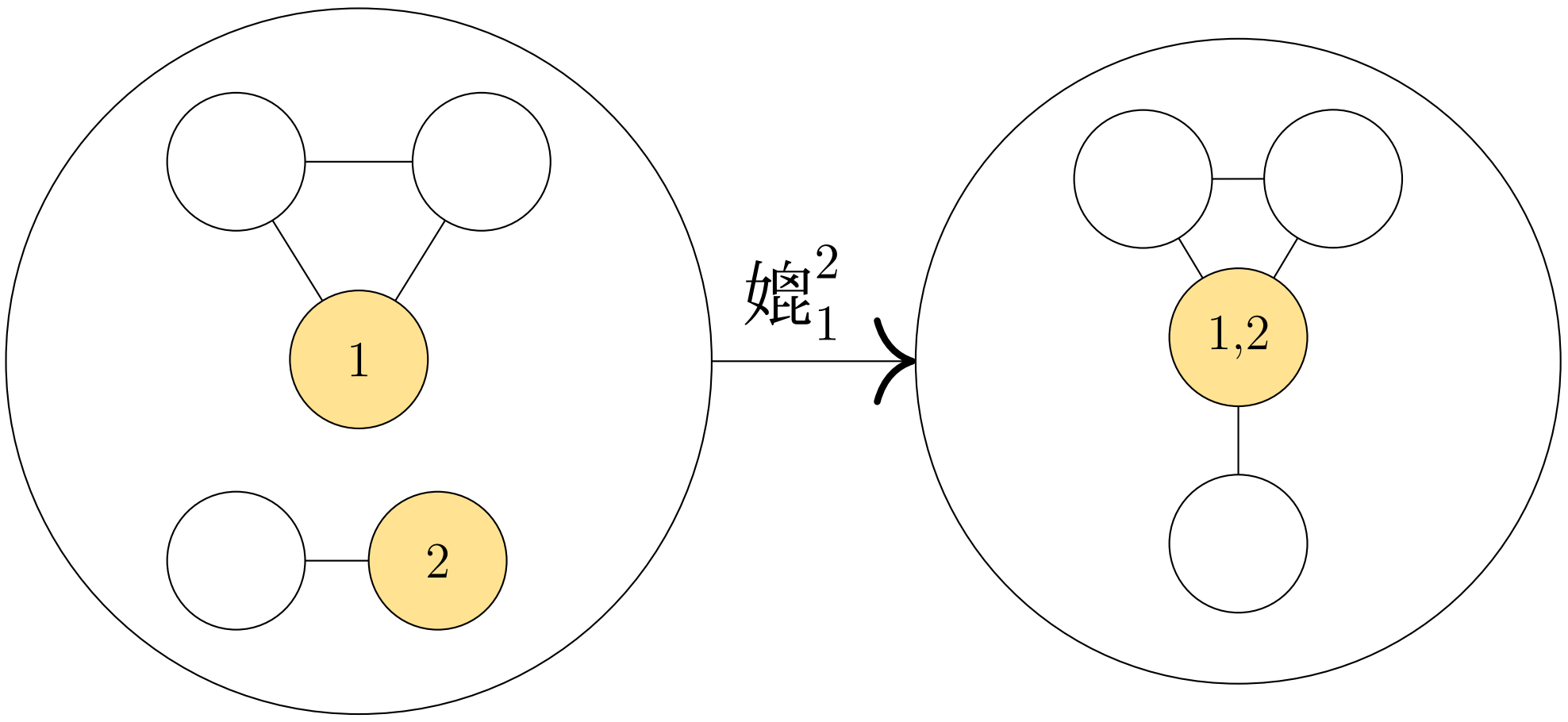


The non-terminal vertices do not change, while the vertices 1 and 2 fuse into one vertex whose degree is the sum of the degrees of the original vertices. Hence, for a graph $G$ with exactly two terminal vertices, we have $\rho_2(\text{媲}_1^2 G) = \top$ if and only if $\rho_2(G)$ and the terminal vertices are the same or their degrees add up to at most two (which is not the case in our picture, whence $\rho_2(\text{媲}_1^2 G) = \bot$).

媲 (tsureai, Japanese for *to marry*) denotes the source fusion. (def. 4.5.5, p. 42)

Therefore, if we restrict ourselves to graphs of type at most 2, we can answer our original question via dynamic programming by breaking down the graph into its constituent parts:

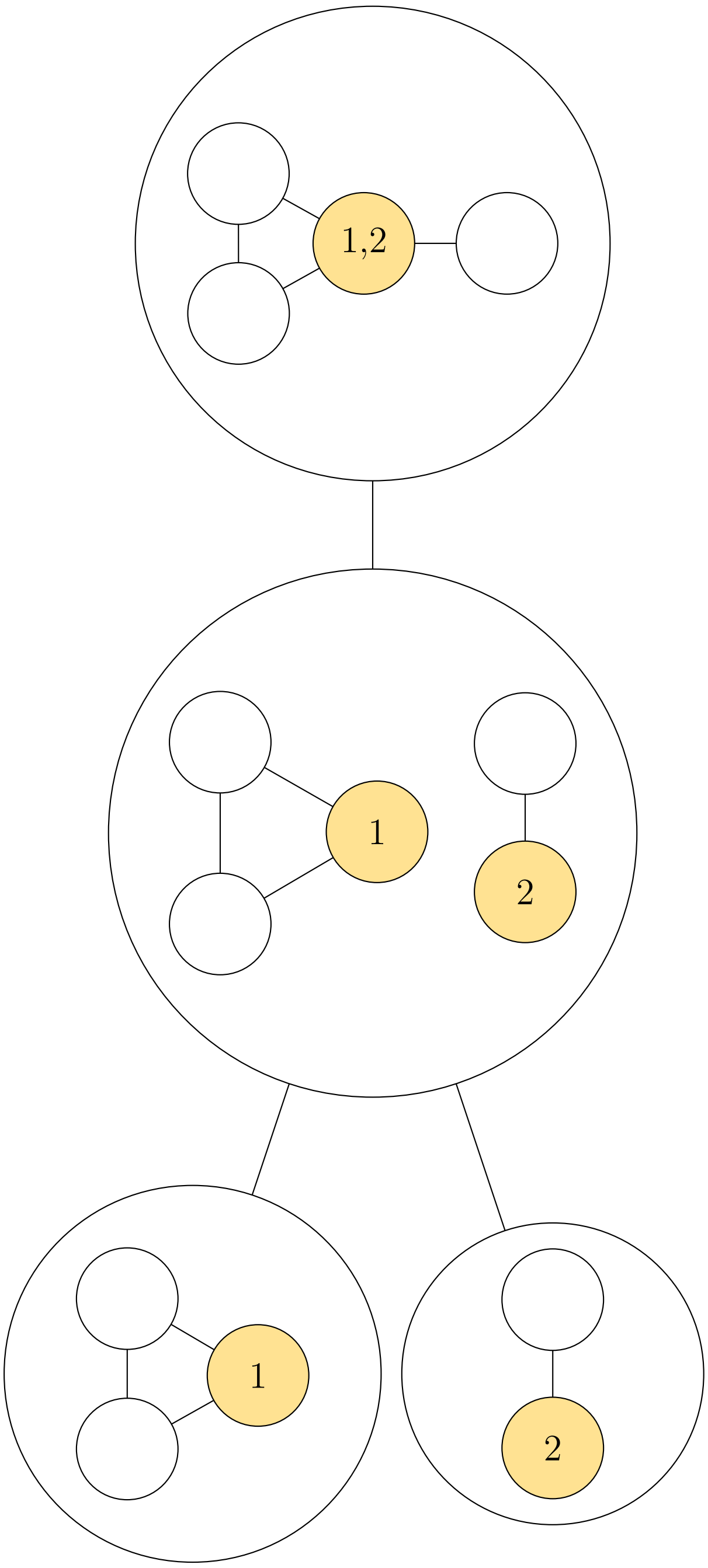


To see whether the graph at the root vertex fulfils $\rho_2$, we simply check

whether the graph at its child vertex fulfils

a. $\mathcal{p}_2$

b. The two terminal vertices are the same or their degrees add up to at most two.

To see whether property a is true, we simply check whether the graphs at the leaves both fulfil $\mathcal{p}_2$. To see whether property b is true, we check the terminals of the graphs at the leaves.

In this way, we have "broken down" a question about the root graph into a series of questions about its constituent paths. This is not very useful in this particular example, but for more complex questions, examining the smaller constituent graphs may be significantly easier than answering the question directly.

Note how, in order to compute the truth value of the predicate $\mathcal{p}_2$, we had to introduce a *second* predicate checking a statement about terminal vertices (even though we did not give it a name). This property of predicates, where we can answer one given question about a graph by computing more and potentially different questions on its constituent parts, is called *inductiveness*.

Due to the need for additional predicates, inductiveness is not a property of a single predicate, but always of a set or family of predicates. More generally, a family of predicates on some $(\mathcal{T}, \mathcal{F}, \langle_\rangle)$-algebra $(\mathcal{C}, \mathcal{O})$ will be called "inductive" if for every function symbol $\mathcal{f} \in \mathcal{F}$ and every predicate $\mathcal{p} \in \mathcal{P}$ that "knows" how to handle that function (that is, where the types are compatible such that evaluating $\mathcal{p}(\mathcal{O}_{\mathcal{f}}(\dots))$ makes sense), instead of evaluating $\mathcal{p}(\mathcal{O}_{\mathcal{f}}(c_1, \dots, c_n))$, we can find a logical formula like $\mathcal{p}_1(c_1) \wedge \mathcal{p}_2(c_2) \vee \mathcal{p}_3(c_3) \dots$ that tells us whether $\mathcal{p}(\mathcal{O}_{\mathcal{f}}(c_1, \dots, c_n))$ is true *without using* $\mathcal{f}$, where$\mathcal{p}_1, \dots$ are other predicates in $\mathcal{P}$. If this formula depends only on $\mathcal{p}$ and $\mathcal{f}$, this allows us to determine the truth of $\mathcal{p}(\mathcal{O}_{\mathcal{f}}(\dots))$ without ever actually having to compute $\mathcal{O}_{\mathcal{f}}$.

## 7.2. Formal Definition

The reader is invited to cross-reference the definition with the example from the previous section to see what is going on.

**Definition 6.7.4**

Let $\mathscr{A} = (\mathscr{C}, \mathscr{O})$ be a $(\mathscr{T}, \mathscr{F}, \langle_\rangle)$-algebra, and let $\mathscr{H} \subseteq \mathscr{F}$. A family of predicates $(\mathscr{P}, \langle_\rangle)$ on $\mathscr{A}$ is called *$\mathscr{H}$-inductive* if for every $f \in \mathscr{H}$ with $\langle f \rangle = (t_1 \dots t_n, t)$ and for every predicate $p \in \mathscr{P}$ with $\langle p \rangle = t$, the following conditions are fulfilled:

There are integers $m_1, \dots, m_n \in \mathbb{N}$ such that:

There exists a formula of propositional logic $\Phi$ with $m := \sum_{i=1}^n m_i$ variable symbols and a sequence of $m$ elements of $\mathscr{P}$, say

$$p_{1,1}, \dots, p_{1,m_1}, p_{2,1}, \dots, p_{n,m_n},$$

such that

$$\forall i \in \{1, \dots, i\}, \forall j \in \{1, \dots, m_i\} : \langle p_{i,j} \rangle = t_i$$

and

$$\begin{aligned} &\forall c_1 \in \mathscr{C}_{t_1}, \dots, \forall c_n \in \mathscr{C}_{t_n} : \\ &p(\mathscr{O}_f(c_1, \dots, c_n)) \\ &\quad = \Phi[p_{1,1}(c_1), \dots, p_{1,m_1}(c_1), p_{2,1}(c_2), \dots, p_{n,m_n}(c_n)]. \end{aligned}$$

The tuple $(\Phi, p_{1,1}, \dots, p_{1,m_1}, p_{2,1}, \dots, p_{n,m_n})$ is called an *inductive decomposition* for $p$ with regard to $f$.

Note that while recognisability is a global property (we either know the algebra morphism generating our set, or we know nothing), inductiveness is a local property – to show that a family of predicates is inductive, we can iterate over all pairs $(p, f)$ of a predicate $p$ and a function symbol $f$ and construct an inductive decomposition for $p$ with regard to $f$.

This is precisely what we shall do to prove Courcelle's Theorem: we want to show that a certain set of graphs is recognisable, but constructing abstract algebra morphisms is hard. Rather than bash our head against this wall, we show inductively that the set of graphs is of the form $\lceil p \rceil$ for some predicate $p$ in an inductive family of predicates.

$\lceil p \rceil$ denotes the set of elements which satisfy $p$. (def. 6.7.3, p. 111)

Theorem 6.8.2 shows that the three concepts we have introduced are in fact

equivalent, whence the approach just described yields the desired result.

# 8. Equivalent Characterisations

We now show that the three concepts above can be reduced to one another.

**Definition 6.8.1**

Let $\mathscr{S} = (\mathscr{T}, \mathscr{F}, \langle_\rangle)$ be a signature, and let $\mathscr{A} = (\mathscr{C}, \mathscr{O})$ be an $\mathscr{S}$-algebra. A family of predicates $(\mathscr{P}, \langle_\rangle)$ on $\mathscr{A}$ is called *locally finite* if for every $t \in \mathscr{T}$, the set $\{\, p \in \mathscr{P} : \langle p \rangle = t \,\}$ is finite.

**Theorem 6.8.2**

Let $\mathscr{S} = (\mathscr{T}, \mathscr{F}, \langle_\rangle)$ be a signature, let $\mathscr{A} = (\mathscr{C}, \mathscr{O})$ be an $\mathscr{S}$-algebra, let $t \in \mathscr{T}$, and let $\mathscr{L} \subseteq \mathscr{C}_t$. Then the following statements are equivalent.

- $\mathscr{L}$ is $\mathscr{A}$-recognisable.
- There exists a locally finite $\mathscr{F}$-inductive family $(\mathscr{P}, \langle_\rangle)$ of predicates on $\mathscr{A}$ and a $p \in \mathscr{P}$ such that $\mathscr{L} = \lceil p \rceil$.
- There exists a locally finite congruence $\sim$ on $\mathscr{A}$ such that $\mathscr{L}$ is saturated with regard to $\sim$.

$\lceil p \rceil$ denotes the set of elements which satisfy $p$. (def. 6.7.3, p. 111)

**Proof.** Fix for the entirety of this proof a signature $\mathscr{S} = (\mathscr{T}, \mathscr{F}, \langle_\rangle)$, an $\mathscr{S}$-algebra $\mathscr{A} = (\mathscr{C}, \mathscr{O})$, a type $t \in \mathscr{T}$, and a set $\mathscr{L} \subseteq \mathscr{C}_t$.

Suppose first that $\mathscr{L}$ is $\mathscr{A}$-recognisable with $\mathscr{L} = h_t^{-1}\beth$ for a locally finite $\mathscr{S}$-algebra $\mathscr{B} = (\mathscr{D}, \mathscr{Q})$, an $\mathscr{S}$-algebra morphism $\{h_t\}_{t \in \mathscr{T}} : \mathscr{A} \to \mathscr{B}$, and some $\beth \subseteq \mathscr{D}_t$.We construct a locally finite $\mathscr{F}$-inductive family of predicates such that $\mathscr{L}$ is of the form $\lceil p \rceil$ for some predicate $p$.

We use the notation $\beth$ (nun) here because this will soon turn out to be (in certain cases) the set of accepting states of a special deterministic bottom-up finite tree automaton.

We fix a set of functions (whose values will be defined in a second)

$$\mathscr{P} := \{\, p^u_d : \mathscr{C}_u \to \{\top, \bot\} : u \in \mathscr{T}, d \in \mathscr{D}_u \,\}$$

as well as one additional function $p : \mathscr{C}_t \to \{\top, \bot\}$ and set

$$\langle_\rangle : \mathscr{P} \cup \{\, p \,\} \to \mathscr{T}, x \mapsto \begin{cases} u & x = p^u_d \text{ for some } d \\ t & x = p, \end{cases}$$

making $\mathscr{P} \cup \{p\}$ into a family of predicates on $\mathscr{A}$.

The predicates' values are defined as

$$\forall u \in \mathscr{T} \forall d \in \mathscr{D}_u : p^u_d : c \mapsto \begin{cases} \top & h_u c = d \\ \bot & \text{otherwise.} \end{cases}$$

In other words, the predicate $p^u_d$ is a sort of indicator function for the set of preimages of $d$. The predicate $p$, then, shall be an indicator function for the preimage of 𐤍:

𐤍 is the Phoenician letter "nun".

$$p : c \mapsto \begin{cases} \top & h_t c \in \text{𐤍} \\ \bot & \text{otherwise.} \end{cases}$$

Thus by definition, we have $\mathscr{L} = \lceil p \rceil$.

$\lceil p \rceil$ denotes the set of elements which satisfy $p$. (def. 6.7.3, p. 111)

Since $\mathscr{B}$ is locally finite, so is $(\mathscr{P} \cup \{p\}, \langle_\rangle)$. It remains to show that this family is $\mathscr{F}$-inductive.

Take a symbol $f \in \mathscr{F}$ with $\langle f \rangle = (u_1 \dots u_n, u)$ and a predicate $p^u_d \in \mathscr{P}$.[7]

We want to find predicates $p_{i,j}$ and a propositional formula $\Phi$ such that we have $p^u_d(\mathscr{O}_f(\dots)) = \Phi[p_{i,j}(\dots)]$.

Fix for now some $c_1 \in \mathscr{C}_{u_1}, \dots, c_n \in \mathscr{C}_{u_n}$. Because $\{h_t\}_{t \in \mathscr{T}}$ is an algebra morphism, by definition we know that

$$h_u \mathscr{O}_f(c_1, \dots, c_n) = \mathscr{Q}_f(h_{u_1} c_1, \dots, h_{u_n} c_n)$$

and thus

$$p^u_d(\mathscr{O}_f(c_1, \dots, c_n)) = \top \Leftrightarrow \mathscr{Q}_f(h_{u_1} c_1, \dots, h_{u_n} c_n) = d.$$

Since each $\mathscr{D}_{u_i}$ is finite, the inputs for which $\mathscr{Q}_f$ equals $d$ can be written down and expressed as a logical disjunction of conjunctions, that is,

$$\begin{aligned} &\mathscr{Q}_f(x_1, \dots, x_n) = d \\ &\Leftrightarrow (x_1 = d^1_1 \wedge \dots \wedge x_n = d^n_1) \vee (x_1 = d^1_2 \wedge \dots \wedge x_n = d^n_2) \vee \dots \end{aligned}$$

Essentially, we build a giant table of the possible inputs of $\mathscr{Q}_f$ (which is a function in $\mathscr{B}$) and express it as a propositional formula.

[7] Apart from potentially $p$, this covers all predicates of compatible type.

Using our family of predicates, we conclude that for $y_1 \in \mathscr{C}_{u_1}, \dots, y_n \in \mathscr{C}_{u_n}$ we have

$$\begin{aligned} &\mathcal{Q}_{f}(h_{u_1} y_1, \dots, h_{u_n} y_n) = d \\ \Leftrightarrow\quad & p^{u_1}_{d^1_1} y_1 = \top \wedge \dots \wedge p^{u_n}_{d^n_1} y_n = \top \\ &\vee p^{u_1}_{d^1_2} y_1 = \top \wedge \dots \wedge p^{u_n}_{d^n_2} y_n = \top \\ &\vee \dots, \end{aligned}$$

whence

$$\begin{aligned} &p^{u}_{d}(\mathcal{O}_{f}(y_1, \dots, y_n)) = \top \\ \Leftrightarrow\quad & p^{u_1}_{d^1_1} y_1 = \top \wedge \dots \wedge p^{u_n}_{d^n_1} y_n = \top \\ &\vee p^{u_1}_{d^1_2} y_1 = \top \wedge \dots \wedge p^{u_n}_{d^n_2} y_n = \top \\ &\vee \dots \end{aligned}$$

The right-hand side depends only on $f$ and $p^u_d$, so we have found an inductive decomposition for $p^u_d$ with regard to $f$.

The only predicate which we have not examined is $p$, but the construction here is analogous, except that the formula is even longer since we have to allow any of the values in the (finite) set 𐤍.

𐤍 is the Phoenician letter “nun”.

We have now constructed an inductive decomposition for every predicate in $\mathscr{P} \cup \{p\}$ with regard to every compatible function symbol, proving that the family $(\mathscr{P} \cup \{p\}, \langle_\rangle)$ is $\mathscr{F}$-inductive, showing the first implication from the theorem.

Let now $(\mathscr{P}, \langle_\rangle)$ be a locally finite $\mathscr{F}$-inductive family of predicates on $\mathscr{A}$ such that $\mathscr{L} = \lceil p \rceil$ for some $p \in \mathscr{P}$. We construct a locally finite congruence $\sim$ on $\mathscr{A}$ such that $\mathscr{L}$ is saturated with regard to $\sim$.

$\lceil p \rceil$ denotes the set of elements which satisfy $p$. (def. 6.7.3, p. 111)

For every $u \in \mathscr{T}$, we define the following equivalence relation on $\mathscr{C}_u$:

$$\forall c, d \in \mathscr{C}_u \colon c \sim_u d :\Leftrightarrow \forall p \in \langle_\rangle^{-1}\{u\} \colon p(c) = p(d),$$

that is, two elements are equivalent if all compatible predicates agree on their truth or falsity.

Due to the inductiveness of $(\mathscr{P}, \langle_\rangle)$, we can evaluate predicates on functions by considering only their inputs and a propositional formula, whence it is

straightforward to see that the set $\{ \sim_u \}_{u \in \mathscr{T}}$ is a congruence on $\mathscr{A}$. The reader is invited to turn back to definition 6.6.1 to convince themselves of this.

Consider now a set $\{ c_1, \dots \} \subseteq \mathscr{C}_u$ of elements of pairwise different equivalence classes under $\sim_u$. By definition of $\sim_u$, this means that every pair $(c_i, c_j)$ disagrees on at least one predicate $p \in \mathscr{P}$. Since there are only finitely many predicates of type $u$, the equivalence relation $\sim_u$ has only finitely many equivalence classes.

$\lceil p \rceil$ denotes the set of elements which satisfy $p$. (def. 6.7.3, p. 111)

Having thus constructed a locally finite congruence on $\mathscr{A}$, it remains to show that $\mathscr{L}$ is saturated with regard to it. But $\mathscr{L} = \lceil p \rceil$ for some predicate $p \in \mathscr{P}$, so if $c \in \mathscr{L}$ and $d \in \mathscr{C}_t$ with $c \sim_t d$, by definition of $\sim_t$ we must have $d \in \mathscr{L}$.

This shows the second implication.

Let finally $\mathscr{L}$ be saturated with regard to some locally finite congruence $\sim = \{ \sim_u \}_{u \in \mathscr{T}}$ on $\mathscr{A}$. We must show that $\mathscr{L}$ is then $\mathscr{A}$-recognisable.

Consider the canonical algebra morphism $\pi \colon \mathscr{A} \to \mathscr{A}/\sim$ mapping each element to its equivalence class. Since $\sim$ is locally finite, the sets $\mathscr{C}_u/\sim_u$ are finite, and hence $\mathscr{A}/\sim$ is locally finite.

𐤍 is the Phoenician letter “nun”.

Set $𐤍 := \pi_t \mathscr{L}$. Because $\mathscr{L}$ is $\sim$-saturated, all elements of $\mathscr{C}_t$ that would be mapped to $\pi_t \mathscr{L}$ must already be in $\mathscr{L}$, whence $\pi_t^{-1} 𐤍 = \mathscr{L}$, proving the final implication.

□

# 9. From Formulas to Predicates

$\mathring{\mathfrak{L}}$ is the circuitous language of graphs. (def. 5.6.1, p. 77)

$\mathring{\Omega}(G)$ is the circuitous induced structure of $G$. (def. 5.6.2, p. 78)

The predicates with which we shall be working are essentially the functions that, for a given sentence $\varphi$ in $\mathring{\mathfrak{L}}$, map a graph $G$ to $\top$ if $\vDash_{\mathring{\Omega}(G)} \varphi$ and to $\bot$ if not. However, in order to apply theorem 6.8.2, we need this family to be locally finite, and we can certainly build infinitely many logical formulas in $\mathring{\mathfrak{L}}$.

We get around this restriction by bounding various parameters of $\varphi$. Recall

that, in the end, we shall only care about *one particular given sentence* (for which we want to build a tree automaton), so these bounds are no real restrictions – once given a sentence, we simply choose our parameters accordingly.

### Definition 6.9.1

Let $w \in \mathbb{N}$. The set of well-formed formulas respectively sentences of $\mathring{\mathfrak{L}}$ that have width at most $w$ is denoted by $|\mathring{\mathfrak{L}}|_w$ respectively by $\|\mathring{\mathfrak{L}}\|_w$.

Let further $l \in \mathbb{N}$. The set of well-formed formulas respectively sentences of $\|\mathring{\mathfrak{L}}\|_w$ that use no variable symbol except $\delta_0, \dots, \delta_l$ is denoted by $|\mathring{\mathfrak{L}}|_w^l$ respectively by $\|\mathring{\mathfrak{L}}\|_w^l$.

Let finally $k \in \mathbb{N}$. The set of well-formed formulas respectively sentences of $\|\mathring{\mathfrak{L}}\|_w^l$ that use no predicate symbol $\lambda_{\text{conn}}^t$ for $t > k$ is denoted by $|\mathring{\mathfrak{L}}|_w^{l,k}$ respectively by $\|\mathring{\mathfrak{L}}\|_w^{l,k}$.

$\mathring{\mathfrak{L}}$ is the circuitous language of graphs. (def. 5.6.1, p. 77)

The width of a formula is the maximum number of nested quantifiers – see definition 5.3.5.

$|\mathring{\mathfrak{L}}|$ is the set of all well-formed formulas over $\mathring{\mathfrak{L}}$. (def. 5.3.5, p. 61)

These restrictions almost suffice to make $\|\mathring{\mathfrak{L}}\|$ into a finite set up to tautological equivalence[8], but $\mathring{\mathfrak{L}}$ has an additional infinite amount of function symbols $\delta_{\text{term}}$ for us to worry about. We get around this by adopting an equivalence relation that is slightly different from tautological equivalence.

$\|\mathring{\mathfrak{L}}\|$ denotes the set of sentences over $\mathring{\mathfrak{L}}$. (def. 5.3.6, p. 62)

### Definition 6.9.2

Let $n \in \mathbb{N}$. Two sentences $\varphi, \psi \in \|\mathring{\mathfrak{L}}\|$ are called *type $n$ equivalent*, written $\varphi \overset{n}{\sim} \psi$, if

$$\forall G \in \mathfrak{G}_n : \vDash_{\mathring{\Omega}(G)} \varphi \Leftrightarrow \vDash_{\mathring{\Omega}(G)} \psi .$$

$\mathfrak{G}_n$ denotes the set of all graphs of type $n$. (def. 4.5.1, p. 39)

In other words, two type $n$ equivalent sentences hold true for exactly the same set of graphs *of type $n$*. The definition says nothing about graphs of other types, which is indeed an important fact since the functions $\delta_{\text{term}}$ can detect terminal vertices and might thus lead to different truth values on the

$\mathring{\Omega}(G)$ is the circuitous induced structure of $G$. (def. 5.6.2, p. 78)

[8] Two sentences being "tautologically equivalent", in layperson's terms, means that one can take the first formula, do a finite sequence of trivial transformations on it such as resolving parentheses or applying de Morgan's laws, and end up with the second formula. This certainly implies (and is often even equivalent to) both formulas yielding the same truth value on all variable assignments. An in-depth discussion can be found in [End72].

same underlying graph depending on how many terminal vertices one adds.

Definition 6.9.2 is actually a special case of the following definition, which will make our proof easier to read.

### Definition 6.9.3

Let $n \in \mathbb{N}$. Two terms $\chi, \xi \in \mathring{\mathfrak{L}}_{\text{Term}}$ are called *type $n$ equivalent*, written $\chi \stackrel{n}{\sim} \xi$, if

- $\overleftarrow{\chi} = \overleftarrow{\xi}$ and
- for every graph $G \in \mathfrak{G}_n$ and for every variable assignment $\tau$ in $\mathring{\Omega}(G)$ which is full for $\chi$, we have

$$\chi[\tau] = \xi[\tau].$$

Two well-formed formulas $\varphi, \psi \in |\mathring{\mathfrak{L}}|$ are called *type $n$ equivalent*, written $\varphi \stackrel{n}{\sim} \psi$, if

- $\overleftarrow{\varphi} = \overleftarrow{\psi}$ and
- for every graph $G \in \mathfrak{G}_n$ and for every variable assignment $\tau$ in $\mathring{\Omega}(G)$ which is full for $\varphi$, we have

$$\varphi[\tau] \leftrightarrow \psi[\tau].$$

$\mathring{\mathfrak{L}}$ is the circuitous language of graphs. (def. 5.6.1, p. 77)

$\mathring{\mathfrak{L}}_{\text{Term}}$ is the set of all terms of $\mathring{\mathfrak{L}}$. (def. 5.3.3, p. 60)

$\overleftarrow{\chi}$ denotes the set of free variables of $\chi$. (def. 5.3.6, p. 62)

$\mathfrak{G}_n$ denotes the set of all graphs of type $n$. (def. 4.5.1, p. 39)

$\mathring{\Omega}(G)$ is the circuitous induced structure of $G$. (def. 5.6.2, p. 78)

$|\mathring{\mathfrak{L}}|$ is the set of all well-formed formulas over $\mathring{\mathfrak{L}}$. (def. 5.3.5, p. 61)

Type $n$ equivalence is manifestly an equivalence relation. We now show that for every $n \in \mathbb{N}$, the number of equivalence classes of sentences under $\stackrel{n}{\sim}$ is finite, yielding the local finiteness we so desperately desire.

First, a small technical observation.

### Lemma 6.9.4

Let $n \in \mathbb{N}$, and let $\varphi \in |\mathring{\mathfrak{L}}|$. Let further $\chi$ be a term occurring in $\varphi$, let $\xi$ be a term with $\xi \stackrel{n}{\sim} \chi$, and let $\varphi'$ denote the formula where every occurrence of $\chi$ is replaced by $\xi$.

Then $\varphi \stackrel{n}{\sim} \varphi'$.

**Proof.** Let $n \in \mathbb{N}$, and let $\varphi \in |\mathring{\mathfrak{L}}|$. Let $\chi$ be a term occurring in $\varphi$, $\xi$ a

term with $\xi \overset{n}{\sim} \chi$, and denote for any formula the same formula with every occurrence of $\chi$ replaced by $\xi$ with a prime.

$\overset{n}{\sim}$ denotes equivalence on graphs of type $n$. (def. 6.9.3, p. 122)

We prove the claim by induction over the cases of definition 5.3.5.

*Case 1: $\varphi$ is atomic, say $\varphi = \rho(\chi_1, \dots, \chi_k)$ for some $k \in \mathbb{N}$.* If $k = 0$, then there is nothing to show.

If $k > 0$, without loss of generality, the term to replace is $\chi_1$. We are given another term $\xi$ with $\xi \overset{n}{\sim} \chi_1$, in particular, $\overleftarrow{\xi} = \overleftarrow{\chi}_1$, whence $\overleftarrow{\varphi}' = \overleftarrow{\varphi}$.

$\overleftarrow{\xi}$ denotes the set of free variables of $\xi$. (def. 5.3.6, p. 62)

Let now $G \in \mathfrak{G}_n$ and let $\tau$ be a variable assignment in $\mathring{\Omega}(G)$ which is full for $\varphi$. Then by definition 5.3.11, we have

$\mathfrak{G}_n$ denotes the set of all graphs of type $n$. (def. 4.5.1, p. 39)

$$\begin{aligned}
\varphi[\tau] \leftrightarrow \top &\Leftrightarrow \begin{cases} (\chi_1[\tau], \chi_2[\tau], \dots, \chi_n[\tau]) \in \mathring{\Omega}(G)(\Lambda) & \text{if } \rho = \Lambda \text{ for some } \Lambda \in {}_n\Lambda_{\mathring{\mathfrak{L}}} \\ (\chi_1[\tau], \chi_2[\tau], \dots, \chi_n[\tau]) \in \tau(\lambda) & \text{if } \rho = \lambda \text{ for some } \lambda \in {}_n\lambda_{\mathring{\mathfrak{L}}} \end{cases} \\
&\Leftrightarrow \begin{cases} (\xi[\tau], \chi_2[\tau], \dots, \chi_n[\tau]) \in \mathring{\Omega}(G)(\Lambda) & \text{if } \rho = \Lambda \text{ for some } \Lambda \in {}_n\Lambda_{\mathring{\mathfrak{L}}} \\ (\xi[\tau], \chi_2[\tau], \dots, \chi_n[\tau]) \in \tau(\lambda) & \text{if } \rho = \lambda \text{ for some } \lambda \in {}_n\lambda_{\mathring{\mathfrak{L}}} \end{cases} \\
&\Leftrightarrow \varphi'[\tau] \leftrightarrow \top,
\end{aligned}$$

$\mathring{\Omega}(G)$ is the circuitous induced structure of $G$. (def. 5.6.2, p. 78)

${}_n\Lambda_{\mathring{\mathfrak{L}}}$ is the set of $n$-place predicate symbols of $\mathring{\mathfrak{L}}$.

$\mathring{\mathfrak{L}}$ is the circuitous language of graphs. (def. 5.6.1, p. 77)

${}_n\lambda_{\mathring{\mathfrak{L}}}$ is the set of $n$-place predicate variables of $\mathring{\mathfrak{L}}$.

whence $\varphi \overset{n}{\sim} \varphi'$.

*Case 2: $\varphi = \neg\psi$ for some $\psi \in |\mathring{\mathfrak{L}}|$.* The claim follows immediately since by induction hypothesis, we have $\psi \overset{n}{\sim} \psi'$.

$|\mathring{\mathfrak{L}}|$ is the set of all well-formed formulas over $\mathring{\mathfrak{L}}$. (def. 5.3.5, p. 61)

*Case 3: $\varphi = \psi \wedge \zeta$ for some $\psi, \zeta \in |\mathring{\mathfrak{L}}|$.* Again, the claim follows immediately since by induction hypothesis, we have $\psi \overset{n}{\sim} \psi'$ and $\zeta \overset{n}{\sim} \zeta'$.

${}_0\delta_{\mathring{\mathfrak{L}}}$ is the set of 0-place function variables of $\mathring{\mathfrak{L}}$.

*Case 4: $\varphi = \forall\delta_0\psi$ for some $\psi \in |\mathring{\mathfrak{L}}|$ and some $\delta_0 \in {}_0\delta_{\mathring{\mathfrak{L}}}$.* By induction hy-

$\overset{n}{\sim}$ denotes equivalence on graphs of type $n$. (def. 6.9.3, p. 122)

$\mathring{\Omega}(G)$ is the circuitous induced structure of $G$. (def. 5.6.2, p. 78)

$\overline{\varphi}$ denotes the set of free variables of $\varphi$. (def. 5.3.6, p. 62)

$\mathring{\mathfrak{L}}$ is the circuitous language of graphs. (def. 5.6.1, p. 77)

pothesis, we have $\psi \overset{n}{\sim} \psi'$, and checking definition 5.3.11 reveals that

$$\begin{aligned}\varphi[\tau] &\Leftrightarrow \psi[\kappa] \leftrightarrow \top \quad \text{for every full variable assignment } \kappa \\ &\qquad\qquad\quad \text{for } \psi \text{ in } \mathring{\Omega}(G) \text{ with } \kappa|_{\overline{\varphi}} = \tau|_{\overline{\varphi}} \\ &\Leftrightarrow \psi'[\kappa] \leftrightarrow \top \quad \text{for every full variable assignment } \kappa \\ &\qquad\qquad\quad \text{for } \psi' \text{ in } \mathring{\Omega}(G) \text{ with } \kappa|_{\overline{\varphi}} = \tau|_{\overline{\varphi}} \\ &\Leftrightarrow \varphi'[\tau],\end{aligned}$$

as desired.

The claim is hence proven for all well-formed formulas of $\mathring{\mathfrak{L}}$.

□

The reader should now refer back to definition 6.9.1 to remind themselves of the restrictions on our formulas.

$|\mathring{\mathfrak{L}}|$ is the set of all well-formed formulas over $\mathring{\mathfrak{L}}$. (def. 5.3.5, p. 61)

**Theorem 6.9.5**

Let $n \in \mathbb{N}$, and let $w, l, k \in \mathbb{N}$. Then the set

$$|\mathring{\mathfrak{L}}|_w^{l,k} \big/ \overset{n}{\sim}$$

is finite.

**Proof.** Let $n \in \mathbb{N}$. Let further $w, l, k \in \mathbb{N}$. We show that the set $|\mathring{\mathfrak{L}}|_w^{l,k}$ admits only finitely many pairwise nonequivalent well-formed formulas for type $n$, implying the claim.

As always, we move along the trail of definitions 5.3.2 to 5.3.5.

We first ascertain that there are only finitely many terms to consider before applying the recursive procedure from definition 5.3.5.

The function symbol $\varnothing$ can only be used to build a single term, $\varnothing()$.

Let now $\chi$ be an elementary term, that is, either a nullary function variable symbol or the term $\varnothing()$. There are exactly $l + 2$ different such terms, a finite number.[9] Let further $K, K' \subseteq \mathbb{N}_{>0}$ be finite. If $\chi$ is a variable

[9] There are variable symbols indexed $0, \dots, l$ plus the symbol $\varnothing$.

symbol and $\tau$ is a variable assignment on a type $n$ graph $G = (V, E, \langle\!\!\langle _ \rangle\!\!\rangle, t)$ with $\tau(\chi) \subseteq V$, then

$$\begin{aligned}\delta_{\text{term}}^{K'}(\delta_{\text{term}}^{K}(\chi))[\tau] &= \tau(\chi) \cup \{\, t(i) : i \in K, i \le n \,\} \cup \{\, t(i) : i \in K', i \le n \,\} \\ &= \tau(\chi) \cup \{\, t(i) : i \in K \cup K', i \le n \,\} \\ &= \delta_{\text{term}}^{K \cup K'}(\chi)[\tau].\end{aligned}$$

If, on the other hand, $\tau(\chi) \subseteq E$, then

$$\delta_{\text{term}}^{K'}(\delta_{\text{term}}^{K}(\chi))[\tau] = \chi[\tau] = \delta_{\text{term}}^{K \cup K'}(\chi)[\tau].$$

For $\chi = \varnothing$, we have $\delta_{\text{term}}^{K'}(\delta_{\text{term}}^{K}(\varnothing)) = \delta_{\text{term}}^{K \cup K'}(\varnothing)$ by the same reasoning as for a vertex variable.

By induction, it is now clear that any term with nested occurrences of $\delta_{\text{term}}$ is type $n$ equivalent to the same term with those nested occurrences replaced by single $\delta_{\text{term}}$ symbols. Since there are no other function symbols in $\mathring{\mathfrak{L}}$, this shows that there are only finitely many type $n$ non-equivalent terms in $\mathring{\mathfrak{L}}$.

$\mathring{\mathfrak{L}}$ is the circuitous language of graphs. (def. 5.6.1, p. 77)

Using lemma 6.9.4, we can now assume without loss of generality that all terms are of the form $\delta_0$, $\varnothing$, $\delta_{\text{term}}^{K}(\delta_0)$, or $\delta_{\text{term}}^{K}(\varnothing)$, where $\delta_0$ is a variable symbol and $K \subseteq \mathbb{N}_{>0}$ is finite, because if the formula we are considering contains a term not of this form, we can replace this term by one of ours and obtain a type $n$ equivalent formula.

It remains to show that we can restrict ourselves to finitely many of the function symbols $\delta_{\text{term}}^{K}$.

For a finite set $K \subseteq \mathbb{N}_{>0}$, denote by $K_{\le n}$ the set $\{\, i \in K : i \le n \,\}$.

Consider a term of the form $\delta_{\text{term}}^{K}(\chi)$, where $\chi$ is an elementary term and $K \subseteq \mathbb{N}_{>0}$ is finite, and consider a type $n$ graph $G = (V, E, \langle\!\!\langle _ \rangle\!\!\rangle, t)$ with a variable assignment $\tau$. If $\chi$ is a variable symbol with $\tau(\chi) \subseteq V$, then

$$\begin{aligned}\delta_{\text{term}}^{K}(\chi)[\tau] &= \tau(\chi) \cup \{\, t(i) : i \in K, i \le n \,\} \\ &= \tau(\chi) \cup \{\, t(i) : i \in K_{\le n} \,\} \\ &= \delta_{\text{term}}^{K_{\le n}}(\chi)[\tau],\end{aligned}$$

and analogously for $\chi = \varnothing$. For $\tau(\chi) \subseteq E$, we get even more easily

$$\delta_{\text{term}}^{K}(\chi)[\tau] = \chi[\tau] = \delta_{\text{term}}^{K_{\le n}}(\chi)[\tau].$$

Since the only thing that went into this reasoning was the type of $G$, we can henceforth (again by virtue of lemma 6.9.4) assume without loss of generality that no function symbol $\delta_{\text{term}}^{K}$ references a set $K$ with any element greater than $n$.

$2^{\{1,\dots,i\}}$denotes the power set of $\{1,\dots,i\}$.

Because $\bigcup_{i=1}^{n} 2^{\{1,\dots,i\}} = 2^{\{1,\dots,n\}}$ is finite and there are no non-nullary function symbols except for the $\delta_{\text{term}}^{K}$, this implies that, if we restrict ourselves to the first $l+1$ variable symbols, there are without loss of generality only finitely many terms that can occur in a formula up to type $n$ equivalence.

We now move on to atomic formulas.

We restrict our predicate symbols to use of the symbols of type $\lambda_{\text{conn}}$ only the first $k+1$ many, $\lambda_{\text{conn}}^{0}, \dots, \lambda_{\text{conn}}^{k}$, leaving only finitely many predicate symbols in total. Because the set of terms is also finite, this yields only finitely many atomic formulas.

We show that these combine to only finitely many non-equivalent well-formed formulas by induction on $w$, that is, on the maximum number of nested quantifiers.

$\overset{n}{\approx}$ denotes equivalence on graphs of type $n$. (def. 6.9.3, p. 122)

We assume without loss of generality that no formula contains the substring $\neg\neg$ since for any formula $\varphi$, we have $\varphi \overset{n}{\approx} \neg\neg\varphi$.

For $w = 0$, formulas cannot contain quantifiers. Hence every formula is of the form $\varphi_1 \wedge \dots \wedge \varphi_r$ for some $r \in \mathbb{N}_{>0}$, where every $\varphi_i$ is either an atomic formula or the negation of an atomic formula.

Since there are only finitely many atomic formulas and two occurrences of the same $\varphi_i$ are type $n$ equivalent to one occurrence of $\varphi_i$ (indeed, $\varphi \wedge \varphi$ always has the same truth value as $\varphi$ under all full variable assignments), there can only be finitely many non-equivalent formulas of this form.

Suppose now the claim has been proven for formulas of width at most $w \in \mathbb{N}$. For a formula of the form $\forall\delta_0\varphi$ of width $w+1$, the formula $\varphi$ has width $w$. In particular, there are only finitely many choices for $\varphi$. There are also only $l+1$ choices for $\delta_0$, meaning that there are only finitely many well-formed formulas of the form $\forall\delta_0\varphi$.

But an arbitrary formula of width $w+1$ must be of the form $\varphi_1 \wedge \dots \wedge \varphi_r$ for some $r \in \mathbb{N}_{>0}$, where each $\varphi_i$ is of the form $\forall\delta_0\varphi$, $\neg\forall\delta_0\varphi$, $\varphi$, or $\neg\varphi$

for some variable symbol $\delta_0$ and some well-formed formula $\varphi$ of width at most $w$. By the same argument as before, there can only be finitely many non-equivalent well-formed formulas of this form.

□

Since every sentence is in particular a well-formed formula, we get an immediate special case.

### Corollary 6.9.6

Let $n \in \mathbb{N}$, and let $w, l, k \in \mathbb{N}$. Then the set

$$\|\mathring{\mathfrak{L}}\|_w^{l,k} \big/ \overset{n}{\approx}$$

is finite.

$\mathring{\mathfrak{L}}$ is the circuitous language of graphs. (def. 5.6.1, p. 77)

$\|\mathring{\mathfrak{L}}\|$ denotes the set of sentences over $\mathring{\mathfrak{L}}$. (def. 5.3.6, p. 62)

$\overset{n}{\approx}$ denotes equivalence on graphs of type $n$. (def. 6.9.3, p. 122)

We now define a family of predicates on the algebra of graphs which will transform a graph property (given as a monadic second-order sentence) into the language of algebras for use in our proofs.

### Definition 6.9.7

Let $w, l, k \in \mathbb{N}$, let $n \in \mathbb{N}$, and let $\overline{\varphi} \in \|\mathring{\mathfrak{L}}\|_w^{l,k} \big/ \overset{n}{\approx}$. We set

$$\mathscr{p}_{\overline{\varphi}}^{n} \colon \mathfrak{G}_n \to \{\top, \bot\}, G \mapsto \begin{cases} \top & \text{if } \vDash_{\mathring{\Omega}(G)} \varphi \text{ for some}^{10}\ \varphi \in \overline{\varphi} \\ \bot & \text{otherwise} \end{cases}$$

and

$$\mathscr{P}_{w,l,k}^{n} := \{ \mathscr{p}_{\overline{\varphi}}^{n} : \overline{\varphi} \in \|\mathring{\mathfrak{L}}\|_w^{l,k} \big/ \overset{n}{\approx} \}.$$

$\mathring{\Omega}(G)$ is the circuitous induced structure of $G$. (def. 5.6.2, p. 78)

$\mathfrak{G}_n$ denotes the set of all graphs of type $n$. (def. 4.5.1, p. 39)

With this definition, for any $w, l, k \in \mathbb{N}$, the pair $\left(\{ \mathscr{P}_{w,l,k}^{n} \}_{n \in \mathbb{N}}, \mathscr{p}_{\overline{\varphi}}^{n} \mapsto n\right)$ is a family of predicates on $\mathfrak{G}$ (see definition 6.7.2). For a given monadic second order sentence $\varphi$, the predicate $\mathscr{p}_{\overline{\varphi}}^{n}$ tells us about every graph of type $n$ whether or not it fulfils the formula $\varphi$. Since the predicates for type $n$ equivalent formulas would be the same anyway, we have immediately identified them.

$\mathfrak{G}$ is the algebra of graphs. (def. 6.2.2, p. 93)

[10] Here, "for some" is equivalent to "for all", since all formulas in $\overline{\varphi}$ are type $n$ equivalent.

We give this family a name for later reference.

**Definition 6.9.8**

Let $w, l, k \in \mathbb{N}$. We set

$$\mathscr{P}_w^{l,k} := \left( \{ \mathscr{P}_{w,l,k}^n \}_{n \in \mathbb{N}}, \langle _ \rangle \colon \mathscr{p}_{\overset{\circ}{\varphi}}^n \mapsto n \right).$$

Remembering why we are doing all this, we want to show that we can apply theorem 6.8.2.

**Lemma 6.9.9**

Let $w, l, k \in \mathbb{N}$. Then $\mathscr{P}_w^{l,k}$ is locally finite.

**Proof.** Let $w, l, k \in \mathbb{N}$. Recall that "locally finite", by definition 6.8.1, means that for every $n \in \mathbb{N}$, the set $\{ \mathscr{p} \in \mathscr{P}_w^{l,k} : \langle \mathscr{p} \rangle = n \}$ should be finite.

But this is precisely the statement of corollary 6.9.6.

□

We are now ready for the final strike: we prove that the family of predicates so defined is inductive. Afterwards, all the obstacles will fall like dominoes before the combined might of our theorems.

We split the proof into three theorems, one for each kind of function symbol in $\mathfrak{G}$.

$\mathfrak{G}$ is the algebra of graphs. (def. 6.2.2, p. 93)

**Theorem 6.9.10**

Let $w, l, k \in \mathbb{N}$. Then $\mathscr{P}_w^{l,k}$ is $\mathfrak{F}_\oplus$-inductive.

$\oplus$ is the disjoint sum. (def. 4.5.3, p. 40)

**Proof.** Let us fix some $w, l, k \in \mathbb{N}$.

Recall first from definition 6.2.1 that $\mathfrak{F}_\oplus = \{ {}_a^b\oplus : a, b \in \mathbb{N} \}$. Checking definition 6.7.4, we discover that we need to find an inductive decomposition for every pair $(\mathscr{p}_{\overset{\circ}{\varphi}}^n, {}_a^b\oplus)$ with $\mathscr{p}_{\overset{\circ}{\varphi}}^n \in \mathscr{P}_w^{l,k}$ and $a + b = n$.

We first observe that in order to decompose a predicate $\mathscr{p}_{\overset{\circ}{\varphi}}^n$ with regard to ${}_a^b\oplus$ with $a + b = n$, we can simply decompose a sentence $\varphi \in \overline{\varphi}$: if we

find a formula $\Phi$ of propositional logic and sentences $\varphi_1^{\mathrm{a}}, \dots, \varphi_r^{\mathrm{a}}, \varphi_1^{\mathrm{b}}, \dots, \varphi_s^{\mathrm{b}}$ such that for all $G \in \mathfrak{G}_a$ and for all $G' \in \mathfrak{G}_b$ we have

$\mathfrak{G}_a$ denotes the set of all graphs of type $a$. (def. 4.5.1, p. 39)

$$\vDash_{\mathring{\Omega}(G \oplus G')} \varphi \Leftrightarrow \Phi[\vDash_{\mathring{\Omega}(G)} \varphi_1^{\mathrm{a}}, \dots, \vDash_{\mathring{\Omega}(G)} \varphi_r^{\mathrm{a}}, \vDash_{\mathring{\Omega}(G)} \varphi_1^{\mathrm{b}}, \dots, \vDash_{\mathring{\Omega}(G)} \varphi_s^{\mathrm{b}}]$$

$\oplus$ is the disjoint sum. (def. 4.5.3, p. 40)

and we can choose $\varphi_1^{\mathrm{a}}, \dots, \varphi_r^{\mathrm{a}}, \varphi_1^{\mathrm{b}}, \dots, \varphi_s^{\mathrm{b}}$ to be in $\|\mathring{\mathfrak{L}}\|_w^{l,k}$, then for all $G \in \mathfrak{G}_a$ and for all $G' \in \mathfrak{G}_b$ we have

$$p_{\varphi}^{n}(G \oplus G') = \Phi[p_{\varphi_1^{\mathrm{a}}}^{a}(G), \dots, p_{\varphi_r^{\mathrm{a}}}^{a}(G), p_{\varphi_1^{\mathrm{b}}}^{b}(G'), \dots, p_{\varphi_s^{\mathrm{b}}}^{b}(G')]$$

$\mathring{\Omega}(G \oplus G')$ is the circuitous induced structure of $G \oplus G'$. (def. 5.6.2, p. 78)

and thus an inductive decomposition for $p_{\varphi}^{n}$ with regard to ${}_a^b\oplus$.

$\mathring{\mathfrak{L}}$ is the circuitous language of graphs. (def. 5.6.1, p. 77)

In order to show that every sentence decomposes like this, we need to once again use the recursive structure of definition 5.3.5 and allow free variables. We proceed as follows: suppose we are given graphs $G \in \mathfrak{G}_a$ and $G' \in \mathfrak{G}_b$ and a variable assignment $\tau$ in $\mathring{\Omega}(G \oplus G')$. We define new variable assignments $\tau_G$ in $\mathring{\Omega}(G)$ and $\tau_{G'}$ in $\mathring{\Omega}(G')$ with each at most as many variables as $\tau$. If we are now given a formula $\varphi$ with width at most $w$, at most $l$ different variable symbols, and predicates constrained by $k$, possibly containing free variables, and a full variable assignment $\tau$ for $\varphi$ in $\mathring{\Omega}(G \oplus G')$, we show that there are formulas $\varphi_1^{\mathrm{a}}, \dots, \varphi_r^{\mathrm{a}}, \varphi_1^{\mathrm{b}}, \dots, \varphi_s^{\mathrm{b}}$ and a propositional formula $\Phi$ such that $\tau_G$ is a full variable assignment for $\varphi_1^{\mathrm{a}}, \dots, \varphi_r^{\mathrm{a}}$, $\tau_{G'}$ is a full variable assignment for $\varphi_1^{\mathrm{b}}, \dots, \varphi_s^{\mathrm{b}}$, we have

$\|\mathring{\mathfrak{L}}\|$ denotes the set of sentences over $\mathring{\mathfrak{L}}$. (def. 5.3.6, p. 62)

$$\varphi[\tau] \leftrightarrow \Phi[\varphi_1^{\mathrm{a}}[\tau_G], \dots \varphi_r^{\mathrm{a}}[\tau_G], \varphi_1^{\mathrm{b}}[\tau_{G'}], \dots \varphi_s^{\mathrm{b}}[\tau_{G'}]],$$

and such that all the formulas $\varphi_1^{\mathrm{a}}, \dots, \varphi_r^{\mathrm{a}}, \varphi_1^{\mathrm{b}}, \dots, \varphi_s^{\mathrm{b}}$ still obey the same restrictions (according to $w, l, k$) as the original formula.

Then in particular, we can decompose every sentence in $\|\mathring{\mathfrak{L}}\|_w^{l,k}$ as discussed above, and the theorem follows.

Note that we need not check that the new formulas obey the type restriction (that is, only use the functions $\delta_{\mathrm{term}}^{K}$ for $K \subseteq \{1, \dots, a\}$ respectively $K \subseteq \{1, \dots, b\}$) because once we find *any* formula that works on $G$, we can find a type $a$ equivalent one which obeys the restriction by the construction in the proof of theorem 6.9.5.

We first construct the new variable assignments. Fix to this end two graphs $G = (V, E, \langle\!\!\!\mid _ \mid\!\!\!\rangle, t) \in \mathfrak{G}_a$ and $G' = (V', E', \langle\!\!\!\mid _ \mid\!\!\!\rangle, t') \in \mathfrak{G}_b$.

$\oplus$ is the disjoint sum. (def. 4.5.3, p. 40)

$\mathring{\Omega}(G \oplus G')$ is the circuitous induced structure of $G \oplus G'$. (def. 5.6.2, p. 78)

$\mathring{\mathfrak{L}}$ is the circuitous language of graphs. (def. 5.6.1, p. 77)

Let now $\tau$ be a variable assignment in $\mathring{\Omega}(G \oplus G')$ with domain $X$.

We set

$$\tau_G : X \to |\mathring{\Omega}(G)|, x \mapsto \tau(x) \cap (V \cup E)\,,$$

recalling that all variables in $\mathring{\mathfrak{L}}$ are set-valued and that all vertices and edges of $G$ are contained in $G \oplus G'$. (As always, we assume that $G$ and $G'$ are disjoint.)

Analogously, we set

$$\tau_{G'} : X \to |\mathring{\Omega}(G')|, x \mapsto \tau(x) \cap (V' \cup E')\,,$$

ending up with variable assignments in $\mathring{\Omega}(G)$ and $\mathring{\Omega}(G')$, respectively.

Let now a well-formed formula $\varphi$ be given, and let $\tau$ be a full variable assignment for $\varphi$ in $\mathring{\Omega}(G \oplus G')$. We show the existence of the claimed formulas by induction.

We introduce the following notation: let $K \subseteq \{1, \dots, a+b\}$. We set

$$K_a := \{ i : i \in K, i \le a \}$$

and

$$K_b := \{ i - a : i \in K, i > a \}.$$

If $K$ is a set of indices of terminal vertices, that means that $K_a$ contains just the terminal vertices in $G$ that correspond to these indices in $G \oplus G'$, and $K_b$ contains those in $G'$. The indices in $K_b$, consequently, must be shifted down to reach the correct correspondence since the terminal vertices of $G \oplus G'$ are those of $G$ concatenated with those of $G'$.

For every term $\chi$, we set

$$\chi^a := \begin{cases} \chi & \text{if } \chi \text{ elementary} \\ \delta^{K_a}_{\text{term}}(\xi) & \text{if } \chi = \delta^{K}_{\text{term}}(\xi) \\ & \text{for some } K \subseteq \{1, \dots, a+b\}, \xi \text{ elementary} \end{cases}$$

and analogously

$$\chi^b := \begin{cases} \chi & \text{if } \chi \text{ elementary} \\ \delta^{K_b}_{\text{term}}(\xi) & \text{if } \chi = \delta^{K}_{\text{term}}(\xi) \\ & \text{for some } K \subseteq \{1, \dots, a+b\}, \xi \text{ elementary.} \end{cases}$$

In the following, let $\chi_0, \dots$ be arbitrary terms. Let $\varphi$ be a well-formed formula in $|\mathring{\mathfrak{L}}|_w^{l,k}$.

$\mathring{\mathfrak{L}}$ is the circuitous language of graphs. (def. 5.6.1, p. 77)

*Case 1: $\varphi$ is atomic.* Recall that all terms are either elementary or of the form $\delta_{\mathrm{term}}^K(\chi)$ for an elementary term $\chi$ and that there are no predicate variables in $\mathring{\mathfrak{L}}$.

$|\mathring{\mathfrak{L}}|$ is the set of all well-formed formulas over $\mathring{\mathfrak{L}}$. (def. 5.3.5, p. 61)

*Case 1.1: $\varphi = \lambda_{\mathrm{sgl}}(\chi_0)$.* We set

$$\begin{aligned}\varphi_1^{\mathrm{a}} &:= \lambda_{\mathrm{sgl}}(\chi_0^a),\\ \varphi_1^{\mathrm{b}} &:= \lambda_{\mathrm{sgl}}(\chi_0^b),\\ \varphi_2^{\mathrm{a}} &:= \chi_0^a \sqsubseteq \varnothing(),\\ \varphi_2^{\mathrm{b}} &:= \chi_0^b \sqsubseteq \varnothing(),\end{aligned}$$

and

$$\Phi := (\varphi_1^{\mathrm{a}} \wedge \varphi_2^{\mathrm{b}}) \vee (\varphi_1^{\mathrm{b}} \wedge \varphi_2^{\mathrm{a}}).$$

Then $\varphi_1^{\mathrm{a}}$, $\varphi_2^{\mathrm{a}}$, $\varphi_1^{\mathrm{b}}$, and $\varphi_2^{\mathrm{b}}$ have the same width, number of variables, and predicate restrictions as $\varphi$, they have the same set of free variables (making $\tau_G$ and $\tau_{G'}$ full), and due to the definition of $\tau_G$ and $\tau_{G'}$, we have

$$\begin{aligned}\varphi[\tau] &\leftrightarrow |\chi_0[\tau]| = 1\\ &\leftrightarrow |\chi_0^a[\tau] \cap |\mathring{\Omega}(G)|| = 1 \wedge \chi_0^b[\tau] \cap |\mathring{\Omega}(G')| = \varnothing\\ &\quad \vee |\chi_0^b[\tau] \cap |\mathring{\Omega}(G')|| = 1 \wedge \chi_0^a[\tau] \cap |\mathring{\Omega}(G)| = \varnothing\\ &\leftrightarrow (\varphi_1^{\mathrm{a}}[\tau_G] \wedge \varphi_2^{\mathrm{b}}[\tau_{G'}]) \vee (\varphi_1^{\mathrm{b}}[\tau_{G'}] \wedge \varphi_2^{\mathrm{a}}[\tau_G])\\ &\leftrightarrow \Phi[\varphi_1^{\mathrm{a}}[\tau_G], \varphi_2^{\mathrm{a}}[\tau_G], \varphi_1^{\mathrm{b}}[\tau_{G'}], \varphi_2^{\mathrm{b}}[\tau_{G'}]],\end{aligned}$$

$\mathring{\Omega}(G)$ is the circuitous induced structure of $G$. (def. 5.6.2, p. 78)

as desired.

In the following cases, we shall omit the step-by-step transformation of $\varphi[\tau]$ and trust that the reader can follow our arguments without this aid (or with a pen and paper).

*Case 1.2: $\varphi = \chi_0 \sqsubseteq \chi_1$.* We set $\varphi^{\mathrm{a}} := \chi_0^a \sqsubseteq \chi_1^a$, analogously $\varphi^{\mathrm{b}} := \chi_0^b \sqsubseteq \chi_1^b$, and $\Phi := \varphi^{\mathrm{a}} \wedge \varphi^{\mathrm{b}}$. Width, number of variables, predicate restrictions, and set of free variables again do not change, and since $G$ and $G'$ partition $G \oplus G'$, we have $\varphi[\tau] \leftrightarrow \varphi^{\mathrm{a}}[\tau_G] \wedge \varphi^{\mathrm{b}}[\tau_{G'}]$.

$\oplus$ is the disjoint sum. (def. 4.5.3, p. 40)

*Case 1.3:* $\varphi = \lambda^i_{\text{conn}}(\chi_0, \chi_1, \dots, \chi_i)$ *for some* $i \in \mathbb{N}$*,* $i \leq k$*.* We set

$$\varphi^{\mathrm{a}} := \lambda^i_{\text{conn}}(\chi^a_0, \chi^a_1, \dots, \chi^a_i),\ \varphi^{\mathrm{b}} := \lambda^i_{\text{conn}}(\chi^b_0, \chi^b_1, \dots, \chi^b_i),$$

$\oplus$ is the disjoint sum. (def. 4.5.3, p. 40)

and $\Phi := \varphi^{\mathrm{a}} \vee \varphi^{\mathrm{b}}$. Width, number of variables, predicate restrictions, and set of free variables again do not change, and since no edge in $G \oplus G'$ can connect a vertex from $G$ to a vertex from $G'$, the claim follows.

The claim is thus shown for all atomic formulas.

*Case 2:* $\varphi$ *is built from formulas for which the claim has been proven.*

*Case 2.1:* $\varphi = \neg\psi$. Let $(\Phi_\psi, \psi^{\mathrm{a}}_1, \dots, \psi^{\mathrm{a}}_r, \psi^{\mathrm{b}}_1, \dots, \psi^{\mathrm{b}}_s)$ be the decomposition of $\psi$. We set

$$\varphi^{\mathrm{a}}_1 := \psi^{\mathrm{a}}_1, \dots, \varphi^{\mathrm{a}}_r := \psi^{\mathrm{a}}_r,\ \varphi^{\mathrm{b}}_1 := \psi^{\mathrm{b}}_1, \dots, \varphi^{\mathrm{b}}_s := \psi^{\mathrm{b}}_s,$$

and $\Phi := \neg\Phi_\psi$. Once more, the restrictions do not change (by induction hypothesis), and of course

$$\varphi[\tau] \leftrightarrow (\neg\psi)[\tau] \leftrightarrow \neg(\psi[\tau]) \leftrightarrow \neg\Phi[\psi^{\mathrm{a}}_1[\tau_G], \dots, \psi^{\mathrm{a}}_r[\tau_G], \psi^{\mathrm{b}}_1[\tau_{G'}], \dots, \psi^{\mathrm{b}}_s[\tau_{G'}]].$$

*Case 2.2:* $\varphi = \psi \wedge \zeta$. Let $(\Phi_\psi, \psi^{\mathrm{a}}_1, \dots, \psi^{\mathrm{a}}_{r_\psi}, \psi^{\mathrm{b}}_1, \dots, \psi^{\mathrm{b}}_{s_\psi})$ be the decomposition of $\psi$ and let $(\Phi_\zeta, \zeta^{\mathrm{a}}_1, \dots, \zeta^{\mathrm{a}}_{s_\zeta}, \zeta^{\mathrm{b}}_1, \dots, \zeta^{\mathrm{b}}_{s_\zeta})$ be the decomposition of $\zeta$. Set

$$\varphi^{\mathrm{a}}_1 := \psi^{\mathrm{a}}_1,\ \dots,\ \varphi^{\mathrm{a}}_{r_\psi} := \psi^{\mathrm{a}}_{r_\psi},\ \varphi^{\mathrm{a}}_{r_\psi+1} := \zeta^{\mathrm{a}}_1,\ \dots,\ \varphi^{\mathrm{a}}_{r_\psi+r_\zeta} := \psi^{\mathrm{a}}_{r_\zeta}$$

and

$$\varphi^{\mathrm{b}}_1 := \psi^{\mathrm{b}}_1,\ \dots,\ \varphi^{\mathrm{b}}_{r_\psi} := \psi^{\mathrm{b}}_{r_\psi},\ \varphi^{\mathrm{b}}_{r_\psi+1} := \zeta^{\mathrm{b}}_1,\ \dots,\ \varphi^{\mathrm{b}}_{r_\psi+r_\zeta} := \psi^{\mathrm{b}}_{r_\zeta},$$

which again preserves all restrictions as per the induction hypothesis. Setting $\Phi := \Phi_\psi \wedge \Phi_\zeta$ yields the desired result by virtue of

$$\begin{aligned}
\varphi[\tau] &\leftrightarrow \psi[\tau] \wedge \zeta[\tau] \\
&\leftrightarrow \quad \Phi_\psi[\psi^{\mathrm{a}}_1[\tau_G], \dots, \psi^{\mathrm{a}}_{r_\psi}[\tau_G], \psi^{\mathrm{b}}_1[\tau_{G'}], \dots, \psi^{\mathrm{b}}_{s_\psi}[\tau_{G'}]] \\
&\qquad \wedge \Phi_\zeta[\zeta^{\mathrm{a}}_1[\tau_G], \dots, \zeta^{\mathrm{a}}_{s_\zeta}[\tau_G], \zeta^{\mathrm{b}}_1[\tau_{G'}], \dots, \zeta^{\mathrm{b}}_{s_\zeta}[\tau_{G'}]] \\
&\leftrightarrow (\Phi_\psi \wedge \Phi_\zeta)[\psi^{\mathrm{a}}_1[\tau_G], \dots, \psi^{\mathrm{a}}_{r_\psi}[\tau_G], \psi^{\mathrm{b}}_1[\tau_{G'}], \dots, \psi^{\mathrm{b}}_{s_\psi}[\tau_{G'}], \\
&\qquad\qquad \zeta^{\mathrm{a}}_1[\tau_G], \dots, \zeta^{\mathrm{a}}_{s_\zeta}[\tau_G], \zeta^{\mathrm{b}}_1[\tau_{G'}], \dots, \zeta^{\mathrm{b}}_{s_\zeta}[\tau_{G'}]].
\end{aligned}$$

*Case 3:* $\varphi = \exists\mu_0\psi$. We use the existential quantifier rather than the universal one simply for convenience. As the reader might remember, it matters not which we use to define our language and which is a shorthand.

Let $(\Phi_\psi, \psi_1^{\mathrm{a}}, \dots, \psi_r^{\mathrm{a}}, \psi_1^{\mathrm{b}}, \dots, \psi_s^{\mathrm{b}})$ be the decomposition of $\psi$.

For a full variable assignment $\tau$ for $\varphi$ whose domain (without loss of generality) does not include $\mu_0$, we have

$$\begin{aligned}
\varphi[\tau] &\leftrightarrow (\exists\mu_0\psi)[\tau] \\
&\leftrightarrow \exists X \subseteq |\mathring{\Omega}(G \oplus G')| \colon \\
&\qquad \psi[\tau, \mu_0 \mapsto X] \\
&\leftrightarrow \exists X \subseteq |\mathring{\Omega}(G \oplus G')| \colon \\
&\qquad \Phi_\psi\big[\psi_1^{\mathrm{a}}[\tau_G, \mu_0 \mapsto X \cap |\mathring{\Omega}(G)|], \dots, \psi_{r_\psi}^{\mathrm{a}}[\tau_G, \mu_0 \mapsto X \cap |\mathring{\Omega}(G)|], \\
&\qquad\quad \psi_1^{\mathrm{b}}[\tau_{G'}, \mu_0 \mapsto X \cap |\mathring{\Omega}(G')|], \dots, \psi_{s_\psi}^{\mathrm{b}}[\tau_{G'}, \mu_0 \mapsto X \cap |\mathring{\Omega}(G')|]\big].
\end{aligned}$$

$\oplus$ is the disjoint sum.
(def. 4.5.3, p. 40)

$\mathring{\Omega}(G \oplus G')$ is the circuitous induced structure of $G \oplus G'$.
(def. 5.6.2, p. 78)

Without loss of generality, we can assume that $\Phi_\psi$ is in disjunctive normal form, that is,

$$\Phi_\psi = \Phi_1 \vee \dots \vee \Phi_c$$

for some $\in \mathbb{N}_{>0}$ such that for all $i \in \{1, \dots, c\}$ we have

$$\Phi_i = \Phi_i^{\mathrm{a},1} \wedge \dots \wedge \Phi_i^{\mathrm{a},t_{\mathrm{a}}} \wedge \Phi_i^{\mathrm{b},1} \wedge \dots \wedge \Phi_i^{\mathrm{b},t_{\mathrm{b}}}$$

for some $t_{\mathrm{a}}, t_{\mathrm{b}} \in \mathbb{N}$ such that for every $j \in \{1, \dots, t_{\mathrm{a}}\}$, the expression $\Phi_i^j$ is either $\psi_u^{\mathrm{a}}$ or $\neg\psi_u^{\mathrm{a}}$ for some $u \in \{1, \dots, r\}$ and for every $j \in \{1, \dots, t_{\mathrm{b}}\}$, the expression $\Phi_i^j$ is either $\psi_u^{\mathrm{b}}$ or $\neg\psi_u^{\mathrm{b}}$ for some $u \in \{1, \dots, s\}$.

Hence

$$\begin{aligned}
\exists X \subseteq |\mathring{\Omega}(G \oplus G')| \colon \Phi_\psi[\dots] &\leftrightarrow \exists X \subseteq |\mathring{\Omega}(G \oplus G')| \colon \Phi_1[\dots] \vee \dots \vee \Phi_c[\dots] \\
&\leftrightarrow \quad \big(\exists X \subseteq |\mathring{\Omega}(G \oplus G')| \colon \Phi_1[\dots]\big) \\
&\qquad \vee \dots \\
&\qquad \vee \big(\exists X \subseteq |\mathring{\Omega}(G \oplus G')| \colon \Phi_c[\dots]\big).
\end{aligned}$$

We pick an $i \in \{1, \dots, c\}$ and notice that, since $|\mathring{\Omega}(G)|$ and $|\mathring{\Omega}(G')|$

$\oplus$ is the disjoint sum. (def. 4.5.3, p. 40)

$\mathring{\Omega}(G \oplus G')$ is the circuitous induced structure of $G \oplus G'$. (def. 5.6.2, p. 78)

partition $|\mathring{\Omega}(G \oplus G')|$, we have

$$\begin{aligned}
& \exists X \subseteq |\mathring{\Omega}(G \oplus G')| \colon \ \Phi_i[...] \\
\leftrightarrow\ & \exists X \subseteq |\mathring{\Omega}(G \oplus G')| \colon \ \Phi_i^{\mathrm{a},1}[\tau_G, \mu_0 \mapsto X \cap |\mathring{\Omega}(G)|] \\
& \qquad \wedge \ldots \\
& \qquad \wedge \Phi_i^{\mathrm{a},t_\mathrm{a}}[\tau_G, \mu_0 \mapsto X \cap |\mathring{\Omega}(G)|] \\
& \qquad \wedge \Phi_i^{\mathrm{b},1}[\tau_{G'}, \mu_0 \mapsto X \cap |\mathring{\Omega}(G')|] \\
& \qquad \wedge \ldots \\
& \qquad \wedge \Phi_i^{\mathrm{b},t_\mathrm{b}}[\tau_{G'}, \mu_0 \mapsto X \cap |\mathring{\Omega}(G')|] \\
\leftrightarrow\ & \exists X_G \subseteq |\mathring{\Omega}(G)|, \exists X_{G'} \subseteq |\mathring{\Omega}(G')| \colon \\
& \qquad \Phi_i^{\mathrm{a},1}[\tau_G, \mu_0 \mapsto X_G] \\
& \qquad \wedge \ldots \\
& \qquad \wedge \Phi_i^{\mathrm{a},t_\mathrm{a}}[\tau_G, \mu_0 \mapsto X_G] \\
& \qquad \wedge \Phi_i^{\mathrm{b},1}[\tau_{G'}, \mu_0 \mapsto X_{G'}] \\
& \qquad \wedge \ldots \\
& \qquad \wedge \Phi_i^{\mathrm{b},t_\mathrm{b}}[\tau_{G'}, \mu_0 \mapsto X_{G'}] \\
\leftrightarrow\ & \exists X_G \subseteq |\mathring{\Omega}(G)| \colon \quad (\Phi_i^{\mathrm{a},1}[\tau_G, \mu_0 \mapsto X_G] \\
& \qquad \wedge \ldots \\
& \qquad \wedge \Phi_i^{\mathrm{a},t_\mathrm{a}}[\tau_G, \mu_0 \mapsto X_G]) \\
& \wedge \exists X_{G'} \subseteq |\mathring{\Omega}(G')| \colon \quad (\Phi_i^{\mathrm{b},1}[\tau_{G'}, \mu_0 \mapsto X_{G'}] \\
& \qquad \wedge \ldots \\
& \qquad \wedge \Phi_i^{\mathrm{b},t_\mathrm{b}}[\tau_{G'}, \mu_0 \mapsto X_{G'}]).
\end{aligned}$$

We set $\varphi_i^\mathrm{a} := \exists \mu_0 (\Phi_i^{\mathrm{a},1} \wedge \ldots \wedge \Phi_i^{\mathrm{a},t_\mathrm{a}})$ and $\varphi_i^\mathrm{b} := \exists \mu_0 (\Phi_i^{\mathrm{b},1} \wedge \ldots \wedge \Phi_i^{\mathrm{b},t_\mathrm{b}})$ and notice that although the height of these formulas is potentially much larger than the height of $\varphi$, their width is exactly the width of $\varphi$. Likewise, the sets of free variables have not grown, and no new variable or predicate symbols have been used, meaning that each of these formulas fulfils the restrictions imposed. We finally set $\Phi := (\varphi_1^\mathrm{a} \wedge \varphi_1^\mathrm{b}) \vee \ldots \vee (\varphi_c^\mathrm{a} \wedge \varphi_c^\mathrm{b})$, yielding

$$\varphi[\tau] \leftrightarrow \Phi\,[\varphi_1^\mathrm{a}[\tau_G], \ldots, \varphi_c^\mathrm{a}[\tau_G], \varphi_1^\mathrm{b}[\tau_{G'}], \ldots, \varphi_c^\mathrm{b}[\tau_{G'}]]\,,$$

as desired.

The claim thus holds for all well-formed formulas, hence in particular for sentences. By the considerations at the beginning of this proof, this implies the statement of the theorem.

□

That $\mathscr{P}^{l,k}_{w}m$ is also inductive with regard to source redefinition should come as no surprise, since the only way in which terminal vertices are even detected in a formula is via the functions $\delta_{\text{term}}$.

**Theorem 6.9.11**

Let $w, l, k \in \mathbb{N}$. Then $\mathscr{P}^{l,k}_{w}$ is $\left(\bigcup_{i\in\mathbb{N}}\bigcup_{j\in\mathbb{N}} {}^{j}_{i}\mathfrak{F}_{\leftrightarrows}\right)$-inductive.

**Proof.** Recall that ${}^{j}_{i}\mathfrak{F}_{\leftrightarrows}$ is the set of all possible source redefinitions between graphs of type $i$ and $j$.

Fix two graph types $n, n' \in \mathbb{N}$ and a function $\sigma\colon \{1, \dots, n\} \to \{1, \dots, n'\}$. For $K \subseteq \{1, \dots, n\}$, set

$$K' := \{\, \sigma(i) : i \in K \,\} \subseteq \{1, \dots, n'\}.$$

Consider now a term of the form $\delta^{K}_{\text{term}}(\chi)$, where $K \subseteq \{1, \dots, n\}$ and $\chi$ is an elementary term, and consider a graph $G = (V, E, ⦉_⦊, t)$ of type $n$ with $G = \leftrightarrows_{\sigma} G'$ for some graph $G' = (V', E', ⦇_⦈, t')$ of type $n'$.

$\leftrightarrows$ denotes the terminal redefinition. (def. 4.5.4, p. 41)

For a full variable assignment $\tau$ for $\delta^{K}_{\text{term}}(\chi)$, we then have

$$\begin{aligned}
\delta^{K}_{\text{term}}(\chi)[\tau] &= \chi[\tau] \cup \{\, t(i) : i \in K \,\} \\
&= \chi[\tau] \cup \{\, t'(\sigma(i)) : i \in K \,\} \\
&= \chi[\tau] \cup \{\, t'(j) : j = \sigma(i) \text{ for some } i \in K \,\} \\
&= \chi[\tau] \cup \{\, t'(i) : i \in K' \,\} \\
&= \delta^{K'}_{\text{term}}(\chi)[\tau].
\end{aligned}$$

Thus by replacing all occurrences of $\delta^{K}_{\text{term}}$ in a sentence $\varphi$ with $\delta^{K'}_{\text{term}}$, we obtain a sentence $\psi$ where all terms evaluate to the same sets on $G'$ as the original terms on $G$.

Since source redefinition does not change the vertices and edges of a graph (only which vertices are terminals with what multiplicity) and the predicates of $\mathring{\mathfrak{L}}$ have no way to detect terminal vertices (that is, they cannot detect the difference between a set $\{v\}$ where $v$ is a terminal vertex and the same set where $v$ is not a terminal vertex), we have

$\mathring{\mathfrak{L}}$ is the circuitous language of graphs. (def. 5.6.1, p. 77)

$$\vDash_{\mathring{\Omega}(G)} \varphi \Leftrightarrow \vDash_{\mathring{\Omega}(\leftrightarrows_\sigma G)} \psi,$$

$\mathring{\Omega}(G)$ is the circuitous induced structure of $G$. (def. 5.6.2, p. 78)

yielding a (trivial) decomposition for $\varphi$.

□

$\leftrightarrows$ denotes the terminal redefinition. (def. 4.5.4, p. 41)

It remains to show that our family plays well with the source fusion.

**Theorem 6.9.12**

Let $w, l, k \in \mathbb{N}$. Then $\mathscr{P}_w^{l;k}$ is $\left(\bigcup_{i\in\mathbb{N}} {}_i\mathfrak{F}_{\text{媲}}\right)$-inductive.

**Proof.** Recall that ${}_i\mathfrak{F}_{\text{媲}}$ is the set of all possible source fusions on graphs of type $i$. We proceed very similarly to the proof of theorem 6.9.10. We fix now until the end of the proof restrictions $w, l, k \in \mathbb{N}$, a graph type $n \in \mathbb{N}$, and two elements $a, b \in \{1, \dots, n\}$.

For any variable assignment $\tau$ on $\mathring{\Omega}(G)$ for some graph $G = (V, E, \langle\!\lvert _ \rvert\!\rangle, t)$ of type $n$ with $G = \text{媲}_a^b G'$ for some graph $G' = (V', E', \langle\!\lvert _ \rvert\!\rangle, t')$ of type $n$, we set

媲 (tsureai, Japanese for *to marry*) denotes the source fusion. (def. 4.5.5, p. 42)

$$\tau_{a,b} : x \mapsto \begin{cases} \tau(x) \subseteq E' = E & \text{if } \tau(x) \subseteq E \\ \tau(x) \subseteq V' & \text{if } \tau(x) \subseteq V \setminus \{t(a)\} \\ \tau(x) \cup \{t'(b)\} \subseteq V' & \text{if } \tau(x) \subseteq V, t(a) \in \tau(x), \end{cases}$$

recalling that due to $G$ being a fusion, $t(a) = t(b)$ while $t'(a)$ and $t'(b)$ need not coincide.

This yields a variable assignment in $G'$. Intuitively, we simply use "the same" assignment, just "expanding" the vertex $t(a)$ into its preimages $t'(a), t'(b)$. That this works might already be apparent to the reader at this point. We prove it nonetheless.

We use the same trick on sets of terminal vertices by setting, for any term $\chi$,

$$\chi_{a,b} := \begin{cases} \chi & \text{if } \chi \text{ elementary} \\ \delta^{K}_{\text{term}}(\xi) & \text{if } \chi = \delta^{K}_{\text{term}}(\xi) \\ & \text{for some } K \subseteq \{1, \dots, n\} \setminus \{a, b\} \\ \delta^{K \cup \{a,b\}}_{\text{term}}(\xi) & \text{if } \chi = \delta^{K}_{\text{term}}(\xi) \\ & \text{for some } K \subseteq \{1, \dots, n\}, a \in K \vee b \in K \end{cases}$$

Given a well-formed formula $\varphi$ with a full variable assignment $\tau$, we construct formulas $\varphi_1, \dots$ and a propositional formula $\Phi$ such that

$$\varphi[\tau] \leftrightarrow \Phi[\varphi_1[\tau_{a,b}], \dots],$$

yielding an inductive decomposition.

We proceed by induction over the structure of $\varphi$.

*Case 1: $\varphi$ is atomic.*

*Case 1.1: $\varphi = \lambda_{\text{sgl}}(\chi)$ for some term $\chi$.* We set

$$\varphi_0 := \lambda_{\text{sgl}}(\chi_{a,b}) \vee \left(\delta^{\{a,b\}}_{\text{term}}(\varnothing) \sqsubseteq \chi_{a,b} \wedge \chi_{a,b} \sqsubseteq \delta^{\{a,b\}}_{\text{term}}(\varnothing)\right).$$

Then $\varphi_0$ has the same width, number of variables, and predicate restrictions as $\varphi$, and they have the same set of free variables, making $\tau_{a,b}$ full. We omit this note in the remaining cases.

We set $\Phi := \varphi_0$.

Now if $\chi_{a,b}[\tau]$ is a singleton that does not contain $t(a)$ (and hence not $t(b)$ either, since $t(a) = t(b)$), then $\chi_{a,b}[\tau_{a,b}] = \chi_{a,b}[\tau]$ and $\varphi_0$ evaluates to $\top$.

If $\chi_{a,b}[\tau]$ is a singleton that *does* contain $t(a)$, then it is equal to $\{t(a)\}$ and we have $\chi_{a,b}[\tau_{a,b}] = \{t'(a), t'(b)\} = \delta^{\{a,b\}}_{\text{term}}(\varnothing)[\tau_{a,b}]$, whence $\varphi_0$ evaluates to $\top$.

If $\chi_{a,b}[\tau]$ is not a singleton in the first place, then its "expanded" version $\chi_{a,b}[\tau_{a,b}]$ can be neither a singleton nor the set $\{a, b\}$, making $\varphi_0$ evaluate to $\bot$.

*Case 1.2: $\varphi = \chi^0 \sqsubseteq \chi^1$ for some terms $\chi^0, \chi^1$.* We set $\varphi_0 := \chi^0_{a,b} \sqsubseteq \chi^1_{a,b}$ and $\Phi := \varphi_0$. By definition, we know that $\chi^0[\tau] \subseteq \chi^1[\tau]$ if and only if $\chi^0_{a,b}[\tau_{a,b}] \subseteq \chi^1_{a,b}[\tau_{a,b}]$.

*Case 1.3:* $\varphi = \lambda^i_{\text{conn}}(\chi^0, \chi^1, \ldots, \chi^i)$ *for some* $i \in \mathbb{N}$ *with* $i < k$ *and terms* $\chi^0, \ldots, \chi^i$. We set

$$\varphi_0 := \lambda^i_{\text{conn}}(\chi^0_{a,b}, \chi^1_{a,b}, \ldots, \chi^i_{a,b})$$

and $\Phi := \varphi_0$. Since fusion changes connections only for the vertices $t'(a)$ and $t'(b)$ and an edge $e$ has $t(a)$ as its $j$-th endpoint in $G$ if and only if it had $t'(a)$ or $t'(b)$ as its $j$-th endpoint in $G'$, we are done.

The claim is thus proven for all atomic formulas.

*Case 2:* $\varphi$ *is built from formulas for which the claim has been proven.*

*Case 2.1:* $\varphi = \neg\psi$. Let $(\Psi, \psi_0, \ldots, \psi_i)$ be the decomposition of $\psi$. We set $\varphi_0 := \psi_0, \ldots, \varphi_i := \psi_i$ and $\Phi := \neg\,\Psi$, trivially fulfilling the requirements.

*Case 2.2:* $\varphi = \psi \wedge \zeta$. Let $\left(\Psi_\psi, \psi_0, \ldots, \psi_i\right)$ be the decomposition of $\psi$, and let $\left(\Psi_\zeta, \zeta_0, \ldots, \zeta_j\right)$ be the decomposition of $\zeta$. We set

$$\varphi_0 := \psi_0,\ \ldots,\ \varphi_i := \psi_i,\ \varphi_{i+1} := \zeta_0,\ \ldots,\ \varphi_{i+j} := \zeta_j$$

and $\Phi := \Psi_\psi \wedge \Psi_\zeta$. This yields a valid decomposition for $\varphi$.

*Case 3:* $\varphi = \exists\mu_0\psi$ *for some variable* $\mu_0$ *and a formula* $\psi$ *for which the claim holds.* We are now in the exact situation of case 3 in the proof of theorem 6.9.10 (page 133), which we shall not reproduce here.

The claim is thus proven for all well-formed formulas, hence in particular for sentences.

□

Theorems 6.9.10 to 6.9.12, combined with the fact that the only remaining symbols in $\mathfrak{F}$, the trivial graphs, are nullary function symbols, yield the inductiveness we desire.

### Corollary 6.9.13

Let $w, l, k \in \mathbb{N}$. Then $\mathscr{P}^{l,k}_w$ is $\mathfrak{F}$-inductive.

# 10. Inheritance

In chapter 9, we shall see that in order to design efficient algorithms to detect graph properties, we need our algebra to be finitely expressible. However, the algebra of all graphs is not so – we shall have to restrict our algorithm to certain subalgebras.

Do we now need to re-prove everything we have just proven, ending up with another ten pages of tedious work?

Luckily, the answer is “no”. We introduce in this section the notion of *inherited* algebra, a type of smaller algebra which we then prove inherits (hence the name) many desirable qualities from its parent algebra; in our case from the algebra of all finite graphs.

**Typing Lessons**

When dealing with inherited signatures, notation becomes tedious quickly.
For convenience, we allow ourselves to drop the typing if it is clear from context, writing $(\mathscr{T}, \mathscr{F}) := (\mathscr{T}, \mathscr{F}, \langle _ \rangle)$. We usually do this in cases where the set $\mathscr{F}$ is simply a set of symbols for functions we already know, such as $\oplus$ or 婚, whose input type (a graph with the correct number of terminal vertices) is obvious.

$\oplus$ is the disjoint sum. (def. 4.5.3, p. 40)

婚 (tsureai, Japanese for *to marry*) denotes the source fusion. (def. 4.5.5, p. 42)

One’s first intuition might be to simply use a *subalgebra*, obtained by

- making the carrier sets of our algebra smaller and
- “forgetting” some function symbols.

The latter may be necessary in order to ensure that the output of all remaining functions lands again in the subalgebra’s (smaller) carrier sets.

### Definition 6.10.1

Let $\mathscr{S} = (\mathscr{T}, \mathscr{F}, \langle _ \rangle)$ be a signature. A *subsignature* of $\mathscr{S}$ is a signature $\mathscr{I} = (\mathscr{U}, \mathscr{G}, \langle _ \rangle)$ with the following properties.

- $\mathscr{U} = \mathscr{T}$.
- $\mathscr{G} \subseteq \mathscr{F}$.
- $\langle _ \rangle = \langle _ \rangle_{|\mathscr{G}}$.

We write $\mathscr{I} \leq \mathscr{S}$.

### Definition 6.10.2

Let $\mathscr{S} = (\mathscr{T}, \mathscr{F}, \langle_\rangle)$ be a signature, and let $\mathscr{A} = (\mathscr{C}, \mathscr{O})$ be an $\mathscr{S}$-algebra. A *subalgebra* of $\mathscr{A}$ is an $\mathscr{I}$-algebra $\mathscr{B} = (\mathscr{D}, \mathscr{Q})$ with the following properties.

- $\mathscr{I} = (\mathscr{U}, \mathscr{G}, \langle_\rangle)$ is a subsignature of $\mathscr{S}$.
- $\forall t \in \mathscr{T} : \mathscr{D}_t \subseteq \mathscr{C}_t$.
- $\forall f \in \mathscr{G} : \langle f \rangle^{\text{in}} = (t_1, \dots, t_n) \Rightarrow \mathscr{Q}_f = \left(\mathscr{O}_f\right)\big|_{(\mathscr{D}_{t_1} \times \dots \times \mathscr{D}_{t_n})}$.

We write $\mathscr{B} \leq \mathscr{A}$.

In other words, a subalgebra has smaller carrier sets, but all the function symbols that it retains are the same functions on those carrier sets. Of course, some function symbols may become "useless" in the smaller algebra, for instance if one of their input carrier sets becomes empty.

Families of predicates restrict to subalgebras in a natural way.

### Definition 6.10.3

Let $\mathscr{S} = (\mathscr{T}, \mathscr{F}, \langle_\rangle)$ be a signature, let $\mathscr{A} = (\mathscr{C}, \mathscr{O})$ be an $\mathscr{S}$-algebra, let $\mathscr{B} = (\mathscr{D}, \mathscr{Q})$ be a subalgebra of $\mathscr{A}$, and let $(\mathscr{P}, \langle_\rangle)$ be a family of predicates on $\mathscr{A}$. We denote by $(\mathscr{P}|_{\mathscr{B}}, \langle_\rangle)$ the family of predicates on $\mathscr{B}$ with

$$\mathscr{P}|_{\mathscr{B}} := \{ p|_{\mathscr{D}_{\langle p \rangle}} : p \in \mathscr{P} \}.$$

With these definitions, we are ready to see that our desired property carries over.

### Lemma 6.10.4

Let $\mathscr{S} = (\mathscr{T}, \mathscr{F}, \langle_\rangle)$ be a signature, let $\mathscr{A} = (\mathscr{C}, \mathscr{O})$ be an $\mathscr{S}$-algebra, let $\mathscr{I} = (\mathscr{U}, \mathscr{G}, \langle_\rangle)$ be a subsignature of $\mathscr{S}$, let $\mathscr{B} = (\mathscr{D}, \mathscr{Q})$ be a $\mathscr{I}$-subalgebra of $\mathscr{A}$, and let $(\mathscr{P}, \langle_\rangle)$ be a family of predicates on $\mathscr{A}$. Let further $\mathscr{H} \subseteq \mathscr{F}$ such that $(\mathscr{P}, \langle_\rangle)$ is $\mathscr{H}$-inductive.

Then $(\mathscr{P}|_{\mathscr{B}}, \langle_\rangle)$ is $\mathscr{H} \cap \mathscr{G}$-inductive.

**Proof.** There is nothing to show – an inductive decomposition for $p \in \mathscr{P}$ is, in particular, an inductive decomposition for $p \in \mathscr{P}|_{\mathscr{B}}$.

□

Sadly, it turns out that this construction is too restrictive. For instance, the algebra of graphs of tree-width at most $k$ (for some $k \in \mathbb{N}$) cannot be obtained as a subalgebra of the algebra of all graphs. However, it still exhibits all the properties we wish for – it is finitely expressible, and the family of predicates we have defined is inductive on it.

We extend the notion of subalgebra by allowing the construction of new functions from existing ones.

First, a quick observation on families of predicates.

### Observation 6.10.5

Let $\mathscr{S} = (\mathscr{T}, \mathscr{F}, \langle_\rangle)$ and $\mathscr{I} = (\mathscr{U}, \mathscr{G}, \langle_\rangle)$ be signatures with $\mathscr{T} \subseteq \mathscr{U}$, let $\mathscr{A} = (\mathscr{C}, \mathscr{O})$ be an $\mathscr{S}$-algebra, and let $\mathscr{B} = (\mathscr{D}, \mathscr{Q})$ be an $\mathscr{I}$-algebra. Let further $(\mathscr{P}, \langle_\rangle)$ be a family of predicates on $\mathscr{A}$.

If for every $t \in \mathscr{T}$ we have $\mathscr{C}_t = \mathscr{D}_t$, then $(\mathscr{P}, \langle_\rangle)$ is also a family of predicates on $\mathscr{B}$.

In particular, as long as we are not talking about inductiveness, the function symbols of an algebra are irrelevant for defining a family of predicates. We now extend a given algebra to include for every element a nullary function which outputs that element.

### Definition 6.10.6

Let $\mathscr{S} = (\mathscr{T}, \mathscr{F}, \langle_\rangle)$ be a signature, and let $\mathscr{A} = (\mathscr{C}, \mathscr{O})$ be an $\mathscr{S}$-algebra.

The *closure* of $\mathscr{S}$ with regard to $\mathscr{A}$ is the signature

$$\overline{\mathscr{S}}^{\mathscr{A}} := \left(\mathscr{T}, \overline{\mathscr{F}}^{\mathscr{A}}, \overline{\langle_\rangle}^{\mathscr{A}}\right)$$

with

$$\overline{\mathscr{F}}^{\mathscr{A}} := \mathscr{F} \cup \{\, f^t_c : t \in \mathscr{T}, c \in \mathscr{C}_t \,\}$$

and

$$\forall f \in \overline{\mathcal{F}}^{\mathcal{A}} : \langle f \rangle := \begin{cases} \langle f \rangle & \text{if } f \in \mathcal{F} \\ (\varepsilon, t) & \text{if } f = f_c^t \text{ for some } t \in \mathcal{T}, c \in \mathcal{C}_t. \end{cases}$$

The *closure* of $\mathcal{A}$ is the $\overline{\mathcal{F}}^{\mathcal{A}}$-algebra $\overline{\mathcal{A}} := \left(\mathcal{C}, \overline{\mathcal{O}}\right)$ with

$$\forall f \in \overline{\mathcal{F}}^{\mathcal{A}} : \overline{\mathcal{O}}_f := \begin{cases} \mathcal{O}_f & \text{if } f \in \mathcal{F} \\ () \mapsto c & \text{if } f = f_c^t \text{ for some } t \in \mathcal{T}, c \in \mathcal{C}_t. \end{cases}$$

We can observe several features right away.

## Observation 6.10.7

Let $\mathcal{A}$ be an algebra. Then the following statements are true.

- $\overline{\overline{\mathcal{A}}} = \overline{\mathcal{A}}$.
- The closure of $\mathcal{A}$ is finitely expressible.
- Every family of predicates on $\mathcal{A}$ is also a family of predicates on $\overline{\mathcal{A}}$.
- Let $\mathcal{H} \subseteq \mathcal{F}$, and let $(\mathcal{P}, \langle _ \rangle)$ be an $\mathcal{H}$-inductive family of predicates on $\mathcal{A}$. Then $(\mathcal{P}, \langle _ \rangle)$ is $\left(\mathcal{H} \cup \overline{\mathcal{F}}^{\mathcal{A}} \setminus \mathcal{F}\right)$-inductive on $\overline{\mathcal{A}}$.

The final property is true because all function symbols in $\overline{\mathcal{F}}^{\mathcal{A}} \setminus \mathcal{F}$ are nullary, yielding a trivial decomposition for any predicate – for a nullary function symbol $f$, we either have $p(\mathcal{O}_f()) = \top$ or $p(\mathcal{O}_f()) = \bot$.

$|\mathcal{A}|$ denotes the set of pre-expressions over $\mathcal{A}$. (def. 6.3.1, p. 98)

We now forget about this construction for a second and look at a different way to extend a given algebra: given a pre-expression $e \in |\mathcal{A}|$, this pre-expression has a type $\langle e \rangle = (t_1 \dots t_n, t)$, and given elements of the correct types, we can plug them into $e$ and obtain a new element of type $t$. In other words, $e$ behaves exactly like a function symbol! It is only natural to extend our signature to include function symbols that arise as pre-expressions.

**Definition 6.10.8**

Let $\mathscr{S} = (\mathscr{T}, \mathscr{F}, \langle_\rangle)$ be a signature. The *flattening* of $\mathscr{S}$ is the signature

$$\mathrm{Fl}(\mathscr{S}) = (\mathscr{T}, \mathrm{Fl}(\mathscr{F}), \mathrm{Fl}\langle_\rangle)$$

with

$$\mathrm{Fl}(\mathscr{F}) := \mathscr{F} \cup \{\, f_e : e \in |\mathscr{S}| \,\}$$

$|\mathscr{S}|$ denotes the set of pre-expressions over $\mathscr{S}$. (def. 6.3.1, p. 98)

and

$$\forall f \in \mathrm{Fl}(\mathscr{F}) \colon \mathrm{Fl}\langle f\rangle := \begin{cases} \langle f\rangle & \text{if } f \in \mathscr{F} \\ \langle e\rangle & \text{if } f = f_e \text{ for some } e \in |\mathscr{S}|. \end{cases}$$

Let now $\mathscr{A}$ be an $\mathscr{S}$-algebra. The *flattening* of $\mathscr{A}$ is the $\mathrm{Fl}(\mathscr{S})$-algebra $\mathrm{Fl}(\mathscr{A}) = (\mathscr{C}, \mathrm{Fl}(\mathscr{O}))$ with $\forall f \in \mathrm{Fl}(\mathscr{F})$:

$$\mathrm{Fl}(\mathscr{O})_f := \begin{cases} \mathscr{O}_f & \text{if } f \in \mathscr{F} \\ (c_1, \dots, c_{|\langle e\rangle^{\mathrm{in}}|}) \mapsto e_{\mathscr{A}}(c_1, \dots, c_{|\langle e\rangle^{\mathrm{in}}|}) & \text{if } f = f_e \text{ for some } e \in |\mathscr{S}|. \end{cases}$$

$e_{\mathscr{A}}(\dots)$ denotes the result of $e$ when evaluated in $\mathscr{A}$. (def. 6.3.2, p. 101)

Another way to look at this construction is to think of the composition of unary functions $g \circ f$, generalised to functions of arbitrary arity. The reason this process is called flattening is that every computation that can be carried out in $\mathscr{A}$ by evaluating an expression tree can be carried out by a single function call in $\mathrm{Fl}(\mathscr{A})$.

While repeated closure is idempotent, repeated flattening always adds new function symbols. Still, it does not make the algebra any “flatter”, as formalised in the following lemma.

**Lemma 6.10.9**

Let $\mathscr{S} = (\mathscr{T}, \mathscr{F}, \langle_\rangle)$ be a signature, and let $\mathscr{A} = (\mathscr{C}, \mathscr{O})$ be an $\mathscr{S}$-algebra. Let further $e \in |\mathrm{Fl}(\mathscr{S})|$ with $\langle e\rangle^{\mathrm{in}} = (t_1, \dots, t_n)$ for some $n \in \mathbb{N}$. Then there is a function symbol $f \in \mathrm{Fl}(\mathscr{F})$ with $\mathrm{Fl}\langle f\rangle = \langle e\rangle$ such that

$$\forall c_1 \in \mathscr{C}_{t_1} \dots \forall c_n \in \mathscr{C}_{t_n} : \mathrm{Fl}(\mathscr{O})_f(c_1, \dots, c_n) = e_{\mathrm{Fl}(\mathscr{A})}(c_1, \dots, c_n).$$

**Proof.** Let $\mathscr{S} = (\mathscr{T}, \mathscr{F}, \langle_\rangle)$ be a signature, and let $\mathscr{A} = (\mathscr{C}, \mathscr{O})$ be an $\mathscr{S}$-algebra. Take a pre-expression

$\mathrm{Fl}(\mathscr{S})$ denotes the flattening of $\mathscr{S}$. (def. 6.10.8, p. 143)

$$e = (V, E, \llparenthesis_\rrparenthesis, \star, \preccurlyeq, \langle_\rangle) \in |\mathrm{Fl}(\mathscr{S})|$$

$|\mathrm{Fl}(\mathscr{S})|$ denotes the set of pre-expressions over $\mathrm{Fl}(\mathscr{S})$. (def. 6.3.1, p. 98)

with $\langle e\rangle = (t_1 \dots t_n, t)$ for some $n \in \mathbb{N}$.

Because this is a pre-expression in $\mathrm{Fl}(\mathscr{S})$, every vertex label $\langle v\rangle, v \in V$ arises in turn from a pre-expression in $\mathscr{S}$ of the same type. Stitching all of these pre-expressions together yields a pre-expression $e' \in |\mathscr{S}|$ with $\langle e'\rangle = \langle e\rangle$ and

$e'_{\mathscr{A}}(\dots)$ denotes the result of $e'$ when evaluated in $\mathscr{A}$. (def. 6.3.2, p. 101)

$$\forall c_1 \in \mathscr{C}_{t_1} \dots \forall c_n \in \mathscr{C}_{t_n} : e'_{\mathscr{A}}(c_1, \dots, c_n) = e_{\mathrm{Fl}(\mathscr{A})}(c_1, \dots, c_n).$$

The following sketch illustrates this point. Input types are indicated at the dashed lines. Primed function symbols are in $\mathscr{F}$, non-primed ones are in $\mathrm{Fl}(\mathscr{F})$. We see first the pre-expression $e \in |\mathrm{Fl}(\mathscr{S})|$, then the transformation into the expression $e' \in |\mathscr{S}|$.

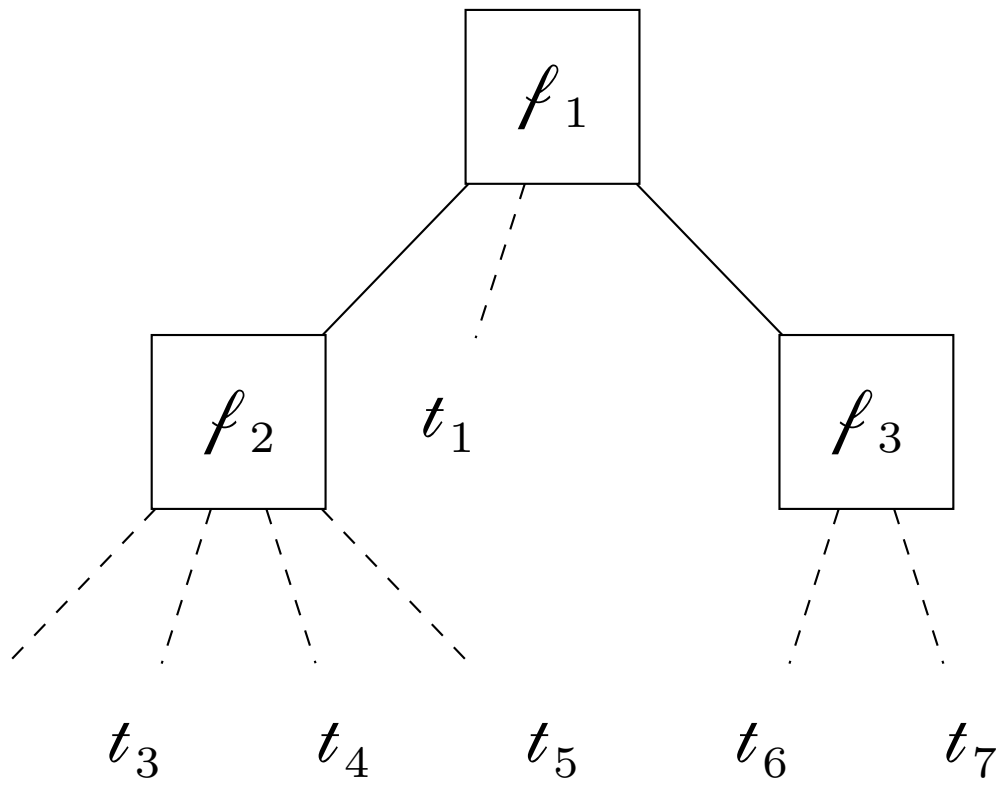

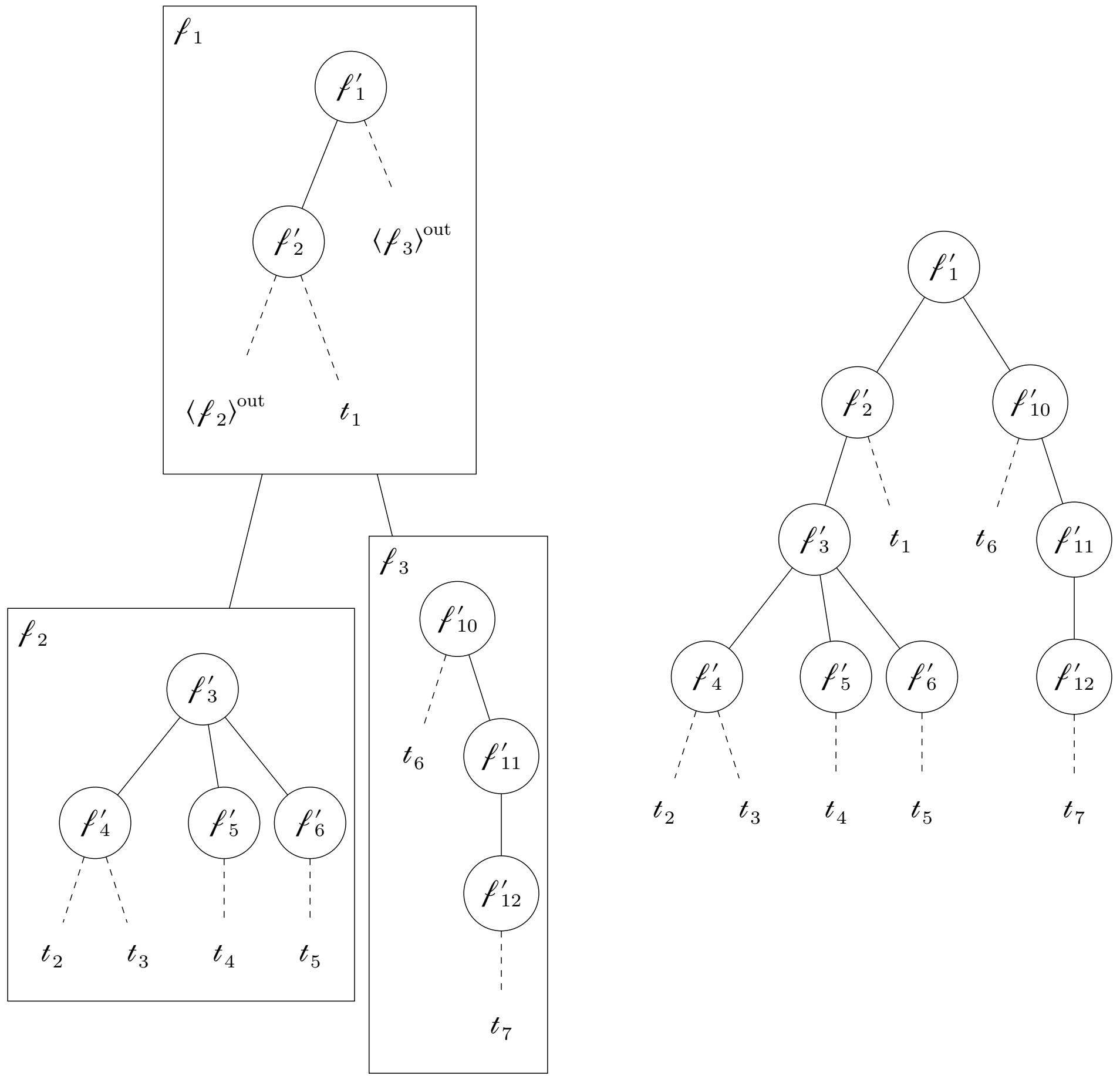

Because $\mathrm{Fl}(\mathscr{A})$ is the flattening of $\mathscr{A}$, by definition there must be a function symbol $\mathscr{g} \in \mathrm{Fl}(\mathscr{F})$ with $\mathrm{Fl}\langle \mathscr{g} \rangle = \langle e' \rangle$

$\mathrm{Fl}(\mathscr{A})$ denotes the flattening of $\mathscr{A}$. (def. 6.10.8, p. 143)

$$\begin{aligned}\forall c_1 \in \mathscr{C}_{t_1} \dots \forall c_n \in \mathscr{C}_{t_n} : \mathrm{Fl}(\mathscr{O})_{\mathscr{g}}(c_1, \dots, c_n) &= e'_{\mathscr{A}}(c_1, \dots, c_n) \\ &= e_{\mathrm{Fl}(\mathscr{A})}(c_1, \dots, c_n),\end{aligned}$$

$e'_{\mathscr{A}}(\dots)$ denotes the result of $e'$ when evaluated in $\mathscr{A}$. (def. 6.3.2, p. 101)

which was the claim.

□

We note some further properties of the flattening for future reference.

### Observation 6.10.10

Let $\mathscr{S} = (\mathscr{T}, \mathscr{F}, \langle_\rangle)$ be a signature, and let $\mathscr{A}$ be an $\mathscr{S}$-algebra. Then

$\mathrm{Fl}(\mathscr{A})$ denotes the flattening of $\mathscr{A}$. (def. 6.10.8, p. 143)

the following are true.

- $\mathrm{Fl}(\mathscr{A})$ is expressible if and only if $\mathscr{A}$ is expressible.
- If $\mathscr{A}$ is finitely expressible, then so is $\mathrm{Fl}(\mathscr{A})$. The converse is not true.
- Let $\mathscr{H} \subseteq \mathscr{F}$, and let $(\mathscr{P}, \langle_\rangle)$ be an $\mathscr{H}$-inductive family of predicates on $\mathscr{A}$. Then $(\mathscr{P}, \langle_\rangle)$ is also $\mathscr{H}$-inductive on $\mathrm{Fl}(\mathscr{A})$.

Of course, what we really want to know is the following stronger result.

### Lemma 6.10.11

Let $\mathscr{S} = (\mathscr{T}, \mathscr{F}, \langle_\rangle)$ be a signature, let $\mathscr{A}$ be an $\mathscr{S}$-algebra, and let $(\mathscr{P}, \langle_\rangle)$ be an $\mathscr{F}$-inductive family of predicates on $\mathscr{A}$. Then $(\mathscr{P}, \langle_\rangle)$ is $\mathrm{Fl}(\mathscr{F})$-inductive on $\mathrm{Fl}(\mathscr{A})$.

**Proof.** Let $\mathscr{S} = (\mathscr{T}, \mathscr{F}, \langle_\rangle)$ be a signature, let $\mathscr{A} = (\mathscr{C}, \mathscr{O})$ be an $\mathscr{S}$-algebra, and let $(\mathscr{P}, \langle_\rangle)$ be an $\mathscr{F}$-inductive family of predicates on $\mathscr{A}$. Let $p \in \mathscr{P}$, and let $f \in \mathrm{Fl}(\mathscr{F})$ with $\mathrm{Fl}\langle f\rangle^{\mathrm{in}} = (t_1, \dots, t_n)$ for some $n \in \mathbb{N}$ and with $\mathrm{Fl}\langle f\rangle^{\mathrm{out}} = \langle p\rangle$. Without loss of generality, we have $f \notin \mathscr{F}$, since otherwise an inductive decomposition already exists by assumption.

$|\mathscr{S}|$ denotes the set of pre-expressions over $\mathscr{S}$. (def. 6.3.1, p. 98)

Hence, we find a pre-expression $e \in |\mathscr{S}|$ such that $\langle e\rangle = \mathrm{Fl}\langle f\rangle$ and

$$\forall c_1 \in \mathscr{C}_{t_1} \dots \forall c_n \in \mathscr{C}_{t_n} : \mathrm{Fl}(\mathscr{O})_f(c_1, \dots, c_n) = e_{\mathscr{A}}(c_1, \dots, c_n).$$

$e_{\mathscr{A}}(\dots)$ denotes the result of $e$ when evaluated in $\mathscr{A}$. (def. 6.3.2, p. 101)

For every function symbol $g$ labelling the vertices of $e$, there is by assumption an inductive decomposition for every $p' \in \mathscr{P}$ with regard to $g$. Therefore, we can find an inductive decomposition for $p$ with regard to $f$ by induction over the vertices of $e$.

□

We are now ready to define inheritance algebras. Intuitively, we want to do the following to our original algebra:

- Make our carrier sets smaller, if desired.
- Introduce any number of new nullary functions, that is, constants (in our new, smaller carrier sets).

- Take any number of pre-expressions in the original algebra, possibly with new nullary functions included, and make them into functions in the canonical way.
- Forget any number of function symbols.
- Restrict the remaining function symbols to our new carrier sets.

The reader should take a moment to convince themselves that the following definition achieves all of these goals.

**Definition 6.10.12**

Let $\mathscr{A}$ be an algebra. We say that an algebra $\mathscr{B}$ is *inherited* from $\mathscr{A}$ if it is a subalgebra of the flattening of the closure of $\mathscr{A}$, that is, if

$$\mathscr{B} \leq \mathrm{Fl}(\overline{\mathscr{A}}).$$

$\overline{\mathscr{A}}$ denotes the closure of $\mathscr{A}$. (def. 6.10.6, p. 141)

Useful graph classes such as the graphs of tree-width bounded by a constant $k \in \mathbb{N}$ will (in chapter 9) turn out to be inheritance algebras of $\mathfrak{G}$, so we want them to retain our favourite property. Of course, we have crafted our definitions in just such a way that the following works out nicely.

$\mathrm{Fl}(\overline{\mathscr{A}})$ denotes the flattening of $\overline{\mathscr{A}}$. (def. 6.10.8, p. 143)

$\mathfrak{G}$ is the algebra of graphs. (def. 6.2.2, p. 93)

**Theorem 6.10.13**

Let $\mathscr{S} = (\mathscr{T}, \mathscr{F}, \langle_\rangle)$ be a signature, let $\mathscr{A}$ be an $\mathscr{S}$-algebra, and let $(\mathscr{P}, \langle_\rangle)$ be an $\mathscr{F}$-inductive family of predicates on $\mathscr{A}$. Let furthermore $\mathscr{I} = (\mathscr{T}, \mathscr{G}, \langle_\rangle)$ be another signature and let $\mathscr{B}$ be an $\mathscr{I}$-algebra inherited from $\mathscr{A}$.[11] Then $\left(\mathscr{P}|_{\mathscr{B}}, \langle_\rangle\right)$ is $\mathscr{G}$-inductive on $\mathscr{B}$.

**Proof.** Let $\mathscr{S} = (\mathscr{T}, \mathscr{F}, \langle_\rangle)$ be a signature, let $\mathscr{A}$ be an $\mathscr{S}$-algebra, and let $(\mathscr{P}, \langle_\rangle)$ be an $\mathscr{F}$-inductive family of predicates on $\mathscr{A}$. Let now $\mathscr{I} = (\mathscr{U}, \mathscr{G}, \langle_\rangle)$ be another signature and let $\mathscr{B}$ be a $\mathscr{I}$-algebra inherited from $\mathscr{A}$, that is, a subalgebra of $\mathrm{Fl}(\overline{\mathscr{A}})$.

We know from our earlier observations that $(\mathscr{P}, \langle_\rangle)$ is a family of predicates on $\overline{\mathscr{A}}$ and hence on $\mathrm{Fl}(\overline{\mathscr{A}})$, and we know from observation 6.10.7 and lemma 6.10.11 that it is $\mathrm{Fl}(\overline{\mathscr{F}}^{\mathscr{A}})$-inductive on $\mathrm{Fl}(\overline{\mathscr{A}})$.

$\overline{\mathscr{F}}^{\mathscr{A}}$ denotes the closure of $\mathscr{F}$ with regard to $\mathscr{A}$. (def. 6.10.6, p. 141)

[11] In particular, we have $\mathscr{U} = \mathscr{T}$ and $\mathscr{G} \subseteq \mathrm{Fl}(\overline{\mathscr{F}}^{\mathscr{A}})$.

The result now follows from lemma 6.10.4.

□

Armed with all this knowledge, we are ready to prove Courcelle's Theorem.

# Chapter 7
## Courcelle's Theorem

We have finally accrued all the machinery necessary to prove Courcelle's Theorem. In chapter 9, we show explicitly some special cases which easily follow from our inheritance technique in a way that has not previously been possible.

## 1. Proof

We remind ourselves what it is that we are trying to prove.

**Theorem 7.1.1**

Let $\varphi$ be a sentence of the monadic second-order language of graphs. Then for every $n \in \mathbb{N}$, the set $\{ G \in \mathfrak{G}_n : \models_{\Omega(G)} \varphi \}$ is $\mathfrak{G}$-recognisable.

$\mathfrak{G}_n$ denotes the set of all graphs of type $n$. (def. 4.5.1, p. 39)

$\mathfrak{G}$ is the algebra of graphs. (def. 6.2.2, p. 93)

$\Omega(G)$ is the induced structure of $G$. (def. 5.4.2, p. 73)

The corollary for the practical-minded graph theorist, that there exists for certain classes of graphs a linear-time algorithm to detect members of this subset, is proven in chapter 9. The extension to "counting" monadic second-order logic is proven in appendix A.

We prove, in fact, the following stronger result.

**Theorem 7.1.2**

Let $\mathscr{A} = (\mathscr{C}, \mathscr{O})$ be an algebra inherited from $\mathfrak{G}$, and let $\varphi$ be a sentence of the monadic second-order language of graphs. Then for every $n \in \mathbb{N}$, the set $\{ G \in \mathscr{C}_n : \models_{\Omega(G)} \varphi \}$ is $\mathscr{A}$-recognisable.

Theorem 7.1.1 follows immediately from theorem 7.1.2 since every algebra inherits from itself.

$\mathfrak{X}$ denotes the direct logical framework of graphs. (def. 5.4.3, p. 73)

**Proof of theorem 7.1.2.** Let $\varphi$ be a sentence of $\mathfrak{X}$, and let $\mathscr{A} = (\mathscr{C}, \mathscr{O})$ be a $(\mathscr{U}, \mathscr{G}, \langle_\rangle)$-algebra inherited from $\mathfrak{G}$.

$\mathfrak{G}$ is the algebra of graphs. (def. 6.2.2, p. 93)

By corollary 5.6.7, we can find a sentence $\mathring{\varphi} \in \mathring{\mathfrak{X}}$ such that

$$\forall n \in \mathbb{N} \colon \forall G \in \mathscr{C}_n \colon \vDash_{\Omega(G)} \varphi \Leftrightarrow \vDash_{\mathring{\Omega}(G)} \mathring{\varphi}\,.$$

$\mathring{\mathfrak{X}}$ denotes the circuitous logical framework of graphs. (def. 5.6.3, p. 79)

This sentence has some width (possibly larger than the width of $\varphi$). Call it $w$. Say further that the highest-indexed variable symbol occurring in $\mathring{\varphi}$ is $\mu_l$ for some $l \in \mathbb{N}$, and that $k \in \mathbb{N}$ is a number such that no predicate symbol $\lambda^t_{\text{conn}}$ occurs in $\mathring{\varphi}$ for $t > k$.

$\Omega(G)$ is the induced structure of $G$. (def. 5.4.2, p. 73)

Let now $n \in \mathbb{N}$.

We set $\mathscr{p} := \mathscr{p}^n_{\mathring{\varphi}} \in \mathscr{P}^n_{w,l,k}$, which we consider as a predicate of the family

$\mathring{\Omega}(G)$ is the circuitous induced structure of $G$. (def. 5.6.2, p. 78)

$$\mathscr{P} := \left( \left( \{\, \mathscr{P}^n_{w,l,k} \,\}_{n\in\mathbb{N}} \right)_{|\mathscr{A}}, \mathscr{p}^n_{\psi} \mapsto n \right)$$

of predicates on $\mathscr{A}$.

We thus have

$$\{\, G \in \mathscr{C}_n : \vDash_{\Omega(G)} \varphi \,\} = \{\, G \in \mathscr{C}_n : \vDash_{\mathring{\Omega}(G)} \mathring{\varphi} \,\} = \lceil \mathscr{p} \rceil.$$

$\lceil \mathscr{p} \rceil$ denotes the set of elements which satisfy $\mathscr{p}$. (def. 6.7.3, p. 111)

By lemma 6.9.9, $\left( \{\, \mathscr{P}^n_{w,l,k} \,\}_{n\in\mathbb{N}}, \mathscr{p}^n_{\psi} \mapsto n \right)$ is locally finite. Since the underlying sets of predicates are the same, so is $\mathscr{P}$.

By corollary 6.9.13, $\left( \{\, \mathscr{P}^n_{w,l,k} \,\}_{n\in\mathbb{N}}, \mathscr{p}^n_{\psi} \mapsto n \right)$ is $\mathfrak{F}$-inductive, and thus (by theorem 6.10.13) $\mathscr{P}$ is $\mathscr{G}$-inductive.

Hence by theorem 6.8.2, the set $\lceil \mathscr{p} \rceil$ is $\mathscr{A}$-recognisable.

□

# 2. A Practical Theorem?

Now that we have proven that graphs modelling a monadic second-order formula are recognisable in $\mathfrak{G}$, can we turn this knowledge into an algorithm?

In other words, can we explicitly compute the algebra morphism whose preimage is our desired set?

Sadly, many of the graph properties expressible in monadic second-order logic are well-known to be **NP**-hard, so even if we could compute said morphism, we could not expect a polynomial-time algorithm.

Does this mean that all of our hard work has been in vain?

Of course not. In the next chapter, we shall show that if the algebra under consideration is *finitely expressible* (definition 6.3.7), then there is a linear-time algorithm for deciding whether an element is in a given recognisable set. In chapter 9, we then show that classes of graphs bounded by constant tree-width and classes of graphs bounded by constant path-width are finitely expressible and give the reader the requisite tools to build their own finitely expressible classes of graphs.

# Chapter 8

# Tree Automata

In order to prove the existence of an algorithm that decides whether a given graph satisfies some monadic second-order formula, we shall construct a so-called deterministic bottom-up finite tree automaton. The reader familiar with tree automata may skip right to section 8.6, where we show how to construct a deterministic bottom-up finite tree automaton for a given recognisable set in a typed algebra.

## 1. Automayton, Automahton

Intuitively, a *finite-state automaton* is a machine that starts in a certain “state” and is fed a finite string of symbols of some alphabet. For each symbol, it transitions to a new state, depending only on its current state and the symbol currently being digested. Once the entire string is digested, the machine halts in whatever state it ended up in.

A certain subset of the machine’s possible states is marked as “accepting”. If the state in which the machine halts is such a state, the string is “accepted”. Otherwise, it is “rejected”.

A very simple example shall illustrate the idea. Suppose we are given the alphabet $\Sigma = \{0, 1\}$ for some binary encoding, and we want for some reason to determine whether the number of times the symbol 0 occurs in a string is divisible by 3. An automaton that accomplishes this might look as follows.

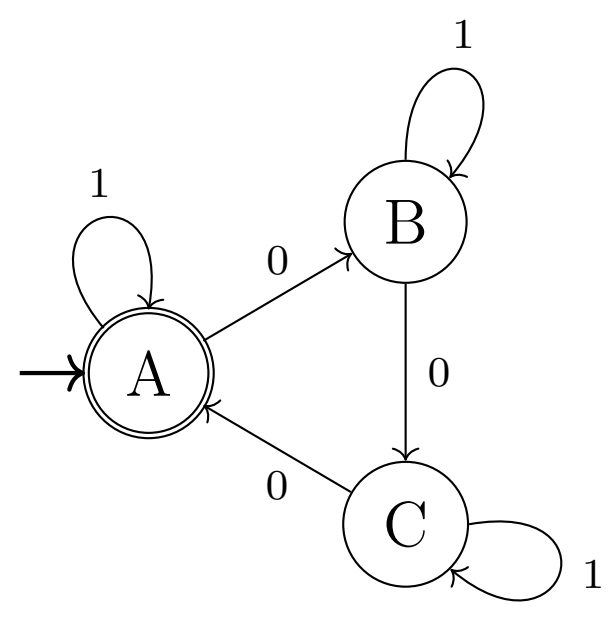


The automaton starts in state $A$ (indicated by the bold arrow) and processes one symbol at a time. State transitions are indicated by edges. The only accepting state is $A$, indicated by the double circle.

Whenever the automaton encounters the symbol 1, its state does not change since the number of ones is irrelevant. Hence all edges labeled 1 are loops.

In its initial state, if the string ends after only ones, the automaton has encountered zero zeroes and accepts. If it encounters a zero, it has now encountered exactly one (which is not divisible by 3) and switches to state $B$, which is not an accepting state. If in state $B$ it encounters a second zero, it switches to state $C$. Two is not divisible by 3, so this is not an accepting state either.

If in state $C$, our automaton encounters another zero, it has now encountered three of them and switches back to state $A$ – if the string ends here, it accepts since three is divisible by 3.

If the string does not end here, the cycle begins anew.

In this manner, the automaton processes any string of zeroes and ones symbol by symbol and finally tells us whether the number of zeroes was divisible by three.

The above is what is called a *deterministic finite-state automaton*, and it is not quite what we are interested in. Rather than a linear string of symbols, the automaton we want to build should process a *tree* of symbols – in the end, we want to feed it an expression over some graph algebra.

We introduce the basic structure our automaton expects before defining the automaton itself.

# 2. Ranked Strings

Consider the alphabet $\Sigma = \{\alpha, \beta, \gamma\}$. Possible words in this alphabet look like, for example, $\alpha\beta\gamma$, $\alpha\alpha\beta\alpha$, or $\alpha\beta\alpha\gamma\alpha$. From the perspective of a finite-state automaton, these strings come as a chain of symbols:

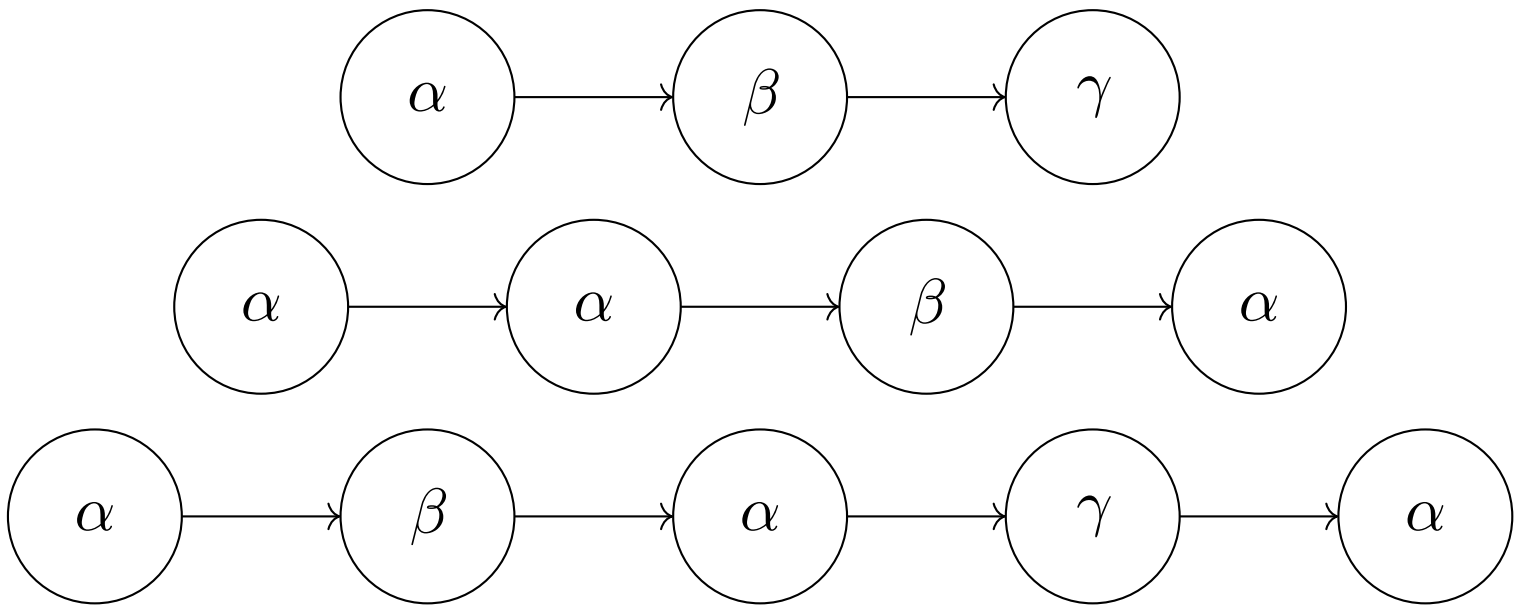


The strings which we need to process, however, arise from repeated application of not necessarily unary functions (think of the disjoint sum), that is, they look something like $f(g(\alpha, \beta), \alpha, f(\alpha, \beta, \gamma))$ for some 3-place function $f$ and some 2-place function $g$. This structure is more intuitively represented as a tree, like we did for expressions already:

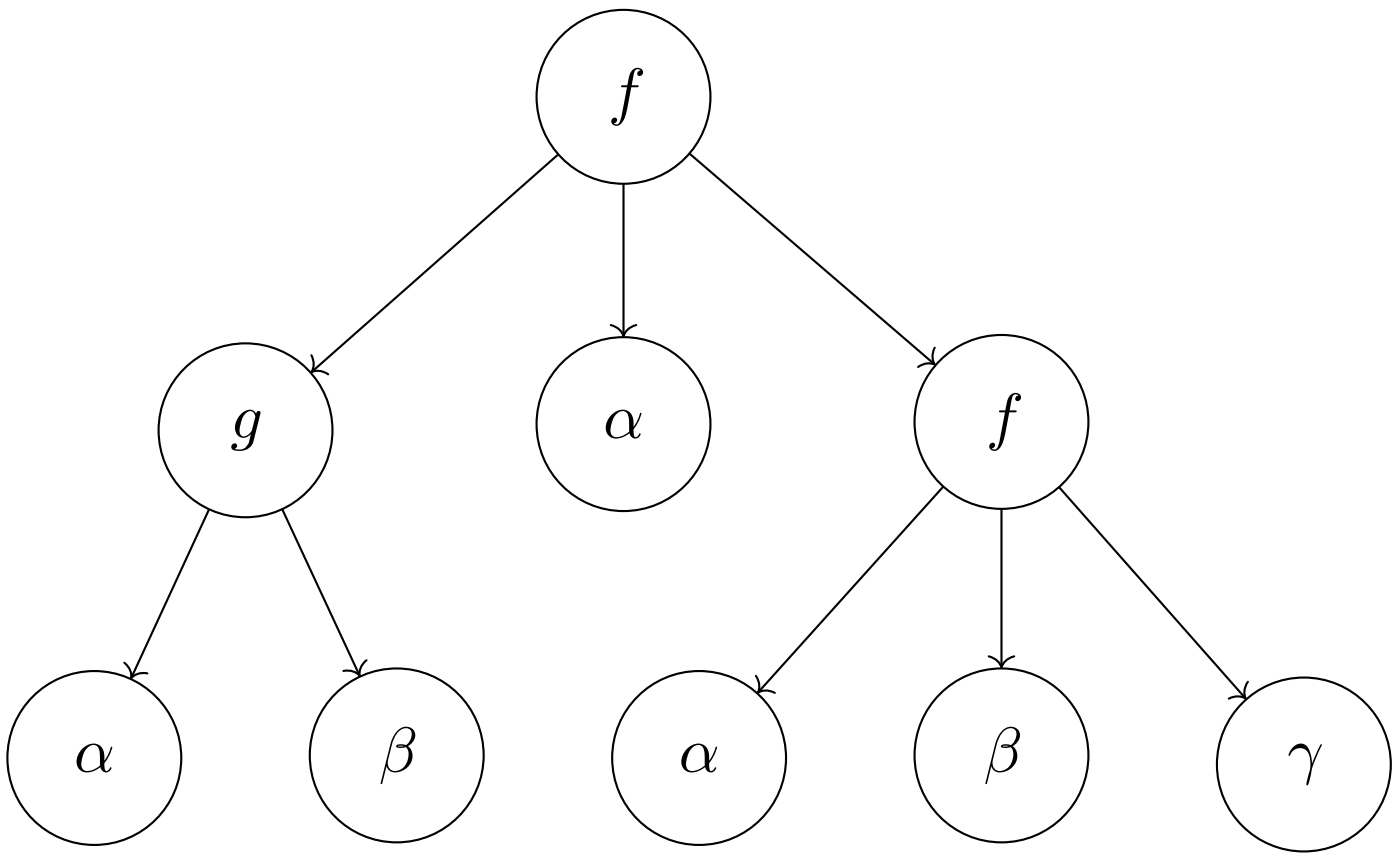


We introduce a variation on alphabets that supports this construction, independently of whether we are actually in an algebra. Expressions will later turn out to be a special case of this.

**Definition 8.2.1**

A *ranked alphabet* is a pair $(\Sigma, |_|)$ fulfilling the following conditions.

- $\Sigma$ is an alphabet.
- $|_|: \Sigma \to \mathbb{N}$ is a function, called the alphabet's *arity function*.

A ranked alphabet $(\Sigma, |_|)$ is called *finite* if $\Sigma$ is.

The arity function tells us how many successors a symbol should have in the string's tree representation. In the example above, $\alpha, \beta, \gamma$ are of arity 0, $f$ is of arity 3, and $g$ is of arity 2.

Nullary elements again correspond to constants. They are the leaves of our tree, like the symbols $\alpha, \beta, \gamma$ in the example above.

Terms built only with unary symbols and constants are the same as strings over a non-ranked alphabet if one interprets the constant as an end-of-string symbol:

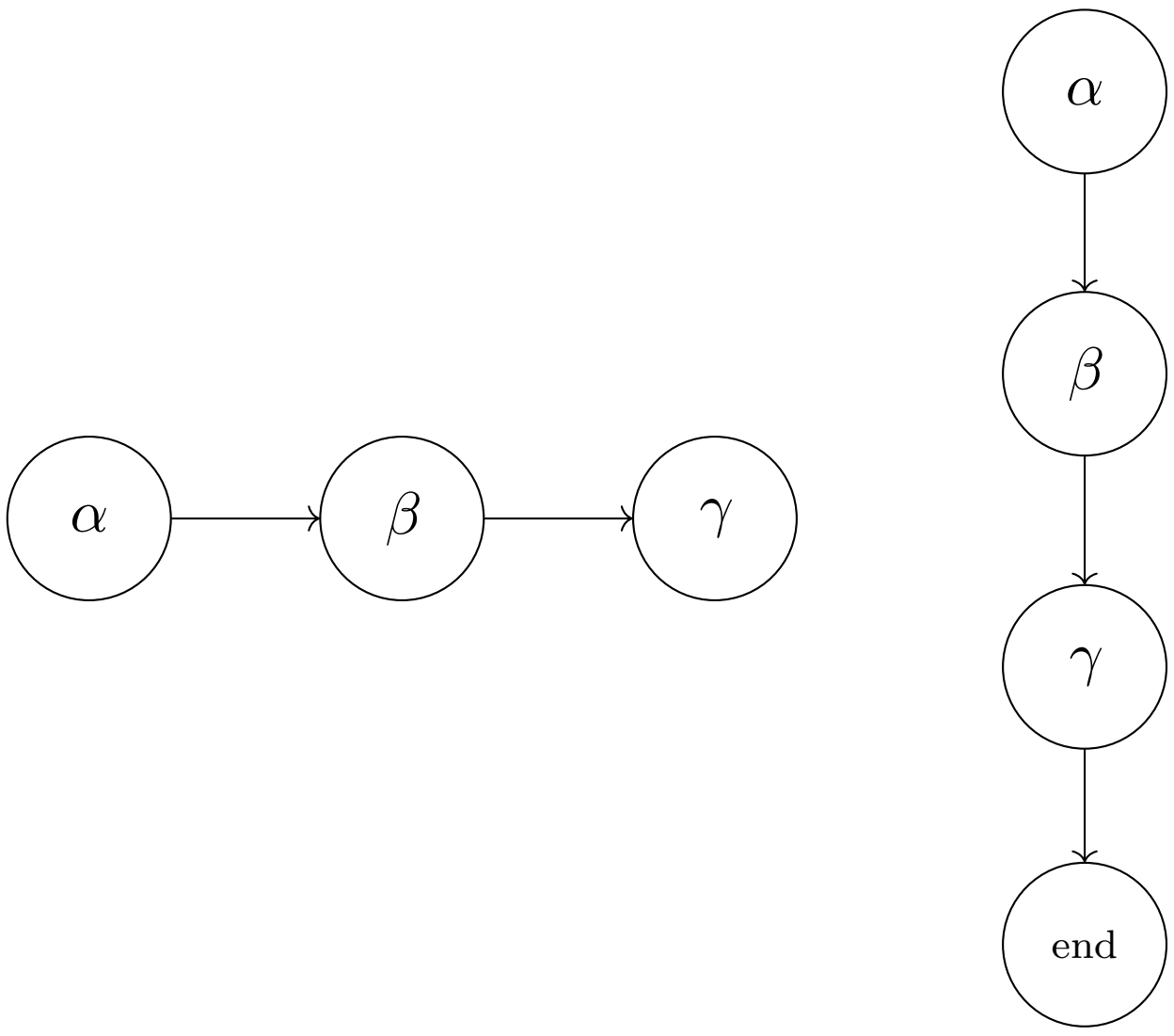


Terms built with symbols of higher arity lead to a tree structure as seen before.

**Definition 8.2.2**

Let $\Pi = (\Sigma, |_|)$ be a ranked alphabet. A *ranked string* in $\Pi$ is a pair $(G, \langle_\rangle)$ fulfilling the following conditions.

- $G = (V, E, \langle\!|_|\!\rangle, \star, \preccurlyeq)$ is a traversal tree.
- $\langle_\rangle\colon V \to \Sigma$ is a function.
- $\forall v \in V\colon \deg^{\mathrm{out}} v = |\langle v\rangle|$.

The set of all ranked strings over $\Pi$ is denoted $||\Pi||$.

We say that two strings $(G, \langle_\rangle)$ and $(G', \langle_\rangle)$ are *isomorphic* if $G$ and $G'$ are isomorphic, say via $(g, h)\colon G \to G'$, and for every vertex $v$ of $G$ we have $\langle gv\rangle = \langle v\rangle$.

Note that since our trees are finite, not every alphabet admits a string – an alphabet without nullary symbols admits no leaves, so it cannot be used to build a nonempty finite tree.[1] Note also that the empty word is *not* a ranked string since we do not allow graphs to be empty.

Note finally that since strings are traversal trees (definition 4.4.5), the isomorphism mentioned is an isomorphism of traversal trees and hence must preserve the ordering of child vertices.

In order to restrict ourselves to those strings that "make sense" (for example, respect type considerations in an algebra), we need to choose a suitable subset of all possible strings. The following definition then ensures that all of our constructions still work.

**Definition 8.2.3**

Let $\Pi$ be a ranked alphabet. A *$\Pi$-vocabulary* is a nonempty set $\Upsilon \subseteq ||\Pi||$ such that

$$\forall (G, \langle_\rangle) \in \Upsilon\colon \forall v \in G\colon (G[v], \langle_\rangle|_{G[v]}) \in \Upsilon.$$

In other words, a vocabulary is a subset of strings that is closed with regard to substrings. This will be helpful when we want to recurse over the structure of a string, as we need not check whether a given substring is

[1] In particular, pre-expressions which expect input are never ranked strings.

contained in our vocabulary – it always is.

Observe that just because a ranked alphabet is finite does not mean that a vocabulary in said alphabet admits only finitely many strings.

# 3. Automata

We know now the structures that will form the input for our automaton. It remains to define how an automaton processes them.

**Definition 8.3.1**

𐤌 is the Phoenician letter "mem".

𐤍 is the Phoenician letter "nun".

𐤎 is the Phoenician letter "semk".

A *deterministic bottom-up finite tree automaton* is a tuple (𐤌, $\Pi$, $\Upsilon$, 𐤍, 𐤎) consisting of the following data:

- a nonempty finite set 𐤌 of *states*,
- a finite ranked alphabet $\Pi = (\Sigma, |_|)$,
- a $\Pi$-vocabulary $\Upsilon$,
- a set 𐤍 ⊆ 𐤌 of *accepting states*,
- and a *transition function* 𐤎: $\Sigma \times$ 𐤌$^* \to$ 𐤌.

Technically it is sufficient if 𐤎 is defined on a subset $A \subseteq (\Sigma \times$ 𐤌$^*)$ that contains all tuples of the form $(f, w)$ with $|w| = |f|$.

𐤐 is the Phoenician letter "pe".

Intuitively, a deterministic bottom-up finite tree automaton starts eating a ranked string from the leaves. Each leaf is assigned a state 𐤐 ∈ 𐤌 by the transition function. The automaton then moves up one level. For a vertex $f$ with successors $a, b, c$, the transition function looks at the states of $a, b, c$ and the symbol $f \in \Sigma$ and assigns on this basis a new state to the vertex $f$. The automaton continues moving up and recursively assigning states to vertices until it reaches the root of the string. The state assigned to the root is the state in which the automaton halts. If it is an accepting state, the entire ranked string is accepted.

As an example, let us build a parser for propositional logic. Consider the alphabet $\Sigma = \{ \top, \bot, \wedge, \neg, \vee, (, ) \}$ where we use only two propositions (one true, one false) for brevity.

Given a string in $\Sigma^*$ that is well-formed, like $(\top \vee \bot) \wedge \neg(\top \wedge (\top \vee \bot))$, it is nontrivial to build a *linear* automaton that recognises whether this formula evaluates to true or to false.

If we consider instead a *ranked* alphabet

$$\Pi = (\{\top, \bot, \wedge, \neg, \vee\}, |_|)$$

with $|\top| = |\bot| = 0$, $|\neg| = 1$ and $|\wedge| = |\vee| = 2$, the same formula is represented as the following tree:

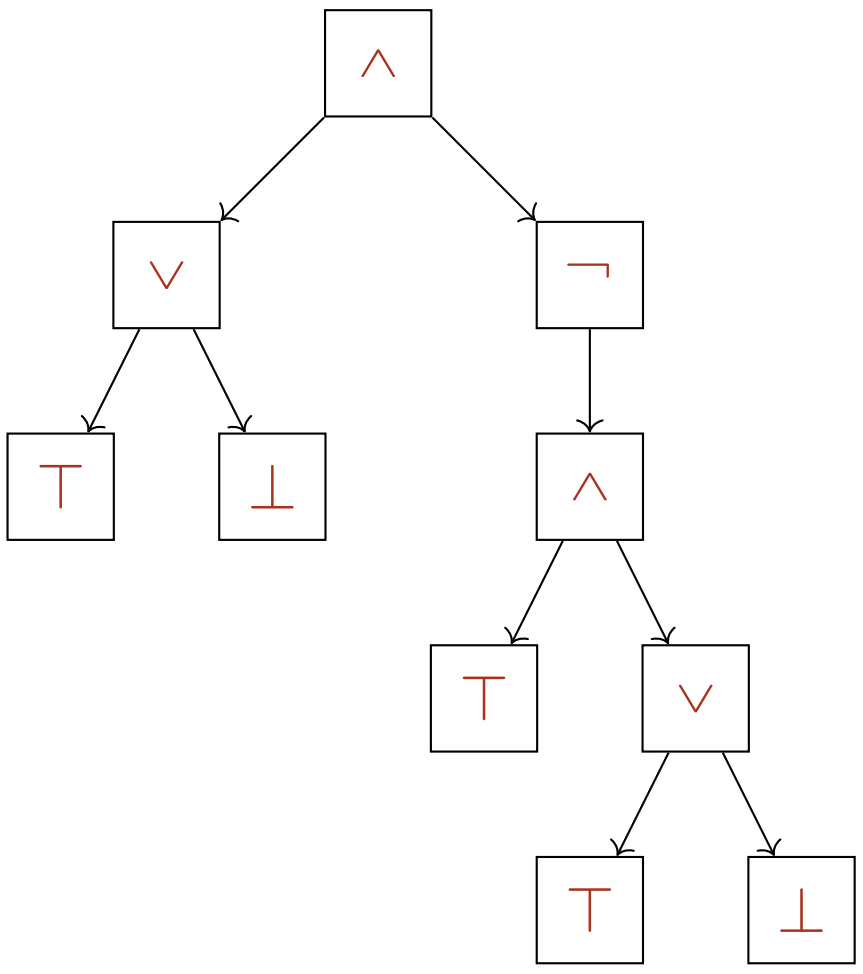


Note how we did not need to include parentheses to make the order of operations well-defined.

Now, in order to know the truth value of the entire formula, we need to know the truth value of the "topmost" expression, which just so happens to be located at the root of our tree. Of course, the value of this $\wedge$ depends on the value of the $\vee$ and the $\neg$ below it, whose values depend on …

The only expressions whose values are known *a priori* are the nullary symbols $\top$ and $\bot$, which correspond to the leaves of our tree. Let us hence assign truth values to these.

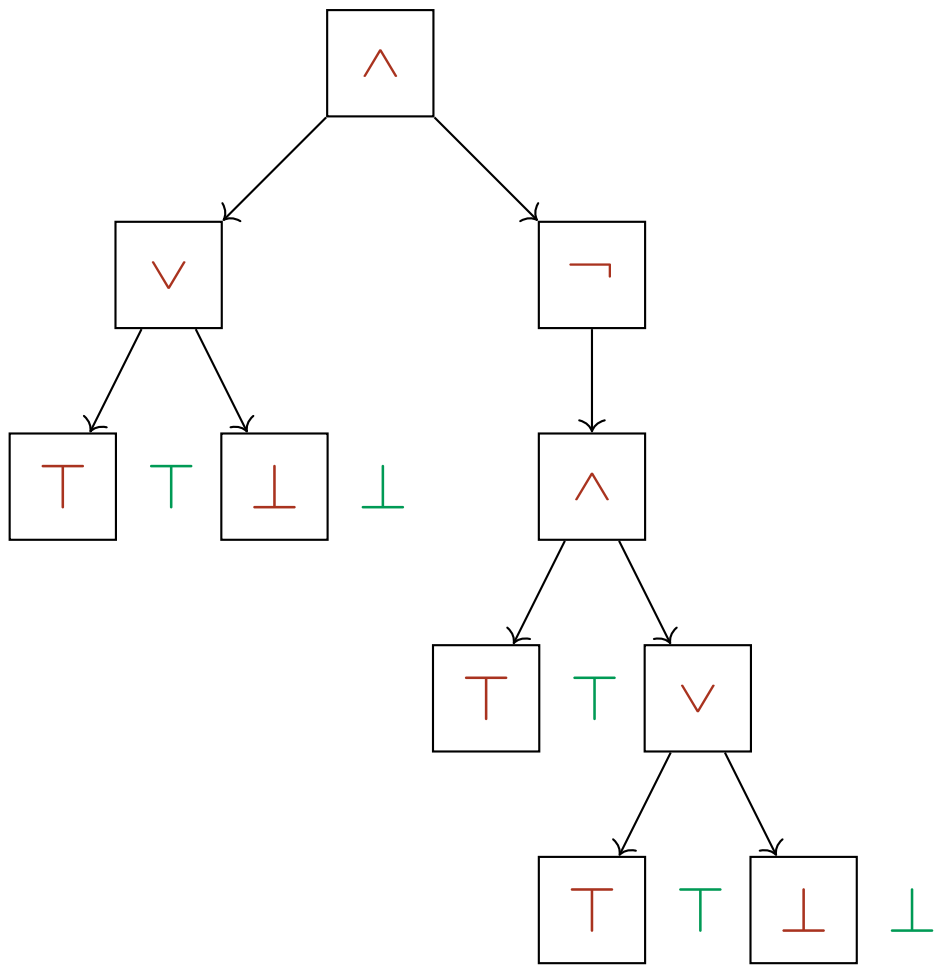


With the leaves translated from syntactical symbols to semantic meaning (even if represented by the same symbol in this particular case), we can now evaluate those vertices of the tree whose only children are leaves.

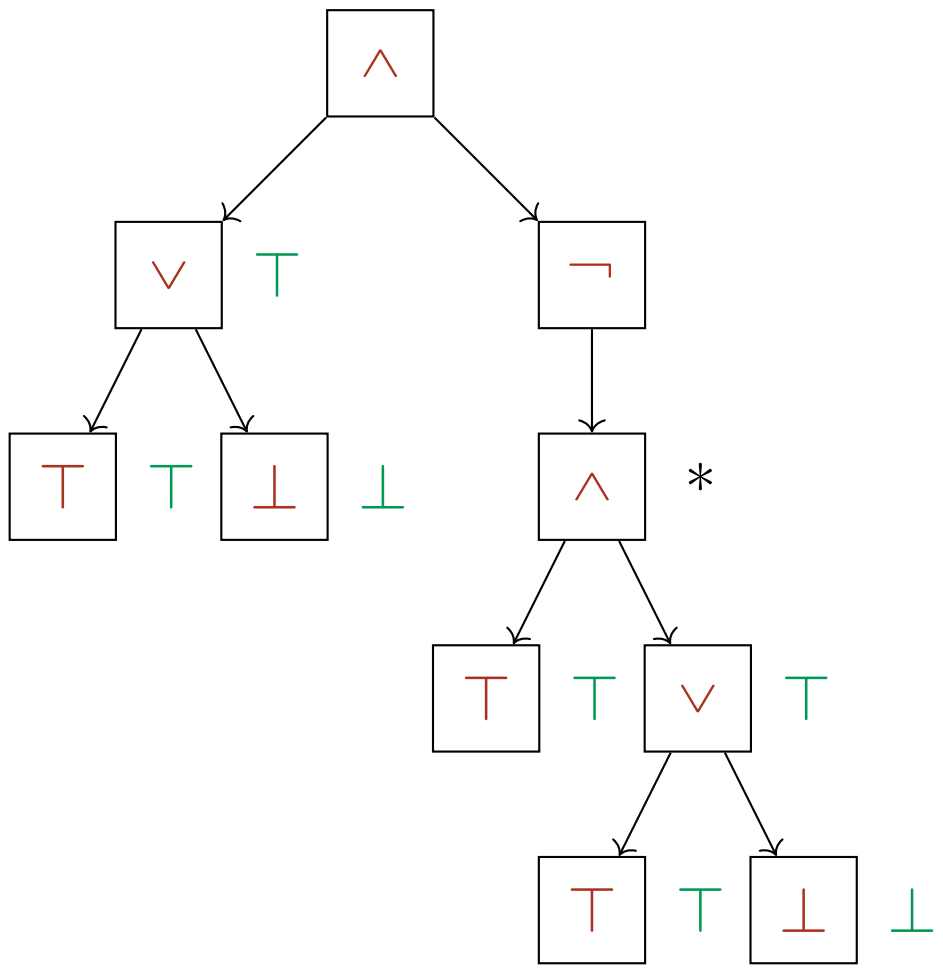


Note how we could not yet evaluate the vertex marked with $*$ because its right successor did not yet have a value assigned to it at the beginning of this step.

We do two more evaluation steps.

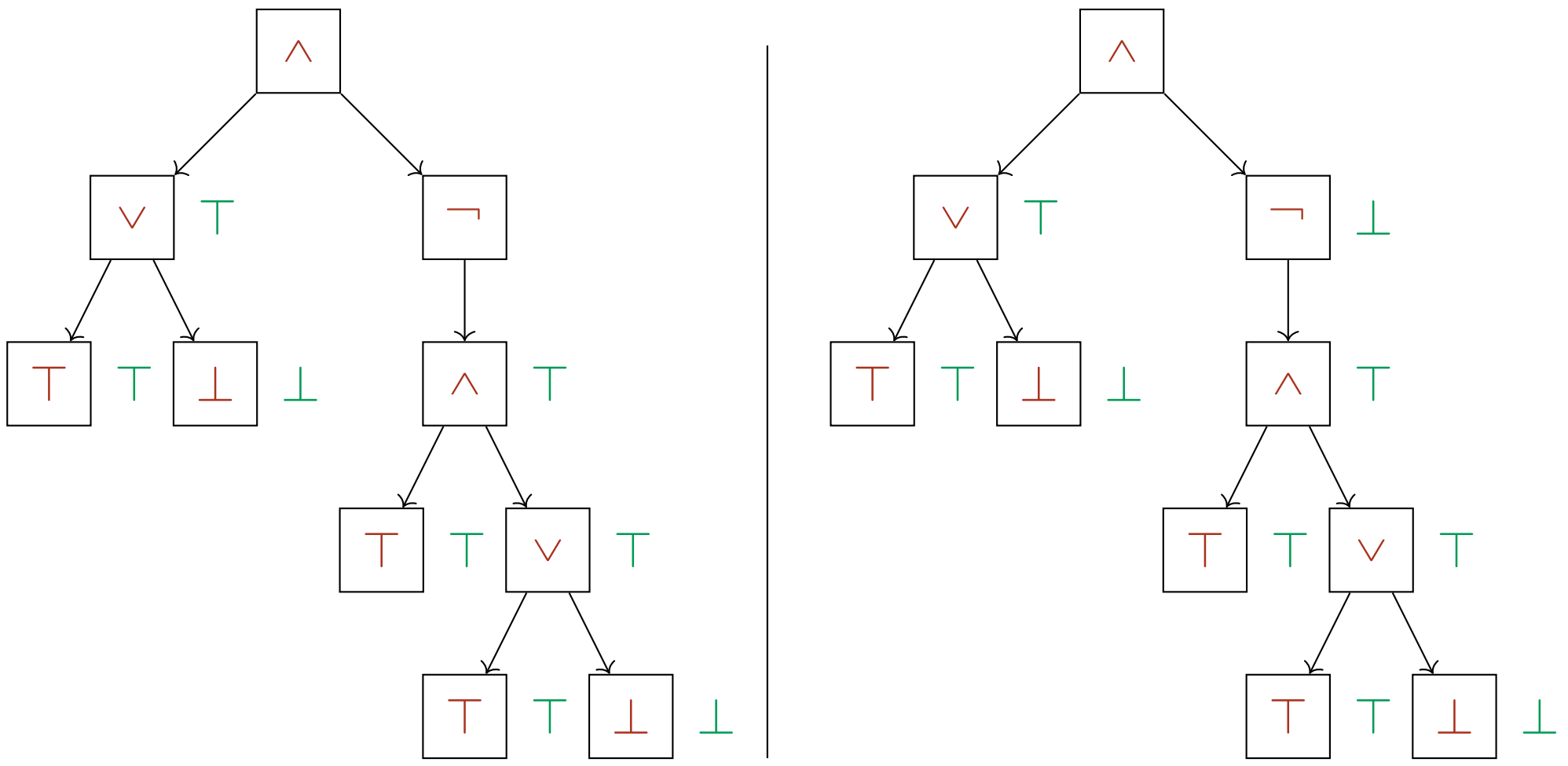

Finally, all successors of the root have been assigned a value, so we can evaluate the root vertex.

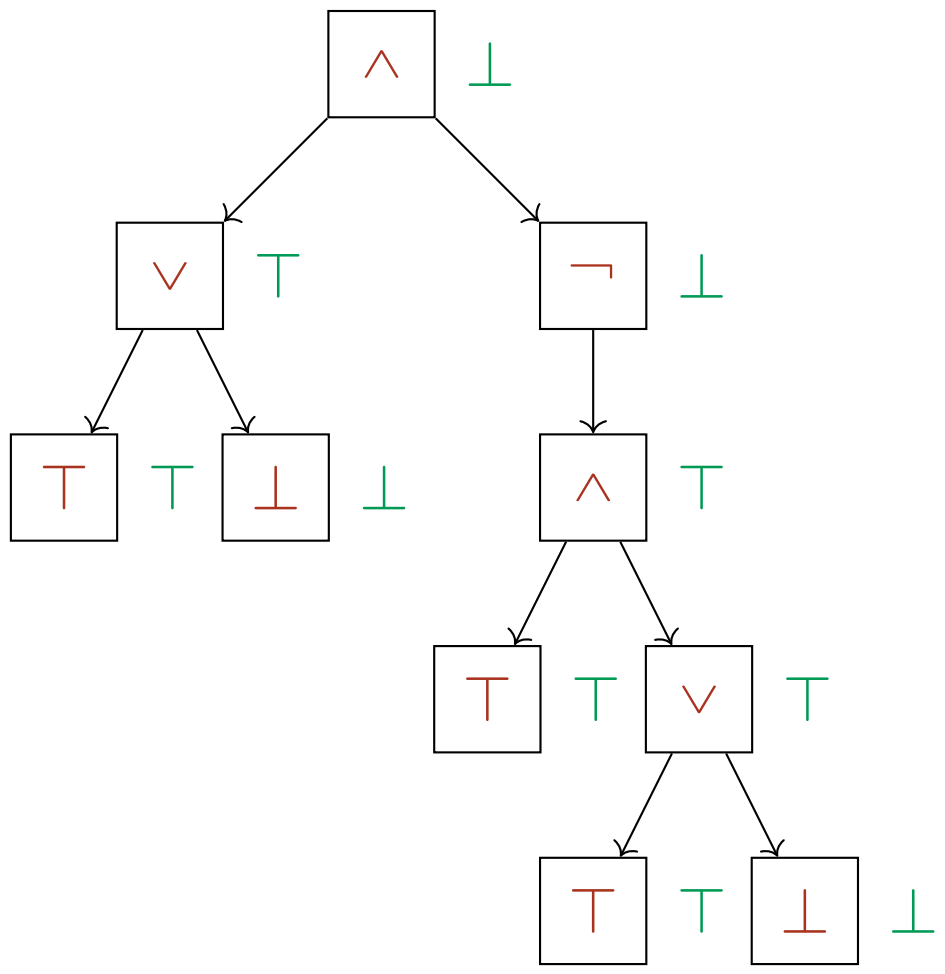

What does this mean? Since the value of the entire formula is equal to the value of the topmost expression, we have discovered that the formula evaluates to "False".

The reader has of course noticed that the procedure above was exactly the same as that of computing the value of an expression in a typed algebra.

Indeed, expressions will be the ranked strings which we input into the deterministic bottom-up finite tree automaton constructed later on.

**Lemma 8.3.2**

> Let $\mathscr{S} = (\mathscr{T}, \mathscr{F}, \langle_\rangle)$ be a signature.
>
> Then $\|\mathscr{S}\|$ is an $(\mathscr{F}, f \mapsto |f|)$-vocabulary.
>
> In particular, every $\mathscr{S}$-expression is a ranked string of the ranked alphabet $(\mathscr{F}, f \mapsto |f|)$.

$\|\mathscr{S}\|$ denotes the set of expressions over $\mathscr{S}$. (def. 6.3.3, p. 103)

**Proof.** Let $\mathscr{S} = (\mathscr{T}, \mathscr{F}, \langle_\rangle)$ be a signature, and let $e \in \|\mathscr{S}\|$. That $e$ is a ranked string over $(\mathscr{F}, f \mapsto |f|)$ follows directly from definitions 6.3.1 and 6.3.3. It remains to show that every subtree of $e$ is again an expression, but this is automatic due to the fact that *every* tree of function symbols which respects input types defines an expression over $\mathscr{S}$.

$\square$

Coming back to our propositional parser, for an automaton that accepts true formulas and rejects false ones, this means that the automaton eating our example tree should end in a non-accepting state. Indeed, we have done here exactly what a deterministic bottom-up finite tree automaton would do: start eating the expression at the leaves, evaluate a vertex as soon as all of its successors are ready, and finally check the value of the root vertex to decide whether the entire tree should be accepted.

Let us formally construct the corresponding automaton. We need a quintuple (𐤌, $\Pi$, $\Upsilon$, 𐤍, 𐤎). We already know the ranked alphabet, and the vocabulary $\Upsilon$ that our automaton can understand should be the set of all well-formed formulas buildable in $\Pi$. In this particular case, this is equal to $||\Pi||$, but it could be smaller for a different use case.

𐤌 is the Phoenician letter "mem".

𐤍 is the Phoenician letter "nun".

𐤎 is the Phoenician letter "semk".

The states of our automaton are simply the labels we attach to vertices during its run. In our case, that means the values $\top$ and $\bot$, so we set

$$𐤌 := \{\top, \bot\}.$$

The automaton should *accept* the tree if the expression evaluates to "True", that is, if the label of the top vertex is $\top$. Formally, we set the set of

accepting states to

$$\text{𐤍} := \{\top\} \subseteq \text{𐤌}.$$

𐤍 is the Phoenician letter "nun".

All that is left is the transition function 𐤎. We remind ourselves how it was defined: 𐤎 should take a symbol from $\Pi$ and a word from $\text{𐤌}^*$ and output a new state. What this means is simply that 𐤎 does precisely what we did above: it looks at a vertex (the symbol from $\Pi$) and at the labels of all of that vertex's successors (which are ordered into a word because our trees are traversal trees) and decides what label the new vertex should get.

𐤌 is the Phoenician letter "mem".

𐤎 is the Phoenician letter "semk".

The relevant parts of 𐤎 are essentially a truth table:

$$\begin{aligned}
\text{𐤎}(\top, \varepsilon) &= \top \\
\text{𐤎}(\bot, \varepsilon) &= \bot \\
\text{𐤎}(\neg, \top) &= \bot \\
\text{𐤎}(\neg, \bot) &= \top \\
\text{𐤎}(\wedge, \bot\bot) &= \bot \\
\text{𐤎}(\wedge, \top\bot) &= \bot \\
\text{𐤎}(\wedge, \bot\top) &= \bot \\
\text{𐤎}(\wedge, \top\top) &= \top \\
\text{𐤎}(\vee, \bot\bot) &= \bot \\
\text{𐤎}(\vee, \top\bot) &= \top \\
\text{𐤎}(\vee, \bot\top) &= \top \\
\text{𐤎}(\vee, \top\top) &= \top
\end{aligned}$$

As per our definition, 𐤎 should assign labels to all other combinations as well, but since these can never be reached by our automaton, it matters not what 𐤎 does with them. We shall for this reason omit unreachable states from definitions of tree automata.

The automaton, just like we did, applies the transition function recursively until every vertex is assigned a label. It then simply checks whether the state of the root vertex is an accepting state.

Every tree in $\Upsilon$ is hence assigned a unique vertex labelling by our deterministic bottom-up finite tree automaton. We give a name to this labelling.

## Definition 8.3.3

𐤄 is the Phoenician letter "he".

𐤌 is the Phoenician letter "mem".

𐤍 is the Phoenician letter "nun".

𐤎 is the Phoenician letter "semk".

Let $𐤄 = (𐤌, \Pi, \Upsilon, 𐤍, 𐤎)$ be a deterministic bottom-up finite tree automaton, and let $T = (V, E, \langle\!|_|\!\rangle, \star, \preccurlyeq, \langle_\rangle) \in \Upsilon$. The unique function $𐤇\colon V \to 𐤌$ such that

$$\forall v \in V\colon\ \mathrm{N}^{\mathrm{out}} v = [v_1, \dots, v_n] \Rightarrow 𐤎(\langle v\rangle, 𐤇v_1 \dots 𐤇v_n) = 𐤇v$$

is called the *interpretation of* $T$ *in* 𐤄, denoted by $𐤇_{𐤄,T}$.

𐤇 (het) is the interpretation function. (def. 8.3.3, p. 164)

That the function $𐤇_{𐤄,T}$ exists and is uniquely determined can be seen by the fact that 𐤎 is a function.

Interpretations make it easy to formally define what an accepted string is.

## Definition 8.3.4

We say that the automaton $𐤄 = (𐤌, \Pi, \Upsilon, 𐤍, 𐤎)$ *accepts* a string $T \in \Upsilon$ if $𐤇_{𐤄,T}\sqrt{T} \in 𐤍$.

Otherwise, we say that 𐤄 *rejects* $T$.

Intuitively, $𐤇_{𐤄,T}$ traces the path the tree $T$ takes through the states of the automaton while being processed. For this reason, we should expect any subtree of $T$ to be interpreted in the same way on its own as it is inside $T$ – in our example above, we did not take into consideration the predecessors of a vertex when deciding upon its label.

## Lemma 8.3.5

Let $𐤄 = (𐤌, \Pi, \Upsilon, 𐤍, 𐤎)$ be a deterministic bottom-up finite tree automaton, $T \in \Upsilon$. Then for any $v \in T$, we have

$$𐤇_{𐤄,T[v]} = (𐤇_{𐤄,T})\big|_{\{\, w \in T[v] \,\}}.$$

**Proof.** Note first that since $\Upsilon$ is a vocabulary, all subtrees involved in the statement lie again in $\Upsilon$.

We proceed by induction on the height $t$ of $v$ in $T$.

For $t = 1$, $v$ has no children, so $𐤇_{𐤄,T[v]} v = 𐤎(\langle v\rangle, \varepsilon) = 𐤇_{𐤄,T} v$.

For $t > 1$, let $\mathrm{N}^{\mathrm{out}}v = [v_1, \ldots, v_n]$. By induction and noticing the fact that $T[v_i] = T[v][v_i]$, we have

$$(𐤇_{𐤄,T})\big|_{\{\,w \in T[v_i]\,\}} = 𐤇_{𐤄,T[v_i]} = 𐤇_{𐤄,T[v][v_i]} = (𐤇_{𐤄,T[v]})\big|_{\{\,w \in T[v][v_i]\,\}}.$$

But now

$$𐤇_{𐤄,T}v = 𐤎(\langle v \rangle, 𐤇_{𐤄,T}v_1 \ldots 𐤇_{𐤄,T}v_n) = 𐤎(\langle v \rangle, 𐤇_{𐤄,T[v]}v_1 \ldots 𐤇_{𐤄,T[v]}v_n) = 𐤇_{𐤄,T[v]}v,$$

which proves the claim.

□

𐤇 (het) is the interpretation function.
(def. 8.3.3, p. 164)

𐤄 is the Phoenician letter “he”.

𐤎 is the Phoenician letter “semk”.

# 4. A Note on Runtime

Runtimes and algorithmical considerations are discussed in more depth in chapter 9, but we want to give a short intuition about how “good” deterministic bottom-up finite tree automata are.

Suppose someone hands us a deterministic bottom-up finite tree automaton $𐤄 = (𐤌, \Pi, \Upsilon, 𐤍, 𐤎)$ where the transition function 𐤎 is given by a constant-time oracle.

Then we can compute the final state of a given input string $T$ as follows.

𐤌 is the Phoenician letter “mem”.

𐤍 is the Phoenician letter “nun”.

**Algorithm 1:** Evaluating a deterministic bottom-up finite tree automaton

```
label each vertex of T with its number of successors
put all vertices with label 0 into a first-in-first-out queue
while the queue is not empty do
    pop the oldest vertex v off the queue
    compute 𐤇_{𐤄,T}(v)
    decrease the label of v's predecessor by 1
    if predecessor's label has become 0 then
        # we know the states of all its successors now
        push predecessor onto queue
```

This runs in time linear in the number of vertices of $T$, a fact which depends only on the constant-time evaluation of 𐤎.

**Lemma 8.4.1**

𐤄 is the Phoenician letter "he".

Let 𐤄 be a deterministic bottom-up finite tree automaton. Then there is an algorithm which, for any input string $T$, decides in time linear in the number of vertices of $T$ whether 𐤄 accepts $T$.

𐤌 is the Phoenician letter "mem".

𐤍 is the Phoenician letter "nun".

𐤎 is the Phoenician letter "semk".

𐤇 (het) is the interpretation function. (def. 8.3.3, p. 164)

**Proof.** Let $𐤄 = (𐤌, \varPi, \varUpsilon, 𐤍, 𐤎)$ be a deterministic bottom-up finite tree automaton. We need only show that in algorithm 1, the computation of $𐤇_{𐤄,T}(v)$ runs in constant time.

By definition of a deterministic bottom-up finite tree automaton, both the alphabet $\varPi$ and the set of states 𐤌 are finite.

Let now

$$n := \max\{\, |\alpha| : \alpha \in \varPi \,\}.$$

Then to compute any transition of 𐤄, it suffices to know 𐤎 on the set

$$\varPi \times \bigcup_{i=0}^{n} 𐤌^i.$$

But this is a finite set of constant size, and this size is independent of $T$, whence we can simply save 𐤎 as a lookup table.

□

In summary, once we have constructed a deterministic bottom-up finite tree automaton, evaluating whether a given ranked string is accepted runs quickly. In order to turn this into an algorithm on finite graphs, two problems remain to solve: construct a deterministic bottom-up finite tree automaton that accepts precisely those graph expressions whose value fulfils a given monadic second-order property, and show that given a graph, we can algorithmically compute a reasonably small $\mathfrak{S}$-expression for it.

$\mathfrak{S}$ is the signature of $\mathfrak{G}$. (def. 6.2.1, p. 92)

We address the former problem first.

# 5. Dictionaries

Suppose we are given for every possible string in a vocabulary a function that tells us which state our automaton should assign to which vertex.

**Definition 8.5.1**

Let $\Pi = (\Sigma, |_|)$ be a finite ranked alphabet, $\Upsilon$ a $\Pi$-vocabulary, and let $\text{𐤌}$ be a finite set. A family

$\text{𐤌}$ is the Phoenician letter "mem".

$$\{ \text{𐤋}_{(V,E,⦉_⦊,\star,\langle_\rangle,\preccurlyeq)} : V \to \text{𐤌} \}_{(V,E,⦉_⦊,\star,\langle_\rangle,\preccurlyeq)\in\Upsilon}$$

of functions which assigns to every string in $\Upsilon$ a vertex labelling in $\text{𐤌}$ is called a $\text{𐤌}$-*dictionary* for $\Upsilon$.

What we expect a dictionary to do is to give us the interpretations of strings in $\Upsilon$ under a deterministic bottom-up finite tree automaton with state set $\text{𐤌}$. This, however, prohibits those dictionaries where the labellings are not compatible with each other.

**Definition 8.5.2**

Let $\text{𐤌}$ be a finite set, $\Pi = (\Sigma, |_|)$ a finite ranked alphabet, and $\Upsilon$ a $\Pi$-vocabulary. We say that a $\text{𐤌}$-dictionary $\{ \text{𐤋}_T \}_{T\in\Upsilon}$ for $\Upsilon$ is *consistent* if it fulfills the following property.

$$\begin{aligned}
&\forall T = (V_T, E_T, ⦉_⦊_T, \star_T, \preccurlyeq_T, \langle_\rangle)) \in \Upsilon: \\
&\forall S = (V_S, E_S, ⦉_⦊_S, \star_S, \preccurlyeq_S, \langle_\rangle) \in \Upsilon: \\
&\forall v \in T: \\
&\forall v' \in S \text{ with } \mathrm{N}^{\mathrm{out}} v = [v_1, \dots, v_n], \mathrm{N}^{\mathrm{out}} v' = [v'_1, \dots, v'_n]: \\
&\quad \langle v \rangle = \langle v' \rangle \\
&\quad \wedge \text{𐤋}_T(v_1) = \text{𐤋}_S(v'_1) \\
&\quad \wedge \dots \\
&\quad \wedge \text{𐤋}_T(v_n) = \text{𐤋}_S(v'_n) \\
&\quad \Rightarrow \text{𐤋}_T(v) = \text{𐤋}_S(v').
\end{aligned}$$

In other words, in a consistent dictionary, the interpretation of a vertex depends only on the vertex's symbol and the interpretation of its children, which coincides with our notion of what a deterministic tree automaton should do.

In particular, a consistent dictionary agrees with itself on subtrees. We note this observation.

**Lemma 8.5.3**

𐤌 is the Phoenician letter “mem”.

Let $\Pi = (\Sigma, |_|)$ be a ranked alphabet, $\Upsilon$ a $\Pi$-vocabulary, 𐤌 a finite set, and let $\{𐤋_T\}_{T \in \Upsilon}$ be a consistent 𐤌-dictionary for $\Upsilon$.

Then for any $T \in \Upsilon$ and any $v \in T$, we have

$$𐤋_{T[v]} = 𐤋_T\big|_{\{\, w \in T[v] \,\}}.$$

**Proof.** This follows by induction from the definition of consistency, similar to the proof of lemma 8.3.5.

□

We now claim that, given a consistent dictionary, we can construct a deterministic bottom-up finite tree automaton that agrees with that dictionary.

**Lemma 8.5.4**

𐤍 is the Phoenician letter “nun”.

𐤄 is the Phoenician letter “he”.

Let $\Pi = (\Sigma, |_|)$ be a finite ranked alphabet, $\Upsilon$ a $\Pi$-vocabulary, 𐤌 a finite set, $𐤍 \subseteq 𐤌$, and let $\{𐤋_T\}_{T \in \Upsilon}$ be a 𐤌-dictionary for $\Upsilon$.

Then if $\{𐤋_T\}_{T \in \Upsilon}$ is consistent, there exists a deterministic bottom-up finite tree automaton $𐤄 = (𐤌, \Pi, \Upsilon, 𐤍, 𐤎)$ such that $\forall T \in \Upsilon$: $𐤇_{𐤄,T} = 𐤋_T$.

𐤎 is the Phoenician letter “semk”.

𐤇 (het) is the interpretation function. (def. 8.3.3, p. 164)

𐤐 is the Phoenician letter “pe”.

**Proof.** We need only construct the transition function 𐤎, and of this only the relevant parts.

For later convenience, we choose once and for all an arbitrary symbol $𐤐 \in 𐤌$.

Let $(f, w) \in \Sigma \times 𐤌^*$ with $|w| = |f|$, $w = w_1 \ldots w_n$.

Case 1: there is a $T = (V, E, ⦉_⦊, \star, \preccurlyeq, \langle_\rangle) \in \Upsilon$ with root $v$ and $\langle v \rangle = f$ that fulfils $\mathrm{N}^{\mathrm{out}} v = [v_1, \ldots, v_n]$ and $\forall i \in \{\, 1, \ldots, n \,\}$: $𐤋_T v_i = w_i$. Then we set $𐤎(f, w) := 𐤋_T v$.

By definition of a consistent dictionary, this is well-defined.

Case 2: there is no such $T \in \Upsilon$. Intuitively, this means that the collection $w$ of states can never be reached anyway, so the transition function for them should not matter (as we shall show in a moment). We set $\text{𐤎}(f, w) := \text{𐤐}$.

We now claim that the function thus constructed makes $\text{𐤄} := (\text{𐤌}, \mathit{\Pi}, \mathit{\Upsilon}, \text{𐤍}, \text{𐤎})$ into a deterministic bottom-up finite tree automaton that satisfies our needs, that is, $\forall T \in \mathit{\Upsilon}$: $\text{𐤇}_{\text{𐤄},T} = \text{𐤂}_T$.

That $\text{𐤄}$ is a deterministic bottom-up finite tree automaton is not a spectacular claim, as any old function $\text{𐤎}$ suffices to make a deterministic bottom-up finite tree automaton. It remains to show that the interpretations coincide with our dictionary.

Let therefore $T = (V, E, ⦉_⦊, \star, \preccurlyeq, \langle_\rangle) \in \mathit{\Upsilon}$, say with height $t$.

For $t = 1$, $T$ has only one vertex $v$ and we have $\text{𐤇}_{\text{𐤄},T} v = \text{𐤎}(\langle v\rangle, \epsilon) = \text{𐤂}_T v$.

For $t > 1$, let $v$ be the root of $T$, $\mathrm{N}^{\mathrm{out}} v = [v_1, \ldots, v_n]$. By induction, the claim holds for $T[v_1], \ldots, T[v_n]$. But now by lemma 8.3.5, we have

$$\begin{aligned}
\text{𐤇}_{\text{𐤄},T} v &= \text{𐤎}(\langle v\rangle, \text{𐤇}_{\text{𐤄},T} v_1, \ldots, \text{𐤇}_{\text{𐤄},T} v_n) \\
&= \text{𐤎}(\langle v\rangle, \text{𐤇}_{\text{𐤄},T[v_1]} v_1, \ldots, \text{𐤇}_{\text{𐤄},T[v_n]} v_n) \\
&= \text{𐤎}(\langle v\rangle, \text{𐤂}_{T[v_1]} v_1, \ldots, \text{𐤂}_{T[v_n]} v_n) \\
&= \text{𐤎}(\langle v\rangle, \text{𐤂}_T v_1, \ldots, \text{𐤂}_T v_n) \\
&= \text{𐤂}_T v.
\end{aligned}$$

Here, the first equality is due to the definition of $\text{𐤇}_{\text{𐤄},T}$, the second is due to lemma 8.3.5, the third is the induction hypothesis, the second-to-last is lemma 8.5.3, and the final equality is the definition of $\text{𐤎}$ according to case 1 since we have found a tree with the required labels.

□

We use this fact in the following section: given a recognisable set $X$ (of graphs, for example), we show the existence of a consistent dictionary for the set of expressions that yield values in $X$, which immediately yields a deterministic bottom-up finite tree automaton.

𐤎 is the Phoenician letter “semk”.

𐤐 is the Phoenician letter “pe”.

𐤄 is the Phoenician letter “he”.

𐤌 is the Phoenician letter “mem”.

𐤍 is the Phoenician letter “nun”.

𐤇 (het) is the interpretation function.
(def. 8.3.3, p. 164)

# 6. Algebras and Tree Automata

It is now time for the reader to recall all that they have forgotten from chapter 6 in order for us to explore the connection between typed algebras and tree automata. There is a whole theory to be developed here, as seen in [GS15] and [Cou89]. We, however, constrain ourselves to the niche pertinent to Courcelle's Theorem and hence focus on finitely expressible algebras.

**Definition 8.6.1**

A signature $(\mathscr{T}, \mathscr{F}, \langle_\rangle)$ is called *weakly locally finite* if for every finite $\mathscr{Y} \subseteq \mathscr{T}$, the set

$$\{ f \in \mathscr{F} : \langle f \rangle \in \mathscr{Y}^* \times \mathscr{Y} \}$$

is finite.

**Theorem 8.6.2**

Let $\mathscr{S} = (\mathscr{T}, \mathscr{F}, \langle_\rangle)$ be a signature, $\mathscr{A} = (\mathscr{C}, \mathscr{O})$ an $\mathscr{S}$-algebra, let $t \in \mathscr{T}$, let $\mathscr{L} \subseteq \mathscr{C}_t$, and let the following conditions be fulfilled:

- $\mathscr{S}$ is weakly locally finite.
- $\mathscr{L}$ is $\mathscr{A}$-recognisable.
- $\mathscr{L}$ is finitely expressible, say $\mathscr{Y}$-expressible for some finite $\mathscr{Y} \subseteq \mathscr{T}$.

Then there exists a deterministic bottom-up finite tree automaton 𐤄 that takes as input all $\mathscr{Y}$-local expressions from $\|\mathscr{S}\|$ and accepts an expression $e$ if and only if $\mathrm{val}_{\mathscr{A}} e \in \mathscr{L}$.

𐤄 is the Phoenician letter "he".

$\|\mathscr{S}\|$ denotes the set of expressions over $\mathscr{S}$. (def. 6.3.3, p. 103)

$\mathrm{val}_{\mathscr{A}} e$ denotes the result of $e$ when evaluated in $\mathscr{A}$. (def. 6.3.3, p. 103)

**Proof.** Let $\mathscr{S} = (\mathscr{T}, \mathscr{F}, \langle_\rangle)$ be a weakly locally finite signature. Let furthermore $\mathscr{A} = (\mathscr{C}, \mathscr{O})$ be an $\mathscr{S}$-algebra, let $t \in \mathscr{T}$, and let $\mathscr{L} \subseteq \mathscr{C}_t$ be $\mathscr{A}$-recognisable and finitely expressible, say $\mathscr{Y}$-expressible for some finite $\mathscr{Y} \subseteq \mathscr{T}$.

We obtain first the ranked alphabet from which our input strings will stem. We set $\Sigma := \{ f \in \mathscr{F} : \langle f \rangle \in \mathscr{Y}^* \times \mathscr{Y} \}$ and $|_| \colon \Sigma \to \mathbb{N}, f \mapsto |f|$. Then $\Pi := (\Sigma, |_|)$ is a ranked alphabet, and the set of all $\mathscr{Y}$-local expressions in $\|\mathscr{S}\|$ is a $\Pi$-vocabulary. Let us call it $\Upsilon$.

Since $\mathscr{Y}$ is finite and $\mathscr{S}$ is weakly locally finite, the alphabet $\Pi$ is finite.

Next, we determine the states of our automaton. Since $\mathscr{L}$ is $\mathscr{A}$-recognisable, by definition we can find a locally finite $\mathscr{S}$-algebra $\mathscr{B} = (\mathscr{D}, \mathscr{Q})$, an $\mathscr{S}$-algebra morphism $\{h_u\}_{u \in \mathscr{T}} \colon \mathscr{A} \to \mathscr{B}$, and a set $\mathscr{M} \subseteq \mathscr{D}_t$ such that $\mathscr{L} = {h_t}^{-1}\mathscr{M}$. We set

$$\text{𐤌} := \bigcup_{u \in \mathscr{Y}} \mathscr{D}_u.$$

𐤌 is the Phoenician letter “mem”.

Because $\mathscr{Y}$ is finite and $\mathscr{B}$ is locally finite, this is a finite (if potentially large) set.

The accepting states shall be $\text{𐤍} := \mathscr{M} \subseteq \text{𐤌}$.

𐤍 is the Phoenician letter “nun”.

We construct a consistent 𐤌-dictionary for $\Upsilon$ in order to apply lemma 8.5.4.

But this is easily done: given a string $T = (V, E, \langle\!\langle _ \rangle\!\rangle, \star, \preccurlyeq, \langle _ \rangle) \in ||\Pi||$, we remember that $T$ is actually also an expression in $\|\mathscr{S}\|$ and set

$$\text{𐤓}_T \colon V \to \text{𐤌}, v \mapsto \mathrm{val}_{\mathscr{B}} T[v].$$

$\|\mathscr{S}\|$ denotes the set of expressions over $\mathscr{S}$.
(def. 6.3.3, p. 103)

The value $\mathrm{val}_{\mathscr{B}} T[v]$ lies in 𐤌 since all function symbols in $\Sigma$ have their output sort in $\mathscr{Y}$. Thus the family $\{\text{𐤓}_T\}_{T \in \Upsilon}$ is a 𐤌-dictionary. Comparing definitions 8.5.2 and 6.3.3 reveals that by construction, this dictionary is consistent. Hence lemma 8.5.4 provides us a function 𐤎 and a deterministic bottom-up finite tree automaton

$\mathrm{val}_{\mathscr{B}} T[v]$ denotes the result of $T[v]$ when evaluated in $\mathscr{B}$.
(def. 6.3.3, p. 103)

$$\text{𐤄} = (\text{𐤌}, \Pi, \Upsilon, \text{𐤍}, \text{𐤎})$$

𐤎 is the Phoenician letter “semk”.

such that

$$\forall T \in \Upsilon \colon \text{𐤇}_{\text{𐤄}, T} = \text{𐤓}_T,$$

𐤄 is the Phoenician letter “he”.

or in other words, such that

$$\forall e \in \|\mathscr{S}\| \text{ such that } e \text{ is } \mathscr{Y}\text{-local} \colon \text{𐤇}_{\text{𐤄}, T} = \mathrm{val}_{\mathscr{B}} e.$$

𐤇 (het) is the interpretation function.
(def. 8.3.3, p. 164)

Thus, 𐤄 accepts an expression $e$ if and only if $\mathrm{val}_{\mathscr{B}} e \in \text{𐤍} = \mathscr{M}$. By lemma 6.3.4, this is equivalent to $h_t(\mathrm{val}_{\mathscr{A}} e) = \mathrm{val}_{\mathscr{B}} e \in \mathscr{M}$, or equivalently, $\mathrm{val}_{\mathscr{A}} e \in {h_t}^{-1}\mathscr{M} = \mathscr{L}$, proving that 𐤄 accepts precisely those expressions that take values in our recognisable set.

□

We state a special case that will be our primary concern in the next chapter.

### Corollary 8.6.3

Let $\mathscr{S} = (\mathscr{T}, \mathscr{F}, \langle_\rangle)$ be a signature, let $\mathscr{A} = (\mathscr{C}, \mathscr{O})$ be a $\mathscr{S}$-algebra, let $t \in \mathscr{T}$, and let $\mathscr{L} \subseteq \mathscr{C}_t$.

Suppose now that the following conditions are fulfilled.

- $\mathscr{S}$ is weakly locally finite.
- $\mathscr{L}$ is $\mathscr{A}$-recognisable.
- $\mathscr{A}$ is finitely expressible.

Then for any finite $\mathscr{Y} \subseteq \mathscr{T}$ such that $\mathscr{L}$ is $\mathscr{Y}$-expressible, there exists a deterministic bottom-up finite tree automaton 𐤄 that takes as input expressions from the set

$$\{\, e \in \|\mathscr{S}\| : e \text{ is } \mathscr{Y}\text{-local and } \mathrm{val}_{\mathscr{A}} e \in \mathscr{C}_t \,\}$$

and accepts an expression $e$ if and only if $\mathrm{val}_{\mathscr{A}} e \in \mathscr{L}$.

𐤄 is the Phoenician letter “he”.

$\|\mathscr{S}\|$ denotes the set of expressions over $\mathscr{S}$. (def. 6.3.3, p. 103)

$\mathrm{val}_{\mathscr{A}} e$ denotes the result of $e$ when evaluated in $\mathscr{A}$. (def. 6.3.3, p. 103)

Combined with lemma 8.4.1, this yields a membership algorithm which, given an expression $e$, is linear in the size of $e$. It remains only to check how large these expressions become and how hard it is to compute them for graphs. The next chapter concerns itself with these questions.

# Chapter 9
# Practical Considerations

For the entirety of this chapter, we assume that our graphs have no loops, that is, no edge is allowed to visit the same vertex more than once. This is done not because the results become false when one introduces loops, but to declutter the proofs. In fact, every result from this chapter holds just as well for graphs with loops. Appendix B shows how to formally arrive at this conclusion. Because of this fact, we do not mention the lack of loops when stating results in this chapter, while at the same time silently assuming all graphs are loop-free in our proofs. The reader may decide whether to restrict themselves to graphs without loops or whether to complement their reading by studying appendix B.

## 1. Applying the Theorem

Using the final result of the previous chapter, we arrive at a more practical version of Courcelle's Theorem.

**Theorem 9.1.1**

Let $\mathscr{S}$ be a signature, let $\mathscr{A} = (\mathscr{C}, \mathscr{O})$ be an $\mathscr{S}$-algebra, and let the following properties be fulfilled:

- $\mathscr{A}$ is inherited from $\mathfrak{G}$.
- $\mathscr{A}$ is finitely expressible.
- $\mathscr{S}$ is weakly locally finite.

$\mathfrak{G}$ is the algebra of graphs. (def. 6.2.2, p. 93)

Let further $\varphi$ be a sentence of the monadic second-order language of graphs, and let $n \in \mathbb{N}$.

Then there exists a deterministic bottom-up finite tree automaton that takes as input a set $\mathscr{K}$ of $\mathscr{S}$-expressions with $\mathscr{A}$-values in $\mathscr{C}_n$ and accepts an expression $e$ if and only if $\models_{\Omega(\mathrm{val}_{\mathscr{A}} e)} \varphi$, and every element of $\mathscr{C}_n$ is the value of some expression in $\mathscr{K}$.

$\mathrm{val}_{\mathscr{A}} e$ denotes the result of $e$ when evaluated in $\mathscr{A}$. (def. 6.3.3, p. 103)

$\Omega(\mathrm{val}_{\mathscr{A}} e)$ is the induced structure of $\mathrm{val}_{\mathscr{A}} e$. (def. 5.4.2, p. 73)

$\mathfrak{G}$ is the algebra of graphs. (def. 6.2.2, p. 93)

$\mathring{\mathfrak{L}}$ is the circuitous language of graphs. (def. 5.6.1, p. 77)

$\|\mathring{\mathfrak{L}}\|$ denotes the set of sentences over $\mathring{\mathfrak{L}}$. (def. 5.3.6, p. 62)

𐤄 is the Phoenician letter "he".

$\|\mathscr{S}\|$ denotes the set of expressions over $\mathscr{S}$. (def. 6.3.3, p. 103)

**Proof.** Suppose we have a weakly locally finite signature $\mathscr{S}$, a finitely expressible $\mathscr{S}$-algebra $\mathscr{A} = (\mathscr{C}, \mathscr{O})$ inherited from $\mathfrak{G}$, and a graph type $n \in \mathbb{N}$. Because $\mathscr{A}$ is finitely expressible, we find a finite set $\mathscr{Y} \subseteq \mathbb{N}$ such that $\mathscr{C}_n$ is $\mathscr{Y}$-expressible.

Let further $\varphi \in \|\mathring{\mathfrak{L}}\|$. By theorem 7.1.2, the set $\{ G \in \mathscr{C}_n : \models_{\Omega(G)} \varphi \}$ is $\mathscr{A}$-recognisable.

Now the assumptions of corollary 8.6.3 are fulfilled, and we get a deterministic bottom-up finite tree automaton 𐤄 that takes as inputs expressions from the set

$$\mathscr{K} := \{ e \in \|\mathscr{S}\| : e \text{ is } \mathscr{Y}\text{-local and } \mathrm{val}_{\mathscr{A}} e \in \mathscr{C}_t \}$$

and accepts an expression $e$ if and only if $\models_{\Omega(\mathrm{val}_{\mathscr{A}} e)} \varphi$.

Because $\mathscr{C}_n$ is $\mathscr{Y}$-expressible, every element of $\mathscr{C}_n$ is the value of some expression in $\mathscr{K}$.

□

Of course, it remains to show that there are nontrivial classes of graphs which fulfil these requirements; otherwise, theorem 9.1.1 would be true, but also purely academic[1].

In section 9.3, we show that, for example, the algebra of all graphs of tree-width at most $k$ for some $k \in \mathbb{N}$, equipped with certain graph operations, fulfils the assumptions of theorem 9.1.1.

Believing for now that such algebras exist, how does the existence of such an automaton translate into an algorithm?

Suppose we are given a finitely expressible algebra $\mathscr{A}$ of weakly locally finite signature inherited from $\mathfrak{G}$, where the type 0 graphs in $\mathscr{A}$ are $\mathscr{Y}$-expressible for some $\mathscr{Y} \subseteq \mathbb{N}$, $|\mathscr{Y}| < \infty$. Suppose we are further given a sentence $\varphi \in \mathfrak{X}$.

[1] A polite term for "useless".

By Courcelle's Theorem, we build a deterministic bottom-up finite tree automaton 𐤄 $=$ (𐤌, $\Pi$, $\Upsilon$, 𐤍, 𐤎) that accepts an $\mathscr{Y}$-local graph expression if and only if its value in $\mathscr{A}$ fulfils $\varphi$.

Then an algorithm to decide $\varphi$-fulfilment would look as follows.

**Algorithm 2:** Applying Courcelle's Theorem

find $\mathscr{Y}$-local expression $e$ such that $\mathrm{val}_{\mathscr{A}}\, e = G$
**if** 𐤇$_{𐤄, e}$ $\in$ 𐤍 **then**
  **return** yes
**else**
  **return** no

Step two, checking whether 𐤄 accepts the ranked string $e$, runs in time linear in the number of vertices of the input tree $e$ (lemma 8.4.1).

All that remains to analyse is how many vertices the expression $e$ has, and how long it takes to construct $e$ in the first place.

$\mathfrak{G}$ is the algebra of graphs. (def. 6.2.2, p. 93)

$\mathfrak{X}$ denotes the direct logical framework of graphs. (def. 5.4.3, p. 73)

𐤄 is the Phoenician letter "he".

𐤌 is the Phoenician letter "mem".

𐤍 is the Phoenician letter "nun".

𐤎 is the Phoenician letter "semk".

$\mathrm{val}_{\mathscr{A}}\, e$ denotes the result of $e$ when evaluated in $\mathscr{A}$. (def. 6.3.3, p. 103)

𐤇 (het) is the interpretation function. (def. 8.3.3, p. 164)

# 2. Computing Graph Expressions

Before we restrict our multiverse to graphs of a certain tree-width, we take another look at our original algebra $\mathfrak{G}$.

In theorem 6.4.2, we showed that every finite hypergraph can be constructed with only disjoint sum, terminal redefinition, and fusion. So why can we not stay in this algebra, if it is already so nice?

While the signature $\mathfrak{S}$ is indeed weakly locally finite, the algebra $\mathfrak{G}$ is not finitely expressible. Nevertheless, we turn the procedure outlined in the proof of theorem 6.4.2 into an algorithm just to see what goes wrong.

$\mathfrak{S}$ is the signature of $\mathfrak{G}$. (def. 6.2.1, p. 92)

The input is a type $n$ graph $G = (V, E, \langle\!|_|\!\rangle, t)$ with $|V| = n$ such that all elements of $V$ are terminal vertices. The output is an expression $e$ with $\mathrm{val}_{\mathfrak{G}}\, e = G$.

Note that the $e$ defined in the algorithm is an *expression*, not a graph –

$\oplus$ is the disjoint sum. (def. 4.5.3, p. 40)

when we write $e := e \oplus G'$, we mean the literal (ranked) string "$e \oplus G'$", not the *value* resulting from this disjoint sum.

$\mathfrak{v}$ is the type 1 graph with one vertex. (def. 4.5.6, p. 43)

$\mathfrak{e}_{|e|}$ is the type $|e|$ graph with $|e|$ vertices and one edge. (def. 4.5.6, p. 43)

嫁 (tsureai, Japanese for *to marry*) denotes the source fusion. (def. 4.5.5, p. 42)

$\leftrightarrows$ denotes the terminal redefinition. (def. 4.5.4, p. 41)

**Algorithm 3:** Constructing a graph expression

$w_{t(1)} := \mathfrak{v}$
$e := w_{t(1)}$
**for** $i \in \{2, \dots, n\}$ **do**
  $w_{t(i)} := \mathfrak{v}$
  $e := e \oplus w_{t(i)}$
**for** $e \in E$ **do**
  $V^e := \mathfrak{e}_{|e|}$
  \# *denote the vertices of* $V^e$ *as* $\{V^e_1, \dots, V^e_{|e|}\}$
  $e := e \oplus V^e$
  **for** $\imath \in \{1, \dots, |e|\}$ **do**
    find $b \in \mathbb{N}$ such that $t'(b) = V^e_\imath$
    set $v$ to be the $\imath$-th vertex in $\langle\!\langle e \rangle\!\rangle$
    find $a \in \mathbb{N}$ such that $t'(a) = w_v$
    $e := 嫁^b_a e$
set $\sigma$ to correct the terminal vertices
**return** $\leftrightarrows_\sigma e$

The correctness of this algorithm has already been shown in the proof of theorem 6.4.2.

$\mathfrak{G}$ is the algebra of graphs. (def. 6.2.2, p. 93)

$\|\mathfrak{G}\|$ denotes the set of expressions over $\mathfrak{G}$. (def. 6.3.3, p. 103)

**Theorem 9.2.1**

Given a type $n \in \mathbb{N}$ graph $G = (V, E, \langle\!\langle _ \rangle\!\rangle, t)$ with $|V| = n$ such that all elements of $V$ are terminal vertices, algorithm 3 outputs an expression $e \in \|\mathfrak{G}\|$ with $\text{val}_{\mathfrak{G}}\, e = G$ in time $\mathcal{O}(|V| + |E| + u|E|)$, where

$$u := \max_{e \in E} |e|.$$

The number of vertices of $e$ is also in $\mathcal{O}(|V| + |E| + u|E|)$.[2]

**Proof.** Only the runtime and the size of the resulting expression remain

[2] It might seem strange at first that we do not simply write $\mathcal{O}(|V| + u|E|)$, but the two differ in the case $u = 0$.

to be shown.

We start with an expression (a ranked string) with one vertex. The first loop in algorithm 3 takes $|V| - 1$ steps and adds $|V| - 1$ vertices to $e$.

The second loop takes $|E|$ steps. Each step adds one vertex to $e$, then goes into the nested loop of at most $|e| \leq u$ steps and adds as many vertices to $e$. Finding $a$ and $b$ runs in linear time, provided we use a reasonable data structure that allows us to look up where in the list of terminals a given vertex lives and we remember at which natural number we "left off" after the previous edge.

This finishes the proof.

□

As seen in lemma 6.4.1, it is a trivial step from here to constructing all possible typed graphs.

The input is an arbitrary type $n$ graph $G = (V, E, \langle\!|_|\!\rangle, t)$. The output is an expression $e$ with $\mathrm{val}_{\mathfrak{G}}\, e = G$.

**Algorithm 4:** Type Schmype

enumerate $V = \{ v_1, \ldots, v_{|V|} \}$
use algorithm 3 to find $e \in \|\mathfrak{S}\|$
with $\mathrm{val}_{\mathfrak{G}}\, e = (V, E, \langle\!|_|\!\rangle, i \mapsto v_i) \in \mathfrak{G}_{|V|}$
set $\sigma\colon v_i \mapsto i$
$e := \leftrightarrows_{\sigma \circ t}\, e$
**return** $e$

This adds a constant number of steps on top of algorithm 3. We restate theorem 9.2.1 in full generality.

**Theorem 9.2.2**

Given a typed graph $G = (V, E, \langle\!|_|\!\rangle, t)$, algorithm 4 outputs an expression $e \in \|\mathfrak{S}\|$ with $\mathrm{val}_{\mathfrak{G}}\, e = G$ in time $\mathcal{O}(|V| + |E| + u|E|)$, where

$$u := \max_{e \in E} |e|.$$

The number of vertices of $e$ is also in $\mathcal{O}(|V| + |E| + u|E|)$.

$\mathrm{val}_{\mathfrak{G}}\, e$ denotes the result of $e$ when evaluated in $\mathfrak{G}$. (def. 6.3.3, p. 103)

$\mathfrak{G}$ is the algebra of graphs. (def. 6.2.2, p. 93)

$\mathfrak{S}$ is the signature of $\mathfrak{G}$. (def. 6.2.1, p. 92)

$\|\mathfrak{S}\|$ denotes the set of expressions over $\mathfrak{S}$. (def. 6.3.3, p. 103)

$\mathfrak{G}_{|V|}$ denotes the set of all graphs of type $|V|$. (def. 4.5.1, p. 39)

$\leftrightarrows$ denotes the terminal redefinition. (def. 4.5.4, p. 41)

For arbitrary graphs, the $u$ defined in theorem 9.2.2 can grow inconveniently large. We state two special cases as corollaries, abusing notation to omit the terminal function – the reader may imagine their favourite number of terminals, for example zero.

### Corollary 9.2.3

For any hypergraph $G = (V, E, \langle\!|_|\!\rangle)$, algorithm 4 outputs an expression $e \in \|\mathfrak{S}\|$ with $\mathrm{val}_{\mathfrak{G}}\, e = G$ in time $\mathcal{O}(|V| + |V| \cdot |E|)$. The number of vertices of $e$ is also in $\mathcal{O}(|V| + |V| \cdot |E|)$.

$\mathfrak{S}$ is the signature of $\mathfrak{G}$. (def. 6.2.1, p. 92)

$\|\mathfrak{S}\|$ denotes the set of expressions over $\mathfrak{S}$. (def. 6.3.3, p. 103)

### Corollary 9.2.4

Let $k \in \mathbb{N}$. For a $k$-uniform hypergraph $G = (V, E, \langle\!|_|\!\rangle)$, algorithm 4 outputs an expression $e \in \|\mathfrak{S}\|$ with $\mathrm{val}_{\mathfrak{G}}\, e = G$ in time $\mathcal{O}(|V| + k|E|)$. The number of vertices of $e$ is also in $\mathcal{O}(|V| + k|E|)$.

$\mathrm{val}_{\mathfrak{G}}\, e$ denotes the result of $e$ when evaluated in $\mathfrak{G}$. (def. 6.3.3, p. 103)

$\mathfrak{G}$ is the algebra of graphs. (def. 6.2.2, p. 93)

So, why can we not find an automaton to detect sets of 2-uniform graphs and solve the minimum vertex cover problem in linear time?

The crux is, of course, the non-finite expressibility. Recall that a deterministic bottom-up finite tree automaton can only deal with a finite alphabet; otherwise the transition table becomes infinite and computing the interpretation function can no longer be accomplished in linear time. But this means that we can not end up with a deterministic bottom-up finite tree automaton detecting a recognisable set $\mathcal{L}$ if $\mathcal{L}$ is not finitely expressible – if there are graphs for which the expressions need arbitrarily large types, our ranked alphabet becomes infinite.

We show that there are such graphs in $\mathfrak{G}$.

### Theorem 9.2.5

Let $k \in \mathbb{N}_{>0}$. Let $e \in \|\mathfrak{S}\|$ such that $\mathrm{val}_{\mathfrak{G}}\, e$ is the complete 2-uniform graph on $k$ vertices with 0 terminal vertices. Then $e$ is not $\{1, \dots, k-1\}$-local.

**Proof.** Fix a positive integer $k \in \mathbb{N}_{>0}$ and denote by $G$ the complete 2-uniform type 0 graph on $k$ vertices. Let further $e = (T, \langle_\rangle) \in \|\mathfrak{S}\|$ such that $\mathrm{val}_{\mathfrak{G}}\, e = G$. Suppose now $e$ were $\{1, \dots, k-1\}$-local.

Consider a subtree $S$ of $T$ such that $G' := \mathrm{val}_{\mathfrak{G}}(S, \langle_\rangle_{|S})$ already contains a vertex of degree $k-1$. We choose $S$ to be minimal with this property, that is, such that for every *proper* subtree $U$ of $S$, the graph $\mathrm{val}_{\mathfrak{G}}(U, \langle_\rangle_{|U})$ has no vertex of degree more than $k-2$.

$\mathrm{val}_{\mathfrak{G}}(S, \langle_\rangle_{|S})$ denotes the result of $(S, \langle_\rangle_{|S})$ when evaluated in $\mathfrak{G}$. (def. 6.3.3, p. 103)

Then $\langle\sqrt{S}\rangle$ must be a fusion symbol: if it were a redefinition, then the subtree rooted at the sole child of $\sqrt{S}$ would result in the same graph except for terminals and would thus contain a vertex of degree $k-1$. If it were a disjoint sum, then the subtree rooted at at least one of its children would contain a vertex of degree $k-1$ because a disjoint sum creates no edges.

$\mathfrak{G}$ is the algebra of graphs. (def. 6.2.2, p. 93)

We can thus assume without loss of generality that $\langle\sqrt{S}\rangle = 婚^2_1$. Denote the first terminal of $G'$ (the result of this fusion) by $v$.

婚 (tsureai, Japanese for *to marry*) denotes the source fusion. (def. 4.5.5, p. 42)

Since $G'$ contains a vertex of degree $k-1$, it must have at least $k$ vertices. Because $\boldsymbol{e}$ is $\{1, \dots, k-1\}$-local, this implies that there exists a vertex $v' \in G'$ which is not a terminal. Since $v$ is a terminal, we have $v' \neq v$.

Because $v'$ is not a terminal, it can never become a terminal in any operation on the path from $\sqrt{S}$ to $\sqrt{T}$. In particular, there can be no fusion involving $v'$ on said path, which means its degree cannot change between $G'$ and $G$. But $G$ contains only vertices of degree $k-1$, whence $v'$ must already have degree $k-1$ in $G'$.

Since the fusion $婚^2_1$ did not change the degree of $v'$, we can conclude that, for the sole child $w$ of $\sqrt{S}$, we must have that in the graph $\mathrm{val}_{\mathfrak{G}}(S[w], \langle_\rangle_{|S[w]})$, the vertex $v'$ already has degree $k-1$.

But $S[w]$ is a proper subtree of $S$, a contradiction.

□

Thus we have proven the following.

### Corollary 9.2.6

The algebra $\mathfrak{G}$ is not finitely expressible.

For this reason, we must restrict ourselves to "smaller" algebras inherited from $\mathfrak{G}$. Maybe the opposite extreme? Can we simply add enough function symbols to make our algebra finitely expressible?

$\mathfrak{G}$ is the algebra of graphs. (def. 6.2.2, p. 93)

The extreme version of this is the closure $\overline{\mathfrak{G}}$. It is trivially finitely expressible, since for any element $\mathfrak{c}$, there exists a nullary function symbol (and hence a very small expression) evaluating to $\mathfrak{c}$.

$\overline{\mathfrak{G}}$ denotes the closure of $\mathfrak{G}$. (def. 6.10.6, p. 141)

Great! We simply use this algebra and all our problems go away.

This is, of course, nonsense. This time, our algebra is finitely expressible, but its signature is no longer weakly locally finite: one nullary function symbol for every finite graph translates into infinitely many nullary function symbols of each type. Looking back at the proof of theorem 8.6.2, we see that if the signature is not weakly locally finite, the ranked alphabet becomes (again!) infinite.

Could there be a sweet spot in between these two extremes that is both finitely expressible and of weakly locally finite signature? Possibly, but it would still not be of much practical use.

### Theorem 9.2.7

$\mathfrak{G}_n$ denotes the set of all graphs of type $n$. (def. 4.5.1, p. 39)

Let $\mathscr{A} = (\mathscr{C}, \mathscr{O})$ be an algebra of weakly locally finite signature inherited from $\mathfrak{G}$ such that there is some $n \in \mathbb{N}$ with $\mathscr{C}_n = \mathfrak{G}_n$, and let $\mathscr{C}_n$ be $\{1, \ldots, k\}$-expressible for some $k \in \mathbb{N}$.

Then unless $\mathbf{P} = \mathbf{NP}$, there does not exist a polynomial-time algorithm which, given a type $n$ graph $G$, computes a $\{1, \ldots, k\}$-local expression $\mathfrak{e}$ with $\mathrm{val}_{\mathscr{A}} \mathfrak{e} = G$.

$\mathrm{val}_{\mathscr{A}} \mathfrak{e}$ denotes the result of $\mathfrak{e}$ when evaluated in $\mathscr{A}$. (def. 6.3.3, p. 103)

**Proof.** Let $\mathscr{A} = (\mathscr{C}, \mathscr{O})$ be an algebra of weakly locally finite signature inherited from $\mathfrak{G}$ such that there is some $n \in \mathbb{N}$ with $\mathscr{C}_n = \mathfrak{G}_n$, and let $\mathscr{C}_n$ be $\{1, \ldots, k\}$-expressible for some $k \in \mathbb{N}$.

If the reader has read the example on page 74, they might remember that the property "this graph is 3-colourable" can be encoded in monadic second-order logic. Let $\varphi$ be the sentence of the monadic second-order language of graphs encoding 3-colourability.

$\Omega(G)$ is the induced structure of $G$. (def. 5.4.2, p. 73)

By theorem 7.1.2, the set $\mathscr{L} := \{G \in \mathscr{C}_n : \vDash_{\Omega(G)} \varphi\}$ is $\mathscr{A}$-recognisable. We set $\mathscr{Y} := \{1, \ldots, k\}$, a finite set. Then by corollary 8.6.3, there exists a deterministic bottom-up finite tree automaton 𐤄 that takes as input all $\mathscr{Y}$-local expressions that result in a graph from $\mathscr{C}_n$ and accepts an expression if and only if its value lies in $\mathscr{L}$.

𐤄 is the Phoenician letter "he".

Suppose now that we know an algorithm that, given a graph $G \in \mathscr{C}_n$, computes a $\{1, \ldots, k\}$-local expression $e$ with $\mathrm{val}_{\mathscr{A}} e = G$ in polynomial time.

$\mathrm{val}_{\mathscr{A}} e$ denotes the result of $e$ when evaluated in $\mathscr{A}$. (def. 6.3.3, p. 103)

Given any (untyped) 2-uniform graph $G$, we can arbitrarily designate one of its vertices to be a terminal with multiplicity $n$, yielding a type $n$ graph, also called $G$. Since $\mathscr{C}_n = \mathfrak{G}_n$, we then have $G \in \mathscr{C}_n$ and can compute in polynomial time a $\{1, \ldots, k\}$-local expression $e$ with $\mathrm{val}_{\mathscr{A}} e = G$. Since the algorithm outputs $e$ in polynomial time, the number of vertices of $e$ must also be polynomial in the size of $G$.

$\mathfrak{G}_n$ denotes the set of all graphs of type $n$. (def. 4.5.1, p. 39)

Because $e$ is $\{1, \ldots, k\}$-local, it is a valid input for 𐤄 and we can decide in time linear in the size of $e$ whether 𐤄 accepts $e$, in other words, whether $G$ is 3-colourable. We have thus constructed an algorithm which, given any graph $G$, decides in polynomial time whether $G$ is 3-colourable.

𐤄 is the Phoenician letter "he".

Since 3-colourability is **NP**-complete, this implies **P** = **NP**.

□

Therefore, when examining a class of graphs, one must do two things:

- Show that the class is a finitely expressible algebra inherited from $\mathfrak{G}$ with weakly locally finite signature and
- provide an algorithm that, given a graph in the class, actually constructs an expression in this inherited algebra in reasonable time.

$\mathfrak{G}$ is the algebra of graphs. (def. 6.2.2, p. 93)

We give two examples of such classes. The first, graphs of tree-width bounded by a constant, is well known.[3] The second one, the narrower case of graphs with path-width bounded by a constant, makes use of the more general formulation of Courcelle's Theorem that we have developed. We include this case as an example of how the reader may construct their own inherited algebra of graphs for use with Courcelle's Theorem.

[3] This result is already mentioned in the 1993 paper [Bod93], which can also point the reader to many more interesting results on graphs of tree-width bounded by a constant.

# 3. Graphs of Bounded Tree-Width

The reader, if not familiar with the concept of tree-decompositions, may want to refer back to section 4.6 for the pertinent definitions.

$\mathfrak{G}$ is the algebra of graphs. (def. 6.2.2, p. 93)

We now construct, for any given $k \in \mathbb{N}$, a graph algebra inherited from $\mathfrak{G}$ where the carrier sets are the graphs of tree-width at most $k$. What remains to construct are the function symbols – we need the range of each function symbol to again lie in our inherited algebra; in other words, the tree-width bound needs to be preserved.

It is easy to see that fusion is "too powerful" for this:

媿 (tsureai, Japanese for *to marry*) denotes the source fusion. (def. 4.5.5, p. 42)

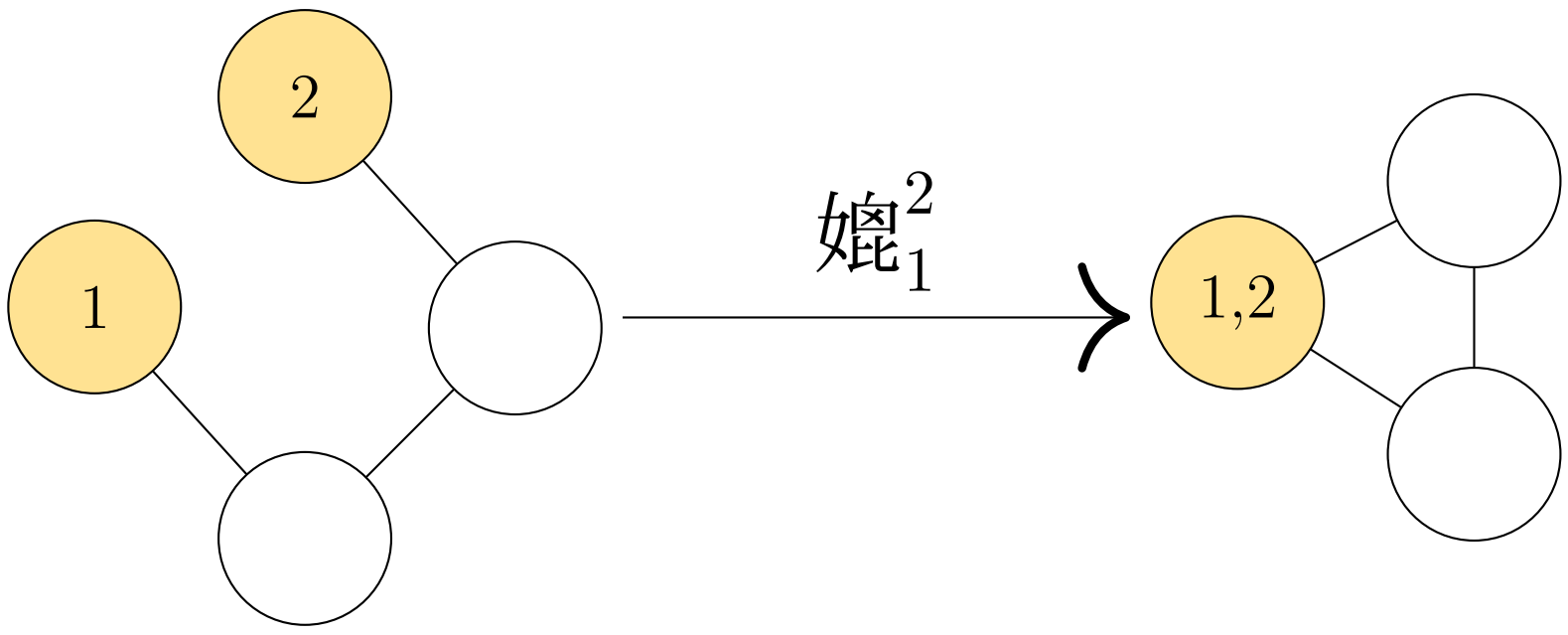

The problem, as the reader can see, is that allowing to fuse arbitrary vertices of a connected graphs can introduce a clique of a larger size into the graph. We hence replace fusion by a construction that is closer (but not equal) to how one constructs partial $k$-trees[4]: we disallow fusing vertices from the same connected component and instead allow to fuse vertices from two different graphs at the same time as we take their direct sum:

[4] The curious reader can learn all about partial $k$-trees in chapter 7 of [Bod98].

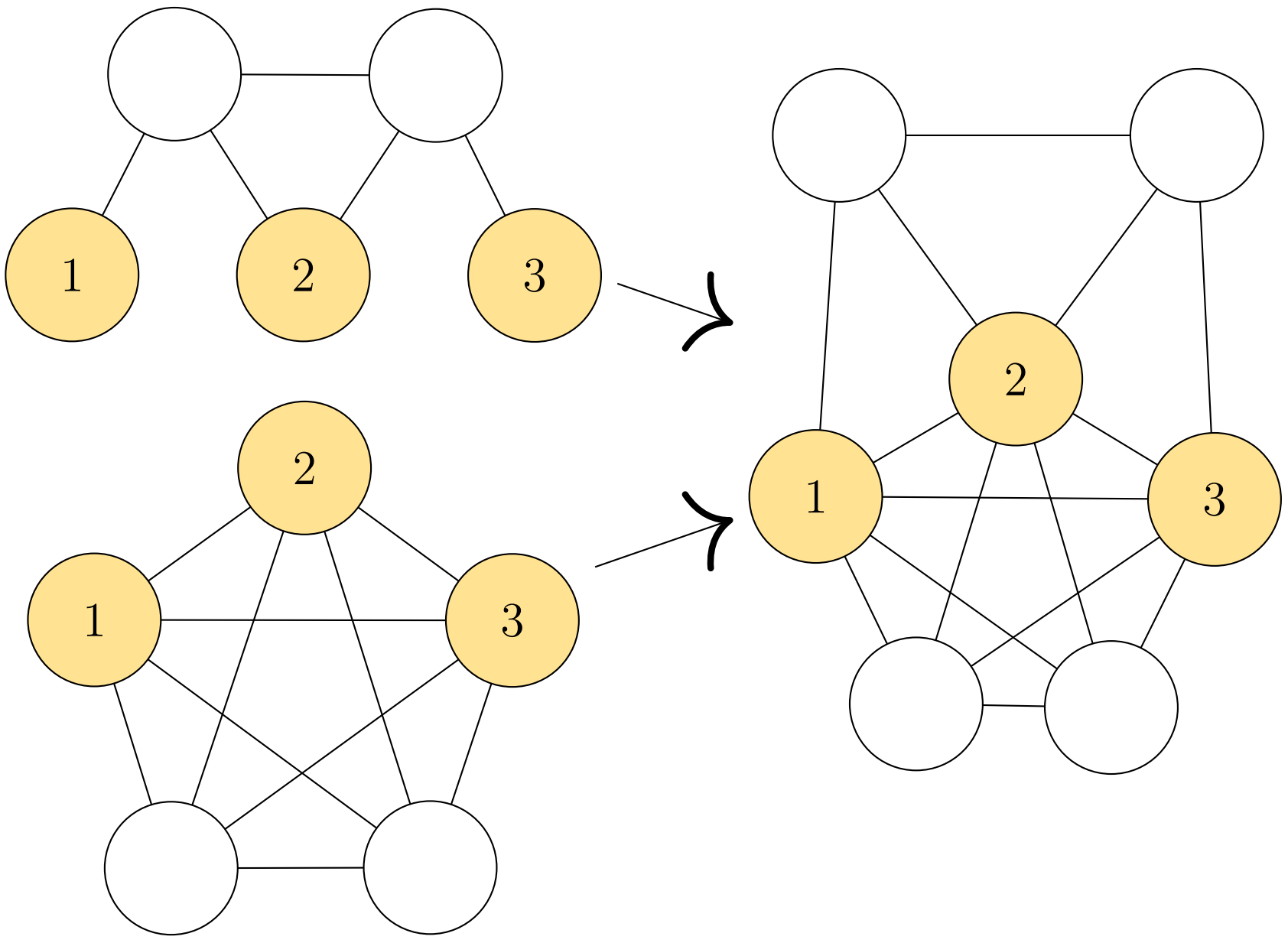

We call this a “3-twine”, or more generally a $k$-twine when we want a graph of type $k$ – the 2-twine simply forgets the third terminal after fusing it:

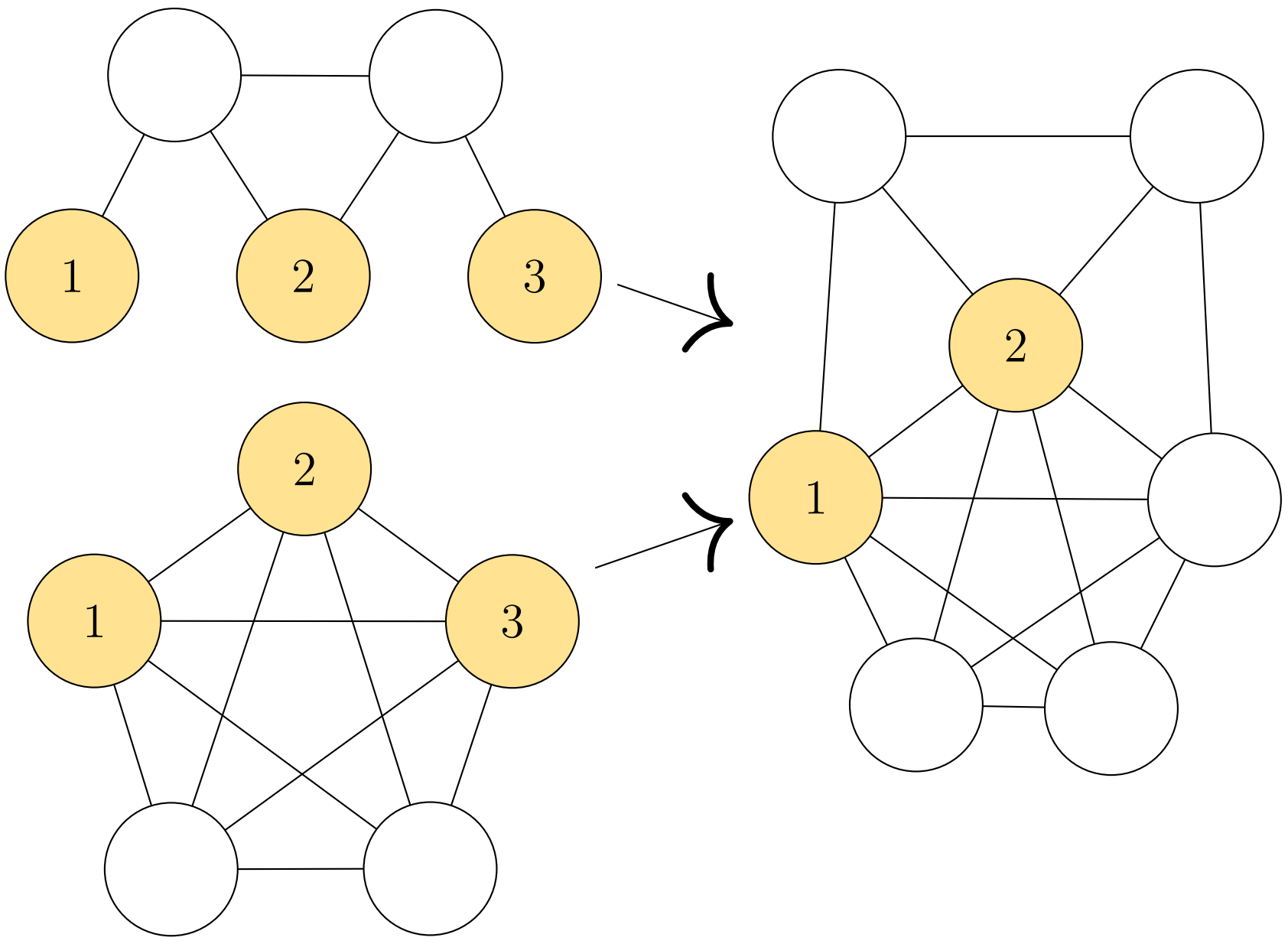

We might sometimes want to fuse only some, not all, of the terminal vertices. To distinguish this, the above is more precisely called the 2-twine *over* $\{1, 2, 3\}$. The 2-twine over $\{1\}$ would look different:

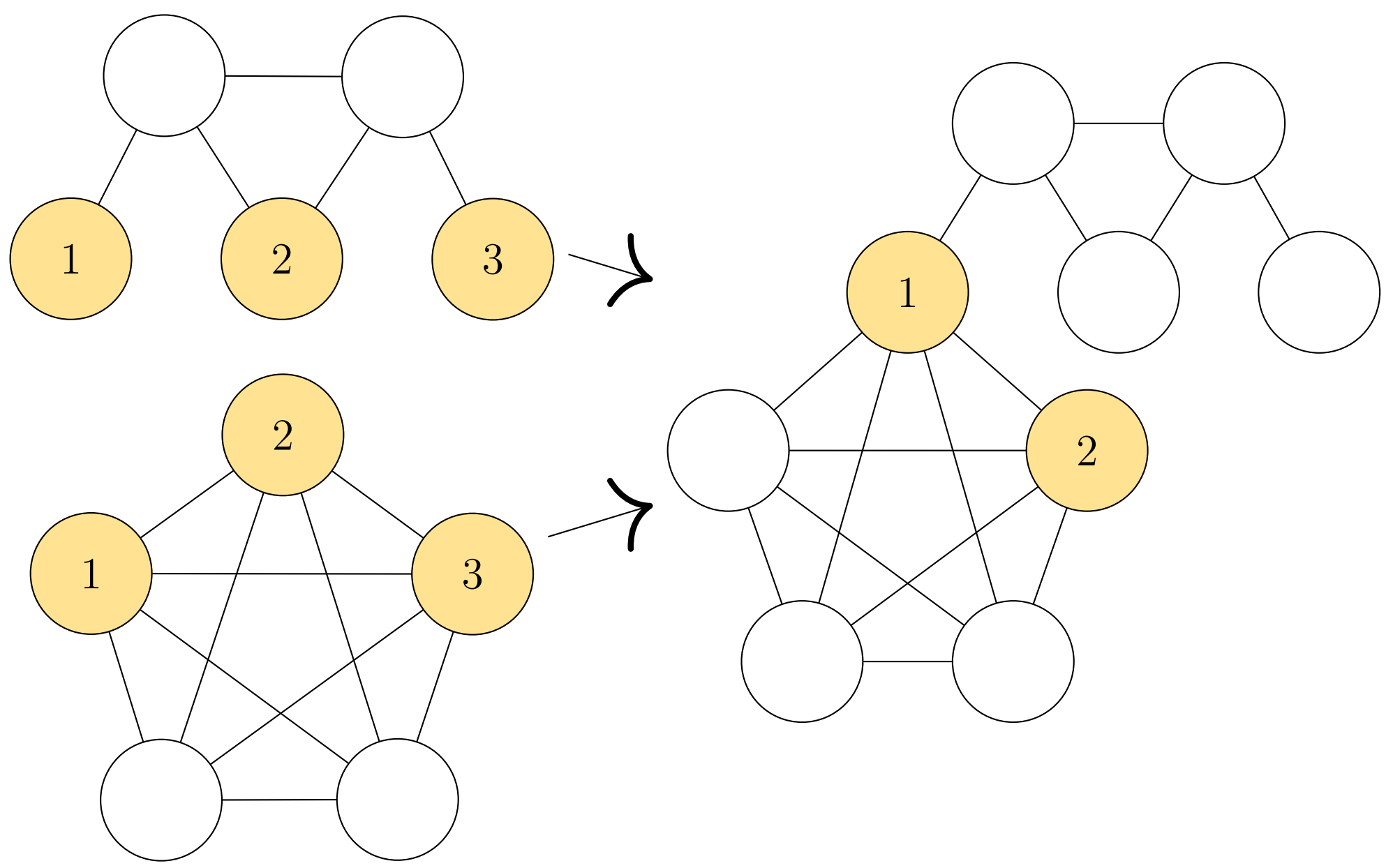

We turn this idea into a formal definition, which the reader may want to check is indeed inherited from $\mathfrak{G}$ (as a composition of existing function symbols).

$\mathfrak{G}$ is the algebra of graphs. (def. 6.2.2, p. 93)

### Definition 9.3.1

$\oplus$ is the disjoint sum. (def. 4.5.3, p. 40)

媲 (tsureai, Japanese for *to marry*) denotes the source fusion. (def. 4.5.5, p. 42)

$\leftrightarrows$ denotes the terminal redefinition. (def. 4.5.4, p. 41)

Let $G$ be a graph of type $n$ and let $G'$ be a graph of type $m$. Let further

$$K = \{ l_1, \dots, l_{|K|} \} \subseteq \{ 1, \dots, \min\{ n, m \} \}$$

and let $k \in \mathbb{N}$, $k \leq n + m - |K|$. We call

$$G \;{}_{m}^{n}\!\otimes_{K}^{k}\, G' := \leftrightarrows_{\sigma} \text{媲}_{l_{|K|}}^{l_{|K|+n}} \text{媲}_{l_{|K|-1}}^{l_{|K|-1+n}} \dots \text{媲}_{l_1}^{l_{1+n}} (G \oplus G')$$

with

$$\sigma \colon \{ 1, \dots, k \} \to \{ 1, \dots, n + m \}, i \mapsto \begin{cases} i & \text{if } i \leq n \\ x_i & \text{otherwise,} \end{cases}$$

where $x_i$ denotes the $i$-th smallest element of

$$\{1+n, \dots, m+n\} \setminus \{x+n : x \in K\},$$

the *$k$-twine* of $G$ and $G'$ *over* $K$.

In words, entwining puts the two graphs side by side and identifies with each other some of the equal-numbered terminals, but without adding to their multiplicity (as ordinary fusion does). If $k < n+m-|K|$, it then throws away the "superfluous" terminal vertices.

Note that entwining does not allow one to create loops.

We now define the carrier sets of our algebra. They will not be *exactly* the sets of graphs of a certain tree-width (the proof of theorem 9.3.6 reveals why), though the latter are a good starting point.

### Definition 9.3.2

Let $k, n \in \mathbb{N}$. We denote

$$\mathfrak{t}_n^k := \{G \in \mathfrak{G}_n : \mathrm{tw}(G) \leq k\}.$$

$\mathfrak{G}_n$ denotes the set of all graphs of type $n$. (def. 4.5.1, p. 39)

The carrier sets shall be as follows.

### Definition 9.3.3

Let $k \in \mathbb{N}$. A tree-decomposition for a graph $G \in \mathfrak{G}$ is called *$k$-verdant* if it has width at most $k$ and there is a bag which contains all terminal vertices of $G$.

$\mathfrak{G}$ is the algebra of graphs. (def. 6.2.2, p. 93)

A graph $G \in \mathfrak{G}$ is called *$k$-verdant* if it admits a $k$-verdant tree-decomposition.

For $n \in \mathbb{N}$, the set of all $k$-verdant graphs of type $n$ is denoted by $\mathfrak{T}_n^k$.

In particular, we have

$$\forall k \in \mathbb{N} \colon \mathfrak{T}_0^k = \mathfrak{t}_0^k,$$

which, if we are honest with ourselves, is all that we *really* care about.

We note the obvious inclusions in one place.

### Observation 9.3.4

$\mathfrak{t}_n^k$ denotes the type $n$ graphs of tree-width at most $k$. (def. 9.3.2, p. 185)

$\mathfrak{T}_n^k$ denotes the $k$-verdant graphs of type $n$. (def. 9.3.3, p. 185)

Let $k, n \in \mathbb{N}$. Then we have the following.

- $\mathfrak{t}_n^k \subseteq \mathfrak{t}_n^{k+1}$.
- $\mathfrak{T}_n^k \subseteq \mathfrak{T}_n^{k+1}$.
- $\mathfrak{T}_n^k \subseteq \mathfrak{t}_n^k$.
- $\mathfrak{T}_0^k = \mathfrak{t}_0^k$.
- $\mathfrak{T}_1^k = \mathfrak{t}_1^k$.

For graphs of type larger than 1, the second and third inclusions are strict – for example, the following type 2 graph of tree-width 1 is 2-verdant but not 1-verdant.

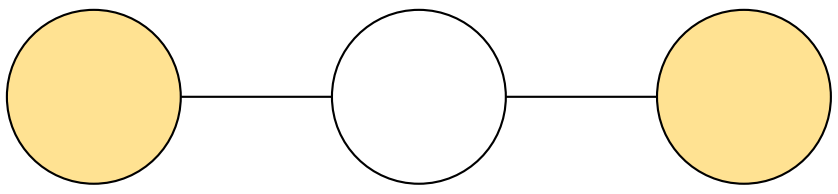

We take the trivial graphs, all entwinements, and all terminal redefinitions, and throw them into an algebra with the $k$-verdant graphs.

### Definition 9.3.5

$\mathfrak{v}$ is the type 1 graph with one vertex. (def. 4.5.6, p. 43)

$\mathfrak{e}_1$ is the type 1 graph with 1 vertices and one edge. (def. 4.5.6, p. 43)

$\otimes$ denotes the graph twine. (def. 9.3.1, p. 184)

$\overline{\mathfrak{D}}$ denotes the closure of $\mathfrak{D}$. (def. 6.10.6, p. 141)

$\mathrm{Fl}(\overline{\mathfrak{D}})$ denotes the flattening of $\overline{\mathfrak{D}}$. (def. 6.10.8, p. 143)

Let $k \in \mathbb{N}$. We set

$$\mathfrak{D}^k := \mathfrak{D}^k_{\text{triv}} \cup \mathfrak{D}^k_{\otimes} \cup \left( \bigcup_{i \in \mathbb{N}} \bigcup_{j \in \mathbb{N}} {}^j_i\mathfrak{F}_{\leftrightarrows} \right)$$

with

$$\mathfrak{D}^k_{\text{triv}} := \{ \mathfrak{v}, \mathfrak{e}_1, \dots, \mathfrak{e}_{k+1} \}$$

and

$$\mathfrak{D}^k_{\otimes} := \{ {}^n_m\otimes^l_K : n, m, l \in \mathbb{N}, K \subseteq \{ 1, \dots, \min\{ n, m \} \}, \\ l \leq n + m - |K|, l \leq k + 1 \},$$

while ${}^j_i\mathfrak{F}_{\leftrightarrows}$ is as in definition 6.2.1. We denote

$$\mathfrak{T}^k := \left( \{\mathfrak{T}^k_n\}_{n \in \mathbb{N}}, \{ \mathrm{Fl}(\overline{\mathfrak{D}})_{\mathfrak{f}} |_{(\bigcup_{n \in \mathbb{N}} \mathfrak{T}^k_n)^*} \}_{\mathfrak{f} \in \mathfrak{D}^k} \right)$$

and call this the *algebra of graphs of tree-width at most $k$ (without loops)*.

We can call this an algebra all we want, but we should probably show that it deserves this name.

**Theorem 9.3.6**

Let $k \in \mathbb{N}$, and let $\mathfrak{f} \in \mathfrak{O}^k$. Then the range of $\mathrm{Fl}(\overline{\mathfrak{O}})_{\mathfrak{f}}|_{(\bigcup_{n\in\mathbb{N}} \mathfrak{T}_n^k)^*}$ lies in $\mathfrak{T}^k_{\langle \mathfrak{f}\rangle^{\text{out}}}$.

$\mathfrak{O}^k$ is the set of function symbols of $\mathfrak{T}^k$. (def. 9.3.5, p. 186)

$\overline{\mathfrak{O}}$ denotes the closure of $\mathfrak{O}$. (def. 6.10.6, p. 141)

$\mathrm{Fl}(\overline{\mathfrak{O}})$ denotes the flattening of $\overline{\mathfrak{O}}$. (def. 6.10.8, p. 143)

$\mathfrak{T}_n^k$ denotes the $k$-verdant graphs of type $n$. (def. 9.3.3, p. 185)

$\leftrightarrows$ denotes the terminal redefinition. (def. 4.5.4, p. 41)

**Proof.** Let $k \in \mathbb{N}$. For each possible function symbol, we verify that plugging in only $k$-verdant graphs yields again a $k$-verdant graph.

Let $\mathfrak{f} \in \mathfrak{O}^k$.

*Case 1: $\mathfrak{f}$ is nullary.* None of the nullary function symbols in $\mathfrak{O}^k$ yields a graph with more than $k+1$ vertices, so a bag consisting of all vertices is of width at most $k$, and it trivially contains all terminal vertices.

*Case 2: $\mathfrak{f} = \leftrightarrows\sigma$ for some function $\sigma$.* Let $G \in \mathfrak{T}_n^k$ for some $n \in \mathbb{N}$ and let $\sigma\colon \{1,\dots,n'\} \to \{1,\dots,n\}$ for some $n' \in \mathbb{N}$. Take the tree-decomposition for $G$ of width at most $k$ which has all $n$ terminals of $G$ in some bag $x$.

Terminal redefinition does not change vertices or edges, so the same tree-decomposition is also a tree-decomposition for $\leftrightarrows_\sigma G$. Terminal redefinition can not promote nonterminal vertices to terminals, so all terminal vertices of $\leftrightarrows_\sigma G$ are again in $x$.

*Case 3: $\mathfrak{f} = {}^{n}_{m}\!\otimes^K_l$ for some set $K$ and some $n, m, l \in \mathbb{N}$.* This is the interesting of the three cases: let $G = (V, E, t) \in \mathfrak{T}_n^k$, let $G' \in \mathfrak{T}_m^k$, and let $K \subseteq \{1,\dots,\min\{n,m\}\}$, $l \le m+n-|K|$ and $l \le k+1$.

$\otimes$ denotes the graph twine. (def. 9.3.1, p. 184)

Let $(T, X, b)$ be a tree-decomposition for $G$ of width at most $k$ with all terminals of $G$ contained in a single bag $x \in X$. Let similarly $(S, Y, c)$ be a tree-decomposition for $G'$ of width at most $k$ with all terminals of $G'$ contained in a single bag $y \in Y$.

Without loss of generality, we assume that $K = \{1,\dots,p\}$ for some $p \in \mathbb{N}$ with $p \le \min\{n,m\}$.

We set $a \coloneqq \min\{\, n, l \,\}$ and

$$z \coloneqq \{\, t(i) : i \leq a \,\} \cup \{\, t'(i+p) : i \in \{\, a+1, \dots, l \,\} \,\},$$

$\otimes$ denotes the graph twine. (def. 9.3.1, p. 184)

in other words, exactly the vertices that will be terminal in $G \otimes_l^K G'$. Note that $|z| \leq l \leq k+1$.[5]

We introduce a new node $r$ and attach to it as subtrees the trees $T$ and $S$ at the nodes $b^{-1}(x)$ and $c^{-1}(y)$. Call this new tree $U$.

We now set

$$d\colon U \to X \cup Y \cup \{z\}, v \mapsto \begin{cases} b(v) & v \in T \\ c(v) & v \in S \\ z & v = r. \end{cases}$$

We claim that this makes $(U, X \cup Y \cup \{z\}, d)$ into a $k$-verdant tree-decomposition for $G \otimes^K G'$.

The terminal vertices of $G \otimes^K G'$ are by definition all in $z$. Hence if $(U, X \cup Y \cup \{z\}, d)$ is a tree-decomposition for $G \otimes^K G'$, it is immediately $k$-verdant (since $z$ contains at most $k+1$ vertices).

Since $(T, X, b)$ and $(S, Y, c)$ were tree-decompositions, all vertices of $G$ and $G'$ are contained in our new decomposition, and hence so are all vertices of $G \otimes^K G'$.

Entwining introduces no edges, so every edge of $G \otimes^K G'$ is still contained in a bag.

Let now $v \in G \otimes^K G'$ and consider the subgraph $U'$ of $U$ induced by $\{\, v' \in U : v \in d(v') \,\}$.

*Case 3.1: $v \notin z$.* Then $U'$ is entirely contained in either $T$ or $S$ and hence connected since they are tree-decompositions.

*Case 3.2: $v \in z$ and $v \notin \{\, t(i) : i \in K \,\}$.* Without loss of generality, $v \in G$. Since $v$ is a terminal of $G$ (by virtue of being in $z$), we have $v \in x$. Since $v$ was not glued, we have $v \notin G'$.

Hence $U' \cap S = \varnothing$ and $U'$ is connected since $T$ was a tree-decomposition.

[5] Strict inequality can occur if either graph has terminals of multiplicity larger than 1.

*Case 3.3:* $v \in \{ t(i) : i \in K \}$. Then $v$ was a terminal of $G$ and of $G'$ and is hence in $x$ and in $y$. It is also in $z$, whence $U'$ is connected.

We have thus shown that $(U, X \cup Y \cup \{z\}, d)$ is indeed a tree-decomposition for $G \otimes^K G'$.

$\otimes$ denotes the graph twine. (def. 9.3.1, p. 184)

This concludes the proof for all allowed function symbols.

□

This is all we need to note the following.

**Corollary 9.3.7**

Let $k \in \mathbb{N}$. Then $\mathfrak{T}^k$ is an inherited $(\mathbb{N}, \mathfrak{O}^k)$-algebra of $\mathfrak{G}$.

$\mathfrak{T}^k$ is the algebra of $k$-verdant graphs. (def. 9.3.3, p. 185)

A short examination of definition 6.2.1 also reveals that there are not too many function symbols of each type.

$\mathfrak{O}^k$ is the set of function symbols of $\mathfrak{T}^k$. (def. 9.3.5, p. 186)

**Observation 9.3.8**

Let $k \in \mathbb{N}$. Then $\mathfrak{T}^k$ is of weakly locally finite signature.

$\mathfrak{G}$ is the algebra of graphs. (def. 6.2.2, p. 93)

It remains to show that this algebra is finitely expressible – otherwise, all of our hard work would have been in vain.

**Theorem 9.3.9**

Let $k, n \in \mathbb{N}$. Then for any $k$-verdant graph $G \in \mathfrak{G}_n$, there exists a $\{ 0, \dots, \max\{ n, k+1 \} \}$-local expression $e \in \|(\mathbb{N}, \mathfrak{O}^k)\|$ with $\mathrm{val}_{\mathfrak{T}^k} e = G$.

$\mathfrak{G}_n$ denotes the set of all graphs of type $n$. (def. 4.5.1, p. 39)

**Proof.** Let $k, n \in \mathbb{N}$, and let $G = (G', t) \in \mathfrak{G}_n$ be $k$-verdant.

$\|(\mathbb{N}, \mathfrak{O}^k)\|$ denotes the set of expressions over $(\mathbb{N}, \mathfrak{O}^k)$. (def. 6.3.3, p. 103)

Since $G$ is $k$-verdant, we can find a tree-decomposition of width at most $k$ and some bag $x$ which contains all terminal vertices of $G$. We root this tree-decomposition at $x$. By theorem 4.6.5, we can then find a nice tree-decomposition $(T, X, b)$ for $G$ of width at most $k$ with the same root bag. In particular, all terminal vertices of $G$ are in $b(\sqrt{T})$.

$\mathrm{val}_{\mathfrak{T}^k} e$ denotes the result of $e$ when evaluated in $\mathfrak{T}^k$. (def. 6.3.3, p. 103)

We denote for $v \in T$ by $G_v$ the induced subgraph $G'[X]$, where

$$X := \bigcup_{v' \in T[v]} bv',$$

that is, the underlying (untyped) graph induced by the vertices that one can find from the bag $v$ “downward”.

Lastly, for $v \in T$ with $b(v) = \{ x_1, \dots, x_l \}$ for some $l \in \mathbb{N}$, we denote by $t_v$ the function $\{ 1, \dots, l \} \to G_v, i \mapsto x_i$.

We now make use of the following property: we call a node $v \in T$ *friendly* if for *every* set $E$ of edges in $G_v$, there exists a $\{ 0, \dots, \max\{ n, k+1 \} \}$-local expression $\boldsymbol{e} \in \|(\mathbb{N}, \mathfrak{O}^k)\|$ with $\mathrm{val}_{\mathfrak{T}^k} \boldsymbol{e} = (G_v - E, t_v)$.

$\mathfrak{O}^k$ is the set of function symbols of $\mathfrak{T}^k$. (def. 9.3.5, p. 186)

In other words, a node is friendly if we can construct the subgraph it induces in the tree-decomposition of $G$ *and* all versions of that subgraph with some or all of its edges removed.

$\|(\mathbb{N}, \mathfrak{O}^k)\|$ denotes the set of expressions over $(\mathbb{N}, \mathfrak{O}^k)$. (def. 6.3.3, p. 103)

$\mathrm{val}_{\mathfrak{T}^k} \boldsymbol{e}$ denotes the result of $\boldsymbol{e}$ when evaluated in $\mathfrak{T}^k$. (def. 6.3.3, p. 103)

$\mathfrak{T}^k$ is the algebra of $k$-verdant graphs. (def. 9.3.3, p. 185)

The reader should first convince themselves that if we can show that the root node $\sqrt{T}$ is friendly, then our proof is complete: we have $G' = G_r$ since every vertex of $G$ must be contained in some bag of the decomposition and all bags can be reached from $\sqrt{T}$. If we can construct the graph $(G', t_{\sqrt{T}})$, then every vertex in $b(\sqrt{T})$ is a terminal vertex of $(G', t_{\sqrt{T}})$. Since all terminal vertices of $G = (G', t)$ are contained in $b(\sqrt{T})$, we can obtain $G$ from $(G', t_{\sqrt{T}})$ by one terminal redefinition. Because $b(\sqrt{T})$ can not contain more than $k+1$ vertices, this terminal redefinition is $\{ 0, \dots, \max\{ n, k+1 \} \}$-local. Hence if we show that $\sqrt{T}$ is friendly, we have shown that $G$ can indeed be obtained as the result of a $\{ 0, \dots, \max\{ n, k+1 \} \}$-local expression in $\mathfrak{T}^k$.

We now show that if all children of a node $v \in T$ are friendly, then so is $v$ itself. By induction, the root $r$ of $T$ must then be friendly.

Let thus $v \in T$.

*Case 1: v is a leaf.* Then $G_v$ contains exactly one vertex and no edges, so it is the result of evaluating the function symbol $\mathfrak{v}$, which is $\{ 0, 1 \}$-local.

$\mathfrak{v}$ is the type 1 graph with one vertex. (def. 4.5.6, p. 43)

*Case 2: v is a forget node.* Let $v'$ be the sole child of $v$. Since $b(v) \subseteq b(v')$, this implies that $G_v = G_{v'}$. The terminal vertices are corrected by a simple terminal redefinition. Since before we correct the root, no terminal vertex we use has multiplicity more than one, this is $\{ 0, \dots, k+1 \}$-local because no bag can contain more than $k+1$ vertices.

*Case 3: v is an introduce node.* Let $v'$ be the sole child of $v$, let the corresponding bag be $b(v') = \{ x_1, \dots, x_l \}$ for some $l \in \mathbb{N}$, and let the bag of $v$ be $b(v) = b(v') \cup \{x_{l+1}\}$.

Note that since $(T, X, b)$ is of width at most $k$, we have $l + 1 \leq k + 1$ (and $l \leq k < k + 1$), meaning that the graph $(G_v, t_v)$, once constructed, is $k$-verdant.

Let now $E$ be any set of edges of $G_v$. As a shorthand, we write $G_E \coloneqq G_v - E$ for the graph $G_v$ with those edges removed.

By the induction hypothesis, we can construct the type $l$ graph

$$(G_E - \{x_{l+1}\}, t_{v'}) = (G_{v'} - E, t_{v'})\,.$$

Consider now the neighbours of $x_{l+1}$ in $G_v$. Take any vertex $y \in G_{v'}$ such that $x_{l+1}$ and $y$ are connected by an edge in $G_v$. Then the tree-decomposition $(T, X, b)$ must contain a node $v''$ such that $x_{l+1}, y \in b(v'')$. This node cannot be in $T[v']$ since $x_{l+1} \notin b(v')$. But we also have $y \in G_{v'}$, whence we must have $y \in b(v)$ by the definition of a tree-decomposition. Because $y \neq x_{l+1}$, this implies $y \in b(v')$.

We know thus that every vertex in $G_v$ that is adjacent to $x_{l+1}$ is already contained in $b(v')$; in particular, it is a terminal vertex of $\big(G_E - \{x_{l+1}\}, t_{v'}\big)$.

We denote $G'_E \coloneqq \big(G_E - \{x_{l+1}\}, t_{v'}\big) \oplus \mathfrak{v}$, which is a $\{\,0, \ldots, k+1\,\}$-local construction since the former has at most $k$ terminal vertices. The vertex added by adjoining $\mathfrak{v}$ is of course the vertex $x_{l+1}$.

$\oplus$ is the disjoint sum. (def. 4.5.3, p. 40)

Now for every edge $e$ incident to $x_{l+1}$ in $G_E$, take the trivial graph $\mathfrak{e}_{|e|}$ and attach it by graph entwining.

$\mathfrak{v}$ is the type 1 graph with one vertex. (def. 4.5.6, p. 43)

*Case 4: v is a join node.* Let us denote the children of $v$ by $v_1$ and $v_2$. By induction, assume that the graphs $G_{v_1}$ and $G_{v_2}$ have been constructed such that the terminals of $G_{v_1}$ are exactly the vertices of $b(v_1)$ and the terminals of $G_{v_2}$ are exactly the vertices of $b(v_2)$.[6]

$\mathfrak{e}_{|e|}$ is the type $|e|$ graph with $|e|$ vertices and one edge. (def. 4.5.6, p. 43)

Suppose we have a set $E$ of edges in $G_v$ and want to construct $G_v - E$. The vertices of $G_v$ are of course the same as those in $G_{v_1} \cup G_{v_2}$. Furthermore, $G_v$ cannot contain edges that are not contained in $G_{v_1}$ or $G_{v_2}$: consider an edge $e$ of $G_v$. All of its end points must be contained in some bag $x$. Either $x$ is contained in the left or right subtree of $v$, in which case $e$ lies in the corresponding graph, or $x$ is a predecessor of $v$, in which case, since $e$ is contained in $G_v$, every end point of $e$ must be in one of the subtrees and hence (because we have a tree-decomposition) also in $b(v)$. But since $b(v) = b(v_1)$, the latter case also implies that $e$ lies in $G_{v_1}$.

[6] Remember that while $b(v_1)$ and $b(v_2)$ contain the same vertices, the graphs $G_{v_1}$ and $G_{v_2}$ are in general not the same.

Consider now an edge $e$ which is in $G_v$ but not in $G_{v_1}$. By the now-familiar tree-decomposition argument, this edge (contained in $G_{v_2}$) cannot be incident to any vertex in $G_{v_1} \setminus b(v_1)$.

Denote by $E'$ the set of all edges contained in $G_{v_1}$ *and* $G_{v_2}$. Per the induction hypothesis, we can construct the graphs $G_{v_1} - E$ and $G_{v_2} - E - E'$ with terminals being the vertices in $b(v)$. By construction, we know that every edge of $G_v$ is contained in exactly one of those graphs and that the only vertices they have in common are the ones in $b(v)$, which are all terminals and of which there are no more than $k+1$. Hence a $\{0, \dots, k+1\}$-local entwinement yields the desired result.

We have thus shown that every node of $T$ is friendly, in particular the root of $T$, concluding the proof.

$\mathfrak{T}^k$ is the algebra of $k$-verdant graphs. (def. 9.3.3, p. 185)

□

$\mathfrak{G}$ is the algebra of graphs. (def. 6.2.2, p. 93)

We restate what we have learnt.

### Corollary 9.3.10

$\mathfrak{O}^k$ is the set of function symbols of $\mathfrak{T}^k$. (def. 9.3.5, p. 186)

Let $k \in \mathbb{N}$. Then $\mathfrak{T}^k$ is a finitely expressible inherited algebra of $\mathfrak{G}$ of weakly locally finite signature.

Hence Courcelle's Theorem becomes for this special case the following.

$\mathfrak{T}^k_n$ denotes the $k$-verdant graphs of type $n$. (def. 9.3.3, p. 185)

### Corollary 9.3.11

$\mathrm{val}_{\mathfrak{T}^k} e$ denotes the result of $e$ when evaluated in $\mathfrak{T}^k$. (def. 6.3.3, p. 103)

Let $\varphi$ be a sentence of the monadic second-order language of graphs, and let $k \in \mathbb{N}$. Then for every $n \in \mathbb{N}$, there exists a deterministic bottom-up finite tree automaton that takes as input a set $\mathscr{H}$ of $(\mathbb{N}, \mathfrak{O}^k)$-expressions with $\mathfrak{T}^k$-values in $\mathfrak{T}^k_n$ and accepts an expression $e$ if and only if $\models_{\Omega(\mathrm{val}_{\mathfrak{T}^k} e)} \varphi$, and every element of $\mathfrak{T}^k_n$ is the value of some expression in $\mathscr{H}$.

$\Omega(\mathrm{val}_{\mathfrak{T}^k} e)$ is the induced structure of $\mathrm{val}_{\mathfrak{T}^k} e$. (def. 5.4.2, p. 73)

We already know from algorithm 2 that once a graph expression is given, applying the theorem can be done in time linear in the number of vertices of said expression. What remains to analyse is how big that number of vertices is, and how long it takes to construct said expression in the first place.

We turn the procedure from the proof of theorem 9.3.9 into an algorithm. Recall how, in the proof, we relied on the fact that for any node $v$ of the nice tree-decomposition, we could build an expression for a certain subgraph of $G$ that we then recursively glue together.

The following subroutine does just that: given a node of the tree-decomposition and a subgraph we want to build, it recurses into the tree-decomposition to build an expression for that subgraph.

**Algorithm 5:** The Building Routine

**Input:** a tree-decomposition, a node in that decomposition, and a typed graph we want to build
**Output:** an expression resulting in that graph

BUILD $(Z, v, (G, t))$:
  **if** $v$ is a leaf **then**
    **return** $\mathfrak{v}$
  **else if** $v$ is a forget node with successor $v'$ **then**
    write $b(v) = \{ x_1, \ldots, x_l \}$, $b(v') = b(v) \cup \{x_{l+1}\}$
    set $\sigma\colon \{ 1, \ldots, l \} \hookrightarrow \{ 1, \ldots, l+1 \}, i \mapsto i$
    set $t'\colon \{ 1, \ldots, l+1 \} \to b(v'), i \mapsto \begin{cases} t(i) & \text{if } i < l+1 \\ x_{l+1} & \text{if } i = l+1 \end{cases}$
    **return** $\leftrightarrows_\sigma \textsc{Build}(Z, v', (G, t'))$
  **else if** $v$ is an introduce node with successor $v'$ **then**
    write $b(v') = \{ x_1, \ldots, x_l \}$, $b(v) = b(v') \cup \{x_{l+1}\}$
    $\boldsymbol{e} := \textsc{Build}(Z, v', (G - \{x_{l+1}\}, i \mapsto x_i))$
    $\boldsymbol{e} := \boldsymbol{e} \oplus \mathfrak{v}$
    **for** every edge $e$ incident to $x_{l+1}$ in $G$ **do**
      for $\langle\!| e |\!\rangle = x_i x_j \ldots$, write $K = \{ i, j, \ldots \}$
      **if** $G$ is directed **then**
        set $\sigma\colon \{ 1, \ldots, |e| \} \to \{ 1, \ldots, l+1 \}, 1 \mapsto i, 2 \mapsto j, \ldots$
        $\boldsymbol{e} := \boldsymbol{e} \otimes_K^{l+1} \leftrightarrows_\sigma \mathfrak{e}_{|e|}$
      **else**
        $\boldsymbol{e} := \boldsymbol{e} \otimes_K^{l+1} \mathfrak{e}_{|e|}$
    set $\sigma\colon \{ 1, \ldots, l+1 \} \to \{ 1, \ldots, l+1 \}, i \mapsto x_{t(i)}$
    **return** $\leftrightarrows_\sigma \boldsymbol{e}$
  **else if** $v$ is a join node with successors $v_1$ and $v_2$ **then**
    $E' :=$ the set of edges contained in both $G_{v_1}$ and $G_{v_2}$
    # *the notation $G_{v_1}$ is as in the proof of theorem 9.3.9:*
    # *the graph induced by all vertices in bags reachable from $v_1$*
    $\boldsymbol{e}_1 := \textsc{Build}(Z, v_1, G_{v_1})$, $\boldsymbol{e}_2 := \textsc{Build}(Z, v_2, G_{v_2} - E)$
    set $\sigma\colon \{ 1, \ldots, l+1 \} \to \{ 1, \ldots, l+1 \}, i \mapsto x_{t(i)}$
    **return** $\leftrightarrows_\sigma(\boldsymbol{e}_1 \otimes_{\{ 1, \ldots, |b(v)| \}}^{|b(v)|} \boldsymbol{e}_2)$

$\mathfrak{v}$ is the type 1 graph with one vertex. (def. 4.5.6, p. 43)

$\oplus$ is the disjoint sum. (def. 4.5.3, p. 40)

$\otimes$ denotes the graph twine. (def. 9.3.1, p. 184)

$\mathfrak{e}_{|e|}$ is the type $|e|$ graph with $|e|$ vertices and one edge. (def. 4.5.6, p. 43)

$\leftrightarrows$ denotes the terminal redefinition. (def. 4.5.4, p. 41)

Once we have this subroutine, building an expression for the original graph is a simple matter of applying it to the root of the tree-decomposition.

**Algorithm 6:** Expressing Yourself

**Input:** a graph $(G, t) \in \mathfrak{T}^k_n$
**Output:** an expression evaluating to $(G, t)$
compute a nice $k$-verdant tree-decomposition $Z = (T, X, b)$ of $G$
write $b(\sqrt{T}) = \{ x_1, \dots, x_l \}$
set $t' \colon \{ 1, \dots, k \} \to b(\sqrt{T}), i \mapsto x_i$
$e :=$ BUILD$(Z, \sqrt{Z}, (G, t'))$
set $\sigma\colon \{ 1, \dots, n \} \to \{ 1, \dots, k \}$ such that $t(i) = \sigma(t'(i))$
**return** $\leftrightarrows_\sigma e$

$\mathfrak{T}^k_n$ denotes the $k$-verdant graphs of type $n$. (def. 9.3.3, p. 185)

The BUILD routine is defined on page 194.

$\leftrightarrows$ denotes the terminal redefinition. (def. 4.5.4, p. 41)

The reader who has understood (and not just read) the proof of theorem 9.3.9 should have no problems with this algorithm.

**Theorem 9.3.12**

Let $n, k \in \mathbb{N}$, and let $G \in \mathfrak{T}^k_n$. Then algorithm 6 terminates and returns an $\|(\mathbb{N}, \mathfrak{O}^k)\|$-expression $e$ with $\mathrm{val}_{\mathfrak{T}^k} e = G$.

$\mathfrak{O}^k$ is the set of function symbols of $\mathfrak{T}^k$. (def. 9.3.5, p. 186)

**Proof.** This is straight-forward from the proof of theorem 9.3.9.

□

$\|(\mathbb{N}, \mathfrak{O}^k)\|$ denotes the set of expressions over $(\mathbb{N}, \mathfrak{O}^k)$. (def. 6.3.3, p. 103)

Of course, what we really want to know is the runtime, which hinges on our ability to quickly compute a verdant tree-decomposition.

Luckily, the following result by Bodlaender helps us out.

$\mathrm{val}_{\mathfrak{T}^k} e$ denotes the result of $e$ when evaluated in $\mathfrak{T}^k$. (def. 6.3.3, p. 103)

**Theorem 9.3.13**

Let $k \in \mathbb{N}$. Then there exists an algorithm which, for any given graph $G = (V, E, \langle\!\langle _ \rangle\!\rangle)$ of tree-width at most $k$, computes a tree-decomposition of width at most $k$ for $G$ in time $\mathcal{O}(|V| + |E|)$. This tree-decomposition has at most $\mathcal{O}(|V| + |E|)$ bags.

$\mathfrak{T}^k$ is the algebra of $k$-verdant graphs. (def. 9.3.3, p. 185)

**Proof.** That such a linear-time algorithm exists is the main result of [Bod96].

Since the algorithm runs in linear time, it cannot output superlinearly many bags.

□

The reader should check the following fact by looking at the definition of a tree-decomposition.

**Observation 9.3.14**

Let $G = (V, E, \langle\!| _ |\!\rangle)$ be a graph, $Z$ a tree-decomposition for $G$, and let $E' \subseteq E$. Then $Z$ is also a tree-decomposition for $G - E'$.

We use this to apply Bodlaender's result to verdant tree-decompositions.

**Theorem 9.3.15**

$\mathfrak{T}^n_k$ denotes the $n$-verdant graphs of type $k$. (def. 9.3.3, p. 185)

Let $n, k \in \mathbb{N}$. Then there exists an algorithm which, for any given graph $G = (V, E, \langle\!| _ |\!\rangle, t) \in \mathfrak{T}^n_k$, computes a $k$-verdant tree-decomposition with at most $\mathcal{O}(|V| + |E|)$ bags for $G$ in time $\mathcal{O}(|V| + |E|)$.

**Proof.** Let $k, n \in \mathbb{N}$, and let $G = (V, E, \langle\!| _ |\!\rangle, t) \in \mathfrak{T}^n_k$. In particular, $G$ is $k$-verdant.

Denote the set of terminal vertices of $G$ as $\{ v_1, \dots, v_l \}$ for some $l \in \mathbb{N}$ (implying that $l \leq k + 1$, $l \leq n$) and by $G'$ the graph $G$ with an additional edge $e$ with $\langle\!| e |\!\rangle = v_1 \dots v_l$.

Consider a $k$-verdant tree-decomposition for $G$. This tree-decomposition has a bag which contains $v_1, \dots, v_l$ and hence all vertices incident to $e$. Thus, it is also a tree-decomposition for $G'$, showing that $G'$ is $k$-verdant as well.

The algorithm by Bodlaender (theorem 9.3.13) enables us to compute in time $\mathcal{O}(|V| + |E \cup \{e\}|) = \mathcal{O}(|V| + |E|)$ a tree-decomposition $Z$ for $G'$ of width at most $k$ with at most $\mathcal{O}(|V| + |E \cup \{e\}|) = \mathcal{O}(|V| + |E|)$ bags. Since $G'$ contains an edge incident to all vertices in $\{ v_1, \dots, v_l \}$, $Z$ must have a bag containing $\{ v_1, \dots, v_l \}$. Therefore, $Z$ is $k$-verdant, and by observation 9.3.14, it is also a tree-decomposition for $G$.

□

**Theorem 9.3.16**

Let $n, k \in \mathbb{N}$. Then there exists an algorithm which, for any given graph $G = (V, E, \langle\!\langle _ \rangle\!\rangle, t) \in \mathfrak{T}^n_k$, computes a nice $k$-verdant tree-decomposition with at most $\mathcal{O}(|V| + |E|)$ bags for $G$ in time $\mathcal{O}(|V| + |E|)$.

$\mathfrak{T}^n_k$ denotes the $n$-verdant graphs of type $k$. (def. 9.3.3, p. 185)

**Proof.** Let $n, k \in \mathbb{N}$, and let $G = (V, E, \langle\!\langle _ \rangle\!\rangle, t) \in \mathfrak{T}^k_n$. Compute a $k$-verdant tree-decomposition for $G$ in linear time with linearly many bags by theorem 9.3.15.

Root this tree-decomposition at a bag containing all terminal vertices and turn it into a nice tree-decomposition in linear time and with (again) linearly many bags by theorem 4.6.5. Since theorem 4.6.5 preserves the root bag, this nice tree-decomposition is again $k$-verdant.

□

We can now analyse the runtime of algorithm 6.

**Theorem 9.3.17**

Let $n, k \in \mathbb{N}$, and let $G = (V, E, \langle\!\langle _ \rangle\!\rangle, t) \in \mathfrak{T}^k_n$. Then algorithm 6 terminates in time $\mathcal{O}(|V| + |E|)$ and returns an expression with $\mathcal{O}(|V| + |E|)$ vertices.

**Proof.** Let $n, k \in \mathbb{N}$, and let $G = (V, E, \langle\!\langle _ \rangle\!\rangle, t) \in \mathfrak{T}^k_n$. Algorithm 6 begins by computing a nice $k$-verdant tree-decomposition for $G$, which by theorem 9.3.16 happens in linear time, and can be done such that the tree-decomposition has at most linearly many bags.

All further steps except the call to the Build subroutine take constant time.

The Build routine is defined on page 194.

The Build subroutine calls itself recursively, but only once for each bag.

It remains to analyse how much time one call to Build takes.

If the node considered is not an introduce node, then Build runs in constant time. If it is an introduce node, it loops over all edges incident to the vertex it introduces. Because $G$ is $k$-verdant, no edge can be incident to more than $k$ vertices, hence no edge can be visited by more than $k$ such Build

The BUILD routine is defined on page 194.

calls (because, of course, no two introduce nodes can introduce the same vertex). Hence, adding up all BUILD calls comes to a runtime of $\mathcal{O}(k \cdot |E| + a)$, where $a$ denotes the number of bags of the nice tree-decomposition, and because there are linearly many bags, we have $\mathcal{O}(k \cdot |E| + a) \subseteq \mathcal{O}(|V| + |E|)$.

This already proves the claim.

□

We can hence prove the main result of Courcelle's Theorem applied to graphs of constantly bounded tree-width.

**Theorem 9.3.18**

Let $\varphi$ be a sentence of the monadic second-order language of graphs, and let $k \in \mathbb{N}$. Then for every $n \in \mathbb{N}$, there exists an algorithm which, given a $k$-verdant type $n$ graph $G = (V, E, \langle\!\!|_|\!\!\rangle, t)$, decides in time $\mathcal{O}(|V| + |E|)$ whether or not we have $\models_{\Omega(G)} \varphi$.

$\Omega(G)$ is the induced structure of $G$. (def. 5.4.2, p. 73)

**Proof.** Let $\varphi$ be a sentence of the monadic second-order language of graphs, and let $k, n \in \mathbb{N}$. By corollary 9.3.11, there exists a deterministic bottom-up finite tree automaton which takes as input graph expressions from $(\mathbb{N}, \mathfrak{O}^k)$ and accepts an expression $e$ if $\models_{\Omega(\mathrm{val}_{\mathfrak{T}^k} e)} \varphi$.

$\mathfrak{O}^k$ is the set of function symbols of $\mathfrak{T}^k$. (def. 9.3.5, p. 186)

Given a graph $G \in \mathfrak{T}^n_k$, we use algorithm 6 to compute an expression $e$ with $\mathrm{val}_{\mathfrak{T}^k} e = G$. By theorem 9.3.17, we do this in linear time and obtain an expression with linearly many vertices.

$\mathrm{val}_{\mathfrak{T}^k} e$ denotes the result of $e$ when evaluated in $\mathfrak{T}^k$. (def. 6.3.3, p. 103)

We can hence implement algorithm 2 to run in linear time by lemma 8.4.1.

□

$\mathfrak{T}^k$ is the algebra of $k$-verdant graphs. (def. 9.3.3, p. 185)

$\mathfrak{T}^n_k$ denotes the $n$-verdant graphs of type $k$. (def. 9.3.3, p. 185)

For "ordinary" (read: untyped) graphs, the statement simplifies even further by observation 9.3.4.

**Corollary 9.3.19**

Let $\varphi$ be a sentence of the monadic second-order language of graphs, and let $k \in \mathbb{N}$. Then there exists an algorithm which, when given a graph $G = (V, E, \langle\!\!|_|\!\!\rangle)$ with $\mathrm{tw}(G) \leq k$, decides in time $\mathcal{O}(|V| + |E|)$ whether or not we have $\models_{\Omega(G)} \varphi$.

Or, in a more citeable version:

**Corollary 9.3.20**

Let $k \in \mathbb{N}$. Then for any graph property that can be expressed in the monadic second-order logic of graphs, there exists an algorithm which, given a graph $G$ of tree-width at most $k$, decides in linear time whether or not $G$ fulfils said property.

# 4. Looking for Graphs

Moving away from the main statement of Courcelle's Theorem, we present an interesting application of some of the intermediary results we have seen. Suppose we are given a class of graphs, described by a logical property. Rather than check an individual graph, we now want to know whether there is some structure that we know can never occur in graphs of the given class, like how no graph of tree-width $k \in \mathbb{N}$ can have a $k+2$-clique.[7]

Courcelle's Theorem, as we know it by now, is not capable of such a global statement – it provides "only" an algorithm that checks whether a given graph fulfils a given formula, so to check for a "forbidden" substructure, we would have to check every graph of our (probably infinitely large) class.

We show now how, despite Courcelle's Theorem not directly providing such a result, we can still easily implement some global checks in finite time. This has found mention in [CM93, pp. 73–74], where several other extensions beyond the scope of this thesis may also be found.

## 4.1. Automata to the Rescue

The centerpiece of our construction is the following fundamental result.

[7] The reader might feel reminded at this point of the theory of graph minors ([Die05, pp. 18–21]). Since minors are encodable by monadic second-order logic, the adventurous reader may work out how to apply some of these results to minors. We have omitted this generalisation since we have not introduced graph minors in this thesis.

**Theorem 9.4.1**

𐤄 is the Phoenician letter "he".
𐤌 is the Phoenician letter "mem".
𐤍 is the Phoenician letter "nun".
𐤎 is the Phoenician letter "semk".

Let 𐤄 = (𐤌, $\Pi$, $\Upsilon$, 𐤍, 𐤎) be a deterministic bottom-up finite tree automaton. Then it can be decided in time linear in the size of 𐤄 whether there exists a string $T \in \Upsilon$ which is accepted by 𐤄.

**Proof.** A proof is found at [Com+08, p. 40].

□

The "size" cited above includes the number of states as well as the size of our transition table, which may be substantial. Consequently, the procedure outlined in the following will *not* be linear (or even polynomial) in the size of some input. We therefore omit runtime considerations.

## 4.2. Peering Into the Emptiness

The application to our case is straightforward.

**Theorem 9.4.2**

$\mathfrak{G}$ is the algebra of graphs. (def. 6.2.2, p. 93)

Let $\varphi$ be a sentence of the monadic second-order language of graphs, and let $\mathscr{A} = (\mathscr{C}, \mathscr{O})$ be a finitely expressible algebra of weakly locally finite signature inherited from $\mathfrak{G}$. Then there exists an algorithm which determines in finite time whether there exists a graph $G \in \mathscr{C}_0$ with $\models_{\Omega(G)} \varphi$.

$\Omega(G)$ is the induced structure of $G$. (def. 5.4.2, p. 73)

$\mathfrak{L}$ denotes the second-order language of graphs. (def. 5.4.1, p. 72)

$\|\mathfrak{L}\|$ denotes the set of sentences over $\mathfrak{L}$. (def. 5.3.6, p. 62)

**Proof.** Let $\varphi \in \|\mathfrak{L}\|$ and let $\mathscr{A} = (\mathscr{C}, \mathscr{O})$ be a finitely expressible algebra of weakly locally finite signature inherited from $\mathfrak{G}$. By theorem 9.1.1, we find a deterministic bottom-up finite tree automaton which handles some set $\mathscr{K}$ of graph expressions capable of expressing every graph in $\mathscr{C}_0$ and accepts an expression if and only if its value fulfils $\varphi$. Hence a graph $G \in \mathscr{C}_0$ fulfilling $\varphi$ exists if and only if there is an expression evaluating to $G$ that gets accepted by our automaton. But whether this is the case can be checked by theorem 9.4.1.

□

Note how this result, contrary to our algorithms on graphs of bounded tree- or path-width, does not require us to be actually able to *compute* an expression for any graph in the algebra. It applies therefore even to algebras of graphs of unbounded tree-width, provided they meet the remaining requirements.

Of course, graphs of bounded tree-width are still interesting.

**Corollary 9.4.3**

Let $\varphi$ be a sentence of the monadic second-order language of graphs, and let $k \in \mathbb{N}$. Then there exists an algorithm which determines in finite time whether or not there exists a graph $G$ of tree-width at most $k$ with $\models_{\Omega(G)} \varphi$.

$\Omega(G)$ is the induced structure of $G$.
(def. 5.4.2, p. 73)

# 5. Graphs of Bounded Path-Width

We now show how Courcelle's Theorem can be adapted to apply to graphs of constantly bounded path-width. Of course, a graph of path-width at most $k \in \mathbb{N}$ has in particular tree-width at most $k$, so this adaptation might at first sight appear to serve no purpose. However, suppose we want to apply theorem 9.4.2 to the class of graphs of path-width at most $k$ for some $k \in \mathbb{N}$, that is, we want to know whether there exists a graph of path-width at most $k$ which fulfils a given monadic second-order sentence $\varphi$. Corollary 9.4.3 provides a way to check whether there exists a graph of *tree*-width at most $k$ fulfilling $\varphi$, but the path-width of that graph might be larger than its tree-width.[8] It is possible to encode the path-width as another monadic second-order sentence $\psi$, whence we can use corollary 9.4.3 to check for a graph fulfilling $\varphi \wedge \psi$, but the construction of $\psi$ makes use of the forbidden minors for path-width $k$, which at the time of this writing are only known up to a value of $k = 2$.

It is hence of interest to consider the class of graphs of path-width at most $k$ (or, indeed, many other subclasses of the class of graphs of tree-width at most $k$) on its own merit.

[8] For example, the graph 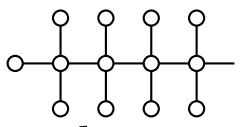 ··· has tree-width 1 (and a maximum node degree of 4), but as one continues the pattern, the path-width tends to infinity.

**Definition 9.5.1**

Let $k, n \in \mathbb{N}$. We denote

$$\mathfrak{p}_n^k := \{\, G \in \mathfrak{G}_n : \mathrm{pw}(G) \le k \,\}.$$

$\mathfrak{G}_n$ denotes the set of all graphs of type $n$. (def. 4.5.1, p. 39)

We introduce almost the same construction as in section 9.3.

**Definition 9.5.2**

Let $k, n \in \mathbb{N}$. We call a graph $G \in \mathfrak{G}_n$ *$k$-verdurous* if it admits a path-decomposition $(T, X, b)$ of width at most $k$ such that there is a node $v \in T$ with $\deg v = 1$ such that $b(v)$ contains all terminal vertices of $G$.

The set of all $k$-verdurous graphs of type $n$ is denoted by $\mathfrak{P}_n^k$.

In particular, we have again

$$\forall k \in \mathbb{N} \colon \mathfrak{P}_0^k = \mathfrak{p}_0^k.$$

Note that this time around, we require that the bag containing the terminal vertices is at one *end* of the path decomposition, not just at any old position. This makes sense when the reader thinks back to the proof of theorem 9.3.6: there we chose as the root of the tree decomposition the bag containing the terminals. For an undirected tree, this makes no difference, but for a path, there are two canonical roots – its end points.

In order to build only graphs of path-width at most $k \in \mathbb{N}$, entwining is too strong – we need to restrict our building techniques further.

Instead of attaching an arbitrary graph at some terminal vertices $1, \ldots, m$, we show that it suffices to be able to attach just one additional edge. For a 2-uniform graph, this looks as follows.

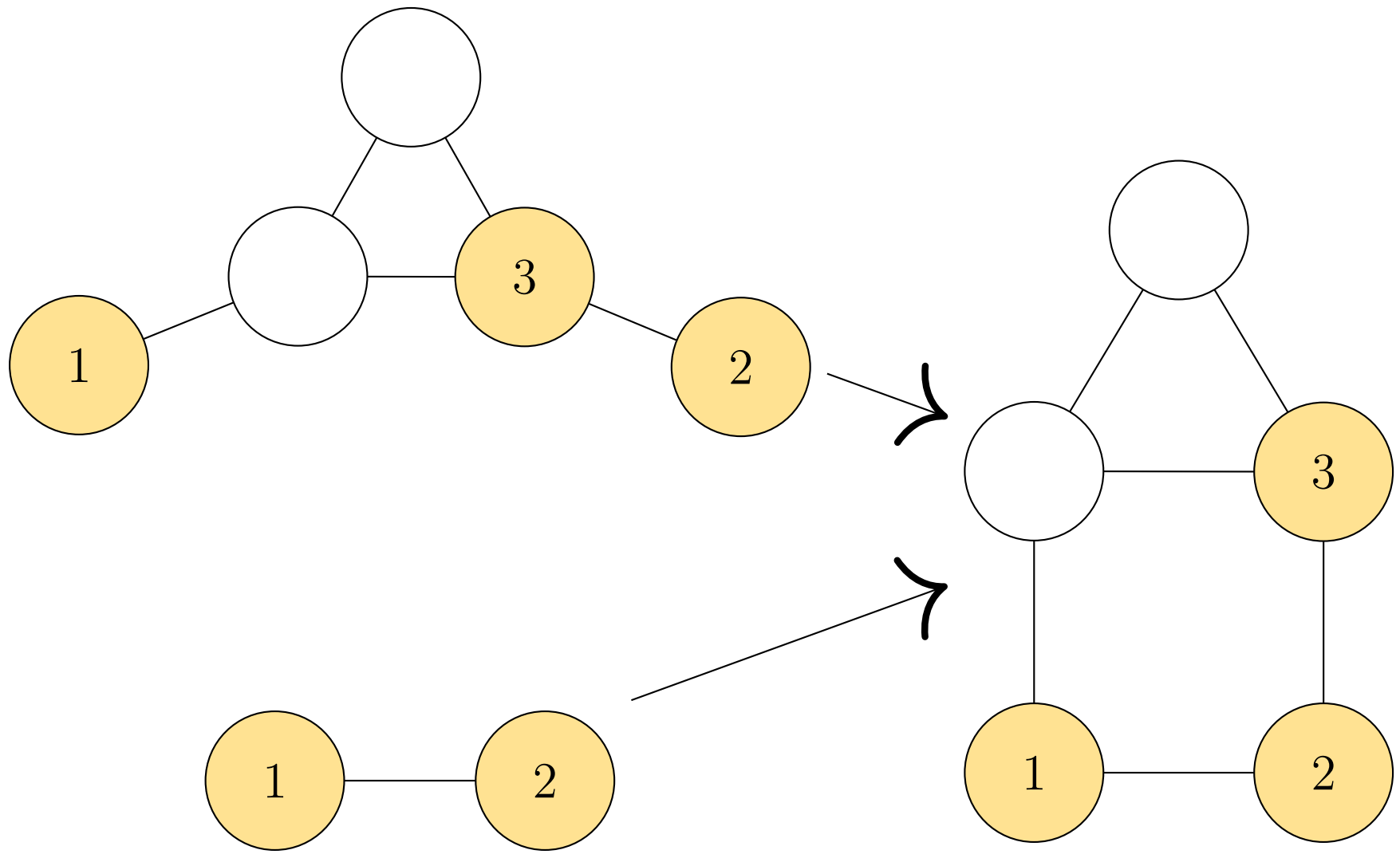

For arbitrary hypergraphs, we can attach edges of higher type.

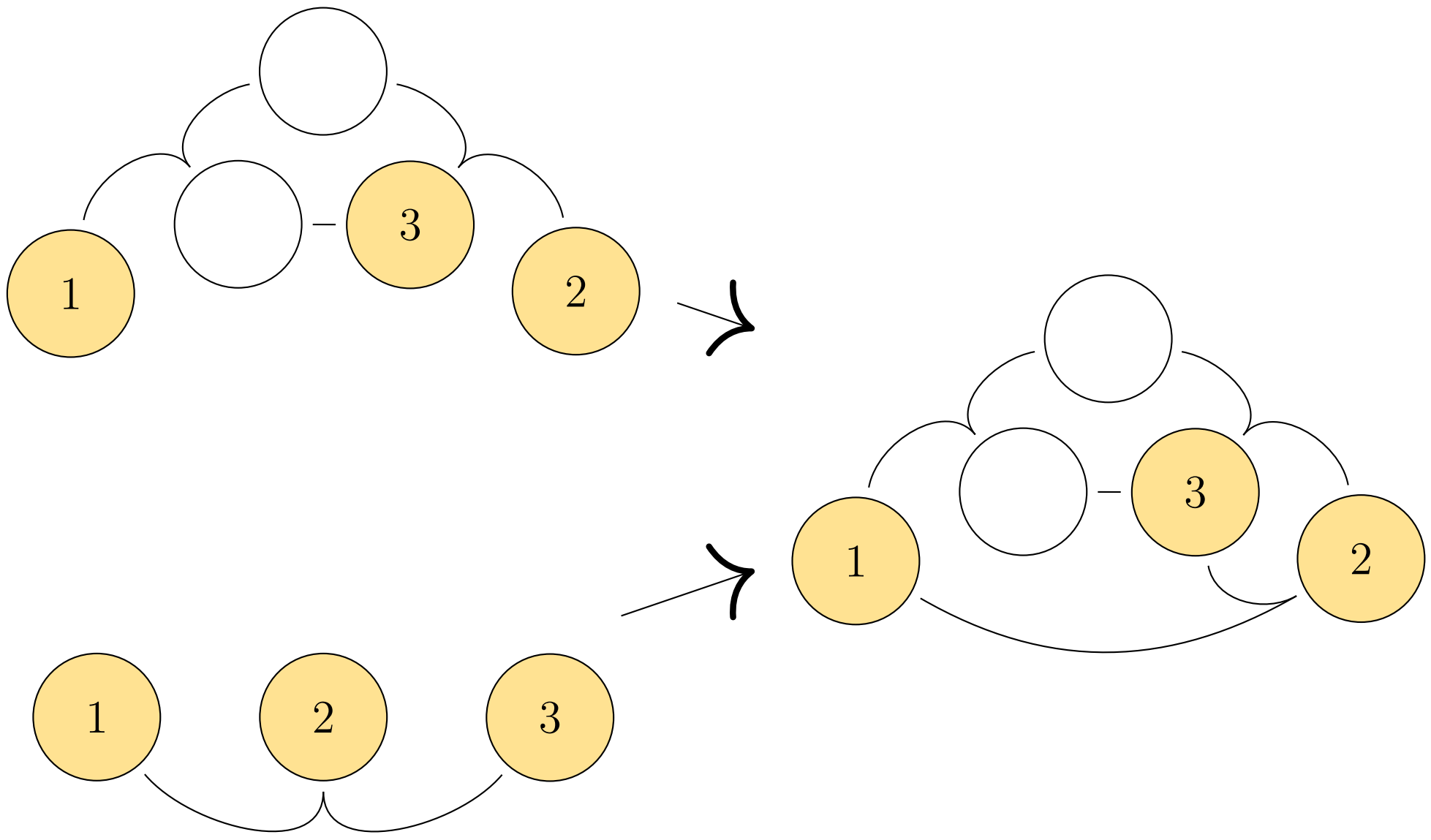

Add in a way to create a new connected component by adding a singular vertex, and we claim that this already suffices to express all verdurous hypergraphs.

### Definition 9.5.3

$\mathfrak{P}^k_n$ denotes the $k$-verdurous graphs of type $n$. (def. 9.5.2, p. 202)

$\oplus$ is the disjoint sum. (def. 4.5.3, p. 40)

$\mathfrak{v}$ is the type 1 graph with one vertex. (def. 4.5.6, p. 43)

$\mathfrak{e}_m$ is the type $m$ graph with $m$ vertices and one edge. (def. 4.5.6, p. 43)

Let $n, k \in \mathbb{N}$ and let $G \in \mathfrak{P}^k_n$.

We call

$$\uparrow G := G \oplus \mathfrak{v}$$

the *sprouting* of $G$.

Let further $m \in \mathbb{N}, m \leq k+1, m \leq n$.

We set

$$\sigma \colon \{1, \dots, n\} \to \{1, \dots, n+m\}, i \mapsto i$$

and call

$$\Lsh_{\mathfrak{e}_m} G := \leftrightarrows_\sigma \text{媲}^{n+m}_m \dots \text{媲}^{n+1}_1 (G \oplus \mathfrak{e}_m)$$

the *$m$-bloom* of $G$.

Proving that these constructions do not leave the realm of verdurous graphs is much easier than for tree-decompositions.

### Lemma 9.5.4

Let $n, k \in \mathbb{N}$ and let $G \in \mathfrak{P}^k_n$. Then if $n \leq k$, we have $\uparrow G \in \mathfrak{P}^k_{n+1}$.

**Proof.** Let $n, k \in \mathbb{N}$ and let $G \in \mathfrak{P}^k_n$. We need only show that for $n \leq k$, the graph $\uparrow G$ is again $k$-verdurous.

Denote the new vertex of $G \oplus \mathfrak{v}$ as $v^+$.

Let $(T, x, b)$ be a $k$-verdurous path-decomposition for $G$ with all terminal vertices of $G$ contained in $b(v)$. We introduce a new node $w$ with $b(w) := \{ v' \in G : v' \text{ is terminal} \}$. Attaching this node to $v$ makes the new path into another path-decomposition for $G$, with the difference that the root bag now contains exactly $n$ vertices.

We introduce another node $x$ with $b(x) := b(w) \cup \{v^+\}$. Attaching $x$ to $w$ then yields a path-decomposition for $G \oplus \mathfrak{v}$, and since $n \leq k$, it is of width at most $k$ and $k$-verdurous.

□

**Lemma 9.5.5**

Let $n, k \in \mathbb{N}$, let $G \in \mathfrak{P}^k_n$, and let $m \in \mathbb{N}$ with $m \leq k+1, m \leq n$. Then $\Game_{\mathfrak{e}_m} G \in \mathfrak{P}^k_n$.

$\mathfrak{P}^k_n$ denotes the $k$-verdurous graphs of type $n$. (def. 9.5.2, p. 202)

**Proof.** Blooming simply attaches an additional edge to the vertices in the root bag of a $k$-verdurous path-decomposition, meaning that if $(T, x, b)$ is a $k$-verdurous path-decomposition for $G$, then it is also a $k$-verdurous path-decomposition for $\Game_{\mathfrak{e}_m} G$.

$\Game$ denotes the graph blooming. (def. 9.5.3, p. 204)

$\mathfrak{e}_m$ is the type $m$ graph with $m$ vertices and one edge. (def. 4.5.6, p. 43)

□

We define a new signature. The reader may want to consult definition 6.2.1 for the symbols used.

**Definition 9.5.6**

Let $k \in \mathbb{N}$. We set

$$\begin{aligned}\mathfrak{Q}^k := & \bigcup_{i \in \mathbb{N}} \bigcup_{j \in \mathbb{N}} {}^j_i\mathfrak{F}_{\leftrightarrows} \\ & \cup \{\mathfrak{v}\} \\ & \cup \{\mathfrak{e}_1, \dots, \mathfrak{e}_{k+1}\} \\ & \cup \{{}_n\!\uparrow : n \in \mathbb{N}, n \leq k\} \\ & \cup \{{}_n\Game_{\mathfrak{e}_m} : m \in \mathbb{N}, m \leq k+1, m \leq n\}\end{aligned}$$

$\mathfrak{v}$ is the type 1 graph with one vertex. (def. 4.5.6, p. 43)

$\uparrow$ denotes the graph sprouting. (def. 9.5.3, p. 204)

and

$$\langle \mathfrak{f} \rangle := \begin{cases} \langle \mathfrak{f} \rangle & \text{if } \mathfrak{f} \in \mathfrak{F} \\ (n, n+1) & \mathfrak{f} = {}_n\!\uparrow \\ (n, n) & \mathfrak{f} = {}_n\Game_{\mathfrak{e}_m}. \end{cases}$$

As is our custom, we shall mostly omit the left index of the new function symbols.

### Definition 9.5.7

Let $k \in \mathbb{N}$. We define

$$\mathfrak{P}^k := \big(\{\, \mathfrak{P}^k_n \,\}_{n\in\mathbb{N}}, \{\mathscr{O}_{\mathfrak{f}}\}_{\mathfrak{f}\in\mathfrak{Q}^k}\big)$$

and let $\mathscr{O}$ assign to every function symbol the corresponding graph construction.

$\mathfrak{P}^k_n$ denotes the $k$-verdurous graphs of type $n$. (def. 9.5.2, p. 202)

$\mathfrak{Q}^k$ is the set of function symbols of $\mathfrak{P}^k$. (def. 9.5.6, p. 205)

$\mathfrak{G}$ is the algebra of graphs. (def. 6.2.2, p. 93)

### Lemma 9.5.8

Let $k \in \mathbb{N}$. Then $\mathfrak{P}^k$ is an inherited algebra of $\mathfrak{G}$.

**Proof.** This is covered by lemmas 9.5.4 and 9.5.5.

□

That the signature $(\mathbb{N}, \mathfrak{Q}^k)$ is weakly locally finite is clear from the definition. It remains to show that these function symbols suffice to construct every $k$-verdurous graph.

### Theorem 9.5.9

Let $k \in \mathbb{N}$. Then $\mathfrak{P}^k$ is finitely expressible.

**Proof.** Let $n, k \in \mathbb{N}$. Pick further a graph $G \in \mathfrak{P}^k_n$. We construct a $\{\, 1, \ldots, \max\{\, k+1, n \,\}\,\}$-local expression that yields $G$.

Note first that without loss of generality, we have $n \leq k+1$, as any $k$-verdurous graph can have at most $k+1$ distinct terminal vertices, and hence graphs of larger type can be constructed from a graph of type at most $k+1$ via $n$-local terminal redefinition. By the same argument, we can assume that the terminal vertices of $G$ are pairwise distinct.

Finally, we assume without loss of generality that, even if $G$ is a proper hypergraph, it contains no edges of type 1 – each edge of type 1 can be attached to its end point when it is introduced by a 1-bloom (which is $(k+1)$-local).

Let now $(T, X, b)$ be a nice $k$-verdurous path-decomposition for $G$. We build $G$ by induction over the nodes of $T$. Denote the nodes of $T$ as $v_1, \ldots, v_r$, where $v_r$ is the root (containing the terminal vertices of $G$) and $v_1$ is the singular leaf.

For $i \in \{1, \dots, r\}$, we write

$$V_i := \bigcup_{j=1}^{i} b(v_j)$$

and denote by $G_i$ the graph $G[V_i]$ with terminals equal to the vertices in $b(v_i)$.

We show by induction that for all $i \in \{1, \dots, r\}$, the graph $G_i$ can be constructed, whence $G = G_r$ can.

For $i = 1$, the graph $G_1$ has one vertex and (without loss of generality, as discussed above) no edges and is hence equal to $\mathfrak{v}$.

$\mathfrak{v}$ is the type 1 graph with one vertex. (def. 4.5.6, p. 43)

For $i > 1$, assume that the graph $G_{i-1}$ has been thus constructed. Since a nice verdurous path-decomposition can have no join nodes, there are only two cases to consider.[9]

*Case 1: $v_i$ is a forget node.* Then $G_i = G_{i-1}$ save for the terminal vertices, and we can erase the superfluous terminal vertex by one terminal redefinition.

*Case 2: $v_i$ is an introduce node.* Say $b(v_i) = b(v_{i-1}) \cup \{v'\}$. Then, in particular, $v' \notin G_{i-1}$, since $v' \in G_{i-1}$ would imply that there is a node $v \in \{v_1, \dots, v_{i-2}\}$ with $v' \in b(v)$, in which case $v'$ would also be contained in every bag between $v$ and $v_i$, in particular in $b(v_{i-1})$.

Now, which edges does $G_i$ contain that $G_{i-1}$ does not? Only edges incident to $v'$, and of those only the ones incident to only vertices in $G_i$. If an edge $e$ is incident to both $v'$ and some $v'' \in G_{i-1}$, then the path-decomposition of $G$ must contain a bag with both $v'$ and $v''$ in it. Since $v' \notin G_{i-1}$, this immediately implies that $v'' \in b(v_j)$ for some $j \geq i$ and hence $v'' \in b(v_{i-1})$ since $v'' \in G_{i-1}$.

Thus, only vertices in $b(v_{i-1})$ can be adjacent to $v'$ in $G_i$.

We take the sprout of $G_{i-1}$, which has at most $k + 1$ terminal vertices since $|b(v_{i-1})| \leq |b(v_i)| - 1 \leq k$. In this graph, all vertices which can be adjacent to $v'$ in $G_i$ are terminal, and we can hence add the desired edges by an appropriate number of blooming operations, which do not increase the number of terminals.

[9] Note that a nice path-decomposition *can* have up to one join node – we need verdurousness to ensure that there are none.

We have thus constructed $G_i$.

By induction, $G = G_r$ can be constructed from a $\{1, \ldots, k+1\}$-local expression, as desired.

□

We collect all these facts.

$\mathfrak{P}^k$ is the algebra of $k$-verdurous graphs. (def. 9.5.2, p. 202)

### Corollary 9.5.10

Let $k \in \mathbb{N}$. Then $\mathfrak{P}^k$ is a finitely expressible algebra of weakly locally finite signature inherited from $\mathfrak{G}$.

$\mathfrak{G}$ is the algebra of graphs. (def. 6.2.2, p. 93)

Hence Courcelle's Theorem becomes for this special case the following.

### Corollary 9.5.11

$\mathfrak{Q}^k$ is the set of function symbols of $\mathfrak{P}^k$. (def. 9.5.6, p. 205)

Let $\varphi$ be a sentence of the monadic second-order language of graphs, and let $k \in \mathbb{N}$. Then for every $n \in \mathbb{N}$, there exists a deterministic bottom-up finite tree automaton that takes as input a set $\mathscr{K}$ of $(\mathbb{N}, \mathfrak{Q}^k)$-expressions with $\mathfrak{P}^k$-values in $\mathfrak{P}^k_n$ and accepts an expression $e$ if and only if $\models_{\Omega(\mathrm{val}_{\mathfrak{P}^k} e)} \varphi$, and every element of $\mathfrak{P}^k_n$ is the value of some expression in $\mathscr{K}$.

We turn the procedure from the proof of theorem 9.5.9 into an algorithm.

---

**Algorithm 7:** The Path to Success

---

**Input:** a $k$-verdurous graph
**Output:** an expression resulting in that graph
Step $(v)$:
  **if** $v$ is a leaf **then**
    **return** $\mathfrak{v}$
  **else if** $v$ is a forget node **then**
    write $v'$ for the successor of $v$
    write $b(v) = \{ x_1, \dots, x_l \}$ and $b(v') = b(v) \cup \{x_{l+1}\}$
    set $\sigma\colon \{ 1, \dots, l \} \to \{ 1, \dots, l+1 \}, i \mapsto i$
    set $t'\colon \{ 1, \dots, l+1 \} \to b(v'), i \mapsto \begin{cases} x_{t(i)} & \text{if } i \leq l \\ x_{l+1} & \text{if } i = l+1 \end{cases}$
    **return** $\leftrightarrows_\sigma$Step$(v')$
  **else if** $v$ is an introduce node **then**
    for $v'$ the successor of $v$,
      write $b(v') = \{ x_1, \dots, x_l \}$, $b(v) = b(v') \cup \{x_{l+1}\}$
    $e :=$ ↑Step$(v')$
    **for** every edge $e$ incident to $x_{l+1}$ in $G_v$ **do**
      **if** $G$ is directed **then**
        set $\sigma$ to correct the terminal vertices for blooming
        $e :=$ ⤴$_{|e|}\leftrightarrows_\sigma e$
      **else**
        $e :=$ ⤴$_{|e|} e$
    **return** $e$

compute a nice $k$-verdurous path-decomposition $Z$ of $G$
$e :=$ Step$(\sqrt{Z})$
set $\sigma\colon \{ 1, \dots, n \} \to \{ 1, \dots, |b(\sqrt{Z})| \}$ such that it corrects the terminals
**return** $\leftrightarrows_\sigma e$

---

$\mathfrak{v}$ is the type 1 graph with one vertex. (def. 4.5.6, p. 43)

↑ denotes the graph sprouting. (def. 9.5.3, p. 204)

$\leftrightarrows$ denotes the terminal redefinition. (def. 4.5.4, p. 41)

⤴ denotes the graph blooming. (def. 9.5.3, p. 204)

The correctness of this algorithm follows directly from the proof of theorem 9.5.9.

### Corollary 9.5.12

Let $k \in \mathbb{N}$, and let $G$ be a $k$-verdurous graph. Then algorithm 7 termi-

$\mathfrak{Q}^k$ is the set of function symbols of $\mathfrak{P}^k$. (def. 9.5.6, p. 205)

nates and returns an $(\mathbb{N}, \mathfrak{Q}^k)$-expression $e$ with $\mathrm{val}_{\mathfrak{P}^k} e = G$.

$\mathrm{val}_{\mathfrak{P}^k} e$ denotes the result of $e$ when evaluated in $\mathfrak{P}^k$. (def. 6.3.3, p. 103)

The runtime, of course, hinges on our ability to efficiently compute a nice $k$-verdurous path-decomposition. We utilise the following result by Bodlaender and Kloks.

$\mathfrak{P}^k$ is the algebra of $k$-verdurous graphs. (def. 9.5.2, p. 202)

### Theorem 9.5.13

Let $k \in \mathbb{N}$. Then there exists an algorithm which, given a graph $G$ of path-width at most $k$ and a tree-decomposition for $G$ of width $k$, computes a path-decomposition of for $G$ of width $k$ in time polynomial in the size of $G$.

**Proof.** This is theorem 6.1 of [BK96] with $l = k$.

□

Since the path-width of a graph is an upper bound for its tree-width, the above combines with theorem 9.3.13 into the following.

### Corollary 9.5.14

Let $k \in \mathbb{N}$. Then there exists an algorithm which, when given a graph $G = (V, E, \langle\!\!\langle _ \rangle\!\!\rangle)$ of path-width at most $k$, computes a path-decomposition for $G$ in time $\mathcal{O}(|V| + |E|)$. This path-decomposition has at most $\mathcal{O}(|V| + |E|)$ bags.

We get the results we desire.

### Theorem 9.5.15

$\mathfrak{P}^n_k$ denotes the $n$-verdurous graphs of type $k$. (def. 9.5.2, p. 202)

Let $n, k \in \mathbb{N}$. Then there exists an algorithm which, for any given graph $G = (V, E, \langle\!\!\langle _ \rangle\!\!\rangle, t) \in \mathfrak{P}^n_k$, computes a $k$-verdurous path-decomposition with at most $\mathcal{O}(|V| + |E|)$ bags for $G$ in time $\mathcal{O}(|V| + |E|)$.

**Proof.** This is entirely analogous to theorem 9.3.15.

□

**Theorem 9.5.16**

Let $n, k \in \mathbb{N}$. Then there exists an algorithm which, for any given graph $G = (V, E, \langle\!\langle _ \rangle\!\rangle, t) \in \mathfrak{P}^n_k$, computes a nice $k$-verdurous path-decomposition with at most $\mathcal{O}(|V|+|E|)$ bags for $G$ in time $\mathcal{O}(|V|+|E|)$.

$\mathfrak{P}^n_k$ denotes the $n$-verdurous graphs of type $k$. (def. 9.5.2, p. 202)

**Proof.** This is again entirely analogous to theorem 9.3.16.

□

We can hence bound the runtime of algorithm 7.

**Theorem 9.5.17**

Let $n, k \in \mathbb{N}$, and let $G = (V, E, \langle\!\langle _ \rangle\!\rangle, t) \in \mathfrak{P}^k_n$. Then algorithm 7 terminates in time $\mathcal{O}(|V| + |E|)$ and returns an expression with $\mathcal{O}(|V| + |E|)$ vertices.

**Proof.** Let $n, k \in \mathbb{N}$, and let $G = (V, E, \langle\!\langle _ \rangle\!\rangle, t) \in \mathfrak{P}^k_n$. Algorithm 7 begins by computing a nice $k$-verdurous path-decomposition for $G$, which by theorem 9.5.16 happens in linear time, and can be done such that the path-decomposition has at most linearly many bags.

The only nontrivial part left to analyse is the Step subroutine, which gets called exactly once for each node. Its runtime is exactly the same as that of the Build subroutine in the tree-decomposition version of the algorithm.

The Step routine is defined on page 209.

The Build routine is defined on page 194.

□

We can hence prove the main result of Courcelle's Theorem applied to graphs of constantly bounded path-width.

**Theorem 9.5.18**

Let $\varphi$ be a sentence of the monadic second-order language of graphs, and let $k \in \mathbb{N}$. Then for every $n \in \mathbb{N}$, there exists an algorithm which, given a $k$-verdurous type $n$ graph $G = (V, E, \langle\!\langle _ \rangle\!\rangle, t)$, decides in time $\mathcal{O}(|V| + |E|)$ whether or not we have $\vDash_{\Omega(G)} \varphi$.

$\Omega(G)$ is the induced structure of $G$. (def. 5.4.2, p. 73)

**Proof.** Let $\varphi$ be a sentence of the monadic second-order language of graphs, and let $k, n \in \mathbb{N}$. By corollary 9.5.11, there exists a deterministic bottom-up finite tree automaton which takes as input graph expressions from $(\mathbb{N}, \mathfrak{Q}^k)$ and accepts an expression $\boldsymbol{e}$ if $\vDash_{\Omega(\mathrm{val}_{\mathfrak{P}^k} \boldsymbol{e})} \varphi$.

Given a graph $G \in \mathfrak{P}^n_k$, we use algorithm 7 to compute an expression $\boldsymbol{e}$ with $\mathrm{val}_{\mathfrak{P}^k} \boldsymbol{e} = G$. By theorem 9.5.17, we do this in linear time and obtain an expression with linearly many vertices.

We can hence implement algorithm 2 to run in linear time by lemma 8.4.1.

□

$\mathfrak{Q}^k$ is the set of function symbols of $\mathfrak{P}^k$. (def. 9.5.6, p. 205)

$\mathrm{val}_{\mathfrak{P}^k} \boldsymbol{e}$ denotes the result of $\boldsymbol{e}$ when evaluated in $\mathfrak{P}^k$. (def. 6.3.3, p. 103)

$\Omega(\mathrm{val}_{\mathfrak{P}^k} \boldsymbol{e})$ is the induced structure of $\mathrm{val}_{\mathfrak{P}^k} \boldsymbol{e}$. (def. 5.4.2, p. 73)

$\mathfrak{P}^k$ is the algebra of $k$-verdurous graphs. (def. 9.5.2, p. 202)

$\mathfrak{P}^n_k$ denotes the $n$-verdurous graphs of type $k$. (def. 9.5.2, p. 202)

For "ordinary" (read: untyped) graphs, the statement simplifies as it did in the section on tree-width.

### Corollary 9.5.19

Let $\varphi$ be a sentence of the monadic second-order language of graphs, and let $k \in \mathbb{N}$. Then there exists an algorithm which, when given a graph $G = (V, E, \langle\!|_|\!\rangle)$ with $\mathrm{pw}(G) \leq k$, decides in time $\mathcal{O}(|V| + |E|)$ whether or not we have $\vDash_{\Omega(G)} \varphi$.

Or, in a more citeable version:

### Corollary 9.5.20

Let $k \in \mathbb{N}$. Then for any graph property that can be expressed in the monadic second-order logic of graphs, there exists an algorithm which, given a graph $G$ of path-width at most $k$, decides in linear time whether or not $G$ fulfils said property.

Thinking back to section 9.4, we get the same corollary as we did for graphs of bounded tree-width, without having to invoke forbidden minors.

### Corollary 9.5.21

Let $\varphi$ be a sentence of the monadic second-order language of graphs, and let $k \in \mathbb{N}$. Then there exists an algorithm which determines in finite time whether or not there exists a graph $G$ of path-width at most $k$ with $\vDash_{\Omega(G)} \varphi$.

# Chapter 10
# Concluding Remarks

Some two hundred pages ago, we set out to introduce the reader to the fascinating world of Courcelle's Theorem. It is our sincere hope that those readers that have stayed with us until this chapter have now an intuitive understanding of not only what Courcelle's Theorem can and cannot do, but also of *why* it can do the things it can. Using the tools developed over the previous chapters, the reader should be able, should the need ever arise, to adapt Courcelle's Theorem to their own peculiar use case.

There has, as ever, more been left unsaid than said, first and foremost the extension of Courcelle's Theorem to cover optimisation problems, that is, to answer not the question "does this graph admit a vertex cover of size 5?", but "what is the smallest $k \in \mathbb{N}$ such that this graph admits a vertex cover of size $k$?" Such an extension is possible and can indeed be found in [CM93]. It does, however, require quite some additional machinery, as [CM93] represents the culmination of an entire series of papers on this topic. The basic building blocks are the same, but there is much additional notation involved – predicates, for example, are replaced by evaluations, which map not to the set $\{ \top, \bot \}$, but to an arbitrary set. The surrounding framework needs to be adjusted accordingly. All this could be explored in a future work.

# Appendix A

# Counting Logic

Courcelle's Theorem, in its original formulation, applies not just to monadic second-order logic, but to so-called *counting* monadic second-order logic. We have opted in this thesis to prove the theorem, first and foremost, for unextended monadic second-order logic, on the one hand to keep notation manageable, on the other hand because most natural examples of graph properties require no counting.

We present in this appendix the missing definitions and extend the proofs from chapters 5 to 9 to properly cover counting logic.

## 1. Return to Logic

We first explain what counting logic entails. For this section, only the concepts from chapter 5 are prerequisite.

### 1.1. Learn to Count Again

We introduce first intuitively the notion of counting monadic second-order logic, which extends the monadic second-order logic the reader knows from chapter 5.

We fix for this section a second-order language $\Gamma$ endowed with a (non-monadic) second-order structure $\Omega$ such that the following conditions are fulfilled.

- The universe $|\Omega|$ is at most countably infinite.
- All 1-place predicates that are true on finitely many constants exist in $(\Gamma, \Omega)$, that is, for every finite subset $X$ of $|\Omega|$, there is a predicate symbol that checks membership in $X$.
- All 1-place functions exist in $(\Gamma, \Omega)$.
- There is a 2-place predicate $\equiv$ that checks equality of constants.

We now present some constructions on this universe that will eventually lead us to the definition of counting monadic second-order logic.

Suppose we are given finite sets $X, Y \subseteq |\Omega|$ and a function $f\colon X \to Y$. We now ask ourselves, as mathematicians sometimes do, "is this function surjective?"

To answer this question, we model it in $(\Gamma, \Omega)$. Modelling $X$ and $Y$ is easy since we have assumed that all finite sets have a predicate symbol representing them, say ${}_1\Lambda_X, {}_1\Lambda_Y \in {}_1\Lambda$ with $X = \Omega({}_1\Lambda_X), Y = \Omega({}_1\Lambda_Y)$.

${}_1\Lambda$ is the set of 1-place predicate symbols.

To model the function $f$, since our structure gives us access to all possible functions, we simply pick a function $\mathring{f}\colon |\Omega| \to |\Omega|$ with $\mathring{f}|_X = f$ and a function symbol ${}_1\Delta_0 \in {}_1\Delta$ with $\Omega({}_1\Delta_0) = \mathring{f}$.

${}_1\Delta$ is the set of 1-place function symbols.

We can then construct the sentence

$$\varphi := \forall \delta_0 : {}_1\Lambda_Y(\delta_0) \to \exists \delta_1 : {}_1\Lambda_X(\delta_1) \wedge \equiv({}_1\Delta_0(\delta_1), \delta_0)$$

("for every element $\delta_0$ of $Y$, there is a variable that lies in $X$ and that is mapped to $\delta_0$ by $f$") such that $\vDash_\Omega \varphi$ if and only if $f$ is surjective.

For the purposes only of this section, we replace this formula by the "pseudo-predicate" symbol **Ep** (for "epic") that takes a function symbol and two predicates (domain and range) and evaluates to $\top$ if the function is surjective from domain onto range and $\bot$ otherwise. It is "pseudo" in the sense that it cannot be a true predicate symbol, since our definitions only allow plugging universe elements into predicates, not function symbols or even other predicates. The symbol **Ep** will simply be a shorthand like $\to$ or $\vee$.

Suppose now we were given a finite set $X \subseteq |\Omega|$ and asked, "is the cardinality of this set even or odd?" Now there is a question! However, armed with our recent epic discovery, we immediately spot a solution: a set of even cardinality can be partitioned into two subsets of equal size, and two finite sets are of equal size if each surjects onto the other. Using the existence of all necessary predicates and functions posited at the beginning of this section, we construct the sentence

$$\begin{aligned}\varphi := {} & \exists_1 \lambda_0 : \exists_1 \lambda_1 : \\ & \quad \forall \delta_0 : {}_1\lambda_0(\delta_0) \to {}_1\Lambda_X(\delta_0) \\ & \wedge \forall \delta_1 : {}_1\lambda_1(\delta_1) \to {}_1\Lambda_X(\delta_1) \\ & \wedge \forall \delta_2 : {}_1\lambda_0(\delta_2) \to \neg {}_1\lambda_1(\delta_2) \\ & \wedge \forall \delta_3 : {}_1\Lambda_X(\delta_3) \to {}_1\lambda_0(\delta_3) \vee {}_1\lambda_1(\delta_3) \\ & \wedge \exists_1 \delta_0 : \mathbf{Ep}({}_1\delta_0, {}_1\lambda_0, {}_1\lambda_1) \\ & \wedge \exists_1 \delta_1 : \mathbf{Ep}({}_1\delta_1, {}_1\lambda_1, {}_1\lambda_0)\end{aligned}$$

with $\vDash_\Omega \varphi$ if and only if $X$ is indeed of even cardinality. To know whether $|X|$ is odd, we simply negate $\varphi$.

An analogous construction reveals whether a set's cardinality is divisible by 3, or indeed by any natural number.

Of course, for the number 2, there are exactly two possibilities: even (that is, 0 modulo 2), or odd (also known as 1 modulo 2). For higher integers, there are more than two options. A simple negation of the sentence for "is a multiple of 3" cannot detect whether the cardinality of $X$ is 1 or 2 modulo 3.

Luckily, this is easily fixed. We give the construction to detect 1-modulo-3-ness; all other cases work analogously.

Denoting the sentence for "$|Y|$ is divisible by 3" as $\mathbf{3}(Y)$, we use the sentence

$$\begin{aligned}\varphi := {} & \exists \delta_0 \exists_1 \lambda_0 : \\ & \quad {}_1\Lambda_X(\delta_0) \\ & \wedge \neg {}_1\lambda_0(\delta_0) \\ & \wedge \forall \delta_1 : {}_1\lambda_0(\delta_1) \to {}_1\Lambda_X(\delta_1) \\ & \wedge \forall \delta_2 : {}_1\Lambda_X(\delta_2) \to {}_1\lambda_0(\delta_2) \vee (\delta_0 = \delta_2) \\ & \wedge \mathbf{3}({}_1\lambda_0)\end{aligned}$$

to determine that we can delete one element from $X$ to end up with a set of cardinality divisible by 3. Hence, $\vDash_\Omega \varphi$ if and only if $|X|$ is 1 modulo 3.

Let us take a step back, breathe deeply, and consider what we have learned. In a second-order logical framework that knows at the least about finite sets and the functions between them, we are for natural number $n, k \in \mathbb{N}_{>0}$ able to detect whether a finite set (given as a 1-place predicate) has cardinality $k$ modulo $n$ with a formula whose length only depends on $n$ and $k$.

Note that this detection is impossible in monadic second-order logical frameworks, as we used quantification over function symbols to ascertain whether two sets have the same cardinality.

We define the following important shorthand for use in the remainder of this chapter.

### Notation A.1.1

In the context of a logical framework as assumed in this section, we denote by $\mathbf{Card}_k^n(\Lambda)$ the formula that determines whether the finite set given by the 1-place predicate $\Lambda$ has cardinality $k$ modulo $n$.

Courcelle's Theorem, in its full generality, works in a subset of second-order logic that is strictly less expressive than full second-order logic, but strictly more expressive than monadic second-order logic. This "counting monadic second-order logic" is most easily imagined by taking a monadic second-order logical framework and adding the "pseudo-predicates" $\mathbf{Card}_k^n$ defined above. Formally defining it, however, is most easily done by taking a larger second-order logical framework satisfying the prerequisites of this section, carrying out the above construction, and then restricting the use of all non-monadic predicate variables and all non-constant function variables except inside that construction. Rather than going through the motions of carrying out this cumbersome and not very enlightening task, we shall trust that the musings above have convinced the reader that such a construction is in theory possible and carry on with the intuitive notion that a counting monadic second-order logical framework is a monadic second-order logical framework with some additional sentences which are always and irretrievably wrapped inside the shorthand **Card**.

Whenever we define a counting monadic second-order logical framework, we define the monadic second-order logical framework just as we have introduced it in definition 5.3.14 and additionally note the subsets of $\mathbb{N}$ from which the indices $n, k$ in $\mathbf{Card}_k^n$ can be chosen.

# 2. You Can Always Count on Graphs

We extend the monadic second-order framework defined in definition 5.4.3.

**Definition A.2.1**

We denote by $\mathfrak{X}^+$ the monadic second-order logical framework obtained by enhancing $\mathfrak{X}$ with the pseudo-predicates **Card** from notation A.1.1.

$\mathfrak{X}$ denotes the direct logical framework of graphs. (def. 5.4.3, p. 73)

This definition is, as discussed in the previous section, not entirely formal. The circuitous logical framework on graphs (as defined in section 5.6) is actually easier to formally change – since its variables are already sets, checking cardinality can be achieved by a simple predicate (and not a pseudo-one).

We reproduce the definitions here, with the new predicate included.

**Definition A.2.2**

We denote by $\mathring{\mathfrak{L}}^+$ the second-order language with

- ${}_1\Lambda_{\mathring{\mathfrak{L}}} := \{ \lambda^0_{\text{conn}}, \lambda_{\text{sgl}}, \lambda_{\text{vert}}, \lambda_{\text{edge}}, \lambda^{n,p}_{\text{card}} : p \in \mathbb{N}_{>0}, n \in \mathbb{N}_{\leq p} \}$,
- ${}_2\Lambda_{\mathring{\mathfrak{L}}} := \{ \lambda^1_{\text{conn}}, \sqsubseteq \}$,
- $\forall n \in \mathbb{N}_{>2} :\ {}_n\Lambda_{\mathring{\mathfrak{L}}} := \{ \lambda^{n-1}_{\text{conn}} \}$,
- ${}_0\Delta_{\mathring{\mathfrak{L}}} := \{ \varnothing \}$,
- ${}_1\Delta_{\mathring{\mathfrak{L}}} := \bigcup_{k \in \mathbb{N}_{>0}} \{ \delta^K_{\text{term}} : K \in \mathbb{2}^{\{1,\ldots,k\}} \}$,
- $\forall n \in \mathbb{N}_{>1} :\ {}_n\Delta_{\mathring{\mathfrak{L}}} := \varnothing$.

${}_1\Lambda_{\mathring{\mathfrak{L}}}$ is the set of 1-place predicate symbols of $\mathring{\mathfrak{L}}$.

$\lambda_{\text{vert}}$ checks whether its argument is a set of vertices.

$\lambda_{\text{edge}}$ checks whether its argument is a set of edges.

$\mathring{\mathfrak{L}}$ is the circuitous language of graphs. (def. 5.6.1, p. 77)

Of course, $\lambda^{n,p}_{\text{card}}(x)$ should be true if and only if $|x| = n \mod p$.

**Definition A.2.3**

${}_0\Delta_{\mathring{\mathfrak{L}}}$ is the set of 0-place function symbols of $\mathring{\mathfrak{L}}$.

Let $G = (V, E, \langle\!\langle _ \rangle\!\rangle, t)$ be a graph of type $k$. The *induced circuitous second-order structure with counting* of $G$, denoted $\mathring{\Omega}^+(G)$, is the following second-order structure on $\mathring{\mathfrak{L}}$.

$\mathbb{2}^{\{1,\ldots,k\}}$ denotes the power set of $\{1, \ldots, k\}$.

$\mathbb{2}^V$ denotes the power set of $V$.

$\lambda_{\text{vert}}$ checks whether its argument is a set of vertices.

$\lambda_{\text{edge}}$ checks whether its argument is a set of edges.

- The universe of $\mathring{\Omega}^+(G)$ is $|\mathring{\Omega}^+(G)| := {}^{\delta}_{0}|\mathring{\Omega}^+(G)| := \mathbb{2}^V \cup \mathbb{2}^E$.
- For $n \in \mathbb{N}_{>0}$, ${}^{\lambda}_{n}|\mathring{\Omega}^+(G)| := \varnothing$.
- For $n \in \mathbb{N}_{>0}$, ${}^{\delta}_{n}|\mathring{\Omega}^+(G)| := \varnothing$.
- $\mathring{\Omega}^+(G)(\lambda_{\text{sgl}}) := \{ x : x \in |\mathring{\Omega}^+(G)|, |x| = 1 \} \subseteq |\mathring{\Omega}^+(G)|$.
- $\mathring{\Omega}^+(G)(\lambda_{\text{vert}}) := \mathbb{2}^V \subseteq |\mathring{\Omega}^+(G)|$.
- $\mathring{\Omega}^+(G)(\lambda_{\text{edge}}) := \mathbb{2}^E \subseteq |\mathring{\Omega}^+(G)|$.
- $\mathring{\Omega}^+(G)(\sqsubseteq) := \{ (x, y) \in |\mathring{\Omega}^+(G)|^2 : x \subseteq y \}$.
- For every $p \in \mathbb{N}_{>0}$ and every $n \in \mathbb{N}_{\leq p}$,

  $$\mathring{\Omega}^+(G)(\lambda^{p,n}_{\text{card}}) := \{ x : x \in |\mathring{\Omega}^+(G)|, |x| = p \mod n \} \subseteq |\mathring{\Omega}^+(G)|^1.$$

- For every $n \in \mathbb{N}$,

  $$\begin{aligned}\mathring{\Omega}^+(G)(\lambda^{n}_{\text{conn}}) :=\{ (E', V_1, \dots, V_n) : \exists e \in E' :&\\ \exists v_1 \in V_1, \dots, \exists v_n \in V_n :&\\ e \in E, \langle\!| e |\!\rangle = v_1 \dots v_n \}&\\ \subseteq |\mathring{\Omega}^+(G)|^{n+1}.\end{aligned}$$

- $\mathring{\Omega}^+(G)(\varnothing)\colon \{ () \} \to |\mathring{\Omega}^+(G)|, () \mapsto \varnothing$.
- For $n \in \mathbb{N}_{>0}$, for $K \subseteq \{ 1, \dots, n \}$,

  $$\mathring{\Omega}^+(G)(\delta^{K}_{\text{term}})\colon |\mathring{\Omega}^+(G)| \to |\mathring{\Omega}^+(G)|,$$
  $$x \mapsto \begin{cases} x \cup \{ t(i) : i \in K, i \leq k \} & \text{if } x \subseteq V \\ x & \text{otherwise.} \end{cases}$$

**Definition A.2.4**

We call

$$\mathring{\mathfrak{M}}^{+} := \bigcup_{n\in\mathbb{N}} \{\, \mathring{\Omega}^{+}(G) : G \in \mathfrak{G}_n \,\}$$

the *circuitous multiverse of finite graphs with counting* and

$$\mathring{\mathfrak{X}}^{+} := \left(\mathring{\mathfrak{L}}^{+}, \mathring{\mathfrak{M}}^{+}\right)$$

the *circuitous logical framework of finite graphs with counting.*

$\mathfrak{G}_n$ denotes the set of all graphs of type $n$. (def. 4.5.1, p. 39)

$\mathring{\mathfrak{L}}^{+}$ is the circuitous language of graphs with counting. (def. A.2.2, p. 219)

It is easy to see that, again, the second language is at least as expressive as the first one.

**Theorem A.2.5**

Let $\varphi$ be a well-formed formula of $\mathfrak{X}^{+}$. Then there exists also a formula $\mathring{\varphi} \in |\mathring{\mathfrak{L}}^{+}|$ with $\overleftarrow{\mathring{\varphi}} = \overleftarrow{\varphi}$ such that for every graph $G \in \mathfrak{G}$ and for every full variable assignment $\tau$ in $\Omega(G)$ for $\varphi$, we have

$$\varphi[\tau] \leftrightarrow \top \Leftrightarrow \mathring{\varphi}[\mathring{\tau}] \leftrightarrow \top.$$

If $\varphi$ did not contain a term of the form $\mathbf{Card}(\,\ldots\,)$, then $\mathring{\varphi}$ can be chosen such that it contains no predicate of the form $\lambda_{\text{card}}$.

$\mathfrak{X}^{+}$ denotes the logical framework of graphs with counting. (def. A.2.1, p. 219)

$|\mathring{\mathfrak{L}}^{+}|$ is the set of all well-formed formulas over $\mathring{\mathfrak{L}}^{+}$. (def. 5.3.5, p. 61)

$\overleftarrow{\varphi}$ denotes the set of free variables of $\varphi$. (def. 5.3.6, p. 62)

$\mathfrak{G}$ is the algebra of graphs. (def. 6.2.2, p. 93)

$\Omega(G)$ is the induced structure of $G$. (def. 5.4.2, p. 73)

**Proof.** We state only the parts of the proof that need to be added to the proof of theorem 5.6.6.

*Case 1: $\varphi$ is an atomic formula.* We add our new case.

*Case 1.1: $\varphi = \mathbf{Card}^{n}_{k}({}_1\lambda_0)$ for some $k \le n \in \mathbb{N}$ and some 1-place predicate variable ${}_1\lambda_0$.* We list this case with the atomic formulas because it serves the same function as a basic building block.

We set $\mathring{\varphi} := \lambda^{k,n}_{\text{card}}({}_1\lambda_0)$, keeping the same set of free variables. Equivalence is immediate by the definition of $\lambda_{\text{card}}$.

The remainder of the proof remains unchanged.

□

$\mathfrak{X}^+$ denotes the logical framework of graphs with counting. (def. A.2.1, p. 219)

$\|\mathfrak{X}^+\|$ denotes the set of sentences over $\mathfrak{X}^+$. (def. 5.3.6, p. 62)

$\mathring{\mathfrak{X}}^+$ denotes the circuitous logical framework of graphs with counting. (def. A.2.4, p. 221)

$\mathfrak{G}$ is the algebra of graphs. (def. 6.2.2, p. 93)

$\Omega(G)$ is the induced structure of $G$. (def. 5.4.2, p. 73)

$\mathring{\Omega}^+(G)$ is the circuitous induced structure with counting of $G$. (def. A.2.3, p. 219)

$\mathring{\mathfrak{L}}^+$ is the circuitous language of graphs with counting. (def. A.2.2, p. 219)

**Corollary A.2.6**

For every sentence $\varphi \in \|\mathfrak{X}^+\|$, there exists a sentence $\mathring{\varphi} \in \|\mathring{\mathfrak{X}}^+\|$ such that

$$\forall G \in \mathfrak{G}\colon \vDash_{\Omega(G)} \varphi \Leftrightarrow \vDash_{\mathring{\Omega}^+(G)} \mathring{\varphi}\,.$$

If $\varphi$ did not contain a term of the form **Card**( … ), then $\mathring{\varphi}$ can be chosen such that it contains no predicate of the form $\lambda_{\text{card}}$.

# 3. Pride and Predicates

Having established an extended language, we need to check that the family of predicates it induces is again locally finite. If we do not restrict the new predicates, there is at least the formula that checks whether $\lambda^{1,p}_{\text{card}}$ holds for the set of all vertices for each $p \in \mathbb{N}_{>0}$, and no two of those formulas are type $n$ equivalent. Consequently, we must add a restriction on the size of $p$. We reproduce definition 6.9.1 with this new restriction.

**Definition A.3.1**

Let $w \in \mathbb{N}$. The set of sentences of $\mathring{\mathfrak{L}}^+$ that have width at most $w$ is denoted by $\|\mathring{\mathfrak{L}}^+\|_w$.

Let further $l \in \mathbb{N}$. The set of sentences of $\|\mathring{\mathfrak{L}}^+\|_w$ that use at most the variable symbols $\mu_0, \dots, \mu_l$ is denoted $\|\mathring{\mathfrak{L}}^+\|^l_w$.

Let further $k \in \mathbb{N}$. The set of sentences of $\|\mathring{\mathfrak{L}}^+\|^l_w$ that use no predicate symbol $\lambda^t_{\text{conn}}$ for $t > k$ is denoted $\|\mathring{\mathfrak{L}}^+\|^{l,k}_w$.

Let finally $m \in \mathbb{N}$. The set of sentences of $\|\mathring{\mathfrak{L}}^+\|^{l,k}_w$ that use no predicate symbol $\lambda^{s,t}_{\text{card}}$ for $t > m$ is denoted $\|\mathring{\mathfrak{L}}^+\|^{l,k,m}_w$.

With this definition, the set of pairwise non-equivalent formulas becomes finite again.

**Theorem A.3.2**

Let $n \in \mathbb{N}$, and let $w, l, k, m \in \mathbb{N}$. Then the set

$$\|\mathring{\mathfrak{L}}^{+}\|_{w}^{l,k,m} \big/ \overset{n}{\sim}$$

is finite.

$\mathring{\mathfrak{L}}^{+}$ is the circuitous language of graphs with counting. (def. A.2.2, p. 219)

$\|\mathring{\mathfrak{L}}^{+}\|$ denotes the set of sentences over $\mathring{\mathfrak{L}}^{+}$. (def. 5.3.6, p. 62)

$\overset{n}{\sim}$ denotes equivalence on graphs of type $n$. (def. 6.9.3, p. 122)

**Proof.** It suffices to notice that in the proof of corollary 6.9.6, if we restrict our predicates to use of $\lambda_{\text{conn}}$ only the predicates $\lambda^{0}_{\text{conn}}, \dots, \lambda^{k}_{\text{conn}}$ and of $\lambda^{a,b}_{\text{card}}$ only those with $b \leq m$ and $a \leq b$, we are still left with only finitely many predicates.

□

Just as before, we turn these formulas into predicates – compare definition 6.9.7.

**Definition A.3.3**

Let $w, l, k, m \in \mathbb{N}$, let $n \in \mathbb{N}$, and let $\mathring{\varphi} \in \|\mathring{\mathfrak{L}}^{+}\|_{w}^{l,k,m} \big/ \overset{n}{\sim}$. We set

$$p^{+n}_{\mathring{\varphi}} : \mathfrak{G}_n \to \{\top, \bot\}, G \mapsto \begin{cases} \top & \models_{\mathring{\Omega}^{+}(G)} \varphi \text{ for some}^{1}\ \varphi \in \mathring{\varphi} \\ \bot & \text{otherwise} \end{cases}$$

and

$$\mathscr{P}^{+n}_{w,l,k,m} := \{\, p^{+n}_{\mathring{\varphi}} : \mathring{\varphi} \in \|\mathring{\mathfrak{L}}^{+}\|_{w}^{l,k,m} \big/ \overset{n}{\sim} \,\}.$$

We then set

$$\mathscr{P}^{+l,k,m}_{w} := \left( \{\, \mathscr{P}^{+n}_{w,l,k,m} \,\}_{n\in\mathbb{N}}, \langle_\rangle : p^{+n}_{\mathring{\varphi}} \mapsto n \right).$$

$\mathring{\Omega}^{+}(G)$ is the circuitous induced structure with counting of $G$. (def. A.2.3, p. 219)

$\mathfrak{G}_n$ denotes the set of all graphs of type $n$. (def. 4.5.1, p. 39)

**Lemma A.3.4**

Let $w, l, k, m \in \mathbb{N}$. Then $\mathscr{P}^{+l,k,m}_{w}$ is locally finite.

[1] Just as in definition 6.9.7, "for some" is equivalent to "for all", since all formulas in $\overline{\varphi}$ are type $n$ equivalent.

**Proof.** There is nothing new to show. Compare lemma 6.9.9.

□

# 4. Cooking with Induction

It remains only to show that our new family of predicates is again inductive. The remainder of the proof of Courcelle's Theorem makes no further mention of the particular family of predicates used, so it can be used verbatim afterwards.

For each of the theorems 6.9.10 to 6.9.12, we amend case 1 of its proof.

**Theorem A.4.1**

$\oplus$ is the disjoint sum. (def. 4.5.3, p. 40)

Let $w, l, k, m \in \mathbb{N}$. Then $\mathscr{P}_w^{+l,k,m}$ is $\mathfrak{F}_\oplus$-inductive.

**Proof.** *Case 1: $\varphi$ is atomic.*

*Case 1.4: $\varphi = \lambda_{\text{card}}^{i,p}(\chi)$ for some $i, p \in \mathbb{N}$, $i \leq p \leq l$.* Recall that while the width of $\varphi^{\text{a}}$ and $\varphi^{\text{b}}$ cannot exceed that of $\varphi$, the height of $\Phi$ suffers no such restriction. We thus set

$$\forall j \in \{0, \dots, i\} \colon \varphi_j^{\text{a}} := \lambda_{\text{card}}^{j,p}(\chi^a), \varphi_j^{\text{b}} := \lambda_{\text{card}}^{i-j,p}(\chi^b).$$

The reader should take a second to convince themselves that $\chi[\tau]$ has cardinality $i \mod p$ if and only if there is a $j$ such that $\chi^a[\tau_G]$ has cardinality $j \mod p$ and $\chi^b[\tau_{G'}]$ has cardinality $i - j \mod p$. We thus choose $\Phi := (\varphi_0^{\text{a}} \wedge \varphi_0^{\text{b}}) \vee \dots \vee (\varphi_i^{\text{a}} \wedge \varphi_i^{\text{b}})$.

As before, the width, number of variables, and set of free variables for the formulas $\varphi_0^{\text{a}}, \dots, \varphi_i^{\text{a}}, \varphi_0^{\text{b}}, \dots, \varphi_i^{\text{b}}$ is the same as for $\varphi$. The predicates used have changed, but all of them still obey the restriction imposed by $m$.

The remainder of the proof remains unchanged.

□

**Theorem A.4.2**

Let $w, l, k, m \in \mathbb{N}$. Then $\mathscr{P}_w^{+l,k,m}$ is $\left(\bigcup_{i\in\mathbb{N}} \bigcup_{j\in\mathbb{N}} {}_i^j\mathfrak{F}_{\leftrightarrows}\right)$-inductive.

**Proof.** The proof of theorem 6.9.11 works for this case without modification.

□

**Theorem A.4.3**

Let $w, l, k, m \in \mathbb{N}$. Then $\mathscr{P}_w^{+l,k,m}$ is $\left(\bigcup_{i\in\mathbb{N}} {}_i\mathfrak{F}_{\text{媲}}\right)$-inductive.

**Proof.** *Case 1: $\varphi$ is atomic.*

*Case 1.4: $\varphi = \lambda_{\text{card}}^{i,p}(\chi)$ for some $i, p \in \mathbb{N}$ with $i < p \leq l$, and some term $\chi$.* We set

$$i_{a,b} := i + 1 \mod p$$

(which is then again less than p) and

$$\begin{aligned}\varphi_0 := &\left(\neg\left(\delta_{\text{term}}^{\{a,b\}}(\varnothing) \sqsubseteq \chi_{a,b}\right) \wedge \lambda_{\text{card}}^{i,p}(\chi_{a,b})\right)\\ &\vee \left(\left(\delta_{\text{term}}^{\{a,b\}}(\varnothing) \sqsubseteq \chi_{a,b}\right) \wedge \lambda_{\text{card}}^{i_{a,b},p}(\chi_{a,b})\right)\end{aligned}$$

and $\Phi := \varphi_0$. If $a \notin \chi[\tau]$, then $\chi_{a,b}[\tau_{a,b}] = \chi[\tau]$ and the cardinality does not change. If $a \in \chi[\tau]$, then the expansion adds exactly one element to $\chi_{a,b}[\tau_{a,b}]$, raising the cardinality by one (modulo $p$).

The remainder of the proof remains unchanged.

□

We obtain the desired result.

**Corollary A.4.4**

Let $w, l, k, m \in \mathbb{N}$. Then $\mathscr{P}_w^{+l,k,m}$ is $\mathfrak{F}$-inductive.

# 5. Courcelle? Je la Connais à Peine!

With all preliminary results carried over to our counting frameworks, it is clear that the proofs for Courcelle's Theorem and its applications apply without modification. We hence state only the results.

**Theorem A.5.1**

Let $\varphi$ be a sentence of the counting monadic second-order language of graphs. Then for every $n \in \mathbb{N}$, the set $\{\, G \in \mathfrak{G}_n : \vDash_{\Omega(G)} \varphi \,\}$ is $\mathfrak{G}$-recognisable.

$\mathfrak{G}_n$ denotes the set of all graphs of type $n$. (def. 4.5.1, p. 39)

$\mathfrak{G}$ is the algebra of graphs. (def. 6.2.2, p. 93)

$\Omega(G)$ is the induced structure of $G$. (def. 5.4.2, p. 73)

**Theorem A.5.2**

Let $\mathscr{A} = (\mathscr{C}, \mathscr{O})$ be an algebra inherited from $\mathfrak{G}$, and let $\varphi$ be a sentence of the counting monadic second-order language of graphs. Then for every $n \in \mathbb{N}$, the set $\{\, G \in \mathscr{C}_n : \vDash_{\Omega(G)} \varphi \,\}$ is $\mathscr{A}$-recognisable.

**Theorem A.5.3**

Let $\mathscr{S}$ be a signature, let $\mathscr{A} = (\mathscr{C}, \mathscr{O})$ be an $\mathscr{S}$-algebra, and let the following properties be fulfilled:

- $\mathscr{A}$ is inherited from $\mathfrak{G}$.
- $\mathscr{A}$ is finitely expressible.
- $\mathscr{S}$ is weakly locally finite.

Let further $\varphi$ be a sentence of the counting monadic second-order language of graphs, and let $n \in \mathbb{N}$.

Then there exists a deterministic bottom-up finite tree automaton that takes as input a set $\mathscr{K}$ of $\mathscr{A}$-expressions with values in $\mathscr{C}_n$ and accepts an expression $e$ if and only if $\vDash_{\Omega(\mathrm{val}_{\mathscr{A}} e)} \varphi$, and every element of $\mathscr{C}_n$ is the value of some expression in $\mathscr{K}$.

$\mathrm{val}_{\mathscr{A}} e$ denotes the result of $e$ when evaluated in $\mathscr{A}$. (def. 6.3.3, p. 103)

### Theorem A.5.4

Let $\varphi$ be a sentence of the counting monadic second-order language of graphs, and let $k \in \mathbb{N}$. Then for every $n \in \mathbb{N}$, there exists a deterministic bottom-up finite tree automaton that takes as input a set $\mathscr{K}$ of $(\mathbb{N}, \mathfrak{O}^k)$-expressions with $\mathfrak{T}^k$-values in $\mathfrak{T}^k_n$ and accepts an expression $e$ if and only if $\models_{\Omega(\mathrm{val}_{\mathfrak{T}^k} e)} \varphi$, and every element of $\mathfrak{T}^k_n$ is the value of some expression in $\mathscr{K}$.

$\mathfrak{O}^k$ is the set of function symbols of $\mathfrak{T}^k$.
(def. 9.3.5, p. 186)

$\mathfrak{T}^k$ is the algebra of $k$-verdant graphs.
(def. 9.3.3, p. 185)

### Theorem A.5.5

Let $\varphi$ be a sentence of the counting monadic second-order language of graphs, and let $k \in \mathbb{N}$. Then for every $n \in \mathbb{N}$, there exists an algorithm which, given a $k$-verdant type $n$ graph $G = (V, E, \langle\!\langle _ \rangle\!\rangle, t)$, decides in time $\mathcal{O}(|V| + |E|)$ whether or not we have $\models_{\Omega(G)} \varphi$.

$\mathfrak{T}^k_n$ denotes the $k$-verdant graphs of type $n$.
(def. 9.3.3, p. 185)

$\mathrm{val}_{\mathfrak{T}^k} e$ denotes the result of $e$ when evaluated in $\mathfrak{T}^k$.
(def. 6.3.3, p. 103)

### Corollary A.5.6

Let $\varphi$ be a sentence of the counting monadic second-order language of graphs, and let $k \in \mathbb{N}$. Then there exists an algorithm which, given a graph $G = (V, E, \langle\!\langle _ \rangle\!\rangle, t)$ of tree-width at most $k$, decides in time $\mathcal{O}(|V| + |E|)$ whether or not we have $\models_{\Omega(G)} \varphi$.

$\Omega(\mathrm{val}_{\mathfrak{T}^k} e)$ is the induced structure of $\mathrm{val}_{\mathfrak{T}^k} e$.
(def. 5.4.2, p. 73)

### Corollary A.5.7

Let $k \in \mathbb{N}$. Then for any graph property that can be expressed in the counting monadic second-order logic of graphs, there exists an algorithm which, given a graph $G$ of tree-width at most $k$, decides in linear time whether or not $G$ fulfils said property.

$\mathfrak{Q}^k$ is the set of function symbols of $\mathfrak{P}^k$.
(def. 9.5.6, p. 205)

### Corollary A.5.8

Let $\varphi$ be a sentence of the counting monadic second-order language of graphs, and let $k \in \mathbb{N}$. Then for every $n \in \mathbb{N}$, there exists a deterministic bottom-up finite tree automaton that takes as input a set $\mathscr{K}$ of $(\mathbb{N}, \mathfrak{Q}^k)$-expressions with $\mathfrak{P}^k$-values in $\mathfrak{P}^k_n$ and accepts an expression $e$ if and only if $\models_{\Omega(\mathrm{val}_{\mathfrak{P}^k} e)} \varphi$, and every element of $\mathfrak{P}^k_n$ is the value of some expression in $\mathscr{K}$.

$\mathfrak{P}^k$ is the algebra of $k$-verdurous graphs.
(def. 9.5.2, p. 202)

$\mathfrak{P}^k_n$ denotes the $k$-verdurous graphs of type $n$.
(def. 9.5.2, p. 202)

**Theorem A.5.9**

Let $\varphi$ be a sentence of the counting monadic second-order language of graphs, and let $k \in \mathbb{N}$. Then for every $n \in \mathbb{N}$, there exists an algorithm which, given a $k$-verdurous type $n$ graph $G = (V, E, \langle\!\lvert _ \rvert\!\rangle, t)$, decides in time $\mathcal{O}(|V| + |E|)$ whether or not we have $\models_{\Omega(G)} \varphi$.

$\Omega(G)$ is the induced structure of $G$. (def. 5.4.2, p. 73)

**Corollary A.5.10**

Let $\varphi$ be a sentence of the counting monadic second-order language of graphs, and let $k \in \mathbb{N}$. Then there exists an algorithm which, given a graph $G = (V, E, \langle\!\lvert _ \rvert\!\rangle, t)$ of path-width at most $k$, decides in time $\mathcal{O}(|V| + |E|)$ whether or not we have $\models_{\Omega(G)} \varphi$.

**Corollary A.5.11**

Let $k \in \mathbb{N}$. Then for any graph property that can be expressed in the counting monadic second-order logic of graphs, there exists an algorithm which, given a graph $G$ of path-width at most $k$, decides in linear time whether or not $G$ fulfils said property.

**Theorem A.5.12**

Let $\varphi$ be a sentence of the counting monadic second-order language of graphs, and let $\mathscr{A} = (\mathscr{C}, \mathscr{O})$ be a finitely expressible algebra of weakly locally finite signature inherited from $\mathfrak{G}$. Then there exists an algorithm which determines in finite time whether or not there exists a graph $G \in \mathscr{C}_0$ with $\models_{\Omega(G)} \varphi$.

$\mathfrak{G}$ is the algebra of graphs. (def. 6.2.2, p. 93)

**Corollary A.5.13**

Let $\varphi$ be a sentence of the counting monadic second-order language of graphs, and let $k \in \mathbb{N}$. Then there exists an algorithm which determines in finite time whether or not there exists a graph $G$ of tree-width at most $k$ with $\models_{\Omega(G)} \varphi$.

**Corollary A.5.14**

Let $\varphi$ be a sentence of the counting monadic second-order language of graphs, and let $k \in \mathbb{N}$. Then there exists an algorithm which determines in finite time whether or not there exists a graph $G$ of path-width at most $k$ with $\models_{\Omega(G)} \varphi$.

$\Omega(G)$ is the induced structure of $G$. (def. 5.4.2, p. 73)

# Appendix B
# Loops

In our practical considerations in chapter 9, we have opted to forgo the case of loops, focusing only on graphs where all end points of an edge must be pairwise distinct. This choice has been made primarily to make the proofs easier to peruse, and the reader who cares not about such strange edges as can connect to the same vertex multiple times may safely skip this appendix.

However, all the versions of Courcelle's Theorem that we have seen still hold true for (hyper-)graphs with loops, so we present in this appendix a formal way to extend our algorithms to include those graphs.

## 1. It's Never Loopus

Why is it that our tree- and path-decomposition constructions do not admit loops in the first place? The problem is that we have restricted our all-powerful graph building tools to only entwining respectively blooming, and neither of those can create a loop. In the original algebra $\mathfrak{G}$, we would have simply taken the graph $\mathfrak{e}_2$ and fused its two vertices to create a loop of type 2.

$\mathfrak{G}$ is the algebra of graphs.
(def. 6.2.2, p. 93)

$\mathfrak{e}_2$ is the type 2 graph with 2 vertices and one edge.
(def. 4.5.6, p. 43)

媲 (tsureai, Japanese for *to marry*) denotes the source fusion. (def. 4.5.5, p. 42)

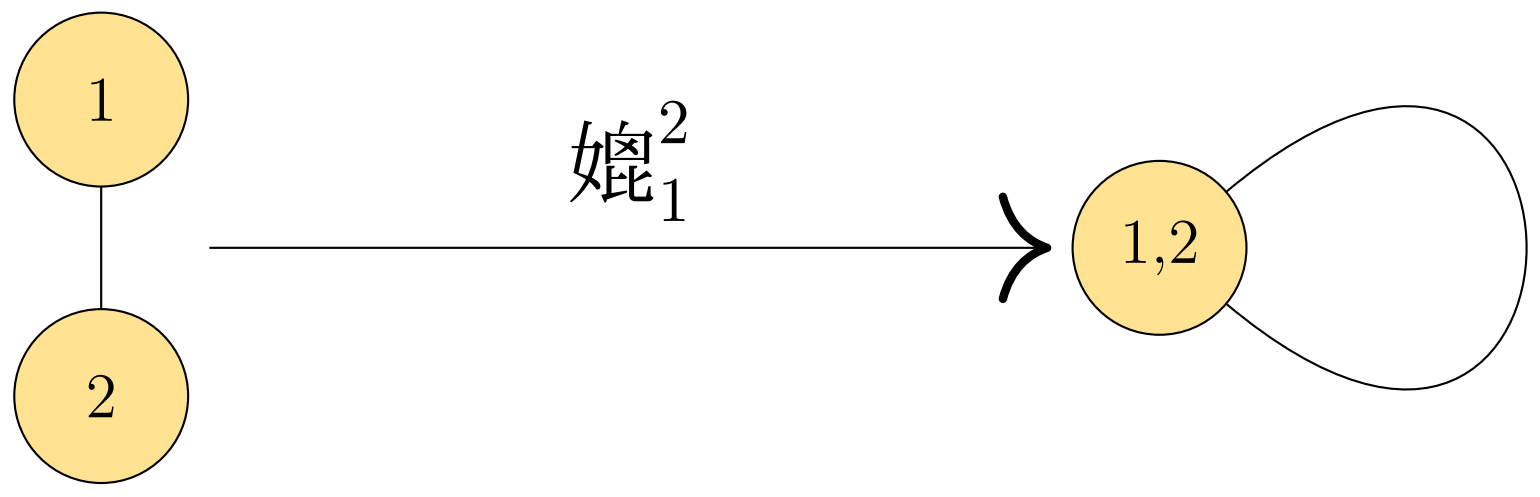

With entwining, we are not allowed to fuse two vertices from the same connected component, which makes the creation of loops impossible. The same holds for blooming.

Of course, there were good reasons for us not to allow this – the example on page 182 showed that fusing two connected vertices can easily increase the tree-width of a graph.

The trouble in said example stems from the fact that the terminals 1 and 2 used to be in different bags of the tree-decomposition, but fusion forces them to be in the same bag (because they are now the same vertex). Conversely, we should be *allowed* to fuse two vertices if they are in the same bag already, since then the tree-decomposition for the original graph is also a tree-decomposition for the new graph (with the missing vertex stricken from all bags):

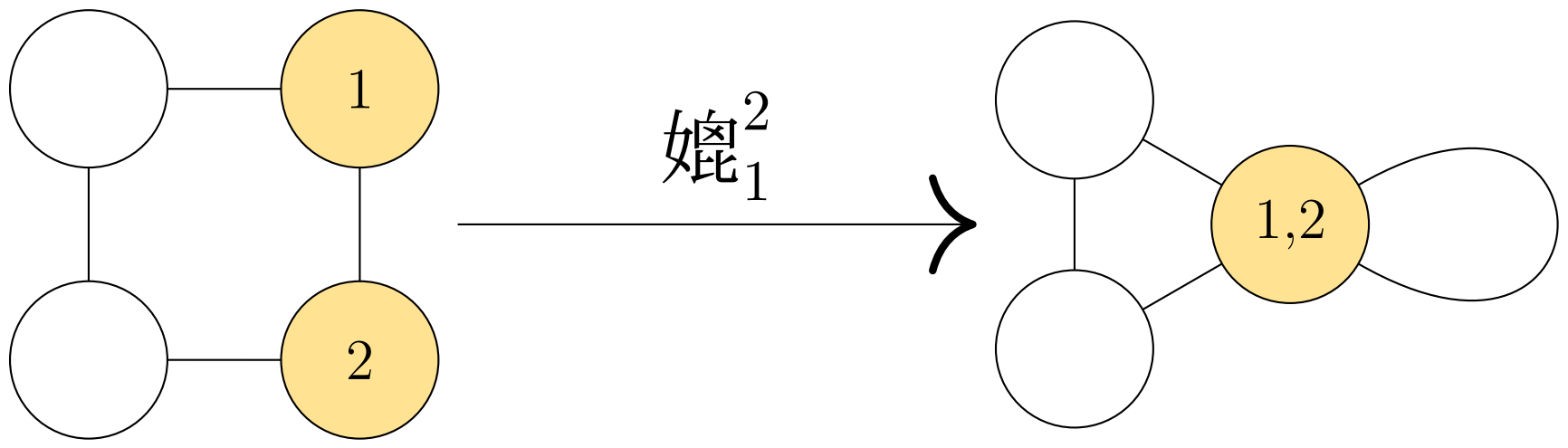

This construction is straightforward to formalise.

### Definition B.1.1

Let $n \in \mathbb{N}$, let $G = (V, E, \langle\!\langle _ \rangle\!\rangle, t) \in \mathfrak{G}_n$, and let $a, b \in \{1, \dots, n\}$. We set

$$\downarrow_{a,b} G := \begin{cases} 嫓_a^b G & \text{if } t(a), t(b) \text{ are adjacent} \\ G & \text{otherwise} \end{cases}$$

and call this the *collapse of $G$ over $a$ and $b$.*

$\mathfrak{G}_n$ denotes the set of all graphs of type $n$. (def. 4.5.1, p. 39)

嫓 (tsureai, Japanese for *to marry*) denotes the source fusion. (def. 4.5.5, p. 42)

And it even preserves verdancy and verdurousness.

### Lemma B.1.2

Let $n \in \mathbb{N}$, let $G \in \mathfrak{G}_n$, and let $a, b \in \{1, \dots, n\}$. Then we always have $\mathrm{tw}(\downarrow_{a,b} G) \leq \mathrm{tw}(G)$ and $\mathrm{pw}(\downarrow_{a,b} G) \leq \mathrm{pw}(G)$. Furthermore, if $G$ is $k$-verdant for some $k \in \mathbb{N}$, then $\downarrow_{a,b} G$ is $k$-verdant, and if $G$ is $k$-verdurous for some $k \in \mathbb{N}$, then $\downarrow_{a,b} G$ is $k$-verdurous.

**Proof.** Let $n \in \mathbb{N}$ and $G = (V, E, \langle\!\langle _ \rangle\!\rangle, t) \in \mathfrak{G}_n$. Let $Z = (T, X, b)$ be a tree-decomposition for $G$, and let $a, b \in \{1, \dots, n\}$. We assume without loss of generality that $t(a) \neq t(b)$, otherwise we have $\downarrow_{a,b} G = G$.

If the reader is familiar with graph minors, the following construction will be familiar to them, and they may safely replace it by their favourite minor argument.

*Case 1: There is no edge connecting $t(a)$ and $t(b)$.* Then by definition, we have $\downarrow_{a,b} G = G$ and there is nothing to show.

*Case 2: There is an edge $e \in E$ such that both $t(a)$ and $t(b)$ occur in $\langle\!\langle e \rangle\!\rangle$.* In this case, there must by definition of a tree-decomposition exist a node $v \in T$ with $t(a), t(b) \in b(v)$. We set

$$f\colon 2^V \to 2^{V \setminus \{t(b)\}}, A \mapsto \begin{cases} A & \text{if } t(b) \notin A \\ \{t(a)\} \cup A \setminus \{t(b)\} & \text{if } t(b) \in A \end{cases}$$

$2^V$ denotes the power set of $V$.

and define

$$Y := \{ f(x) : x \in X \}$$

and

$$c\colon T \to Y, v \mapsto f(b(v))$$

↓ denotes the graph collapse. (def. B.1.1, p. 233)

and claim that $(T, Y, c)$ is then a tree-decomposition for $\downarrow_{a,b} G$. Indeed: take a vertex $w \in \downarrow_{a,b} G$ which is contained in $c(v) \cap c(v')$ for some nodes $v, v' \in T$. If $w \neq t(a)$, a path from $v$ to $v'$ containing $w$ in each bag already existed since $Z$ was a tree-decomposition. If, on the other hand, $w = t(a)$, then in the original tree-decomposition $Z$ we either have $t(a) \in b(v) \cap b(v')$, $t(b) \in b(v) \cap b(v')$, or (without loss of generality) $t(a) \in b(v)$ and $t(b) \in b(v')$. In the first two cases, a path from $v$ to $v'$ containing $t(a)$ respectively $t(b)$ exists. In the third case, because $t(a)$ and $t(b)$ are adjacent, we can find a node $v'' \in T$ with $t(a), t(b) \in b(v'')$. Hence, there exists a path from $v$ to $v''$ containing $t(a)$ in each bag and a path from $v''$ to $v'$ containing $t(b)$ in each bag. Since in 媲$_a^b G$ we have $t(a) = t(b)$, we are done.

媲 (tsureai, Japanese for *to marry*) denotes the source fusion. (def. 4.5.5, p. 42)

We have thus shown that in each case, a tree-decomposition of at most the same width can be found, and since we did not change the bags except for identifying terminals, we know that if all terminals were in one bag of the tree-decomposition for $G$, the tree-decomposition for $\downarrow_{a,b} G$ we constructed again has all terminal vertices in one bag.

The proof for path-decompositions is entirely analogous.

□

However, as the reader has perhaps noticed, there is one small, but significant problem with our construction: it discriminates between graphs based on whether two vertices are adjacent, which neither fusion nor redefinition nor disjoint sum can do – and as such, it is not an inherited operation! Consequently, we cannot use it for Courcelle's Theorem.

# 2. Staying in the Loop

We solve our problem by taking a page out of the playbook of 2-uniform graphs. These admit only one "kind" of loop.

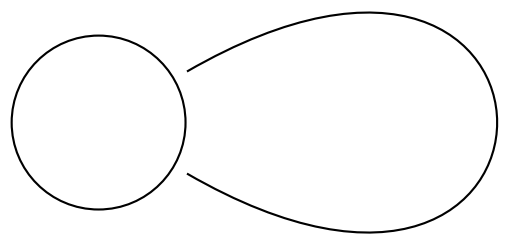

If our graph has loops, it suffices to add one trivial graph to our list (which previously contained only $\mathfrak{v}$ and $\mathfrak{e}_2$).

$\mathfrak{v}$ is the type 1 graph with one vertex.
(def. 4.5.6, p. 43)

$\mathfrak{e}_2$ is the type 2 graph with 2 vertices and one edge.
(def. 4.5.6, p. 43)

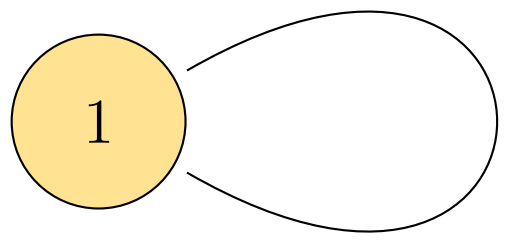

Whenever a graph contains loops, we can just fuse this new trivial graph to the correct vertices as many times as needed. This introduces only one new function symbol, namely the nullary one evaluating to this trivial graphs, and everything works out as expected.

More care is needed for proper hypergraphs, which admit many different "kinds" of loops.

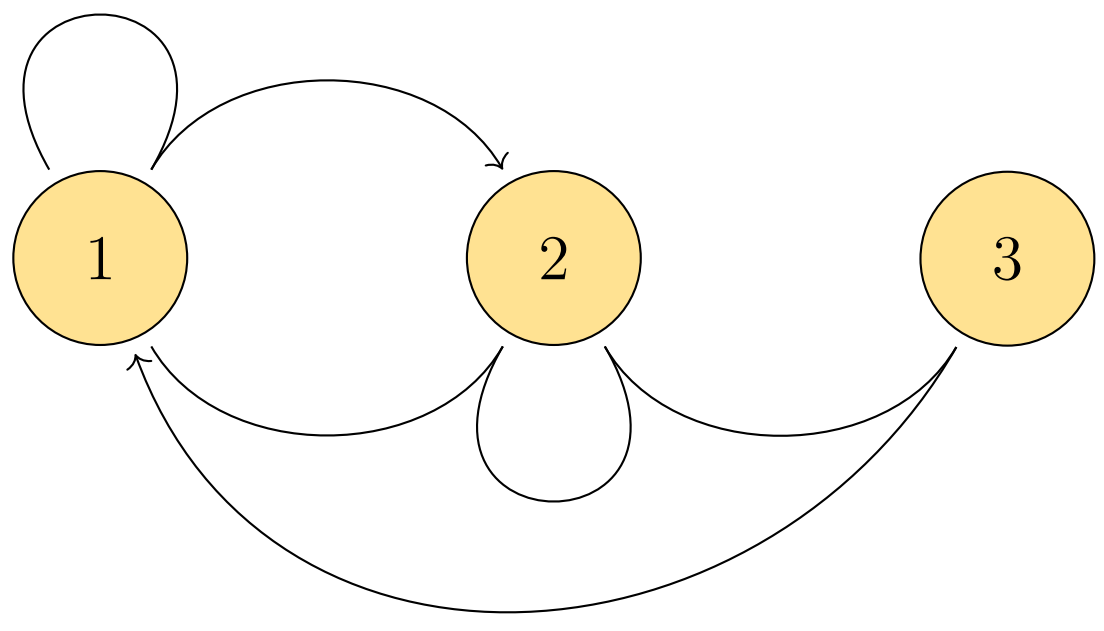

We use the brute-force method and add a nullary function symbol for each different kind of loop.

**Definition B.2.1**

> Let $n \in \mathbb{N}_{>0}$, and let $\omega \in \{1, \dots, n\}^*$ be a word that contains each number from 1 to $n$ at least once. We set $\mathfrak{l}_\omega$ to be the type $n$ graph
>
> $$(\{v_1, \dots, v_n\}, \{e\}, \langle\!\langle _ \rangle\!\rangle, t)$$
>
> with terminal function $t\colon i \mapsto v_i$ and $\langle\!\langle e \rangle\!\rangle = t^*(\omega)$. The set of all these graphs is denoted $\mathfrak{L}_n$. The set of all these graphs where $|\omega| \le k$ for some $k \in \mathbb{N}$ is denoted $\mathfrak{L}_n^k$.

For example, the (directed) graph above is the result of the expression

$$\mathfrak{l}_{112}\ {}^{2}_{3}\!\otimes^{3}_{\{1,2\}}\ \mathfrak{l}_{12231}.$$

$\otimes$ denotes the graph twine. (def. 9.3.1, p. 184)

$\mathfrak{O}^k$ is the set of function symbols of $\mathfrak{T}^k$. (def. 9.3.5, p. 186)

Notice that, for any $n \in \mathbb{N}_{>0}$, the set $\mathfrak{L}_n$ is countably infinite – adding these nullary function symbols to, for example, $(\mathbb{N}, \mathfrak{O}^k)$ would therefore not preserve weak local finiteness. The set $\mathfrak{L}_n^k$, however, is finite for any $k, n \in \mathbb{N}$. We must therefore restrict ourselves to a certain finite subset of edge types.[1]

# 3. Applying Some Loop Makes It Easy

As promised, we add the new trivial graphs to our signatures.

**Definition B.3.1**

> Let $k, n \in \mathbb{N}$, and let $c \in \mathbb{N}$ with $c \ge k + 1$. We write $\mathring{\mathfrak{T}}_n^{k,c}$ for the set of all $k$-verdant type $n$ graphs, possibly with loops, where no edge has type larger than $c$.

[1] The loop-free version $\mathfrak{T}^k$ cleverly avoids this problem since without loops, an edge of type $k$ always induces a bag of width at least $k - 1$, so a graph of tree-width $k$ cannot have an edge of type $k + 2$ there. With loops, an edge of type $k + 2$ might still only connect to one vertex, thus fitting into the tree-decomposition, but not into our finitely expressible algebra.

We require $c$ to be larger than the size of the largest bag in order to get $\mathfrak{T}_n^k \subseteq \mathring{\mathfrak{T}}_n^{k,c}$.

$\mathfrak{T}_n^k$ denotes the $k$-verdant graphs of type $n$. (def. 9.3.3, p. 185)

### Definition B.3.2

Let $k \in \mathbb{N}$, and let $c \in \mathbb{N}$ with $c \geq k+1$. We set

$$\mathring{\mathfrak{O}}^{k,c} := \mathfrak{O}^k \cup \bigcup_{i=1}^{k} \mathfrak{L}_i^c.$$

We denote

$$\mathring{\mathfrak{T}}^{k,c} := \left( \{ \mathring{\mathfrak{T}}_n^{k,c} \}_{n \in \mathbb{N}}, \{ \mathrm{Fl}(\overline{\mathfrak{O}})_{\mathfrak{f}} \big|_{\left(\bigcup_{n \in \mathbb{N}} \mathring{\mathfrak{T}}_n^{k,c}\right)^*} \}_{\mathfrak{f} \in \mathring{\mathfrak{O}}^{k,c}} \right).$$

$\mathfrak{O}^k$ is the set of function symbols of $\mathfrak{T}^k$. (def. 9.3.5, p. 186)

$\mathfrak{L}_i^c$ is the collection of all loops with $i$ vertices and edge type at most $c$. (def. B.2.1, p. 236)

$\overline{\mathfrak{O}}$ denotes the closure of $\mathfrak{O}$. (def. 6.10.6, p. 141)

Since we have added only nullary function symbols, and of those only finitely many per type, we immediately get the following.

$\mathrm{Fl}(\overline{\mathfrak{O}})$ denotes the flattening of $\overline{\mathfrak{O}}$. (def. 6.10.8, p. 143)

### Corollary B.3.3

Let $k \in \mathbb{N}$, and let $c \in \mathbb{N}$ with $c \geq k+1$. Then $\mathring{\mathfrak{T}}^{k,c}$ is an inherited $(\mathbb{N}, \mathring{\mathfrak{O}}^{k,c})$-algebra of $\mathfrak{G}$, and $(\mathbb{N}, \mathring{\mathfrak{O}}^{k,c})$ is weakly locally finite.

$\mathfrak{G}$ is the algebra of graphs. (def. 6.2.2, p. 93)

What remains to show is that this algebra is finitely expressible.

### Theorem B.3.4

Let $k, n \in \mathbb{N}$, and let $c \in \mathbb{N}$ with $c \geq k+1$. Let now $G \in \mathfrak{G}_n$ be a $k$-verdant graph that contains no edge of type larger than $c$. Then there exists a $\{0, \ldots, \max\{n, k+1\}\}$-local expression $e \in \|(\mathbb{N}, \mathring{\mathfrak{O}}^{k,c})\|$ with $\mathrm{val}_{\mathring{\mathfrak{T}}^{k,c}} e = G$.

$\mathfrak{G}_n$ denotes the set of all graphs of type $n$. (def. 4.5.1, p. 39)

$\|(\mathbb{N}, \mathring{\mathfrak{O}}^{k,c})\|$ denotes the set of expressions over $(\mathbb{N}, \mathring{\mathfrak{O}}^{k,c})$. (def. 6.3.3, p. 103)

**Proof.** The proof of theorem 9.3.9 can easily be adapted by adding, whenever a new vertex is introduced, all loops that can be added at this point. □

$\mathrm{val}_{\mathring{\mathfrak{T}}^{k,c}} e$ denotes the result of $e$ when evaluated in $\mathring{\mathfrak{T}}^{k,c}$. (def. 6.3.3, p. 103)

It is straightforward to add the loops to algorithm 6 – whenever a new vertex is introduced, we first add to it any and all loops that can be added at this stage.

Hence the results from section 9.3 carry over immediately.

**Theorem B.3.5**

Let $\varphi$ be a sentence of the monadic second-order language of graphs, and let $k, c \in \mathbb{N}$. Then for every $n \in \mathbb{N}$, there exists an algorithm which, given a $k$-verdant type $n$ graph $G = (V, E, \langle\!\!|_|\!\!\rangle, t)$ with no edge of type larger than $c$, decides in time $\mathcal{O}(|V| + |E|)$ whether or not we have $\vDash_{\Omega(G)} \varphi$.

$\Omega(G)$ is the induced structure of $G$. (def. 5.4.2, p. 73)

**Corollary B.3.6**

Let $\varphi$ be a sentence of the monadic second-order language of graphs, and let $k, c \in \mathbb{N}$. Then there exists an algorithm which, given a graph $G = (V, E, \langle\!\!|_|\!\!\rangle, t)$ of tree-width at most $k$ with no edge of type larger than $c$, decides in time $\mathcal{O}(|V| + |E|)$ whether or not we have $\vDash_{\Omega(G)} \varphi$.

**Corollary B.3.7**

Let $k, c \in \mathbb{N}$. Then for any graph property that can be expressed in the monadic second-order logic of graphs, there exists an algorithm which, given a graph $G$ of tree-width at most $k$ with no edge of type larger than $c$, decides in linear time whether or not $G$ fulfils said property.

Generalising the results about paths of constantly-bounded path-width should be straightforward using this appendix as a guide. We state only the end result.

**Theorem B.3.8**

Let $\varphi$ be a sentence of the monadic second-order language of graphs, and let $k, c \in \mathbb{N}$. Then for every $n \in \mathbb{N}$, there exists an algorithm which, given a $k$-verdurous type $n$ graph $G = (V, E, \langle\!\!|_|\!\!\rangle, t)$ with no edge of type larger than $c$, decides in time $\mathcal{O}(|V| + |E|)$ whether or not we have $\vDash_{\Omega(G)} \varphi$.

**Corollary B.3.9**

Let $\varphi$ be a sentence of the monadic second-order language of graphs, and let $k, c \in \mathbb{N}$. Then there exists an algorithm which, given a graph $G = (V, E, \langle\!|_|\!\rangle, t)$ of path-width at most $k$ with no edge of type larger than $c$, decides in time $\mathcal{O}(|V| + |E|)$ whether or not we have $\models_{\Omega(G)} \varphi$.

**Corollary B.3.10**

Let $k, c \in \mathbb{N}$. Then for any graph property that can be expressed in the monadic second-order logic of graphs, there exists an algorithm which, given a graph $G$ of path-width at most $k$ with no edge of type larger than $c$, decides in linear time whether or not $G$ fulfils said property.

# Index

The page where a term is defined is printed in **bold**. References to sidenotes or footnotes are shown in *italics*.